\documentclass[11pt]{amsart}         
\usepackage{amscd}                   
\usepackage{amssymb} 

\newtheorem{thm}{Theorem}[section]
\newtheorem{lem}[thm]{Lemma}

\newtheorem{cor}[thm]{Corollary}
\newtheorem{claim}[thm]{Claim}
\newtheorem{prop}[thm]{Proposition}

\theoremstyle{definition}

\renewcommand{\thecase}{}

\newtheorem{conj}[thm]{Conjecture}

\newtheorem{defn}[thm]{Definition}

\newtheorem{rmk}[thm]{Remark} 
\renewcommand{\thestep}{}
\newtheorem{prob}[thm]{Problem}

\theoremstyle{remark}

\makeatletter
\def\alphenumi{
  \def\theenumi{\alph{enumi}}
  \def\p@enumi{\theenumi}
  \def\labelenumi{(\@alph\c@enumi)}}
\makeatother

\makeatletter
\def\thecase{\@arabic\c@case}
\makeatother

\numberwithin{equation}{section}

\makeatletter
\def\thestep{\@arabic\c@step}
\makeatother

\newenvironment{pf}{\begin{proof}[\proofname]}{\end{proof}}
\newenvironment{pf*}[1]{\begin{proof}[#1]}{\end{proof}}

\newcommand\barB{{\bar{B}}}

\newcommand\barM{{\bar{M}}}

\newcommand\barOm{{\bar\Omega}}

\newcommand\CC{\mathbb{C}}

\newcommand\PP{\mathbb{P}}

\newcommand\RR{\mathbb{R}}

\newcommand\ZZ{\mathbb{Z}}

\newcommand\bdelta{{\boldsymbol{\delta}}}

\newcommand\bfeta{{\boldsymbol{\eta}}}
\newcommand\bga{{\boldsymbol{\gamma}}}
\newcommand\bgamma{{\boldsymbol{\gamma}}}
\newcommand\bchi{{\boldsymbol{\chi}}}

\newcommand\blambda{{\boldsymbol{\lambda}}}

\newcommand\bvarphi{{\boldsymbol{\varphi}}}

\newcommand\bPhi{{\boldsymbol{\Phi}}}
\newcommand\bpsi{{\boldsymbol{\psi}}}

\newcommand\bxi{{\boldsymbol{\xi}}}
\newcommand\bXi{{\boldsymbol{\Xi}}}

\newcommand\bA{{\mathbf{A}}}

\newcommand\bD{{\mathbf{D}}}

\newcommand\bFr{{\mathbf{Fr}}}

\newcommand\bG{{\mathbf{G}}}
\newcommand\bGl{{\mathbf{Gl}}}

\newcommand\bx{{\mathbf{x}}}

\newcommand\bZ{{\mathbf{Z}}}

\newcommand{\cov}{\nabla}

\newcommand{\rd}{\partial}

\newcommand\thalf{{\textstyle{\frac{1}{2}}}}

\newcommand\tthreehalf{{\textstyle{\frac{3}{2}}}}
\newcommand\half{{{\frac{1}{2}}}}

\newcommand\quarter{{{\frac{1}{4}}}}

\newcommand\eighth{{{\frac{1}{8}}}}

\newcommand\fB{{\mathfrak{B}}}

\newcommand\fg{{\mathfrak{g}}}

\newcommand\fs{{\mathfrak{s}}}
\newcommand\fS{{\mathfrak{S}}}
\newcommand\ft{{\mathfrak{t}}}

\newcommand\fV{{\mathfrak{V}}}

\newcommand\be{\beta}

\newcommand\de{\delta}
\newcommand\eps{\varepsilon}

\newcommand\Ga{\Gamma}
\newcommand\la{\lambda}
\newcommand\La{\Lambda}

\newcommand\Om{\Omega}
\newcommand\si{\sigma}

\newcommand\so{{\mathfrak{s}\mathfrak{o}}}

\newcommand\su{{\mathfrak{s}\mathfrak{u}}}
\newcommand\fu{{\mathfrak{u}}}

\newcommand\GL{\operatorname{GL}}

\newcommand\PU{\operatorname{PU}}

\newcommand\SO{\operatorname{SO}}

\newcommand\SU{\operatorname{SU}}
\newcommand\U{\operatorname{U}}

\newcommand\less{\setminus}

\newcommand{\8}{\infty}

\newcommand\asd{{\operatorname{asd}}}

\newcommand\Cl{\operatorname{Cl}}

\newcommand\dist{\operatorname{dist}}

\newcommand\End{\operatorname{End}}

\newcommand\Fr{\operatorname{Fr}}

\newcommand\Hom{\operatorname{Hom}}

\newcommand\Ind{\operatorname{Index}}

\newcommand\Ker{\operatorname{Ker}}

\newcommand\loc{\operatorname{loc}}

\newcommand\PD{\operatorname{PD}}

\newcommand\Ric{\operatorname{Ric}}

\newcommand\Ran{\operatorname{Ran}}

\newcommand\Spec{\operatorname{Spec}}

\newcommand\Stab{\operatorname{Stab}}

\newcommand\supp{\operatorname{supp}}
\newcommand\sw{{\operatorname{sw}}}
\newcommand\SW{SW}
\newcommand\Sym{\operatorname{Sym}}
\newcommand\Tor{\operatorname{Tor}}

\newcommand\vol{\operatorname{vol}}
\newcommand\Vol{\operatorname{Vol}}

\newcommand\ad{{\mathrm{ad}\,}}

\newcommand\id{{\mathrm{id}}}
\newcommand\LC{{\mathrm{LC}}}

\newcommand\Rm{{\mathrm{Rm}}}

\newcommand\spinc{\text{$\text{spin}^c$ }}
\newcommand\spinu{\text{$\text{spin}^u$ }}
\newcommand\Spinc{\text{$\mathrm{Spin}^c$}}

\newcommand\sA{{\mathcal{A}}}
\newcommand\sB{{\mathcal{B}}}
\newcommand\sC{{\mathcal{C}}}
\newcommand\sD{{\mathcal{D}}}

\newcommand\sG{{\mathcal{G}}}
\newcommand\sH{{\mathcal{H}}}

\newcommand\sK{{\mathcal{K}}}
\newcommand\sL{{\mathcal{L}}}

\newcommand\sR{{\mathcal{R}}}

\newcommand\sU{{\mathcal{U}}}

\newcommand\sW{{\mathcal{W}}}

\newcommand\tsC{{\tilde\sC}}

\newcommand\tg{{\tilde g}}

\newcommand\tM{{\tilde M}}

\begin{document}
\title[$\PU(2)$ monopoles. III]
{\boldmath{$\PU(2)$} monopoles. III: Existence of gluing and obstruction maps}
\author[Paul M. N. Feehan]{Paul M. N. Feehan}
\address{Department of Mathematics\\
Ohio State University\\
Columbus, OH 43210}
\email{feehan@math.ohio-state.edu; http://www.math.ohio-state.edu/$\sim$feehan}
\curraddr{School of Mathematics\\
Institute for Advanced Study\\
Olden Lane\\
Princeton, NJ 08540}
\email{feehan@math.ias.edu; http://www.math.ias.edu/$\sim$feehan}
\author[Thomas G. Leness]{Thomas G. Leness}
\address{Department of Mathematics\\
Florida International University\\
Miami, FL 33199}
\email{lenesst@fiu.edu}
\dedicatory{}
\thanks{The first author was supported in part by an NSF Mathematical 
Sciences Postdoctoral Fellowship under grant DMS 9306061 and by NSF grants
DMS 9704174 and, through the Institute for Advanced Study, by DMS 9729992}
\date{This version: July 20, 1999. First version: July 15,
1999. math.DG/9907107.} 
\begin{abstract}
This is the third installment
in our series of articles on the application of the PU(2)
monopole equations to prove Witten's conjecture concerning the relation
between the Donaldson and Seiberg-Witten invariants of smooth
four-manifolds. The moduli space of solutions to the PU(2) monopole
equations provides a noncompact cobordism between links of compact moduli
spaces of U(1) monopoles of Seiberg-Witten type and the moduli space of
anti-self-dual SO(3) connections, which appear as singularities in this
larger moduli
 space. In this paper we prove the first part of a general gluing theorem
for PU(2) monopoles. The ultimate purpose of the gluing theorem is to
provide topological models for neighborhoods of ideal Seiberg-Witten moduli
spaces appearing in lower levels of the Uhlenbeck compactification of the
moduli space of PU(2) monopoles and thus permit calculations of their
contributions to Donaldson invariants using the PU(2)-monopole cobordism.
\end{abstract}
\maketitle
%\pagenumbering{roman}
\setcounter{tocdepth}{1}
\tableofcontents

%\pagenumbering{arabic}

%end of file
%file:  intro3.tex 

%\newpage 
%\pagenumbering{arabic}
\section{Introduction} 
\label{subsec:Introduction}
The Pidstrigatch-Tyurin program \cite{PTLectures}, \cite{PTLocal} to relate
the Donaldson and Seiberg-Witten invariants, as considered in
\cite{FeehanGenericMetric}, \cite{FL1}, \cite{FLGeorgia}, \cite{FL2}
requires a gluing theory, for $\PU(2)$ monopoles, possessing the properties
we establish here and in the sequel \cite{FL4}. The purpose of gluing, as
we outlined in \cite{FLGeorgia}, is to provide topological models for
neighborhoods of ideal Seiberg-Witten moduli spaces appearing in lower
levels of the Uhlenbeck compactification of the moduli space of $\PU(2)$
monopoles and thus enable calculations of their contributions to Donaldson
invariants using the $\PU(2)$-monopole cobordism \cite{FLConj}.

\subsection{Statement of results} 
The main purpose of the present article and its companion \cite{FL4} is to
prove a general gluing theorem for $\PU(2)$ monopoles, adequate for the
topological calculations mentioned in \cite{FLGeorgia} and carried out in
detail in \cite{FLLevelOne}, \cite{FLConj}.

Throughout our work we let $(X,g)$ be a closed, connected, oriented, $C^\8$
Riemannian four-manifold with $b^+(X)>0$ and \spinc structure
$\fs_0=(\rho,W^+,W^-)$. We require that the perturbation parameters
defining the $\PU(2)$ monopole equations \eqref{eq:PT} be generic, so the
transversality results of \cite{DK}, \cite{FeehanGenericMetric}, and
\cite{FU} ensure that the moduli spaces of $\PU(2)$ monopoles,
Seiberg-Witten $\U(1)$ monopoles, and anti-self-dual $\SO(3)$ connections
possess the usual smoothness properties and have the expected dimension.
 
While we shall direct the reader to \S \ref{sec:Prelim},
\S \ref{sec:Splicing}, and \S \ref{sec:Existence} 
for detailed definitions and notation, we note that $M^{\sw}(\fs_0)$ is the
moduli space of Seiberg-Witten monopoles on $X$, that $M(\ft)$ is the
moduli space of $\PU(2)$ monopoles on $X$ for a \spinu structure
$\ft=(\rho,W^+,W^-,E)$, where $E$ is a Hermitian, rank-two bundle over $X$, and
that $M(\ft,\mu)$ is a finite-dimensional thickened or virtual moduli space
of $\PU(2)$ monopoles which contains $M(\ft)$, as a submanifold away from
singularities, with $\mu>0$ a small-eigenvalue bound arising in the
extended $\PU(2)$ monopole equations which define $M(\ft,\mu)$.  We let
$w\in H^2(X;\ZZ)$ be a class such that no $\SO(3)$ bundle over $X$ with
second Stiefel-Whitney class equal to $w\pmod{2}$ admits a flat connection.
The open subset $M^{*,0}(\ft)\subset M(\ft)$ is represented by pairs such
that the $\SO(3)$ connection is not reducible and the spinor is not
identically zero.

\begin{thm} 
\label{thm:GluingTheorem1} 
Let $E$ be a rank-two, Hermitian vector bundle over $X$ with $c_1(E)=w$ and
$c_2(E)-\quarter c_1(E)^2 = \kappa\geq 1$.  Let
$0<\ell<\lfloor\kappa\rfloor\leq\kappa$ be an integer, let $E_{\ell}$ be
the rank-two, Hermitian vector bundle over $X$ with $\det E_{\ell} = \det
E$ and $c_2(E_{\ell}) = c_2(E)-\ell$. Let $\ft = (\rho,W^+,W^-,E)$ and
$\ft_{\ell} = (\rho,W^+,W^-,E_{\ell})$. Let $\Sigma\subset\Sym^\ell(X)$ be a
smooth stratum and let $\sU_{\ell,\mu}\Subset M(\ft_{\ell},\mu)$ be a
finite-dimensional, precompact, open, $S^1$-invariant subset. The space
$\sU_{\ell,\mu}$ is a tubular neighborhood of $\sU_{\ell,\mu}\cap
M^{*,0}(\ft_{\ell})$, of size $\eps_\ell>0$ via the tubular distance function
\eqref{eq:DeformCondn}. Then, for a small enough 
positive constants $\eps_\ell$ and $\lambda_0$, there are
\begin{itemize}
\item 
A compact Lie group $\bG(\Sigma)$, an $S^1$-equivariant $C^\8$ principle
$\bG(\Sigma)$-bundle $\bFr(\sU_{\ell,\mu},\Sigma)\to
\sU_{\ell,\mu}\times\Sigma$, a universal topological $\bG(\Sigma)$-space
$\bZ(\Sigma)$, and an $S^1$-equivariant fiber bundle
$$
\bGl(\sU_{\ell,\mu},\Sigma,\lambda_0) 
= 
\bFr(\sU_{\ell,\mu},\Sigma)\times_{\bG(\Sigma)}\bZ(\Sigma)
\to
\sU_{\ell,\mu}\times\Sigma. 
$$
Here, $\lambda_0$ is a sufficiently small positive constant, depending at
most on $(g,\kappa,\sU_{\ell,\mu})$.
\item 
An $S^1$-equivariant $C^\8$ map
$\bgamma_{\mu,\Sigma}:\bGl(\sU_{\ell,\mu},\Sigma,\lambda_0)
\to M^{*,0}(\ft,\mu)$, where $\bGl^+(\sU_{\ell,\mu},\Sigma,\lambda_0)
\subset \bGl(\sU_{\ell,\mu},\Sigma,\lambda_0)$ is the open subset
\eqref{eq:ConstrainedGluingDataBundle}.
\item 
An $S^1$-equivariant $C^\8$ section $\bchi_{\mu,\Sigma}$ of an
$S^1$-equivariant  $C^\8$ vector bundle $\Xi_{\mu}$ over
$\bGl^+(\sU_{\ell,\mu},\Sigma,\lambda_0)$. The bundle
$\Xi_{\mu}$ has real rank 
equal to $2\ell$ plus the real codimension of $M^{*,0}(\ft)\subset
M^{*,0}(\ft,\mu)$. 
\end{itemize}
Together, the map $\bgamma_{\mu,\Sigma}$ and section $\bchi_{\mu,\Sigma}$ have
the property that
$$
\bgamma_{\mu,\Sigma}(\bchi_{\mu,\Sigma}^{-1}(0)) \subset M^{*,0}(\ft).
$$
\end{thm}

We call $\bGl(\sU_{\ell,\mu},\Sigma,\lambda_0)$ the {\em bundle of local
gluing data}, $\bgamma_{\mu,\Sigma}$ the {\em local gluing map\/},
$\Xi_{\mu}$ the {\em local obstruction bundle\/}, and
$\bchi_{\mu,\Sigma}$ the {\em local obstruction section\/}. The gluing maps
and obstruction sections are defined in \S \ref{sec:Splicing} and \S
\ref{sec:Existence}. Theorem
\ref{thm:GluingTheorem1}, along with many other results, was announced in
\cite{FeehanStrasbourgTalk}. An earlier version of this result, leaning
more toward the gluing method of Donaldson \cite{DonConn}, \cite{DS},
\cite{DK} --- which is less satisfactory for our application (see the
explanations below) --- was announced in \cite{FeehanOberwolfachTalk} and
the corresponding preprint \cite{FL0} (containing a detailed proof) was
written during 1995 and distributed to a small audience during Spring 1996.
The solution to the gluing problem described here 
and in \cite{FL4} follows more in the
tradition of Taubes' approach \cite{TauSelfDual}, \cite{TauIndef},
\cite{TauStable} to the problem of gluing anti-self-dual connections.
Though better suited to our application \cite{FLConj},
it has long appeared that Taubes' method was much more difficult to adapt
to the case of $\PU(2)$ monopoles than that of Donaldson, due primarily to
the technical Problem \ref{prob:FA+FA-BochnerProblem} discussed in \S 
\ref{subsec:TechDiffForGluingPUMonopoles}
below and, indeed, its solution
eluded us for some time despite considerable effort.

Naturally, many statements of gluing theorems for anti-self-dual
connections have appeared before in the literature
\cite{DonConn}, \cite{DK}, \cite{DS}, \cite{MrowkaThesis}, 
\cite{TauSelfDual}, \cite{TauIndef}, \cite{TauFrame}, \cite{TauStable}, so
the reader may reasonably ask what is new and different about 
$\PU(2)$ monopoles. The gluing-theorem statements that are closest to those
needed for a proof of the Witten conjecture are given by $\PU(2)$-monopole
analogues of the gluing
Theorems 3.4.10 and 3.4.17 (for anti-self-dual connections)
in \cite{FrM}: both assertions are based on
work of Taubes, principally \cite{TauSelfDual}, \cite{TauIndef},
\cite{TauFrame}, \cite{TauStable}. In any event, the
difficulties which are peculiar to gluing $\PU(2)$ monopoles are discussed
in \S \ref{subsec:TechDiffForGluingPUMonopoles}.
Furthermore, there are several important properties which are required of a
useful gluing theory --- in any context, but particularly for applications
to Witten's conjecture --- which are {\em not\/} asserted by Theorem
\ref{thm:GluingTheorem1}. These extensions are discussed in \S
\ref{subsec:Sequels} and are derived in the companion article \cite{FL4}.

\subsection{PU(2) monopoles and Witten's conjecture}
\label{subsec:WittensConjecture}
Recall that a closed, smooth four-manifold $X$ has Kronheimer-Mrowka simple
type provided the Donaldson invariants corresponding to products $z$ of
homology classes in $H_\bullet(X;\ZZ)$ and a generator $x\in H_0(X;\ZZ)$
are related by $D_X^w(x^2z)=4D_X^w(z)$. Kronheimer and Mrowka
\cite{KMStructure} showed that the Donaldson series of a four-manifold of
Kronheimer-Mrowka simple type with $b^1(X)=0$ and odd $b^+(X)\ge 3$ is
given by
\begin{equation}
\bD^w_X = e^{Q_X/2}\sum_{r=1}^s(-1)^{(w^2+w\cdot K_r)/2}a_r e^{K_r},
\label{eq:KMFormula}
\end{equation}
where $w\in H_2(X;\ZZ)$, $Q_X$ is the intersection form on
$H_2(X;\ZZ)$, the coefficients $a_r$ are non-zero rational numbers, and the
$K_r\in H^2(X;\ZZ)$ are the Kronheimer-Mrowka basic classes. They also
conjectured the shape of the general structure theorem for four-manifolds
of possibly non-simple type \cite{KMGeneralType}.

\begin{conj}[Witten]\cite{Witten}
The four-manifold $X$ has Kronheimer-Mrowka simple type if and only if it has
Seiberg-Witten simple type. If $X$ has simple type, then the
Kronheimer-Mrowka basic classes coincide with the Seiberg-Witten basic
classes and
\begin{equation}
\bD^w_X
=2^{2-c(X)}e^{Q_X/2}
\sum_{\fs\in\Spinc(X)}(-1)^{(w^2+w\cdot c_1(\fs))/2}\SW_X(\fs)e^{c_1(\fs)}.
\label{eq:WittenFormula}
\end{equation}
\end{conj}

The quantum field theory argument giving the relation 
\eqref{eq:WittenFormula} when $b^+(X)\ge
3$ has been extended by Moore and Witten \cite{MooreWitten} to
allow $b^+(X) \ge 1$, $b^1(X)\ge 0$, and four-manifolds $X$ of possibly
non-simple type.

The idea of Pidstrigach and Tyurin \cite{PTLocal} for proving Witten's
conjecture
\cite{MooreWitten}, \cite{Witten} has two principal steps. The first step is
to try to show that the Donaldson invariants are given by a sum over the
Seiberg-Witten basic classes of terms given by universal polynomials in the
intersection form $Q_X$ and the classes $c_1(\fs)$ and $\Lambda$, with
coefficients depending only on products of the classes $c_1(\fs)$, $\Lambda$,
and $w$, the Donaldson-invariant degree $\delta$, the Seiberg-Witten
invariant $\SW_X(\fs)$, and the signature $\sigma(X)$ and Euler
characteristic $\chi(X)$ (see Conjectures 4.1 and 4.2 in
\cite{FLGeorgia}).  While this ``homotopy version'' of Witten's conjecture 
would not immediately yield the formula \eqref{eq:WittenFormula}, it does
imply that the Donaldson invariants are completely determined by and
computable (in principle) from the Seiberg-Witten invariants.

The existence of a universal formula is expected
because the moduli space of $\PU(2)$ monopoles contains both the Donaldson
moduli space of anti-self-dual connections and moduli spaces of $\U(1)$
monopoles of Seiberg-Witten type and so provides a cobordism between
suitably defined, codimension-one links of these moduli spaces.  The second
step is to try to determine these universal formulas explicitly via known
examples and recursion relations, including those obtained from blow-up
formulas for Donaldson and Seiberg-Witten invariants \cite{FSTurkish},
\cite{FSBlowUp}.  The work of G\"ottsche
\cite{Goettsche} suggests that the second part of this strategy should
ultimately succeed as, for example, he was able to compute the
wall-crossing formula for the Donaldson invariants of simply-connected
four-manifolds with $b^+(X)=1$ and $b^1(X)=0$, under the assumption that
the Kotschick-Morgan conjecture
\cite{KotschickMorgan} holds for such four-manifolds; that conjecture
asserts that for four-manifolds with $b^+(X)=1$, the Donaldson invariants
computed using metrics lying in two different chambers of the positive cone
in $H^2(X;\RR)$ differ by a universal polynomial whose coefficients depend
only on homotopy data \cite{KotschickMorgan}.

\subsection{Some technical difficulties associated with gluing PU(2) monopoles}
\label{subsec:TechDiffForGluingPUMonopoles}
Aside from certain aspects of gluing theory for instantons which were not
completely addressed prior to the Seiberg-Witten revolution, there are a
number of features of the $\PU(2)$ monopole equations which makes the
gluing theory especially challenging. Since it will take a sequel
\cite{FL4} to complete the proofs of the additional properties required
(the embedding property, surjectivity, and Uhlenbeck continuity) for our
application \cite{FLConj}, it seems worthwhile to motivate the lengthy
development by commenting on a few of the technical difficulties at the
outset. 

The splicing (or ``pre-gluing'' or ``cut-and-paste'') construction yields
families of approximate $\PU(2)$ monopoles, $(A,\Phi)$,
$$
\fS(A,\Phi) \approx 0,
$$
where the map $\fS(\cdot,\cdot)$ defines the $\PU(2)$ monopole equations
\eqref{eq:PT}, $A$ is a connection on the $\SO(3)$ bundle $\fg_E$, and
$\Phi$ a section of $V^+=W^+\otimes E$. Given an approximate $\PU(2)$
monopole $(A,\Phi)$, we would like to solve for a deformation $(a,\phi)$
such that
\begin{equation}
\label{eq:IntroBasicPTDef}
\fS(A+a,\Phi+\phi) 
= 0.
\end{equation}
The preceding equation can be rewritten as 
\begin{equation}
\label{eq:IntroFirstPTDef}
d_{A,\Phi}^1(a,\phi) + \{(a,\phi),(a,\phi)\} 
= -\fS(A,\Phi),
\end{equation}
where $d_{A,\Phi}^1 = (D\fS)_{A,\Phi}$. We now seek a solution
$(a,\phi)= d_{A,\Phi}^{1,*}(v,\psi) \in (\Ker d_{A,\Phi}^1)^\perp$, so
\begin{equation}
\label{eq:IntroSecondPTDef}
d_{A,\Phi}^1d_{A,\Phi}^{1,*}(v,\psi) 
+ 
\{d_{A,\Phi}^{1,*}(v,\psi),d_{A,\Phi}^{1,*}(v,\psi)\} 
= 
-\fS(A,\Phi).
\end{equation}
This brings us to a familiar difficulty in gauge theory when attempting to
solve non-linear partial differential equations by splicing solutions:

\begin{prob}
\label{prob:SmallEigenvalues}
The Laplacian $d_{A,\Phi}^1d_{A,\Phi}^{1,*}$ is {\em not} uniformly
invertible. It has {\em small eigenvalues}, tending to {\em zero\/} as
$(A,\Phi)$ bubbles and $\fS(A,\Phi)$ converges to zero.
\end{prob}

Small eigenvalues occur here because the adjoint Dirac operators
$D_{A_i}^*$ over $S^4$ have nontrivial kernels, as does the adjoint
linearization $d_{A_0,\Phi_0}^{1,*}$ over $X$, where $(A_0,\Phi_0)$ is a
background pair over $X$ and $A_i$ is an instanton over $S^4$.  Instead,
motivated by \cite{TauIndef}, we should try to solve an {\em extended
$\PU(2)$ monopole equation},
$$
\Pi_{A,\Phi,\mu}^\perp\fS(A+a,\Phi+\phi) = 0,
$$
or equivalently,
\begin{equation}
\label{eq:IntroLongSecOrderExtPUMonEqnForvpsi}
d_{A,\Phi}^1d_{A,\Phi}^{1,*}(v,\psi) 
+ \Pi_{A,\Phi,\mu}^\perp
\{d_{A,\Phi}^{1,*}(v,\psi), d_{A,\Phi}^{1,*}(v,\psi)\}
= 
- \Pi_{A,\Phi,\mu}^\perp\fS(A,\Phi),
\end{equation}
for deformations $(a,\phi)= d_{A,\Phi}^{1,*}(v,\psi)$ such that
$\Pi_{A,\Phi,\mu}(v,\psi)=0$, where $\Pi_{A,\Phi,\mu}$ is the
$L^2$-orthogonal projection onto the eigenspaces of
$d_{A,\Phi}^1d_{A,\Phi}^{1,*}$ with eigenvalues less than $\mu>0$.  To
solve equation \eqref{eq:LongSecOrderExtPUMonEqnForvpsi} we need a uniform
estimate for a right inverse $P_{A,\Phi,\mu}$ of the operator
$d_{A,\Phi}^1\circ\Pi_{A,\Phi,\mu}^\perp$. For this purpose we try to
adapt Taubes' method in \cite{TauSelfDual}, which relies on a clever use of
the Bochner formula for $d_{A}^+d_{A}^{+,*}$
and choice of Banach spaces (namely $L^{\sharp}(X)$, as defined in \S
\ref{sec:Decay}) to achieve the desired
estimate in the case of the anti-self-dual equation.  It is at this point
that we encounter one of the more persistent problems peculiar to
$\PU(2)$-monopole gluing:

\begin{prob}
\label{prob:FA+FA-BochnerProblem}
The Laplacian $d_{A,\Phi}^1d_{A,\Phi}^{1,*}$ has the shape
$$
d_{A,\Phi}^1d_{A,\Phi}^{1,*}
=
\begin{pmatrix}
d_{A}^+d_{A}^{+,*} & 0
\\
0 & D_{A}D_{A}^*
\end{pmatrix}
+
\text{non-diagonal lower-order terms}
$$
have the shape
\begin{align*}
2d_{A}^+d_{A}^{+,*}
&= 
\cov_{A}^*\cov_{A} + \{F_{A}^+,\cdot\} + \text{bounded zeroth-order terms},
\quad\text{on }\Gamma(\Lambda^+\otimes\fg_E),
\\
D_{A}D_{A}^* 
&= 
\cov_{A}^*\cov_{A} + \rho(F_{A}^-) + \text{bounded zeroth-order terms},
\quad\text{on }\Gamma(V^-),
\end{align*}
where $F_{A}^+$ is $L^{\sharp,2}$ {\em small\/}, but $F_{A}^-$ is
$L^{\sharp,2}$ {\em large} (as, roughly speaking, it is the component of
$F_A$ which ``bubbles''), with the remaining zeroth-order terms given by
curvatures of connections which remain fixed throughout. {\em A priori}
estimates \cite{FL1} give a uniform $L^\8$ bound on $F_A^+$ but at best an
$L^2$ bound on $F_A^-$, for $\PU(2)$ monopoles $(A,\Phi)$. 
\end{prob}

The practical effect of the last problem is that it becomes exceedingly
difficult to obtain an estimate for $P_{A,\Phi,\mu}$, that is a bound of
the form
\begin{equation}
\label{eq:MainEstimate}
\|d_{A,\Phi}^{1,*}(v,\psi)\|_{L^2_{1,A}(X)}
\leq
C\|d_{A,\Phi}^1d_{A,\Phi}^{1,*}(v,\psi)\|_{L^{\sharp,2}(X)}, 
\end{equation}
with constant $C$ which is uniform with respect to the bubbling pair
$(A,\Phi)$. Naively, everywhere one sees a lower-order term $F_A^+$ in
Taubes' arguments in \cite[\S 5]{TauStable}, one should replace this by a
term $F_A$, in view of the Bochner formula for $(2d_{A}^+d_{A}^{+,*},
D_{A}D_{A}^*)$. In so doing, one quickly realizes that Taubes' argument
breaks down (bear in mind that Taubes is solving the self-dual rather than
the anti-self-dual equation and that we have switched conventions to the
now customary anti-self-dual equation). For example, in passing from
equation (5.17) to (5.18) in
\cite{TauStable} Taubes uses the fact that $F_A^+$ is $L^{\sharp,2}$-small
to permit rearrangement: the condition on $F_A^+$ is embodied in the first
equation of \cite[p. 193]{TauStable} and the hypothesis of his Lemma
5.6. In addition to the preceding difficulty, one needs to derive
suitable small-eigenvalue bounds in order to address Problem
\ref{prob:SmallEigenvalues} using the above fairly weak Bochner
formulas. Thus, to say the least, it is a challenging task to reproduce
analogues of Taubes' results in the case of $\PU(2)$ monopoles. Consequently,
one of the main results of this paper is to show that in fact this can be
done, despite these obstacles, and derive the key estimate 
\eqref{eq:MainEstimate} (see Corollary \ref{cor:L21AEstPAaphi}). One of the
major obstacles is the problem of bounding negative spinors uniformly in
$L^q_{1,A}$, $q>2$, with
respect to $F_A$ and the radii of small balls surrounding curvature bubbles
when all one has is the weak Bochner formula for $D_AD_A^*$ given above
(this problem is addressed in \S \ref{sec:Global} and \S \ref{sec:Dirichlet},
using decay estimates for eigenspinors in \cite{FeehanKato}).

Instead of trying to adapt Taubes' method for constructing a right-inverse,
one could alternatively try to adapt the Donaldson-Kronheimer method for
writing the extended anti-self-dual equations and constructing a right
inverse with the appropriate bound \cite[\S 7.2]{DK}. Indeed, this is the
approach we initially took in \cite{FL0}: one has to modify their method,
as the $\PU(2)$ monopole equations are not conformally invariant \cite[\S
4.2]{FL1}. The principal difference is that, instead of constructing an
``intrinsic'' right inverse $P_{A,\Phi,\mu}$ based on spectral bounds, one can
modify the traditional method of constructing a parametrix for a
pseudodifferential operator \cite{Gilkey}, \cite{Hormander} and patch
together right inverses corresponding to
the background pairs $(A_0,\Phi_0)$ away from a
a collection of small balls and the instantons $(A_i,0)$ spliced in from
$S^4$ onto small balls in $X$, keeping the metric $g$ on $TX$ fixed rather
than allowing it to vary conformally as in \cite{DK}.  
In this situation, the main technical problem 
addressed in \cite{FL0} was the need to allow the instantons $A_i$ to vary
over their entire moduli spaces rather than precompact subsets (as in
\cite{DK}). A third possibility, not considered here, might be to use the
long-cylinder model for gluing as in \cite{MST} for the Seiberg-Witten
equations, following the example of \cite{MrowkaThesis}. However, in this
situation, the metric is varying (with multiple tube lengths) and appears less
appropriate for our purpose, for reasons discussed in the paragraph below.

With either of these two gluing methods, one needs to fit together the
resulting topological models as the stratum $\Sigma$ varies through the
symmetric product $\Sym^\ell(X)$ so as to form a global model 
(as in \cite{FLConj}) and it turns
out that the method of \cite[\S 7.2]{DK} appears not well suited to this
task. It is for this reason that we ultimately chose the technically more
demanding path initiated in this paper and concluded in
\cite{FL4}. As we shall explain in \cite{FLConj}, Taubes' method does adapt
well to the requirement that one can fit together the local gluing-data bundles
$\bGl(\sU_{\ell,\mu},\Sigma,\lambda_0)$ in such a way as to construct a
global gluing-data space $\bGl(\sU_\ell,\lambda_0)$.

We note in passing that if $(A,\Phi)+d_{A,\Phi}^{1,*}(v,\psi)$ solves
the extended $\PU(2)$ monopole equation, then it is a true 
$\PU(2)$ monopole if and only if $(v,\psi)$ also solves the obstruction
equation
$$
\bchi_\mu(A,\Phi)
\equiv
\Pi_{A,\Phi,\mu}
\left(\{d_{A,\Phi}^{1,*}(v,\psi), d_{A,\Phi}^{1,*}(v,\psi)\}
+\fS(A,\Phi)\right) = 0.
$$
However, for our topological calculations, it is enough to only solve the
extended equation.  The map $(A,\Phi) \mapsto \bchi_\mu(A,\Phi)$
defines a section of a local obstruction bundle $\bXi_\mu$ over
$\bGl(\sU_{\ell,\mu},\Sigma,\mu)$, while the assignment $(A,\Phi) \mapsto
\bgamma_{\mu,\Sigma}(A,\Phi) = (A,\Phi) + d_{A,\Phi}^{1,*}(v,\psi)$ 
defines the local gluing map.

There is one remaining difficulty which is most convenient to address in
\cite{FLConj}, when actually carrying out the requisite topological
calculations: 

\begin{prob}
\label{prob:SpectralFlow}
If $\dim M_\fs^{\sw} > 0$, one cannot necessarily fix a single, uniform
positive upper bound $\mu$ for the small eigenvalues of
$d_{A,\Phi}^1d_{A,\Phi}^{1,*}$, due to spectral flow as the point
$[A,\Phi]$ varies in a space $\sU_\ell\subset\sC(\ft_\ell)$ containing
$M_\fs^{\sw}$. 
\end{prob}

If it were not for spectral flow, one could define thickened moduli spaces
of $\PU(2)$ monopoles using the {\em intrinsic, extended $\PU(2)$ monopole
equations\/}, $\Pi_{A,\Phi,\mu}^\perp\fS(A,\Phi) = 0$, the equations being
``intrinsic'' because they are defined on the configuration space $\sC(\ft)$
of all pairs rather than only on the space of gluing data, as in the case
of Taubes' version of the extended equations. In \cite{FL2a} we constructed
a thickened moduli space of $\PU(2)$ monopoles by adapting the
Atiyah-Singer stabilization method, but --- unlike Taubes' method --- this
technique does not mesh particularly well with the non-compactness due to
bubbling and with gluing. We remark that it would be preferable, for the
purposes of fitting together the spaces of local gluing data indexed by the
strata $\Sigma \subset\Sym^\ell(X)$, to exclusively use the intrinsic
extended equations rather than the extended equations. However, while
perhaps not impossible, the gluing theory for the intrinsic, extended
$\PU(2)$ monopole equations appears considerably more difficult; the
problems are explained in a little more detail in \S \ref{sec:Existence}.

\subsection{Outline}
\label{subsec:Outline}
We now summarize the contents of the remainder of this article.  

We recall the definition of the $\PU(2)$ monopole equations and the basic
transversality and compactness properties of the $\PU(2)$ monopole moduli
space in \S \ref{sec:Prelim}.

In \S \ref{sec:Splicing} we define the {\em local splicing\/} or {\em
pre-gluing map\/} $\bgamma_{\mu,\Sigma}'$ for monopoles by adapting the
cut-and-paste techniques of \cite{TauSelfDual}, \cite{TauIndef},
\cite{TauFrame}, and \cite{TauStable}.  The important technical issue here
is to keep track of the constants defining the splicing map itself. The
splicing maps yield continuous maps (smooth away from singularities) from
bundles of gluing data, one bundle per smooth stratum
$\Sigma\subset\Sym^\ell(X)$, into the configuration space $\sC(\ft)$ of
pairs $[A,\Phi]$, whose image is Uhlenbeck-close to the space
$M(\ft_{\ell})\times\Sigma$. The question of exactly how ``close'' is
addressed in and \S \ref{sec:Decay}. Later, in \cite{FLConj}, we
return to the issue of how the bundles of local gluing data,
$\bGl(\sU_{\mu,\ell},\Sigma,\lambda_0)\to\sU_{\mu,\ell}\times\Sigma$, and 
local splicing maps, $\bgamma_{\mu,\Sigma}'$, may be assembled into a space
of {\em global gluing data\/} with projection map,
$\bGl(\sU_{\mu,\ell},\lambda_0)\to\sU_{\mu,\ell}\times\Sym^\ell(X)$, and 
{\em global splicing map}, $\bgamma_\mu'$.  

In \S \ref{sec:Regularity} we recall the regularity theory for the $\PU(2)$
monopole equations \cite[\S 3]{FL1} which, indeed, applies to any
quasi-linear first-order elliptic equation with a quadratic non-linearity
on a four-manifold. The results are used repeatedly throughout the present
article and its sequels, so precise statements are given. The general
pattern is that one typically has uniform, global $L^2_1$ estimates for
solutions to the $\PU(2)$ monopole equations (obtained, for example, by
deforming an approximate, spliced solution), but when we need estimates
with respect to stronger norms --- uniform with respect to a non-compact
parameter --- we must get by with only local estimates. The regularity
section concludes with a proof that an $L^2_1$ gluing solution $(a,\phi)$
to the extended $\PU(2)$ monopole equations is in $L^2_k$ (with $k\geq 4$)
if the spliced pair $(A,\Phi)$ is in $L^2_k$ (Proposition
\ref{prop:RegularityTaubesSolution}).  

Section \ref{sec:Decay} gives an estimate for $\fS(A,\Phi)$, when
$(A,\Phi)$ is an approximate solution to the $\PU(2)$ monopole equations
\eqref{eq:PT} produced by splicing. A point that first arises here --- but
will recur in sections \ref{sec:Global}, \ref{sec:Dirichlet}, 
\ref{sec:Eigenvalue}, and \ref{sec:Existence} and which needs to be
continually borne in mind --- is that we shall almost always want estimates
for pairs where the constants depend at most on the global $L^{\sharp,2}$
norms of the curvatures of the potentially bubbling $\SO(3)$ connections or
their $L^\8$ norms when restricted to suitable open sets. There are
different ways these estimates can be achieved, depending on the geometry
of the region in $X$ where estimates are required. The cleanest approach
is to use decay estimates for Yang-Mills connections (which we apply to
the instantons we splice in from $S^4$),
due in varying degrees of generality to Donaldson, Groisser-Parker,
Morgan-Mrowka-Ruberman, or R\aa de; in sections \ref{sec:Dirichlet},
\ref{sec:Eigenvalue}, and \ref{sec:Existence} we also use the 
decay estimates for eigenspinors derived in \cite{FeehanKato}.

Section \ref{sec:Global} collects the basic, global, uniform elliptic
estimates for the anti-self-dual and Dirac operators which are used to
estimate small eigenvalues of the Laplacian in \S \ref{sec:Eigenvalue} and
to estimate the partial right inverse $P_{A,\Phi,\mu}$ in \S
\ref{sec:Existence}. 

At an important crux step (see inequalities \eqref{eq:LpBoundWorstTerm} and
\eqref{eq:Lp1psiU}) in our derivation of a uniform bound for the partial
right inverse $P_{A,\Phi,\mu}$ in \S \ref{sec:Existence} we shall need a
uniform $L^q_{1,A}(U)$ estimate for a negative spinor $\psi$ in terms of
the $\|D_{A}D_{A}^*\psi\|_{L^p(X)}$ when $1< p \leq 2$ and $4/3 < q \leq
4$ are defined by $1/p=1/q+1/4$, with $q$ at least slightly greater than
$2$ and $U\subset X$ the complement of small balls $B_i$ in $X$ where we
splice in instantons from $S^4$ (and the spinor $\Phi$ is zero). This
crucial estimate is provided by Theorem \ref{thm:Lp2SpinorAnnulusEst} and
the purpose of \S \ref{sec:Dirichlet} is to prove this result, one of the
more delicate ones in this paper. It is important to recognize at the
outset that such an estimate (with the needed uniformity property) does
{\em not\/} follow from a naive application of a global elliptic
estimate over all of $X$ to a cutoff spinor $\chi\psi$, where $d\chi$ is
supported on small annuli surrounding the balls $B_i$ and $\chi$ is zero on
the balls $B_i$ (where the curvature of $A$ is large).

Our technique hinges on the fact that
\cite{FeehanKato} provides us with good decay estimates for eigenspinors
away from the region of $X$ where the curvature of $A$ is large, namely
the balls $B_i$. In \S \ref{sec:Global} we use these decay estimates to
obtain a uniform $L^2_{1,A}(U)$ elliptic bound for $\psi$ by splitting
$\psi$ into a ``small-eigenvalue'' component and its $L^2$-orthogonal
complement, as the latter component can be estimated separately (over all
of $X$) using results of \cite{FeehanSlice}.  Unfortunately, the method of
\S \ref{sec:Global} does not give $L^q_{1,A}(U)$ estimates when
$q>2$. Instead, in order to obtain an $L^q_{1,A}(U)$ estimate for
$\psi|_U$ which is uniform with respect to the radii of the balls $B_i$, we
would like to apply the Calderon-Zygmund theorem to $\psi|_U$. The latter
application would require that $\psi = 0$ on the boundary $\partial U$ and,
as we cannot assume this, we instead apply the Calderon-Zygmund theorem to
$(\psi-h)|_U$, where $h$ obeys $D_{A}D_{A}^*h=0$ and $h=\psi$ on
$\partial U$. The requisite $L^2_{1,A}(U)$ estimate for $\psi$ is provided
by Proposition \ref{prop:LinftyL22EstphiCompOfBubbles} (see inequality
\eqref{eq:LinftyL22EstEigenphiCompOfBubbles}). 

Splicing produces a family, $(A,\Phi)$, of approximate $\PU(2)$ monopoles
on $X$ and thus a family of Laplacians, $d_{A,\Phi}^1d_{A,\Phi}^{1,*}$. We
shall need upper and lower bounds for the small eigenvalues of this
Laplacian, along with estimates for its eigensections. These and related
estimates are derived in \S \ref{sec:Eigenvalue}.

Section \ref{sec:Existence} assembles the proof of the existence of a
partial right inverse to the linearization, with useful bounds, and then
deduces the existence of solutions to the {\em extended $\PU(2)$ monopole
equations\/}. Viewed globally, this step yields the local gluing and
obstruction maps and allows us to conclude the proof of Theorem
\ref{thm:GluingTheorem1}.

\subsection{Sequels and related articles}
\label{subsec:Sequels}
Some remarks are in order concerning the contents of the companion
articles \cite{FL4} and \cite{FLConj}.

\subsubsection{Properties of gluing maps}
The reader will note that the main result proved here is at most the first
half of a desired ``gluing theorem''. To be of use for parameterizing
neighborhoods of ideal, reducible $\PU(2)$ monopoles, we also need 
the following properties:
\begin{enumerate}
\item {\em Continuity.\/}
The local gluing map $\bgamma_{\mu,\Sigma}$ extends to a continuous map on
the Uhlenbeck closure of the gluing data,
$\bar{\bGl}^+(\sU_{\ell,\mu},\Sigma,\lambda_0)$.
\item {\em Embedding property.\/}
The map $\bgamma_{\mu,\Sigma}$ is a smooth embedding of
$\bGl^+(\sU_{\ell,\mu},\Sigma,\lambda_0)$ and a topological embedding of
the Uhlenbeck closure of the gluing data,
$\bar{\bGl}(\sU_{\ell,\mu},\Sigma,\lambda_0)$.
\item {\em Surjectivity.\/}
The image of
$\bar{\bGl}^+(\sU_{\ell,\mu},\Sigma,\lambda_0)\cap\bchi_{\mu,\Sigma}^{-1}(0)$
under $\bgamma_{\mu,\Sigma}$ is an open subset of $\barM(\ft)$ and
the space $\barM(\ft)$ has a finite covering by such open subsets.
\end{enumerate}
These remaining properties (1), (2), (3) are proved in \cite{FL4}
and comprise the second half of the gluing theorem. 
It is important to note that preceding three
gluing-map properties are not simple consequences of the proof of existence
of solutions to the (extended) $\PU(2)$ monopole equations
and their justification constitutes the more
difficult half of the proof of the full gluing theorem. In particular,
despite an extensive literature on gluing theory for anti-self-dual
connections, the existing accounts (see, for example,
\cite{DonConn}, \cite{DK}, \cite{MrowkaThesis}, 
\cite{TauSelfDual}, \cite{TauIndef}, \cite{TauFrame}, \cite{TauStable})
even there only address somewhat special cases which do not capture all of
the difficulties one encounters when attempting to solve the complete
gluing problem for anti-self-dual connections.

\subsubsection{Intersection theory and completion of the proof of the
homotopy Witten conjecture} In work in preparation, \cite{FLConj}, we
complete the proof of the ``homotopy version'' of the Witten conjecture.
This is carried out by first altering the domains of the gluing maps in a
manner allowing a description of their overlaps which depends only on the
homotopy type of $X$, but does not change the gluing (deformation) theory
of the present article and its companion \cite{FL4}.  A topological
description of the obstruction bundle is then given which, together with a
technique for describing cohomology classes with compact support, allows us
to evaluate the integrals over the links in terms of Seiberg-Witten
invariants and homotopy data for the manifold $X$.

\subsection{Acknowledgments}
The first author is grateful to Tom Mrowka and Jeff McNeal for helpful
discussions. He warmly acknowledges the hospitality and support of the many
institutions where various parts of the present article and its companions
were written, including Harvard University, Ohio State University
(Columbus), the Max Planck Institut f\"ur Mathematik (Bonn), together with
the Institute des Hautes Etudes Scientifiques (Bures-sur-Yvette) and the
Institute for Advanced Study (Princeton), where the article was finally
completed.  He would also like to thank the National Science Foundation for
their support during the preparation of this article. The second author
also expresses his gratitude for the hospitality of Harvard University,
Ohio State University, and the Institute for Advanced Studies during his
visits and for the support of Florida International University.

%end of file
%file: prelim3.tex

\section{Preliminaries}
\label{sec:Prelim}
In this section we recall the framework for gauge theory for $\PU(2)$
monopoles established in 
\cite{FeehanGenericMetric}, \cite{FL1}.  In \S
\ref{subsec:NonAbelianMonopoles} we 
describe the $\PU(2)$ monopole equations while in \S
\ref{subsec:CompactnessTransversality} we recall our Uhlenbeck compactness
and transversality results from
\cite{FeehanGenericMetric}, \cite{FL1}.  

\subsection{PU(2) monopoles}
\label{subsec:NonAbelianMonopoles}
Throughout this article, $X$ denotes a closed, connected, oriented, smooth
four-manifold. We shall briefly recall the description of the moduli space
of $\PU(2)$ monopoles from
\cite{FeehanGenericMetric}, \cite{FL1}, \cite{FLGeorgia}. We give $X$ a
Riemannian metric $g$ and consider Hermitian two-plane bundles $E$ over $X$
whose determinant line bundles $\det E$ are isomorphic to a fixed Hermitian
line bundle endowed with a fixed $C^\8$, unitary connection $A_e$.  Let
$\fs_0 = (\rho,W^+,W^-)$ be a
\spinc structure on $X$, where $\rho:T^*X\to\Hom_\CC(W^+,W^-)$ is a
Clifford map compatible with $g$ and so obeys the properties of
\eqref{eq:CliffordMap}, and the Hermitian four-plane bundle $W=W^+\oplus
W^-$ is endowed with a $C^\8$, unitary connection. Clifford multiplication
extends to an injective, linear map
$\rho:\Lambda^\bullet(T^*X)\to\End_\CC(W)$ in the usual way. The unitary
connection on $W$ uniquely determines a Riemannian connection on $T^*X$, by
requiring that it act as a derivation with respect to the Clifford map
$\rho$, and determines a unitary connection on $\det W^+$; conversely, a
choice of Riemannian connection on $T^*X$ and unitary connection on $\det
W^+$ induce a unitary connection on $W$. We shall require the connection on
$W$ to be \spinc, that is, it induces the Levi-Civita connection on $T^*X$
for the given Riemannian metric. In order to take advantage of decay
estimates for Yang-Mills connections \cite{DK}, \cite{GroisserParkerDecay},
\cite{Rade}, we can assume without loss that the fixed connections $A_d$
and $A_e$ are Yang-Mills when convenient.

Let $k\ge 4$ be an integer and let $\sA_E$ be the space of $L^2_k$
connections $A$ on the $\U(2)$ bundle $E$ all inducing the fixed
determinant connection $A_e$ on $\det E$.  Equivalently, following \cite[\S
2(i)]{KMStructure}, we may view $\sA_E$ as the space of $L^2_k$ connections
$A$ on the $\PU(2)=\SO(3)$ bundle $\su(E)$.  We shall pass back and forth
between these viewpoints, via the fixed connection on $\det E$, relying on
the context to make the distinction clear.  Given a 
unitary connection $A$ on $E$
with curvature $F_A\in L^2_{k-1}(\La^2\otimes\fu(E))$, then $(F_A^+)_0 \in
L^2_{k-1}(\La^+\otimes\su(E))$ denotes the traceless part of its self-dual
component. Equivalently, if $A$ is a connection on $\su(E)$ with curvature
$F_A\in L^2_{k-1}(\La^2\otimes\so(\su(E)))$, then $\ad^{-1}(F_A^+) \in
L^2_{k-1}(\La^+\otimes\su(E))$ is its self-dual component, viewed as a
section of $\La^+\otimes\su(E)$ via the isomorphism
$\ad:\su(E)\to\so(\su(E))$. When no confusion can arise, the isomorphism
$\ad:\su(E)\to\so(\su(E))$ will be implicit and so we regard $F_A$ as a
section of $\La^+\otimes\su(E)$ when $A$ is a connection on $\su(E)$.
The {\em instanton number\/} for $E$ is defined by $\kappa = -\quarter
p_1(\su(E)) = c_2(E)-\quarter c_1(E)^2$. 

It will be convenient to define Hermitian vector bundles
$V = W\otimes E$ and $V=V^+\oplus V^-$, where $V^\pm = W^\pm\otimes E$. 
For an $L^2_k$ section $\Phi$ of $W^+\otimes E$, let $\Phi^* =
\langle\cdot,\Phi\rangle$ be its pointwise Hermitian dual and let
$(\Phi\otimes\Phi^*)_{00}$ be the component of the Hermitian endomorphism
$\Phi\otimes\Phi^*$ of $W^+\otimes E$ which lies in
$\su(W^+)\otimes\su(E)$. The Clifford multiplication $\rho$ defines an
isomorphism $\rho:\La^+\to\su(W^+)$ and thus an isomorphism
$\rho=\rho\otimes\id_{\su(E)}$ of $\La^+\otimes\su(E)$ with
$\su(W^+)\otimes\su(E)$. We shall generally write $\fg_E =
\su(E)$ henceforth, for convenience.

Unless $\Tor_2 H^2(X;\ZZ) = 0$, our monopole moduli spaces are labeled in
part by a choice of \spinc structure on $(X,g)$ rather than just a choice
of characteristic class in $H^2(X;\ZZ)$.
\begin{itemize}
\item
Given $\ft = (\rho,W^+,W^-,E)$, define $(w,p,\Lambda)$ by $w=c_1(E)$, $p =
p_1(\fg_E)$, and $\Lambda = \half c_1(V^+)$, or
\item
Given $(w,p,\Lambda)$ with $p=w^2\pmod{4}$ and $\Lambda-w\equiv
w_2(X)\pmod{2}$, we can choose $\ft = (\rho,W^+,W^-,E)$.
\end{itemize}
We call $\ft = (\rho,W^+,W^-,E)$ a {\em \spinu structure\/} on
$(X,g)$. Generalized \spinc structures, or spin-$G$ structures, and
associated $G$-monopole equations were described by Witten in
\cite{WittenCambridgeNewtonTalk}; these ideas were further developed by Teleman
in \cite{TelemanNonAbelian}, \cite{TelemanMonopole}, where the generalized
\spinc structure underlying the $\PU(2)$ monopole equations is called a
spin-$\U(2)$ structure and is equivalent to a \spinu structure.  More
recently, an elegant repackaging of the $\PU(2)$ monopole equations has
been described by Mrowka \cite{MrowkaPrincetonMorseTalk}.

While our choices of orientations of the $\PU(2)$-monopole moduli spaces
may depend on the classes $c_1(W^+)$ or $c_1(E)$, rather than just $\Lambda
= c_1(W^+)+c_1(E)$, the moduli spaces of solutions to equations
\eqref{eq:PT} for \spinu structures $\ft' = (\rho,W^+\otimes L,W^-\otimes
L,E\otimes L^{-1})$ and $\ft = (\rho,W^+,W^-,E)$ are otherwise identical,
for any complex line bundle $L$. Thus, aside from the issue of
orientations, we need not distinguish between such pairs of moduli spaces.

Let $\sG_E$ be the Hilbert Lie group of $L^2_{k+1}$ unitary gauge
transformations of $E$ with determinant one, with Lie algebra
$L^2_{k+1}(\fg_E)$. 
Our pre-configuration space of pairs on $(\fg_E,V^+)$ is given by
$$
\tsC(\ft) := \sA_E\times L^2_k(V^+),
$$
with tangent spaces $L^2_k(\fg_E)\oplus L^2_k(V^+)$. We call a
pair $(A,\Phi) \in \tsC(\ft)$ a {\em $\PU(2)$ monopole\/} if 
\begin{equation}
\label{eq:PT}
\fS(A,\Phi)
=
\begin{pmatrix}
F_A^+ - \tau\rho^{-1}(\Phi\otimes\Phi^*)_{00},
\\
D_A\Phi + \rho(\vartheta)\Phi
\end{pmatrix}
=
0,
\end{equation}
where $D_A=\rho\circ\cov_A:L^2_k(V^+)\to L^2_{k-1}(V^-)$
is the Dirac operator, while $\tau \in L^2_k(X,\GL(\Lambda^+))$ and
$\vartheta \in L^2_k(\Lambda^1\otimes\CC)$ are perturbation
parameters. We let
$$
M(\ft) := \{[A,\Phi] \in \sC(\ft):\text{ $(A,\Phi)$ satisfies
\eqref{eq:PT}}\},
$$
be the moduli space of solutions to \eqref{eq:PT} cut out of the
configuration space,
$$
\sC(\ft) := \tsC(\ft)/\sG_E,
$$ 
where $u\in\sG_E$ acts by $u(A,\Phi) :=
(u_*A,u\Phi)$. The linearization of the map $\sG_E\to \tsC(\ft)$,
$u\mapsto u(A,\Phi)$ at $\id_E$, is given by
$$
\zeta 
\mapsto 
-d^0_{A,\Phi}(\zeta) := (-d_A\zeta,\zeta\Phi),
$$
with $L^2$-adjoint $d^{0,*}_{A,\Phi}$.

\begin{rmk}
\label{rmk:GaugeGroupConvention}
Note that we break here with our former convention
(\cite{FL1}, \cite{FLGeorgia}, \cite{FL2a}, \cite{FL2b}) of denoting
$M(\ft)$ and $\sC(\ft)$ as quotients by
$S_Z^1\times_{\{\pm\id_E\}}\sG_E$. Our new convention reduces notational
clutter and, when we wish to consider quotients by
$S_Z^1\times_{\{\pm\id_E\}}\sG_E$, we shall simply write $M(\ft)/S^1$ and
$\sC(\ft)/S^1$, noting that the spaces $M(\ft)$ and $\sC(\ft)$ carry a
standard action of $S_Z^1\subset\U(2)$ \cite{FL1}.
\end{rmk}

As customary, we say that an $\SO(3)$ connection $A$ on $\fg_E$ is {\em
irreducible\/} if its stabilizer in $\sG_E$ is $\{\pm\id_E\}$,
corresponding to the center of $\SU(2)$, and {\em reducible\/} otherwise;
we say that a pair $(A,\Phi)$ on $(\fg_E,V^+)$ is irreducible
(respectively, reducible) if the connection $A$ is irreducible
(respectively, reducible). We let $\sC^*(\ft)\subset \sC(\ft)$ be the open
subspace of gauge-equivalence classes of irreducible pairs. If $\Phi\equiv
0$ on $X$, we call $(A,\Phi)$ a {\em zero-section\/} pair.  We let
$\sC^0(\ft)\subset \sC(\ft)$ be the open subspace of gauge-equivalence
classes of non-zero-section pairs and recall that $\sC^{*,0}(\ft) :=
\sC^*(\ft)\cap\sC^0(\ft)$ is a Hausdorff, Hilbert manifold
\cite[Proposition 2.12]{FL1} represented by pairs with stabilizer
$\{\id_E\}$ in $\sG_E$. Let $M^{*,0}(\ft) = M(\ft)\cap\sC^{*,0}(\ft)$ be
the open subspace of the moduli space $M(\ft)$ represented by irreducible,
non-zero-section $\PU(2)$ monopoles.

Let $\sB_E(X) := \sA_E(X)/\sG_E$ be the quotient of the space of $\SO(3)$
connections on $\fg_E$ by the induced action of $\sG_E$ on $\fg_E$ and
let $\sA^*_E(X)$ and $\sB_E^*(X)$ be the subspace of irreducible $L^2_k$
connections and its quotient. As before, we may equivalently view
$\sB_E(X)$ as the quotient space of $L^2_k$
connections on $E$ which induce the fixed connection on $\det
E$. Finally, we may view $\sA_E\subset\tsC(\ft)$ and
$\sB_E\subset\sC(\ft)$ as the subspace of zero-section pairs and its
quotient. Hence, the moduli space of anti-self-dual connections,
$$
M_\kappa^w = M_E^{\asd} := \{[A]\in\sB_E: F_A^+ = 0\},
$$
is identified with the subspace of $M(\ft)$ given by zero-section
solutions to \eqref{eq:PT}. 

As in the hypothesis of Theorem \ref{thm:GluingTheorem1},
we further constrain the class $w$ by requiring that no $\SO(3)$ bundle $P$
over $X$ with $w_2(P)\equiv w\pmod{2}$ admit a flat, reducible connection:
a sufficient condition for this to hold is that there exist a spherical
class $\beta \in H^2(X;\ZZ)$ for which $\langle w,\beta\rangle\neq
0\pmod{2}$ \cite[p. 226]{MorganMrowkaPoly}. For any closed four-manifold
$X$ and integral class $w$, the pair $(\hat X, \hat w)$ obeys the
Morgan-Mrowka criterion, where $X=X\#\bar{\CC\PP}^2$ and $\hat w = w +
\PD(e)$, with $e$ the class of the exceptional curve in $\hat X$. Here, the
Morgan-Mrowka criterion ensures that the Uhlenbeck compactification
$\barM(\ft)$ contains no reducible, zero-section pairs if $b^+(X)>0$ and
the metric $g$ is generic (see Propositions 2.9(1) and 3.1(3) and Lemma
3.3 in \cite{FL2a}). The Morgan-Mrowka condition can always be satisfied
by blowing up if necessary and so, as far as the computation of Donaldson
or Seiberg-Witten invariants is concerned, there is no loss in generality
(because of \cite[Theorem 6.9]{KotschickBPlus1}, \cite[Theorem
1.4]{FSTurkish}).

\subsection{Uhlenbeck compactness and transversality for PU(2) monopoles}
\label{subsec:CompactnessTransversality}
We briefly recall our Uhlenbeck compactness and transversality results 
\cite{FeehanGenericMetric}, \cite{FL1} for the moduli space of $\PU(2)$
monopoles 
with the perturbations discussed in the preceding section. The moduli space
of $\PU(2)$ monopoles is non-compact but has an Uhlenbeck closure analogous
to that of the moduli space of anti-self-dual connections on $\fg_E$
\cite[\S 4.4]{DK}.

We say that a sequence
of points $[A_\alpha,\Phi_\alpha]$ in $\sC(\ft)$ {\em converges\/} to a
point $[A,\Phi,\bx]$ in $\sC(\ft_{\ell})\times\Sym^\ell(X)$, where
$\ft_{\ell} = (\rho,W^+,W^-,E_{\ell})$ and
$E_{\ell}$ is a Hermitian two-plane bundle over $X$ such that 
$$
\det(E_{\ell}) = \det E
\quad\text{and}\quad 
c_2(E_{\ell}) = c_2(E)-\ell, 
\quad\text{with}\quad 
\ell\in\ZZ_{\ge 0},
$$
if the following hold:
\begin{itemize}
\item There is a sequence of $L^2_{k+1,\loc}$ determinant-one, unitary
bundle isomorphisms $u_\alpha:E|_{X\less\bx}\to
E_{\ell}|_{X\less\bx}$ such that the sequence of monopoles
$u_\alpha(A_\alpha,\Phi_\alpha)$ converges to
$(A,\Phi)$ in $L^2_{k,\loc}$ over $X\less\bx$, and 
\item The sequence of measures  
$|F_{A_\beta}|^2$ converges
in the weak-* topology on measures to $|F_A|^2 +
8\pi^2\sum_{x\in\bx}\delta(x)$.
\end{itemize}
Given a Hermitian two-plane bundle $E$, we call the intersection of
$\barM(\ft)$ with $M(\ft_{\ell})\times\Sym^\ell(X)$ a {\em
lower-level\/} of the compactification $\barM(\ft)$ if $\ell>0$ and call
$M(\ft)$ the {\em top\/} or {\em highest level\/}.

\begin{thm}
\label{thm:Compactness}
\cite{FL1}
Let $X$ be a closed, oriented, smooth four-manifold with $C^\8$ Riemannian
metric, \spinc structure $(\rho,W^+,W^-)$ with \spinc connection on $W =
W^+\oplus W^-$, and a Hermitian two-plane bundle $E$ with unitary
connection on $\det E$.  Then there is a positive integer $N_p$, depending
at most on the curvatures of the fixed connections on $W$ and $\det E$
together with $c_2(E)$, such that for all $N\ge N_p$ the topological space
$\barM(\ft)$ is second countable, Hausdorff, compact, and given by the
closure of $M(\ft)$ in the space of {\em ideal $\PU(2)$ monopoles\/}
$\sqcup_{\ell=0}^{N}(M(\ft_{\ell})\times\Sym^\ell(X)$ with respect to the
Uhlenbeck topology, where $E_{\ell}$ is a Hermitian two-plane bundle over
$X$ with $\det E_{\ell} = \det E$ and $c_2(E_{\ell})=c_2(E)-\ell$
for each integer $\ell\ge 0$.
\end{thm}

Theorem \ref{thm:Compactness} is a special case of the more general result
proved in \cite{FL1} for the moduli space of solutions to the $\PU(2)$
monopole equations in the presence of holonomy perturbations.
The existence of an Uhlenbeck compactification for the moduli space of
solutions to the unperturbed $\PU(2)$ monopole equations
\eqref{eq:PT} was announced by Pidstrigach
\cite{PTLectures} and an argument was outlined in \cite{PTLocal}.
A similar argument for the equations 
\eqref{eq:PT} (without perturbations)
was outlined by Okonek and Teleman in
\cite{OTQuaternion}. An independent proof of Uhlenbeck compactness for
\eqref{eq:PT} and other perturbations of these equations is
also given in \cite{TelemanGenericMetric}.

We recall from \cite[Equation (2.25)]{FL1} that the elliptic deformation
complex for the moduli space $M(\ft)$ is given by
\begin{equation}
\begin{CD}
L^2_{k+1}(\fg_E)
@>{d_{A,\Phi}^0}>> 
\begin{matrix}
L^2_k(\Lambda^1\otimes\fg_E) \\
\oplus \\
L^2_k(V^+)
\end{matrix}
@>{d_{A,\Phi}^1}>> 
\begin{matrix}
L^2_{k-1}(\Lambda^+\otimes\fg_E) \\
\oplus \\
L^2_{k-1}(V^-)
\end{matrix}
\end{CD} 
\label{eq:ConstDetDefComplex}
\end{equation}
with elliptic deformation operator
\begin{equation}
\sD_{A,\Phi} := d_{A,\Phi}^{0,*} + d^1_{A,\Phi}
\label{eq:ConstDetDefOperator}
\end{equation}
and cohomology $H_{A,\Phi}^\bullet$. Here, $d^1_{A,\Phi}$ is the
linearization at the pair $(A,\Phi)$ of the gauge-equivariant map $\fS$
defined by the equations \eqref{eq:PT}, so
$$
d^1_{A,\Phi}(a,\phi)
:=
(D\fS)_{A,\Phi}(a,\phi)
=
\begin{pmatrix}
d_A^+a-\half\tau\rho^{-1}(\phi\otimes\Phi^*+\Phi\otimes\phi^*)_{00} \\
(D_A+\rho(\vartheta))\phi+\rho(a)\Phi
\end{pmatrix}.
$$
The space $H_{A,\Phi}^0 = \Ker d_{A,\Phi}^0$ is the Lie algebra
of the stabilizer $\Stab_{A,\Phi}$ in $\sG_E$ of a pair
$(A,\Phi)$ and $H^1_{A,\Phi}$ is the Zariski tangent space to $M(\ft)$ at
a point $[A,\Phi]$.  If $H_{A,\Phi}^2=0$, then $[A,\Phi]$ is a regular point of
the zero locus of the $\PU(2)$ monopole equations \eqref{eq:PT} on
$\sC(\ft)$.

We now turn to the question of transversality. Recall that a real-linear
map $\rho:T^*X\to\End_\CC(W)$ defines a Clifford-algebra representation
$\rho:\Cl_\CC(T^*X)\to\End_\CC(W)$ if and only if \cite{LM},
\cite{SalamonSWBook} 
\begin{equation}
\rho(\alpha)^\dagger = - \rho(\alpha)
\quad\text{and}\quad
\rho(\alpha)^\dagger\rho(\alpha) = g(\alpha,\alpha)\id_W,
\qquad \alpha\in C^\8(T^*X),
\label{eq:CliffordMap}
\end{equation}
where $g$ denotes the Riemannian metric on $T^*X$.

For the Riemannian metric $g$ on $T^*X$, let $\cov^g$ be an $\SO(4)$
connection on $T^*X$. A unitary connection $\cov$ on $W$ is called 
{\em spinorial with respect to $\cov^g$\/} 
if it induces the $\SO(4)$ connection $\cov^g$ on $T^*X$,
or, equivalently, if
\begin{equation}
[\cov_\eta,\rho(\alpha)] = \rho(\cov_\eta^g\alpha),
\label{eq:RhoCompatibility}
\end{equation}
for all $\eta \in C^\8(TX)$ and $\alpha \in \Omega^1(X,\RR)$. The unitary
connection $\cov$ on $W$ uniquely determines a unitary connection on $\det W^+
\simeq \det W^-$. Any two unitary
connections on $W$, which are both spinorial with respect to $\cov^g$,
differ by an element of $\Omega^1(X,i\RR)$. Conversely, a unitary
connection $\cov$ on a Hermitian four-plane bundle $W = W^+\oplus W^-$
over an oriented four-manifold $X$ is uniquely determined by
\begin{itemize}
\item A Clifford map $\rho:T^*X\to\Hom_\CC(W^+,W^-)$ satisfying
\eqref{eq:CliffordMap} for the Riemannian metric $g$ on $T^*X$,
\item An $\SO(4)$ connection $\cov^g$ on $T^*X$ for the metric $g$, which
need not be torsion free, and, 
\item A unitary connection on $\det W^+$.
\end{itemize}
The resulting connection $\cov$ on $W$ is then spinorial with respect to
$\cov^g$; we call $\cov$ a {\em \spinc\/} connection if $\cov^g$ is the
Levi-Civita connection. The Dirac operator $D_A$ is not self-adjoint on
$C^\8(W\otimes E)$ unless the $\SO(4)$ connection $\cov^g$ is torsion-free, so
$\cov^g$ is the Levi-Civita connection, as one can see from examples.

Given a Riemannian metric $g$, Clifford map $\rho$, and unitary connection
on $\det W^+$, the Dirac operators $D_A^g$ and $D_A^{\LC}$ defined by any
$\SO(4)$ connection $\cov^g$ and by the Levi-Civita connection $\cov^{\LC}$
on $T^*X$ differ by an element $\rho(\vartheta') \in \Hom_\CC(W^+,W^-)$
\cite[Lemma 3.1]{FeehanGenericMetric}, where $\vartheta' \in
\Omega^1(X,\CC)$. Even though a fixed unitary connection on $W$ will not
necessarily induce a torsion-free connection on $T^*X$ for generic
pairs $(g,\rho)$ of Riemannian metrics and compatible Clifford maps, we can
assume that the Dirac operator $D_A$ in 
\eqref{eq:PT} is defined using the Levi-Civita connection for the metric
$g$ by absorbing the difference term $\rho(\vartheta')$ into the
perturbation term $\rho(\vartheta)$. Given any fixed pair $(g_0,\rho_0)$
satisfying \eqref{eq:CliffordMap} and automorphism $f\in
C^\8(\GL(T^*X))$, then $(g,\rho) := (f^*g_0,f^*\rho_0)$ is again a compatible
pair; the pair $(g,\rho)$ is {\em generic\/} if $f$ is generic.  

\begin{thm}
\label{thm:Transversality}
\cite{FeehanGenericMetric}
Let $X$ be a closed, oriented, smooth four-manifold with Hermitian
two-plane bundles $(W^+,W^-)$ with unitary connection on 
$\det W^+$, and a Hermitian line bundle $\det E$ with unitary
connection. Then for a generic, $C^\8$ pair $(g,\rho)$ satisfying
\eqref{eq:CliffordMap} and generic, $C^\8$ parameters
$(\tau,\vartheta)$, the
moduli space $M^{*,0}(\ft,g,\rho,\tau,\vartheta)$ of $\PU(2)$ monopoles is a
smooth manifold of the expected dimension,
\begin{align*} 
\dim M^{*,0}(\ft,g,\rho,\tau,\vartheta)
&=
\dim M^{*,\asd}_E(g) + 2\Ind_\CC D_A
\\
&= 
-2p_1(\fg_E)-3(1-b^1(X)+b^+(X))
\\
&\quad + \thalf p_1(\fg_E)+\thalf(\Lambda^2-b^+(X)+b^-(X)),
\end{align*}
where $\Lambda := \half c_1(V^+) = c_1(W^+)+c_1(E)$.
\end{thm}

Theorem \ref{thm:Transversality} was proved independently, using a somewhat
different method, by A. Teleman \cite{TelemanGenericMetric}.

Note that $\dim M^{*,0}(\ft)$ depends only on $\chi$, $\sigma$, $\Lambda$,
and $p_1(\fg_E)$: in particular, it does not depend on $c_1(W^+)$,
$w=c_1(E)$, or even $w_2(\fg_E) = w\pmod{2}$, but only on $\Lambda =
c_1(W^+) + c_1(E) = \half c_1(V^+)$.

\begin{rmk}
\label{rmk:DiracConvention}
In Theorem \ref{thm:Transversality}, the Dirac operator $D_A$ of
\eqref{eq:PT} is defined by the unitary connection on $\det W^+$,
Levi-Civita connection on $T^*X$ for the metric $g$, and the Clifford map
$\rho$, together with the unitary connection on $E$ induced by that on
$\det E$ and by an $\SO(3)$ connection $A$ on $\fg_E$. To reduce notational
clutter and when it makes no difference to the argument at hand, we shall
often write $D_{A,\vartheta}$ or even simply $D_A$ for the perturbation
Dirac operator $D_A + \rho(\vartheta)$ on $\Gamma(V^+)$, regarding the
perturbation $\vartheta\in\Omega^1(X,\CC)$ as a component of the
connection $A_d$ on $\det W^+$, now simply viewed as a complex rather than
a unitary connection.
\end{rmk}

For the remainder of this article, we shall write the expected dimensions
appearing in the statement of Theorem \ref{thm:Transversality} as
\begin{align}
d_a(\fg_E) 
:= 
\dim M^{\asd}_E
&= 
-2p_1(\fg_E)-3(1-b^1(X)+b^+(X))
\notag\\
&= 
-2p_1(\fg_E)-\tthreehalf(\chi(X)+\sigma(X)), 
\\
2n_a(\fg_E,\Lambda) 
:=
2\Ind_\CC D_A
&= 
\thalf p_1(\fg_E)+\thalf(\Lambda^2-\sigma(X)),
\notag
\end{align}
where $\chi(X)$ and $\sigma(X)$ are the Euler class and signature of $X$,
respectively.

%end of file
%file: splicing.tex

\section{Splicing extended PU(2) monopoles}
\label{sec:Splicing}
We construct a {\em splicing map\/} $\bga'$ which will give a homeomorphism
from the total space of a gluing data bundle into an open neighborhood in
$\sC(\ft)$ of a level $M(\ft)\times \Sigma$, where
$\Sigma\subset\Sym^\ell(X)$ is a smooth stratum.  When defining our spaces
of connections and gauge transformations, we shall always choose to work
$L^2_k$ connections modulo $L^2_{k+1}$ gauge transformations with $k\geq
4$: though only $k\geq 2$ is strictly necessary, we shall occasionally need
to consider estimates whose constants depend on the $C^0$ norm of $F_A$ (on
a subset of $X$ where $A$ does not bubble)
or $C^1$ norm of $\Phi$.  Thus, {\em for the remainder of this article and
its sequels\/} we shall assume $k\geq 4$ in order to make use of the
Sobolev embedding $L^2_4\subset C^1$.

\subsection{Gluing data}
\label{subsec:GluingData}
In this section we describe the gluing-data bundles we shall need to
parameterize an open neighborhood in $\barM(\ft)$ of a level
$(M(\ft_{\ell})\times\Sigma)\cap \barM(\ft)$.

\subsubsection{Connections over the four-sphere}
\label{subsubsec:ConnectionsOnFourSphere}
We begin by defining the required gluing data associated with the
four-spheres: this will, essentially, define the fibers of our
gluing-data bundles.  Let $\kappa\ge 1$ be an integer and let $S^4$
have its standard round metric of radius one. Let $E$ be a
Hermitian, rank-two vector bundle over $S^4$ with $c_2(E) =
\kappa$ and let $\fg_E = \su(E)$ be the corresponding $\SO(3)$ 
bundle with $-\quarter p_1(\fg_E) = \kappa$.

A choice of frame $f$ in the principal
$\SO(4)$ frame bundle $\Fr(TS^4)$ for
$TS^4$, over the north pole $n\in S^4$, defines a conformal
diffeomorphism,
\begin{equation}
\label{eq:ConformalDiffeo}
\varphi_n:\RR^4 \to S^4\less\{s\},
\end{equation}
that is inverse to a stereographic projection from the south pole
$s\in S^4\subset \RR^5$; let $y(\,\cdot\,): S^4\less\{s\}\to\RR^4$ be
the corresponding coordinate chart.

\begin{defn}
\label{defn:CenterScale}
\cite[pp. 343--344]{TauFrame}
Continue the notation of the preceding paragraph.  The {\em center\/}
$z=z([A],f)\in\RR^4$ and the {\em scale\/} $\lambda=\lambda([A],f)$ of a
point $[A]\in\sB_E(S^4)=\sB_\kappa(S^4)$ are defined by
\begin{align}
z([A],f) 
&:= \frac{1}{8\pi^2\kappa}\int_{{\Bbb R}^4}y|F_A|^2\,d^4y, 
\label{eq:MassCenter}\\
\lambda([A],f)^2 
&:= \frac{1}{8\pi^2\kappa}\int_{{\Bbb R}^4} |y-z([A],f)|^2|F_A|^2\,d^4y. 
\label{eq:Scale}
\end{align}
A point $[A]$ is {\em mass-centered\/} if $z([A],f)=0$ and {\em centered\/}
if it also obeys $\lambda([A],f) = 1$.  An {\em ideal point\/} $[A,\bx]\in
\sB_{\kappa-\ell}(S^4)\times\Sym^\ell(S^4)$ is {\em mass-centered\/} if
$z([A,\bx],f)=0$, where the center and scale are defined by replacing
$|F_A|^2$ with $|F_A|^2+8\pi^2\sum_{x\in\bx}\delta_x$ in equations
\eqref{eq:MassCenter}, with the analogous definition holding for {\em
centered\/} ideal points.
\end{defn}

For any $(\lambda,z)\in\RR^+\times\RR^4$, we define a conformal
diffeomorphism of $\RR^4$ by
\begin{equation}
\label{eq:CenterScaleDiffeo}
f_{\lambda,z}:\RR^4\to\RR^4,\qquad y\mapsto (y-z)/\lambda.
\end{equation}
It is easy to see that the point $[(f_{\lambda,z}^{-1})^*A]$ is
centered if $z = z([A],f)$ and $\lambda = \lambda([A],f)$.  Via the
standard action of $\SO(4)\subset\GL(\RR^4)$ on $\RR^4$ and the
conformal diffeomorphism $\varphi_n^{-1}:S^4\less\{s\}\to\RR^4$, we
obtain a homomorphism from $\SO(4)$ into the group of conformal
diffeomorphisms of $S^4$ fixing the south pole $s$. The group $\SO(4)$
acts on the quotient space of connections based at the south pole,
$\sB_\kappa^s(S^4) \cong \sA_\kappa(S^4)\times_{\sG_E}\Fr(\fg_E)|_s$, by
\begin{equation}
\label{eq:BasedSphereModSpaceSO(4)Action}
\SO(4)\times \sB_\kappa^s(S^4) \to \sB_\kappa^s(S^4),
\quad
(u,[A,q]) \mapsto u\cdot[A,q] = [(u^{-1})^*A,q].
\end{equation}
The assignment $f\mapsto z([A],f)$ is $\SO(4)$-equivariant, that is,
$z([A],uf) = uz([A],f)$ for $u\in \SO(4)$, and so $\lambda[A] =
\lambda([A],f)$ is independent of the choice of $f \in \Fr(TS^4)|_n$. The
center and scale functions define a smooth map
$(z,\lambda):\sB_\kappa(S^4)\to\RR^4\times\RR^+$ and the preimage of $(0,1)\in
\RR^4\times\RR^+$, namely the codimension-five submanifold of {\em
centered connections\/}
\begin{equation}
\label{eq:QuotientSpaceCenteredConn}
\sB_\kappa^{\diamond}(S^4)
:=
\{[A]: z[A] = 0 \text{ and } \lambda[A] = 1\}
\subset
\sB_\kappa(S^4), 
\end{equation}
serves as a canonical slice for the action of $\RR^4\times\RR^+$ on
$\sB_\kappa(S^4)$; the preimage of $0\in
\RR^4$, namely the codimension-four submanifold of {\em
mass-centered connections\/}
\begin{equation}
\label{eq:QuotientSpaceMassCenteredConn}
\sB_\kappa^{\natural}(S^4)
:=
\{[A]: z[A] = 0\}
\subset
\sB_\kappa(S^4), 
\end{equation}
serves as a canonical slice for the action of $\RR^4$ on
$\sB_\kappa(S^4)$. We simply write $z[A]$ when the frame $f$ for $TS^4|_n$
is understood, employ the decoration $\natural$ to indicate a subspace of
mass-centered connections, and employ the decoration $\diamond$ to indicate
a subspace of centered connections.  Since the action of $\SO(4)$
preserves the north and south poles of $S^4$, it acts on
$\sB_\kappa^{\natural,s}(S^4)$, and so $\sB_\kappa^{\natural,s}(S^4)$
admits an action of $\SO(3)\times\SO(4)$, with $\SO(3)$ varying the frame
of $\fg_E|_s$; the same holds for the spaces of centered connections.

A simple {\em a priori\/} decay estimate for the curvature of
connections is provided by the Chebychev inequality:

\begin{lem}\label{lem:Chebychev}
Let $[A] \in \sB_\kappa(S^4)$ and suppose $z([A],f) = 0$. If
  $\lambda = \lambda[A]$ and $R>0$, then
$$
\frac{1}{8\pi^2\kappa} \int_{|y|\ge R\lambda}|F_A|^2\,d^4y \le R^{-2}.
$$
\end{lem}

\begin{pf}
From equation \eqref{eq:Scale} we see that
$$
\lambda^2 \ge \frac{1}{8\pi^2\kappa}\int_{|y|\ge R\lambda}|y|^2|F_A|^2\,d^4y
\ge 
\frac{R^2\lambda^2}{8\pi^2\kappa}\int_{|y|\ge R\lambda}|F_A|^2\,d^4y
$$
and this gives the required bound.
\end{pf}

An immediate consequence of the Chebychev inequality is that if $A$ is
a centered anti-self-dual connection over $S^4$, with second Chern
number $\kappa$ so its energy is equal to $8\pi^2\kappa$, then $A$
cannot bubble outside the ball $B(0,1/2\pi)\subset\RR^4$, as it has
most energy $4\pi^2$ in this region by Lemma \ref{lem:Chebychev}.

We let $M^\diamond_\kappa(S^4)\subset \sB_\kappa^\diamond(S^4)$ denote the
subspace of anti-self-dual connections, centered at the north pole, while
$M^s_\kappa(S^4)\subset \sB_\kappa^s(S^4)$ is the subspace of
anti-self-dual connections, framed at the south pole. Finally, we let
$M^{\diamond,s}_\kappa(S^4)\subset \sB_\kappa^{\diamond,s}(S^4)$ denote the
subspace of anti-self-dual connections, which framed at the south pole and
centered at the north pole. The definition of Uhlenbeck convergence
(see, for example, \S \ref{subsec:CompactnessTransversality}) 
for sequences of points $[A_\alpha]$
in $\sB_\kappa(S^4)$ generalizes easily to the case of sequences
$[A_\alpha,q_\alpha]$ in $\sB_\kappa^s(S^4)$, when $\Fr(\fg_E|_s)\cong
\SO(3)$ is given its usual topology. Since Lemma \ref{lem:Chebychev}
implies that points in $M^{\diamond,s}_\kappa(S^4)$ cannot bubble near the
south pole, the proof of \cite[Theorem 4.4.3]{DK} adapts with no change to
show that $M_\kappa^{\diamond,s}(S^4)$ has compact closure
$\barM^{\diamond,s}_\kappa(S^4)$ in the space of framed, ideal
anti-self-dual connections,
$$
\bigsqcup_{\ell=0}^\kappa
\left(M^s_{\kappa-\ell}(S^4) \times \Sym^\ell(S^4) \right).
$$
Note that it is the ideal points $[A_0,\bx]$ in
$\barM^{\diamond,s}_\kappa(S^4)$ which are centered in the sense of
Definition \ref{defn:CenterScale} and not the background connections $A_0$.
Note also that the points in $S^4=\RR^4\cup\{\8\}$ representing the
multisets $\bx$ will be constrained to lie in the ball $\bar B(0,1/2\pi)$
by Lemma \ref{lem:Chebychev}.

Observe that since the bundle map $u$ in the definition of an Uhlenbeck
neighborhood respects the $\SO(3)$ action on frames, there is a global
action of $\SO(3)$ on the space $\barM^{\diamond,s}_\kappa(S^4)$.  This
action is free on all strata except $[\Theta]\times\Sym^\kappa(S^4)$, where
the background connection is trivial.  In addition, the $\SO(4)$ action on
$M^{\diamond,s}_\kappa(S^4)$ extends over $\barM^{\diamond,s}_\kappa(S^4)$.

\subsubsection{Gluing data bundles}
\label{subsubsec:GluingDataBundles}
We now define bundles of gluing data, $\bGl(\sU_\ell,\Sigma,\lambda_0)$,
associated with a choice of \spinu structure $\ft$ on $(X,g)$, smooth
stratum $\Sigma\subset\Sym^\ell(X)$, precompact open submanifold
$\sU_{\ell}\subset \sC(\ft_{\ell})$, and positive constant
$\lambda_0$. Our definition is motivated by those of Friedman-Morgan
\cite[\S\S 3.4.2--3.4.4]{FrM}, Kotschick-Morgan
\cite[\S 4]{KotschickMorgan}, and Taubes \cite{TauFrame}, \cite{TauStable}.

If $\ft = (\rho,W^+,W^-,E)$, we denote $\kappa = -\quarter p_1(\su(E))$. The
integer $\ell$ obeys $1\le\ell\le \lfloor\kappa\rfloor$ and we let
$\Sigma\subset\Sym^\ell(X)$ be a smooth stratum, defined by an integer
partition $\kappa_1\ge\cdots\ge\kappa_m\ge 1$ of $\ell$ such that
$\kappa_1+\cdots+\kappa_m=\ell$. Let $E_{\ell}$ be the Hermitian, rank-two
vector bundle over $X$ with $\det E_{\ell} = \det E$ and $c_2(E_{\ell}) =
c_2(E)-\ell$, so that $\ft_\ell = (\rho,W^+,W^-,E_\ell)$.  Let $\Fr(TX)$ denote
the principal $\SO(4)$ bundle of oriented, $g$-orthonormal frames for $TX$.
Suppose $\tilde\sU_\ell\subset
\tsC(\ft_\ell)$ projects to the submanifold $\sU_\ell$. Observe that the
fibered product 
\begin{equation}
\tilde\sU_\ell\times_{\sG(\ft_{\ell})}\prod_{i=1}^m 
\left(\Fr(\fg_{E_{\ell}})\times_X\Fr(TX)\right)
\to \sU_\ell\times \prod_{i=1}^m X
\label{eq:PreFiberedProductBundle}
\end{equation}
is a principal $\prod_{i=1}^m(\SO(3)\times\SO(4))$ bundle. Let
$\Delta\subset\prod_{i=1}^m X$ be the full diagonal, that is, the subset of
$m$-tuples of points in $X$ where some point occurs with multiplicity at
least two. Let $\fS = \fS_\Sigma = \fS(\kappa_1,\dots,\kappa_m)$ be the
largest subgroup of the symmetric group on $m$ symbols which preserves the
multiplicities.  The group $\fS(\kappa_1,\dots,\kappa_m)$ acts freely on
$\prod_{i=1}^m X-\Delta$ and we have an identification
$$
\left(\prod_{i=1}^m X-\Delta\right)/\fS(\kappa_1,\dots,\kappa_m)
=
\Sigma\subset\Sym^\ell(X).
$$
The group $\fS(\kappa_1,\dots,\kappa_m)$ acts on the fibered product bundle
\eqref{eq:PreFiberedProductBundle} to give a principal bundle
\begin{equation}
\bFr(\sU_\ell,\Sigma) \to \sU_\ell\times \Sigma,
\end{equation}
with total space 
\begin{equation}
\bFr(\sU_\ell,\Sigma)
:= 
\left(\tilde\sU_\ell\times_{\sG(\ft_{\ell})}\prod_{i=1}^m 
\left(\Fr(\fg_{E_{\ell}})\times_X\Fr(TX)\right)\right)|_{\prod_{i=1}^m X-\Delta}
\end{equation}
and structure group
\begin{equation}
\bG(\Sigma) 
= 
\bG(\kappa_1,\dots,\kappa_m) 
:= 
\left(\prod_{i=1}^m (\SO(3)\times\SO(4))\right)
\rtimes \fS(\kappa_1,\dots,\kappa_m). 
\end{equation}
Let $S^4$ have its standard round metric of radius one, let
$E^1,\dots,E^m$ be Hermitian, rank-two vector bundles over $S^4$ with
$c_2(E^j) = \kappa_j$ for $j=1,\dots,m$, and let $\fg_{E^j} \equiv
\su(E^j)$ be the corresponding $\SO(3)$ bundles with $-\quarter
p_1(\fg_{E^j}) = \kappa_j$.  Denote $\RR^+ := (0,\8)$. The smooth
manifold
\begin{equation}
\label{eq:ConeBundleFiber}
\bZ(\Sigma) 
= 
\bZ(\kappa_1,\dots,\kappa_m) 
:=
\prod_{i=1}^m \left(\tM_{E^i}^{\diamond}\times_{\sG_{E^i}}
\Fr(\fg_{E^i})|_s\times\RR^+\right)
\end{equation}
carries a natural action of $\bG(\kappa_1,\dots,\kappa_m)$ coming from
the action of $\SO(4)$ rotating $S^4$ and the action of $\SO(3)$ on
the choice of frame for $\fg_{E^i}|_s$ in each factor, and the action
of $\fS(\kappa_1,\dots,\kappa_m)$ permuting the factors whose
underlying bundles are isomorphic.  The bundle $\Fr(TX)$ has its usual
right $\SO(4)$ action given by
\begin{equation}
\label{eq:TangentFrameBundleAction}
\SO(4)\times\Fr(TX)\to \Fr(TX), \quad (u,f) \mapsto fu^{-1}.
\end{equation}
We form the fibered product 
\begin{equation}
\label{eq:FiberedProduct}
\bGl(\sU_\ell,\Sigma) 
:= 
\bFr(\sU_\ell,\Sigma)\times_{\bG(\Sigma)} \bZ(\Sigma),
\end{equation}
together with the natural projection map
\begin{equation}
\label{eq:ConeBundleOverStratum}
\pi:\bGl(\sU_\ell,\Sigma) \to \sU_\ell\times \Sigma.
\end{equation}
This is a locally trivial fiber bundle with structure group
$\bG(\kappa_1,\dots,\kappa_m)$. We let
$\bZ(\Sigma,\lambda_0)\subset\bZ(\Sigma)$ denote the open subspace
obtained by replacing $\RR^+$ by $(0,\lambda_0)$ in the definition
\eqref{eq:ConeBundleFiber} and let 
$\bGl(\sU_\ell,\Sigma,\lambda_0)$ denote the corresponding gluing data bundle.

\subsection{Splicing pairs over $X$} 
\label{subsec:SplicingConstruction}
In this section we describe the {\em splicing construction\/} for
pairs on $(\fg_E,V^+)$ via a partition of unity on $X$ and
define the {\em splicing maps} $\bga'$ for pairs.  Our
description is largely motivated by those of Taubes in 
\cite{TauSelfDual}, \cite{TauPath}, \cite{TauIndef}, \cite[\S
4]{TauFrame}, \cite{TauStable}. We keep the notation of \S
\ref{subsec:GluingData}.

\subsubsection{Cutting and splicing pairs}
\label{subsubsec:CutAndSpliceConnections}
Fix a smooth stratum $\Sigma\subset\Sym^\ell(X)$, let $\bx\in\Sigma$, and
choose $r_0=r_0(\bx)>0$ to be one half the smaller of
\begin{itemize}
\item The minimum geodesic distance between distinct points $x_i$,
  $x_j$ of the representative $(x_1,\dots,x_m)$ of $\bx$, and
\item 
The injectivity radius of the Riemannian manifold $(X,g)$.
\end{itemize}
A choice of $L^2_k$ connection $A_0$ on $\fg_{E_{\ell}}$ and $\SO(3)$
frames $p_i\in \Fr(\fg_{E_{\ell}})|_{x_i}$ define $L^2_{k+1}$ local sections of
$\Fr(\fg_{E_\ell})$, 
$$
\sigma_i = \sigma_i(A_0,p_i)
:B(x_i,r_0)\to \Fr(\fg_{E_\ell}), \quad i=1,\dots,m, 
$$
using parallel transport along radial geodesics emanating from $x_i$
together with smoothing for sections of $\fg_{E_\ell}$ using the heat kernel
$\exp(-td_{A_0}^*d_{A_0})$ for small $t$, as described in \cite[\S
A.1]{FL1}; see also the remarks in \cite[p. 177]{TauStable}. We shall
ultimately confine our attention to the case where $(A_0,\Phi_0)$ is a
solution to an (extended) $\PU(2)$ monopole equation, in which case the
$L^2_k$ pair $(A_0,\Phi_0)$ is $L^2_{k+1}$-gauge equivalent to a $C^\8$
pair by Proposition 3.7 in \cite{FL1} (which, though only stated for the
$\PU(2)$ monopole equations, also holds for extended equations). Hence, when
we are splicing background pairs which are (extended) $\PU(2)$ monopoles,
then the heat-kernel smoothing described above is not required as all
pairs, connections, and trivializations can be assumed to be $C^\8$. Let
\begin{equation}
\fg_{E_\ell}|_{B(x_i,r_0)}\simeq B(x_i,r_0)\times \SO(3)
\label{eq:BackroundBundleTriv}
\end{equation}
be the $L^2_{k+1}$ local trivializations defined by the sections $\sigma_i$. 

Define a smooth cutoff function $\chi_{x_0,\eps}:X\to [0,1]$ by setting
\begin{equation}
\label{eq:ChiCutoffFunctionDefn}
\chi_{x_0,\eps}(x) 
:= 
\chi(\dist(x,x_0)/\eps)
\end{equation}
where $\chi:\RR\to [0,1]$ is a smooth function such that $\chi(t)=1$
for $t\ge 1$ and $\chi(t)=0$ for $t\le 1/2$.  Thus, we have
$$
\chi_{x_0,\eps}(x) 
=
\begin{cases}
1 &\text{for } x\in X - B(x_0,\eps),
\\
0 &\text{for } x\in B(x_0,\eps/2).
\end{cases}
$$
We define a cut-off connection on the bundle $\fg_{E_\ell}$ over $X$ by
setting
\begin{equation}
\label{eq:CutOffBackgroundConnection}
\chi_{\bx,4\sqrt{\blambda}}A_0
:=
\begin{cases}
A_0 
&\text{over $X - \cup_{i=1}^m B(x_i,4\sqrt{\lambda_i})$}, 
\\
\Gamma + \chi_{x_i,4\sqrt{\lambda_i}}\sigma_i^*A_0
&\text{over $\Omega(x_i,2\sqrt{\lambda_i},4\sqrt{\lambda_i})$,} 
\\
\Gamma 
&\text{over $B(x_i,2\sqrt{\lambda_i})$,} 
\end{cases} 
\end{equation}
where $\Gamma$ denotes the product connection on $B(x_i,r_0)\times \SO(3)$,
while $\bx = (x_1,\dots,x_m)$ and $\blambda =
(\lambda_1,\dots,\lambda_m)$. 

Similarly, we define a cut-off spinor on the bundles $V^+_\ell$ and $V^+$
over $X$ by setting
\begin{equation}
\label{eq:SplicedSpinor}
\Phi \equiv \chi_{\bx,8\blambda^{1/3}}\Phi_0
:=
\begin{cases}
\Phi_0 
&\text{over $X - \cup_{i=1}^m B(x_i,8\lambda_i^{1/3})$}, 
\\
\chi_{x_i,8\lambda_i^{1/3}}\sigma_i^*A_0
&\text{over $\Omega(x_i,4\lambda_i^{1/3},8\lambda_i^{1/3})$,} 
\\
0
&\text{over $B(x_i,4\lambda_i^{1/3})$,} 
\end{cases} 
\end{equation}
with the same notation. The reason for the different choice of annulus
radii will be explained in \S \ref{sec:Eigenvalue}.

Let $A_i$ be centered $L^2_k$ connections on the bundles $\fg_{E^i}$
over $S^4$. Let 
$$
\tau_i
=
\tau_i(A_i,q_i)
:S^4\less\{n\} \to \Fr(\fg_{E^i}), \quad i=1,\dots,m, 
$$
be the analogously-defined $L^2_{k+1}$ local sections induced by the
connections $A_i \in \sA_{E^i}$ and $\SO(3)$ frames $q_i \in
\Fr(\fg_{E^i})|_s$. Let
\begin{equation}
\fg_{E^i}|_{S^4\less\{n\}}\simeq S^4\less\{n\}\times \SO(3)
\label{eq:SphereBundleTriv}
\end{equation}
be the $L^2_{k+1}$ local trivializations defined by the $\tau_i$. Let
$\delta_\lambda:S^4\to S^4$ be the conformal diffeomorphism of $S^4$
defined by
\begin{equation}
\delta_\lambda := \varphi_n\circ\lambda\circ \varphi_n^{-1},
\end{equation}
where $\lambda:\RR^4\to\RR^4$ is the dilation given by $y\mapsto
y/\lambda$, so $\delta_\lambda$ is a conformal diffeomorphism of $S^4$
which fixes the north and south poles.  We define cut-off, rescaled
connections on the bundles $\fg_{E^i}$ over $S^4$, with mass centers at the
north pole, by setting
\begin{equation}
(1-\chi_{x_i,\sqrt{\lambda_i}/2})\delta_{\lambda_i}^*A_i
:=
\begin{cases}
\delta_{\lambda_i}^*A_i
&\text{over $\varphi_n(B(0,\quarter\sqrt{\lambda_i}))$}, 
\\
\Gamma + (1-\chi_{x_i,\sqrt{\lambda_i}/2})\tau_i^*\delta_{\lambda_i}^*A_i
&\text{over $\varphi_n(\Omega(0,\quarter\sqrt{\lambda_i},
\half\sqrt{\lambda_i}))$,} 
\\
\Gamma 
&\text{over $S^4 - \varphi_n(B(0,\half\sqrt{\lambda_i}))$,} 
\end{cases} 
\end{equation}
where $\Gamma$ denotes the product connection on $(S^4-\{n\})\times \SO(3)$.
Note that $\tau_i^*\delta_{\lambda_i}^*A_i
= \delta_{\lambda_i}^*\tau_i^*A_i$, since the section $\tau_i$ is
defined by parallel translation from the south pole via the connection $A_i$.

A choice of $\SO(4)$ frames $f_i\in\Fr(TX)|_{x_i}$ defines coordinate charts
$$
\varphi_i^{-1}:B(x_i,r_0)\subset X\to\RR^4, \quad i=1,\dots,m, 
$$ 
via the exponential maps $\exp_{f_i}:B(0,r_0)\subset TX|_{x_i}\to X$.  
The orientation-preserving diffeomorphism $\varphi_i\circ \varphi_n^{-1}$
identifies the annulus
$\varphi_n(\Omega(0,\half\sqrt{\lambda_i},2\sqrt{\lambda_i}))$ in $S^4$, 
$$
\Omega\left(0,\thalf\sqrt{\lambda_i},2\sqrt{\lambda_i}\right)
:=
\varphi_n
\left(\left\{x\in\RR^4:\thalf\sqrt{\lambda_i}<|x|<2\sqrt{\lambda_i}
\right\}\right)
\subset \RR^4
$$
with the annulus in $X$,
$$
\Omega\left(x_i,\thalf\sqrt{\lambda_i},2\sqrt{\lambda_i}\right)
:=
\left\{x\in X: \thalf\sqrt{\lambda_i} < \dist_g(x,x_i) < 2\sqrt{\lambda_i} 
\right\}
\subset X.
$$
We define a glued-up $\SO(3)$ bundle $\fg_E$ over $X$ by setting
\begin{equation}
\fg_E 
= 
\begin{cases}
\fg_{E_\ell} &\text{over $X - \cup_{i=1}^m B(x_i,\half\sqrt{\lambda_i})$,} \\
\fg_{E^i} &\text{over $B(x_i,2\sqrt{\lambda_i})$.} 
\end{cases}
\label{eq:DefnGluedUpSO(3)Bundle}
\end{equation}
The bundles $\fg_{E_\ell}$ and $\fg_{E^i}$ are identified over the annuli
$\Omega(x_i,\half\sqrt{\lambda_i},2\sqrt{\lambda_i}\})$ in $X$ via the
isomorphisms of $\SO(3)$ bundles defined by the 
orientation-preserving diffeomorphisms $\varphi_i\circ \varphi_n^{-1}$,
identifying the annuli
$\Omega(x_i,\half\sqrt{\lambda_i},2\sqrt{\lambda_i}\})$ with the
corresponding annuli in $S^4$ and the $\SO(3)$ bundle maps defined by the
trivializations \eqref{eq:BackroundBundleTriv} and
\eqref{eq:SphereBundleTriv}.   
Define a spliced $L^2_k$ connection $A$ on $\fg_E$ by setting
\begin{equation}
\label{eq:SplicedConnection}
A
:= 
\begin{cases}
  A_0 &\text{over $X - \cup_{i=1}^m B(x_i,4\sqrt{\lambda_i})$},
  \\
  \Gamma +
  \chi_{x_i,4\sqrt{\lambda_i}}A_0 +
  (1-\chi_{x_i,\sqrt{\lambda_i}/2})\delta_{\lambda_i}^*A_i &\text{over
    $\Omega(x_i,\quarter\sqrt{\lambda_i},4\sqrt{\lambda_i})$},
  \\
  \delta_{\lambda_i}^*A_i &\text{over $B(x_i,\quarter\sqrt{\lambda_i})$,
    $i=1,\dots,m$,}
\end{cases} 
\end{equation}
where the cut-off connections are defined as above; the bundle and
annulus identifications are understood. Also implicit in definitions
\eqref{eq:SplicedSpinor} and \eqref{eq:SplicedConnection} 
are the following constraints on the
scales, given an $m$-tuple of distinct points, $(x_1,\dots,x_m)$,
\begin{equation}
\label{eq:ConstraintScales}
8\sqrt{\lambda_i} + 8\sqrt{\lambda_j} < \dist_g(x_i,x_j),
\quad i\neq j,
\end{equation}
representing a multiset $\bx\in\Sigma$.  The {\em splicing map\/},
$\bga'$, is defined on the open subset of the gluing-data bundle,
$\bGl(\sU_\ell,\Sigma)$ in definition
\eqref{eq:FiberedProduct}, defined by
\begin{equation}
\label{eq:ConstrainedGluingDataBundle}
\bGl^+(\sU_\ell,\Sigma)
:=
\{\pi^{-1}([A_0,\Phi_0],\bx)\in \bGl(\sU_\ell,\Sigma):
\text{$(\blambda,\bx)$ obeys \eqref{eq:ConstraintScales}} \}. 
\end{equation} 
We shall normally require that the scales $\lambda_i=\lambda[A_i]$ be
``sufficiently small'' (the precise requirement will emerge in \S
\ref{sec:Existence}) and so we shall write
$$
\bGl(\sU_\ell,\Sigma,\lambda_0)
\subset
\bGl^+(\sU_\ell,\Sigma),
$$
for the open subset defined by the requirement that $\lambda[A_i] <
\lambda_0$ for all $i$, for some positive constant $\lambda_0$, when this
bound is significant.

It remains to summarize what the splicing construction has achieved.
The construction, thus far, yields a splicing map giving a smooth
embedding
\begin{equation}
\label{eq:SplicingMap}
\bga':\bGl(\sU_\ell,\Sigma,\lambda_0)\to \sC^{*,0}(\ft)
\end{equation}
The fact that the map is a smooth embedding is a fairly straightforward
consequence of the definitions; it will follow from the stronger result in
\cite{FL4} that the gluing map $\bga$ gives a smooth embedding of the
gluing data. Similarly, the fact that the image of $\bga'$ contains only
irreducible, non-zero-section pairs is also a fairly easy consequence of
the definitions; it is proved in \cite{FL4}.

%end of file
%file: regularity.tex

\section{Regularity of solutions to extended PU(2) monopole equations}
\label{sec:Regularity}
In this section we recall the basic regularity results for solutions to the
$\PU(2)$ monopole equations \eqref{eq:PT}: these were proved in \cite{FL1},
so we just summarize the conclusions here, as well as proving regularity
for (gluing) solutions to the {\em extended $\PU(2)$ monopole equations\/}
(see \S \ref{subsec:MonoEqns}).  In common with the anti-self-dual
equations considered by Taubes in
\cite{TauIndef}, \cite{TauStable}, there are obstructions to solving the
full $\PU(2)$ monopole equations \eqref{eq:PT} and so we shall only prove
existence of solutions to the weaker `extended $\PU(2)$ monopole
equations', described in \S \ref{subsec:MonoEqns}. 

To encompass the regularity theory for these extended monopole equations,
it will be convenient to consider a more general, inhomogeneous version of
the equations \eqref{eq:PT}. Therefore, as in \cite[\S 3]{FL1}, we examine
a quasi-linear, inhomogeneous elliptic system consisting of a
generalization of the equations \eqref{eq:PT} and Coulomb gauge equation
for a pair $(A,\Phi)+(a,\phi)\in \tsC(\ft)$ relative to $(A,\Phi)$.  As
before, we choose $\ell\geq 4$ (so that $L^2_\ell \subset C^1$, with
$\tsC(\ft)$ the pre-configuration space of $L^2_\ell$ pairs acted on by
$L^2_{\ell+1}$ gauge transformations, $\sG_E$. Combining \eqref{eq:PT}
with the Coulomb gauge equation, and allowing inhomogeneous terms, we obtain
an elliptic system of equations for a pair $(a,\phi)$ in
$L^2_\ell(\Lambda^1\otimes\fg_E)\oplus L^2_\ell(V^+)$,
\begin{equation}
\label{eq:Coulomb+PT}
\begin{aligned}
d_{A,\Phi}^{0,*}(a,\phi) &= \zeta, 
\\
\fS(A+a,\Phi+\phi) &= (w_0,s_0).
\end{aligned}
\end{equation}
Considering $A$ to be a connection on $\fg_E$ and using the isomorphism
$\ad:\fg_E\to\so(\fg_E)$ to view $F_A$ as a section of $\La^2\otimes\fg_E$,
we write \eqref{eq:Coulomb+PT} as
\begin{align*}
d_{A}^*a - (\cdot\Phi)^*\phi &= \zeta, 
\\
F_{A+a}^+ 
- \tau\rho^{-1}((\Phi+\phi)\otimes(\Phi+\phi)^*)_{00} &=  w_0, 
\\ 
(D_{A+a}+\rho(\vartheta))(\Phi+\phi) &= s_0,
\end{align*}
where $(w_0,s_0)\in L^2_{\ell-1}(\Lambda^+\otimes\fg_E)\oplus
L^2_{\ell-1}(V^-)$. Recalling that $d_{A,\Phi}^0$ and
$d_{A,\Phi}^1$ are the differential operators in the elliptic
deformation complex \cite[Equation (2.37)]{FL1} for the $\PU(2)$ monopole
equations \eqref{eq:PT}, the above system may be rewritten in the form
\begin{align*}
d_{A,\Phi}^{0,*}(a,\phi) &= \zeta, 
\\
d_{A,\Phi}^1(a,\phi) + \{(a,\phi),(a,\phi)\} &= 
-\fS(A,\Phi) + (w_0,s_0) =: (w,s),
\end{align*}
where the differentials $d_{A,\Phi}^{0,*}$ and $d_{A,\Phi}^1$ are
given by \cite[Equations (2.38) \& (2.36)]{FL1}, respectively.
It will be convenient to view the
quadratic term $\{(a,\phi),(a,\phi)\}$ as being defined via the following
bilinear form, 
$$
\{(a,\phi),(b,\varphi)\}
:= \left(\begin{matrix}
(a\wedge b)^+ - \tau\rho^{-1}(\phi\otimes\varphi^*)_{00} \\
\rho(a)\varphi
\end{matrix}\right), 
$$
with $(b,\varphi)$ in $L^2_\ell(\Lambda^1\otimes\fg_E)\oplus
L^2_\ell(V^+)$. Our elliptic system \eqref{eq:Coulomb+PT} then takes the
simple shape
\begin{equation}
\sD_{A,\Phi}(a,\phi) + \{(a,\phi),(a,\phi)\} = (\zeta,w,s),
\label{eq:PTEllReg}
\end{equation}
recalling from \cite[Equation (2.39)]{FL1} that $\sD_{A,\Phi} =
d_{A,\Phi}^{0,*} + d_{A,\Phi}^1$.  This is the form of the
(inhomogeneous) Coulomb gauge and $\PU(2)$ monopole equations we will
use for the majority of the basic regularity arguments. 

By analogy with \cite{TauSelfDual}, \cite{TauStable}, we use
$(a,\phi)=d_{A,\Phi}^{1,*}(v,\psi)$, where $(v,\psi)\in
L^2_{\ell-1}(\Lambda^+\otimes\fg_E)\oplus L^2_{\ell-1}(V^-)$, to rewrite
the first-order quasi-linear equation \eqref{eq:PTEllReg}
for $(a,\phi)$ in terms of $(w,s)$,
\begin{equation}
\label{eq:QuasiSecondPUMonCoulomb}
d_A^+d_{A,\Phi}^{1,*}(v,\psi) 
+ \{d_{A,\Phi}^{1,*}(v,\psi),d_{A,\Phi}^{1,*}(v,\psi)\} = (w,s),
\end{equation}
to give a second-order, quasi-linear elliptic equation for $(v,\psi)$ in
terms of $(w,s)$; the Coulomb-gauge equation reduces to
$d_{A,\Phi}^{0,*}d_{A,\Phi}^{1,*}(v,\psi) = (\fS(A,\Phi)\cdot)^*(v,\psi)$.

The precise dependence of the constants $C$ below on the reference pair
$(A,\Phi)$ in this section is not something we track carefully, as we are
primarily interested in local $L^2_k$ estimates of solutions $(a,\phi)$ to
equations \eqref{eq:PTEllReg} or of solutions $(v,\psi)$ to
equations \eqref{eq:QuasiSecondPUMonCoulomb} when $k\geq 2$ or $k\geq 3$,
respectively. For this range of the Sobolev indices, the constants no
longer depend on simply $\|F_A\|_{L^2(X)}$ or $\|F_A^+\|_{L^\8(X)}$, for
example. Rather, they depend on higher-order covariant derivatives of the
pair $(A,\Phi)$. In general, throughout this section,
constants depending on $(A,\Phi)$ can be gauge-invariantly given as constants
depending on $L^2_{k,A}$ norms of $F_A$ or $L^2_{k+1,A}$ norms of $\Phi$,
following the method of \cite[Lemma 2.3.11]{DK}.

\subsection{Regularity for \boldmath{$L^2_1$} solutions to the inhomogeneous
Coulomb gauge and PU(2) monopole equations}
\label{subsec:L2_1InhomoGlobal}
We recall in this subsection that an $L^2_1$ solution $(a,\phi)$ to the
$\PU(2)$ monopole and Coulomb-gauge equations \eqref{eq:PTEllReg}, with an
$L^2_k$ inhomogeneous term (with $k \ge 1$) is in $L^2_{k+1}$. Thus, if the
inhomogeneous term is in $C^\8$ then $(a,\phi)$ is in $C^\8$. 

Given any $L^2_\ell$ orthogonal connection $A$ on $\fg_E$, fixed
unitary connections on $\det E$, and fixed \spinc connection on $W$
our Sobolev norms are defined in the usual way: for example, if $a\in
C^\8(\Lambda^1\otimes\fg_E)$, we write
$$
\|a\|_{L^p_{k,A}(X)} 
:= \left(\sum_{j=0}^k\|\cov_{A}^j a\|_{L^p(X)}^p\right)^{1/p},
$$
and if $(a,\phi)\in C^\8(\Lambda^1\otimes\fg_E)\oplus C^\8(V^+)$, we write 
$$
\|(a,\phi)\|_{L^p_{k,A}(X)} 
:= \left(\|a\|_{L^p_{k,A}(X)}^p 
+ \|\phi\|_{L^p_{k,A}(X)}^p\right)^{1/p}, 
$$
for any $1\le p\le \8$ and integer $k\ge 0$ for which $L^2_\ell\subset
L^p_k$ (see \cite[Theorem 5.4]{Adams}). Let
$L^p_k(\Lambda^1\otimes\fg_E)$ be the completion of
$C^\8(\Lambda^1\otimes\fg_E)$ with respect to this norm, defining
Sobolev norms and Banach-space completions of the other function spaces we
shall need in the obvious way. Occasionally, it will be useful to consider
different Sobolev exponents on one-forms and spinors when defining Banach
spaces of pairs, such as
$$
\|(a,\phi)\|_{L^{p;p'}(X)} 
=
\|a\|_{L^p(X)} 
+
\|\phi\|_{L^{p'}(X)}, 
$$
with other Sobolev norms defined in the analogous way.

\begin{prop}
\label{prop:L2_1InhomoReg}
\cite[Proposition 3.2]{FL1}
Let $X$ be a closed, oriented, Riemannian four-manifold with metric $g$,
\spinc structure $(\rho,W^+,W^-)$, and let $E$ be a Hermitian two-plane bundle
over $X$.  Let $(A,\Phi)$ be a $C^\8$ pair on the bundles
$(\fg_E,V^+)$ over $X$ and let $2\le p < 4$.  Then there are
positive constants $\eps=\eps(A,\Phi,p)$ and $C=C(A,\Phi,p)$ with
the following significance. Suppose that $(a,\phi)\in
L^2_1(X,\La^1\otimes\fg_E)\oplus L^2_1(X,V^+)$ is an $L^2_1$
solution on $(\fg_E,V^+)$ to the elliptic system
\eqref{eq:PTEllReg} over $X$, where
$(\zeta,w,s)$ is in $L^p$.  If $\|(a,\phi)\|_{L^4(X)}<\eps$ then
$(a,\phi)$ is in $L^p_1$ and
$$
\|(a,\phi)\|_{L^p_{1,A}(X)} \le C\left(\|(\zeta,w,s)\|_{L^p(X)}
+ \|(a,\phi)\|_{L^2(X)}\right).
$$
\end{prop}

\begin{prop}
\label{prop:Lp_1InhomoReg}
\cite[Proposition 3.3]{FL1}
Continue the notation of Proposition \ref{prop:L2_1InhomoReg}.  Let $k\ge
1$ be an integer and let $2<p<\8$. Let $(A,\Phi)$ be a $C^\8$ pair
on the bundles $(\fg_E,V^+)$ over $X$. Suppose that
$(a,\phi)\in L^p_1(X,\La^1\otimes\fg_E)\oplus L^p_1(X,V^+)$ 
is a solution on $(\fg_E,V^+)$ to the elliptic system
\eqref{eq:PTEllReg} over $X$, where $(\zeta,w,s)$ is in $L^2_k$. Then
$(a,\phi)$ is in $L^2_{k+1}$ and there is a universal polynomial
$Q_k(x,y)$, with positive real coefficients, depending at most on
$(A,\Phi),k$, such that $Q_k(0,0)=0$ and
$$
\|(a,\phi)\|_{L^2_{k+1,A}(X)} 
\le Q_k\left(\|(\zeta,w,s)\|_{L^2_{k,A}(X)},
\|(a,\phi)\|_{L^p_{1,A}(X)}\right).
$$
In particular, if $(\zeta,w,s)$ is in $C^\8$ then $(a,\phi)$ is in
$C^\8$ and if $(\zeta,w,s)=0$, then
$$
\|(a,\phi)\|_{L^2_{k+1,A}(X)} \le C\|(a,\phi)\|_{L^p_{1,A}(X)}.
$$
\end{prop}

By combining Propositions \ref{prop:L2_1InhomoReg} and \ref{prop:Lp_1InhomoReg}
we obtain the desired regularity result for
$L^2_1$ solutions to the inhomogeneous Coulomb gauge and $\PU(2)$
monopole equations:

\begin{cor}
\label{cor:L2_1InhomoReg}
\cite[Corollary 3.4]{FL1}
Continue the notation of Proposition \ref{prop:L2_1InhomoReg}.  Let
$(A,\Phi)$ be a $C^\8$ pair on the bundles $(\fg_E,V^+)$ over
$X$.  Then there is a positive constant $\eps=\eps(A,\Phi)$ such that
the following hold.  Suppose that $(a,\phi)\in
L^2_1(X,\La^1\otimes\fg_E)\oplus L^2_1(X,V^+)$ is a solution on
$(\fg_E,V^+)$ to the elliptic system
\eqref{eq:PTEllReg} over $X$, where $(\zeta,w,s)$ is in $L^2_k$ and 
$\|(a,\phi)\|_{L^4(X)}<\eps$ and $k\ge 0$ is an integer.
Then $(a,\phi)$ is in $L^2_{k+1}$ and 
there is a universal polynomial $Q_k(x,y)$, with positive real coefficients,
depending at most on $(A,\Phi),k$, such that $Q_k(0,0)=0$ and
$$
\|(a,\phi)\|_{L^2_{k+1,A}(X)} 
\le Q_k\left(\|(\zeta,w,s)\|_{L^2_{k,A}(X)},
\|(a,\phi)\|_{L^2(X)}\right).
$$
In particular, if $(\zeta,w,s)$ is in $C^\8$ then $(a,\phi)$ is in
$C^\8$ and if $(\zeta,w,s)=0$, then
$$
\|(a,\phi)\|_{L^2_{k+1,A}(X)} \le C\|(a,\phi)\|_{L^2(X)}.
$$
\end{cor}

\subsection{Regularity for \boldmath{$L^2_2$} solutions to the inhomogeneous,
second-order quasi-linear equation}
\label{subsec:L2_2SecondOrderInhomoGlobal}
It remains to consider the regularity properties of the solution $(v,\psi)$
to the second-order equation \eqref{eq:QuasiSecondPUMonCoulomb}, where
$(a,\phi)=d_{A,\Phi}^{1,*}(v,\psi)$ and thus
$$
d_{A,\Phi}^1d_{A,\Phi}^{1,*}(v,\psi) = d_{A,\Phi}^1(a,\phi).
$$
The Laplacian $d_{A,\Phi}^1d_{A,\Phi}^{1,*}$ is elliptic, with $C^\8$
coefficients, and thus standard regularity theory implies that if
$(a,\phi)$ is in $L^p_k$ then $d_{A,\Phi}^1(a,\phi)$ is in $L^p_{k-1}$ and
so $(v,\psi)$ is in $L^p_{k+1}$.

In order to simultaneously treat regularity for solutions to the (extended)
$\PU(2)$ monopole equations and their differentials, we consider the
following generalization of the second-order system
\eqref{eq:QuasiSecondPUMonCoulomb}:
\begin{equation}
\label{eq:GeneralSecondOrderPUMon}
d_{A,\Phi}^1d_{A,\Phi}^{1,*}(v,\psi) 
+ \{d_{A,\Phi}^{1,*}(v,\psi), d_{A,\Phi}^{1,*}(v,\psi)\}
+ \langle(\alpha,t),d_{A,\Phi}^{1,*}(v,\psi)\rangle = (w,s),
\end{equation}
where $(v,\psi)\in L^2_2(X,\Lambda^+\otimes\fg_E)\oplus L^2_2(X,V^-)$,
$(w,s)\in L^p_k(X,\Lambda^+\otimes\fg_E)\oplus L^p_k(X,V^+)$, and
$(\alpha,t) \in C^\8(X,\Lambda^1\otimes\fg_E)\oplus C^\8(X,V^+)$. Here,
$$
\langle\, ,\, \rangle:\Gamma(\Lambda^1\otimes\fg_E)\oplus\Gamma(V^+)
\times
\Gamma(\Lambda^1\otimes\fg_E)\oplus\Gamma(V^+)
\to 
\Gamma(\Lambda^+\otimes\fg_E)\oplus\Gamma(V^-),
$$
denotes a bilinear map. Then the proofs of the regularity results of
\S \ref{subsec:L2_1InhomoGlobal} adapt with essentially no change to give
the following regularity results for 
equation \eqref{eq:GeneralSecondOrderPUMon}.

\begin{prop}
\label{prop:L2_2SecondOrderInhomoReg}
Let $X$ be a closed, oriented, Riemannian four-manifold with metric $g$,
\spinc structure $(\rho,W^+,W^-)$, and let $E$ be a Hermitian two-plane bundle
over $X$.  Let $(A,\Phi)$ be a $C^\8$ pair on the bundles
$(\fg_E,V^+)$ over $X$ and let $2\le p < 4$.  Suppose $(\alpha,t)
\in C^\8(X,\Lambda^1\otimes\fg_E)\oplus C^\8(X,V^+)$.  Then there are
positive constants $\eps=\eps(A,\Phi,\alpha,t,p)$ and
$C=C(A,\Phi,\alpha,t,p)$ with the following significance. Suppose that
$(v,\psi)\in L^2_2(X,\La^1\otimes\fg_E)\oplus L^2_2(X,V^-)$ is an
$L^2_2$ solution to the elliptic system
\eqref{eq:GeneralSecondOrderPUMon} over $X$, where
$(w,s)$ is in $L^p$.  If $\|d_{A,\Phi}^{1,*}(v,\psi)\|_{L^4(X)}<\eps$ then
$(v,\psi)$ is in $L^p_2$ and
$$
\|(v,\psi)\|_{L^p_{2,A}(X)} \le C\left(\|(w,s)\|_{L^p(X)}
+ \|(v,\psi)\|_{L^2(X)}\right).
$$
\end{prop}

\begin{prop}
\label{prop:Lp_2SecondOrderInhomoReg}
Continue the notation of Proposition \ref{prop:L2_1InhomoReg}.  Let $k\ge
1$ be an integer and let $2<p<\8$. Let $(A,\Phi)$ be a $C^\8$ pair
on the bundles $(\fg_E,V^+)$ over $X$. Suppose $(\alpha,t)
\in C^\8(X,\Lambda^1\otimes\fg_E)\oplus C^\8(X,V^+)$.  Suppose that
$(v,\psi)\in L^p_1(X,\La^1\otimes\fg_E)\oplus L^p_1(X,V^+)$ 
is a solution to the elliptic system
\eqref{eq:GeneralSecondOrderPUMon} over $X$, where $(w,s)$ is in
$L^2_k$. Then 
$(v,\psi)$ is in $L^2_{k+2}$ and there is a universal polynomial
$Q_k(x,y)$, with positive real coefficients, depending at most on
$(A,\Phi),(\alpha,t),k$, such that $Q_k(0,0)=0$ and
$$
\|(v,\psi)\|_{L^2_{k+2,A}(X)} 
\le Q_k\left(\|(w,s)\|_{L^2_{k,A}(X)},
\|(v,\psi)\|_{L^p_{2,A}(X)}\right).
$$
In particular, if $(w,s)$ is in $C^\8$ then $(v,\psi)$ is in
$C^\8$ and if $(w,s)=0$, then
$$
\|(v,\psi)\|_{L^2_{k+2,A}(X)} \le C\|(v,\psi)\|_{L^p_{2,A}(X)}.
$$
\end{prop}

By combining Propositions \ref{prop:L2_1InhomoReg} and \ref{prop:Lp_1InhomoReg}
we obtain the desired regularity result for
$L^2_2$ solutions to the inhomogeneous, general second-order $\PU(2)$
monopole equations:

\begin{cor}
\label{cor:L2_2SecondOrderInhomoReg}
Continue the notation of Proposition \ref{prop:L2_1InhomoReg}.  Let
$(A,\Phi)$ be a $C^\8$ pair on the bundles $(\fg_E,V^+)$ over
$X$.  Suppose $(\alpha,t)
\in C^\8(X,\Lambda^1\otimes\fg_E)\oplus C^\8(X,V^+)$.  
Then there is a positive constant $\eps=\eps(A,\Phi,\alpha,t)$ such
that the following hold.  Suppose that $(v,\psi)\in
L^2_1(X,\La^1\otimes\fg_E)\oplus L^2_1(X,V^+)$ is a solution
to the elliptic system
\eqref{eq:GeneralSecondOrderPUMon} over $X$, where $(w,s)$ is in
$L^2_k$ and $\|d_{A,\Phi}^{1,*}(v,\psi)\|_{L^4(X)}<\eps$ and $k\ge 0$
is an integer.  Then $(v,\psi)$ is in $L^2_{k+2}$ and there is a universal
polynomial $Q_k(x,y)$, with positive real coefficients, depending at most
on $(A,\Phi),(\alpha,t),k$, such that $Q_k(0,0)=0$ and
$$
\|(v,\psi)\|_{L^2_{k+2,A}(X)} 
\le Q_k\left(\|(w,s)\|_{L^2_{k,A}(X)},
\|(v,\psi)\|_{L^2(X)}\right).
$$
In particular, if $(w,s)$ is in $C^\8$ then $(v,\psi)$ is in
$C^\8$ and if $(w,s)=0$, then
$$
\|(v,\psi)\|_{L^2_{k+2,A}(X)} \le C\|(v,\psi)\|_{L^2(X)}.
$$
\end{cor}

\subsection{Local regularity and interior estimates for \boldmath{$L^2_1$}
solutions to the inhomogeneous Coulomb gauge and PU(2) monopole equations}
\label{subsec:L2_1InhomoLocal} 
In this section we specialize the results of \S 
\ref{subsec:L2_1InhomoGlobal} to the case where the reference pair is a
trivial $\PU(2)$ monopole, so $(A,\Phi)=(\Ga,0)$ on the bundles
$(\fg_E,V^+)$, over an open subset $\Om\subset X$, where $\Ga$ is
a flat connection.

We continue to assume that $X$ is a closed, oriented four-manifold with
metric $g$, \spinc bundle $W$, and Hermitian two-plane bundle $E$ extending
those on $\Om\subset X$.  We use these inhomogeneous estimates and
regularity results to show that a global $L^2_1$ gluing solution $(a,\phi)$
to the extended $PU(2)$ monopole equations is actually $C^\8$.

We have the following local versions of Propositions
\ref{prop:L2_1InhomoReg} and \ref{prop:Lp_1InhomoReg} and
Corollary \ref{cor:L2_1InhomoReg}:

\begin{prop}
\label{prop:L2_1InhomoRegLocal}
\cite[Proposition 3.9]{FL1}
Continue the notation of the preceding paragraph. Let $\Om'\Subset\Om$ be a
precompact open subset and let $2\le p < 4$.  Then there are positive
constants $\eps=\eps(\Om,p)$ and $C=C(\Om',\Om,p)$ with the following
significance. Suppose that $(a,\phi)$ is an $L^2_1(\Om)$ solution to the
elliptic system \eqref{eq:PTEllReg} over $\Om$, with $(A,\Phi)=(\Ga,0)$
and where $(\zeta,w,s)$ is in $L^p(\Om)$. If $\|(a,\phi)\|_{L^4(\Om)}<\eps$
then $(a,\phi)$ is in $L^p_1(\Om')$ and
$$
\|(a,\phi)\|_{L^p_{1,\Ga}(\Om')} 
\le C\left(\|(\zeta,w,s)\|_{L^p(\Om)} 
+ \|(a,\phi)\|_{L^2(\Om)}\right).
$$
\end{prop}

\begin{prop}
\label{prop:Lp_1InhomoRegLocal}
\cite[Proposition 3.10]{FL1}
Continue the notation of Proposition \ref{prop:L2_1InhomoRegLocal}.  Let
$k\ge 1$ be an integer, and let $2<p<\8$. Suppose that $(a,\phi)$ is an
$L^p_1(\Om)$ solution to the elliptic system
\eqref{eq:PTEllReg} over $\Om$ with $(A,\Phi)=(\Ga,0)$, 
where $(\zeta,w,s)$ is in $L^2_k(\Om)$. Then $(a,\phi)$ is in
$L^2_{k+1}(\Om')$ and there is a universal polynomial $Q_k(x,y)$, with
positive real coefficients, depending at most on
$k$, $\Om'$, $\Om$, such that
$Q_k(0,0)=0$ and
$$
\|(a,\phi)\|_{L^2_{k+1,\Ga}(\Om')} 
\le Q_k\left(\|(\zeta,w,s)\|_{L^2_{k,\Ga}(\Om)},
\|(a,\phi)\|_{L^p_{1,\Ga}(\Om)}\right).
$$
If $(\zeta,w,s)$ is in $C^\8(\Om)$ then $(a,\phi)$ is in $C^\8(\Om')$
and if $(\zeta,w,s)=0$, then
$$
\|(a,\phi)\|_{L^2_{k+1,\Ga}(\Om')} 
\le C\|(a,\phi)\|_{L^p_{1,\Ga}(\Om)}.
$$
\end{prop}

\begin{cor}
\label{cor:L2_1InhomoRegLocal}
\cite[Corollary 3.11]{FL1}
Continue the notation of Proposition \ref{prop:Lp_1InhomoRegLocal}.  Then
there is a positive constant $\eps=\eps(\Om)$ with the following
significance. Suppose that $(a,\phi)$ is an $L^2_1(\Om)$ solution to the
elliptic system \eqref{eq:PTEllReg} over $\Om$ with $(A,\Phi)=(\Ga,0)$,
where $(\zeta,w,s)$ is in $L^2_k(\Om)$ and $\|(a,\phi)\|_{L^4(\Om)}<\eps$.
Then $(a,\phi)$ is in $L^2_{k+1}(\Om')$ and there is a universal polynomial
$Q_k(x,y)$, with positive real coefficients, depending at most on $k$,
$\Om'$, $\Om$, such that $Q_k(0,0)=0$ and
$$
\|(a,\phi)\|_{L^2_{k+1,\Ga}(\Om')} 
\le Q_k\left(\|(\zeta,w,s)\|_{L^2_{k,\Ga}(\Om)},
\|(a,\phi)\|_{L^2(\Om)}\right).
$$
If $(\zeta,w,s)$ is in 
$C^\8(\Om)$ then $(a,\phi)$ is in $C^\8(\Om')$ and if
$(\zeta,w,s)=0$, then
$$
\|(a,\phi)\|_{L^2_{k+1,\Ga}(\Om')} \le C\|(a,\phi)\|_{L^2(\Om)}.
$$
\end{cor}

Corollary \ref{cor:L2_1InhomoRegLocal} thus yields a sharp local
elliptic regularity result for $\PU(2)$ monopoles $(A,\Phi)$ in $L^2_1$ 
which are given to us in Coulomb gauge relative to $(\Ga,0)$.

\begin{prop}
\label{prop:L2_1CoulMonoRegLocal}
\cite[Proposition 3.12]{FL1}
Continue the notation of Corollary \ref{cor:L2_1InhomoRegLocal}.  Then
there is a positive constant $\eps=\eps(\Om)$ and, if $k\ge 1$ is an
integer, there is a positive constant $C=C(\Om',\Om,k)$ with the following
significance. Suppose that $(A,\Phi)$ is an $L^2_1$ solution to the
$\PU(2)$ monopole equations \eqref{eq:PT} over $\Om$, which is in Coulomb
gauge over $\Om$ relative to $(\Ga,0)$, so $d_\Ga^*(A-\Ga) = 0$, and obeys
$\|(A-\Ga,\Phi)\|_{L^4(\Om)}<\eps$.  Then $(A-\Ga,\Phi)$ is in $C^\8(\Om')$
and for any $k\ge 1$,
$$
\|(A-\Ga,\Phi)\|_{L^2_{k,\Ga}(\Om')} 
\le C\|(A-\Ga,\Phi)\|_{L^2(\Om)}.
$$
\end{prop}

\subsection{Estimates for PU(2) monopoles in a good local
gauge}\label{subsec:Uhlenbeck} 
It remains to combine the local regularity results and estimates of \S
\ref{subsec:L2_1InhomoLocal}, for $\PU(2)$ monopoles $(A,\Phi)$ where
the connection $A$ is assumed to be in Coulomb gauge relative to the
product $\SO(3)$ connection $\Ga$, with Uhlenbeck's local, Coulomb
gauge-fixing theorem. We then obtain regularity results and estimates for
$\PU(2)$ monopoles $(A,\Phi)$ with small curvature $F_A$, parallel to those
of Theorem 2.3.8 and Proposition 4.4.10 in \cite{DK} for anti-self-dual
connections.

In order to apply Corollary \ref{cor:L2_1InhomoRegLocal} we need
Uhlenbeck's Coulomb gauge-fixing result
\cite[Theorem 2.1 \& Corollary 2.2]{UhlLp}). Let $B$ (respectively,
$\barB$) be the open (respectively, closed) unit ball centered at
the origin in $\RR^4$ and let $G$ be a compact
Lie group. In order to provide universal constants we assume $\RR^4$ has
its standard metric, though the results of this subsection naturally hold
for any $C^\8$ Riemannian metric, with comparable constants for metrics
which are suitably close. 

\begin{thm}
\label{thm:CoulombBallGauge}
There are positive constants $c$ and $\eps$ with the following
significance. If $2\le p < 4$ is a constant and $A\in
L^p_1(B,\La^1\otimes\fg)\cap L^p_1(\rd B,\La^1\otimes\fg)$ is a
connection matrix whose curvature satisfies $\|F_A\|_{L^p(B)} < \eps$,
then there is a gauge transformation 
$u\in L^p_2(B, G)\cap  L^p_2(\rd B, G)$ such that $u(A)
:= uAu^{-1} - (du)u^{-1}$ satisfies
\begin{align}
d^*u(A) &= 0\text{ on }B, \tag{1}\\
{\textstyle{\frac{\rd}{\rd r}}}\lrcorner u(A) &= 0\text{ on }\rd B, \tag{2}\\
\|u(A)\|_{L^p_1(B)} &\le c\|F_A\|_{L^p(B)}.\tag{3}
\end{align}
If $A$ is in $L^p_k(B)$, for $k\ge 2$, then $u$ is in $L^p_{k+1}(B)$. The gauge
transformation $u$ is unique up to multiplication by a constant element of $G$.
\end{thm}

\begin{rmk}
If $G$ is Abelian then the requirement that $\|F_A\|_{L^p(B)} < \eps$ can be
omitted. 
\end{rmk} 

It is often useful to rephrase Theorem \ref{thm:CoulombBallGauge} in two
other slightly different ways. Suppose $A$ is an $L^2_k$ connection on a
$C^\8$ principal $G$ bundle $P$ over $B$ with $k\ge 2$ and $\|F_A\|_{L^2(B)} <
\eps$. Then the assertions of Theorem \ref{thm:CoulombBallGauge} are 
equivalent to each of the following:
\begin{itemize}
\item There is an $L^2_{k+1}$ trivialization $\tau:P\to B\times G$ such
that (i) $d_\Gamma^*(\tau(A)-\Gamma)=0$, where $\Gamma$ is the product
connection on $B\times G$, (ii) $\frac{\rd}{\rd r}
\lrcorner (\tau(A)-\Gamma) = 0$, and (iii)
$\|(\tau(A)-\Gamma)\|_{L^2_1(B)} \le c\|F_A\|_{L^2(B)}$.

\item There is an $L^2_{k+1}$ flat connection $\Gamma$ on $P$ such
that (i) $d_\Gamma^*(A-\Gamma)=0$, (ii) $\frac{\rd}{\rd r}
\lrcorner (A-\Gamma) = 0$, and (iii)
$\|(A-\Gamma)\|_{L^2_1(B)} \le c\|F_A\|_{L^2(B)}$, and
an $L^2_{k+1}$ trivialization $P|_B\simeq B\times G$
taking $\Gamma$ to the product connection.
\end{itemize}

We can now combine Theorem \ref{thm:CoulombBallGauge} with Proposition
\ref{prop:L2_1CoulMonoRegLocal} to give the following analogue of Theorem
2.3.8 in \cite{DK} --- the interior estimate for anti-self-dual connections
with $L^2$-small curvature.  

\begin{cor}
\label{cor:PTLocalReg}
Let $B\subset\RR^4$ be the open unit ball with center at the origin with
\spinc structure $(\rho,W^+,W^-)$, let $U\Subset B$ be an open subset, and let
$\Ga$ be the product connection on $B\times\SO(3)$.  Then there is a
positive constant $\eps$ and if $\ell\ge 1$ is an integer, there is a
positive constant $C(\ell,U)$ with the following significance.  Suppose
that $(A,\Phi)$ is an $L^2_1$ solution to the $\PU(2)$ monopole equations
\eqref{eq:PT} over $B$ and that the curvature of the $\SO(3)$ connection
matrix $A$ obeys $\|F_A\|_{L^2(B)} < \eps$.  Then there is an $L^2_{k+1}$
gauge transformation $u:B\to\SU(2)$ such that $(u(A)-\Ga,u\Phi)$ is in
$C^\8(B)$ with $d^*(u(A)-\Ga)=0$ over $B$ and
$$
\|(u(A)-\Ga,u\Phi)\|_{L^2_{\ell,\Gamma}(U)} 
\le C\|F_A\|_{L^2(B)}^{1/2}.
$$
\end{cor}

Again, it is often useful to rephrase Corollary \ref{cor:PTLocalReg} in the two
other slightly different ways. Suppose $\ell\ge 1$ and that $(A,\Phi)$ is a
$\PU(2)$ monopole in $L^2_k$ on $(\fg_E,V^+)$ over the unit ball
$B\subset\RR^4$ with $k\ge \max\{2,\ell\}$, 
$\|F_A\|_{L^2(B)}<\eps$ and $U\Subset B$. Then the
assertions of Corollary \ref{cor:PTLocalReg} are  
equivalent to each of the following: 
\begin{itemize}
\item There is a $C^\8$ trivialization $\tau:E|_B\to
B\times \CC^2$ and a $L^2_{k+1}$
determinant-one, unitary bundle automorphism $u$ of
$E|_B$ such that, with respect to the product connection $\Gamma$ on
$B\times \su(2)$, we have
(i) $d_\Gamma^*(\tau u(A)-\Gamma)=0$, and (ii)
$\|(\tau u(A)-\Gamma,\tau u\Phi)\|_{L^2_{\ell,\Gamma}(U)} 
\le C\|F_A\|_{L^2(B)}^{1/2}$.

\item There is an $L^2_{k+1}$ flat connection $\Gamma$ on $\fg_E|_B$ such
that (i) $d_\Gamma^*(A-\Gamma)=0$, and (ii)
$\|(A-\Gamma,\Phi)\|_{L^2_{\ell,\Gamma}(U)} \le c\|F_A\|_{L^2(B)}^{1/2}$,
and an $L^2_{k+1}$ trivialization $\fg_E|_B\simeq B\times \su(2)$
taking $\Gamma$ to the product connection.
\end{itemize}

We will also need interior estimates for $\PU(2)$ monopoles in a
good local gauge over more general simply-connected
regions than the open balls
considered in Corollary \ref{cor:PTLocalReg}. Specifically, recall
that a domain $\Om\subset X$ is {\em strongly simply-connected\/} if it has
an open covering by balls $D_1,\dots,D_m$ (not necessarily geodesic) such
that for $1\le r\le m$ the intersection $D_r\cap (D_1\cup\cdots\cup
D_{r-1})$ is connected. We recall (see \cite[Proposition 2.2.3]{DK} or
\cite[Proposition I.2.6]{Kobayashi}):

\begin{prop}
\label{prop:FlatProduct}
If $\Gamma$ is a $C^\8$ flat connection on a principal $G$ bundle
$P$ over a simply-connected manifold $\Om$, then there is a $C^\8$
isomorphism $P\simeq \Om\times G$ taking $\Gamma$ to
the product connection on $\Om\times G$. 
\end{prop}

More generally, if $A$ is $C^\8$ connection on a $G$ bundle $P$ over 
a simply-connected manifold-with-boundary $\barOm=\Om\cup\rd\Om$
with $L^p$-small curvature (with $p>2$), then Uhlenbeck's theorem implies
that $A$ is $L^p_2$-gauge equivalent to a connection which is $L^p_1$-close
to an $L^p_1$ flat connection on $P$ (see \cite[Corollary 4.3]{UhlChern} or
\cite[p. 163]{DK}). The following {\em a priori\/}
interior estimate is a straightforward generalization of
\cite[Proposition 4.4.10]{DK}.

\begin{prop}
\label{prop:MonopoleGoodGaugeIntEst}
Let $X$ be a closed, oriented, Riemannian four-manifold with \spinc
structure $(\rho,W^+,W^-)$ and let $\Om\subset X$ be a strongly simply-connected
open subset. Then there is a positive constant $\eps(\Om)$ with the
following significance. For $\Om'\Subset \Om$ a precompact open subset and
an integer $\ell\ge 1$, there is a constant $C(\ell,\Om',\Om)$ such that
the following holds. Suppose $(A,\Phi)$ is a $\PU(2)$ monopole in $L^2_k$
on $(\fg_E,V^+)$ over $\Om$ with $k\ge \max\{2,\ell\}$ such that
$$
\|F_A\|_{L^2(\Om)} < \eps.
$$
Then there is an $L^2_{k+1}$ flat connection $\Gamma$ on $\fg_E|_{\Om'}$ such
that 
$$
\|(A-\Ga,\Phi)\|_{L^2_{\ell,\Ga}(\Om')} \le C\|F_A\|_{L^2(\Om)}^{1/2},
$$
and an $L^2_{k+1}$ trivialization $\fg_E|_{\Om'}\simeq \Om'\times \su(2)$
taking $\Gamma$ to the product connection.
\end{prop}

\subsection{Regularity of gluing solutions to the
extended PU(2) monopole equations}
\label{subsec:TaubesRegularity}
Finally, we come to one of the main results of this section.

\begin{prop}
\label{prop:RegularityTaubesSolution}
Let $X$ be a closed, oriented, $C^\8$ four-manifold with 
$C^\8$ Riemannian metric $g$, \spinc structure $(\rho,W^+,W^-)$,
and let $E$ be a Hermitian, rank-two bundle
over $X$.  Let $(A,\Phi)$ be a $C^\8$ pair on the bundles
$(\fg_E,V^+)$ over $X$. Suppose that $(a,\phi)\in
L^2_1(X,\La^1\otimes\fg_E)\oplus L^2_1(X,V^+)$ is a
solution to an extended $\PU(2)$ monopole equation over $X$,
$$
\fS(A+a,\Phi+\phi) = (w,s),
$$
where $(a,\phi)=d_{A,\Phi}^{1,*}(v,\psi)$, for some $(v,\psi)\in (C^0\cap
L^2_2)(X,\La^+\otimes\fg_E)\oplus (C^0\cap L^2_2)(X,V^-)$, and $(w,s)\in
C^\8(X,\La^+\otimes\fg_E)\oplus C^\8(X,V^-)$. Then $(a,\phi)$ is in
$C^\8(X,\La^1\otimes\fg_E)\oplus C^\8(X,V^+)$ and 
$(v,\psi)\in C^\8(X,\La^+\otimes\fg_E)\oplus C^\8(X,V^+)$.
\end{prop}

\begin{rmk}
\label{rmk:AnalysisForTopologists}
It is important to note that the hypotheses of Proposition
\ref{prop:RegularityTaubesSolution} are strictly weaker than those of
Corollary \ref{cor:L2_1InhomoReg} because we not assume that $(a,\phi)$ is
$L^4$ small relative to some constant $\eps(A,\Phi)$ which depends, in a
possibly unfavorable way, on the pair $(A,\Phi)$.

Otherwise, the reader may wonder why we could not simply apply Corollary
\ref{cor:L2_1InhomoReg} to the elliptic system
\eqref{eq:GluingExtPUMonforRegularityProof} and deduce straightaway that
the solution $(a,\phi)$ is in $C^\8$.  The difficulty, of course, is that
the hypotheses of Corollary
\ref{cor:L2_1InhomoReg} are not necessarily satisfied due to the possibly
unfavorable dependence of the constant $\eps(A,\Phi)$ on the pair
$(A,\Phi)$. For example, by combining results of \S \ref{sec:Decay} and \S
\ref{sec:Existence}, we can show that $\|(a,\phi)\|_{L^4(X)} \leq C\lambda$
when $(A_\lambda,\Phi_\lambda)$ is a spliced pair (produced by the
algorithm of \S
\ref{sec:Splicing}), $\lambda$ is a parameter which can be made arbitrarily
small, and $C$ is an essentially universal constant, depending at most on
$\|F_{A_\lambda}^+\|_{L^{\sharp,2}(X)}$, $\|F_{A_{\lambda}}\|_{L^2}$, and
$\|\Phi_\lambda\|_{L^2_{2,A_\lambda}(X)}$: in particular, $C$ is independent
of $\lambda$. Now, as we make $\lambda$ smaller to try to satisfy the
constraint of Corollary
\ref{cor:L2_1InhomoReg}, we will in general find that
$\eps(A_\lambda,\Phi_\lambda)$ also becomes smaller (as it will, for
example, if $\eps(A_\lambda,\Phi_\lambda)$ depends on an $L^p_k$ norm of
$F_{A_\lambda}$ with $p>2$, $k=0$ or $k\geq 1$, $p=2$). Thus, as we do not
necessarily have a lower bound for $\eps(A,\Phi)$ which is uniform with
respect to $(A,\Phi)$, we are left ``chasing our tail'' and Corollary
\ref{cor:L2_1InhomoReg} will not apply.
\end{rmk}

\begin{proof}[Proof of Proposition \ref{prop:RegularityTaubesSolution}]
Since $a = d_{A,\Phi}^{1,*}(v,\psi)$, the extended $\PU(2)$ monopole and
Coulomb-gauge equations become
\begin{equation}
\begin{aligned}
\label{eq:GluingExtPUMonforRegularityProof}
d_{A,\Phi}^1(a,\phi) + \{(a,\phi),(a,\phi)\} &= (w,s)-\fS(A,\Phi),
\\
d_{A,\Phi}^{0,*}(a,\phi) &= (\fS(A,\Phi)\cdot)^*(v,\psi),
\end{aligned}
\end{equation}
where $\fS(A,\Phi) \in C^\8(\La^+\otimes\fg_E)\oplus C^\8(V^-)$ and
$(\fS(A,\Phi)\cdot)^*(v,\psi) \in (C^0\cap L^2_2)(\fg_E)$. Observe that
$$
(\zeta,w',s') := \left((\fS(A,\Phi)\cdot)^*(v,\psi),(w,s)-\fS(A,\Phi)\right)
\in 
L^2_2(\fg_E)\oplus L^2_2(\La^+\otimes\fg_E)\oplus L^2_2(V^-).
$$
{}From Remark \ref{rmk:AnalysisForTopologists} we note that we {\em
cannot\/} take the apparently obvious route and apply Corollary
\ref{cor:L2_1InhomoReg}. 
Instead, we fix a point $x\in X$ and a geodesic ball $B(x,\delta)$ centered
at $x$, fix a $C^\8$ trivialization $\fg_E|_{B(x,\delta)}\simeq
B(x,\delta)\times\su(2)$, and let $\Gamma$ be the resulting product
connection. Then
\begin{equation}
\label{eq:LocalPUMonRegularityEqn}
\begin{aligned}
d^1_{\Gamma,0}(b,\varphi) + \{(b,\varphi),(b,\varphi)\}
&= (w,s)-\fS(\Gamma,0) = (w,s),
\\
d^{0,*}_{\Gamma,0}(b,\varphi) 
&= (\fS(A,\Phi)\cdot)^*(v,\psi)
\end{aligned}
\end{equation}
where $b := A-\Gamma+a \in L^2_1(B(x,\delta),\La^1\otimes\fg_E)$, $\varphi
:= \Phi+\phi \in L^2_1(B(x,\delta),V^+)$, and the first equation is the
expansion of $\fS(\Gamma+b,0+\Phi+\phi)=(w,s)$. In particular, the right-hand
side of equation \eqref{eq:LocalPUMonRegularityEqn} is in
$L^2_1(B(x,\delta))$.

In Corollary \ref{cor:L2_1InhomoRegLocal}, let $\Omega = B\subset \RR^4$ be
the unit ball, let $\eps = \eps(\Omega)$ be the corresponding positive
constant, and let $\Omega' = B'$ be a slightly smaller concentric ball.
Since $(a,\phi)\in L^4(X,\La^1\otimes\fg_E)$, we know that
$\|(a,\phi)\|_{L^4(U)}\to 0$ as $\Vol(U)\to 0$, where $U$ is any measurable
subset of $X$. In particular, for a small enough ball $B(x,\delta)$, we
have $\|(b,\varphi)\|_{L^4(B(x,\delta))} < \eps$. We now rescale the metric
$g$ on the ball $B(x,\delta)$ to give a ball $\tilde B$ with unit radius,
center $x$, and metric which is approximately flat. Note that the equations
\eqref{eq:LocalPUMonRegularityEqn} are scale-equivariant (see \cite[\S
4.2]{FL1}), when $\varphi$ is 
replaced by $\tilde\varphi = \delta\varphi$ and 
the metric $g$ on $TX$ by $\tilde g =
\delta^{-2}g$, and that
$$
\|(b,\tilde\varphi)\|_{L^4(\tilde B,\tilde g)} 
=
\|(b,\varphi)\|_{L^4(B(x,\delta),g)} < \eps,
$$
by the scale-invariance of the $L^4$ norm on one-forms and rescaled
spinors. (Bear in mind that we write ``$g$'' in \cite[\S 4.2]{FL1} for the
metric $g^*$ on $T^*X$ and so our rescaling rule for monopoles agrees with
\cite[\S 4.1]{TelemanGenericMetric} and \cite[\S 8.1]{SalamonSWBook}, with
$\tilde g^* = \delta^2 g^*$.)  Corollary \ref{cor:L2_1InhomoRegLocal} now
implies that $(b,\tilde\varphi)$ --- and thus $(a,\tilde\phi)$ --- is in
$L^2_2$ on the slightly smaller concentric rescaled ball $\tilde B'\Subset
\tilde B$. Therefore, $(a,\phi)$ is in $L^2_2$ on the slightly smaller ball
$B'\Subset B$. Since $d_{A,\Phi}^1d_{A,\Phi}^{1,*}(v,\psi) =
d_{A,\Phi}^1(a,\phi)$ is now in $L^2_1(B')$, we have $(v,\psi)\in
L^2_3(B'')$ by standard elliptic regularity theory for a Laplacian,
$d_{A,\Phi}^1d_{A,\Phi}^{1,*}$, with $C^\8$ coefficients \cite{Hormander},
where $B''\Subset B'$ is again a slightly smaller concentric ball.  Hence,
the right-hand side of equation
\eqref{eq:LocalPUMonRegularityEqn} is in $L^2_2(B'')$.  Repeating the
above argument shows that $(a,\phi)$ is in $C^\8$ on $B(x,\delta/2)$ and
so, as the point $x\in X$ was arbitrary, we see that both $(a,\phi)$ and
$(v,\psi)$ are $C^\8$ over all of $X$.
\end{proof}

%end of file
%file: decay.tex

\section{Estimates for approximate PU(2) monopoles}
\label{sec:Decay}
The main purpose of this section is to derive an $L^{\sharp,2;2}$ estimate
(see Proposition \ref{prop:CutoffMonoPairEst}) for $\fS(A,\Phi)$, when
$(A,\Phi)$ is a pair --- an approximate solution to the extended $\PU(2)$
monopole equations --- produced by the splicing construction of \S
\ref{sec:Splicing}.

\subsection{Technical preliminaries}
\label{sec:DecayTechPrelim}
In this subsection we gather a few of the technical estimates we shall need
for the proof of Proposition \ref{prop:CutoffMonoPairEst} and, moreover,
when we derive estimates for eigenvectors and eigenvalues of the Laplacian
$d_{A,\Phi}^1 d_{A,\Phi}^{1,*}$ in \S \ref{sec:Eigenvalue}. 

The pointwise bounds we require are due to J. R\aa de \cite{Rade} in the
case of a Yang-Mills connection defined on an annulus in $\RR^4$ with its
Euclidean metric and to D. Groisser and T. Parker in the case of an annulus
in a four-manifold $X$ with an arbitrary Riemannian metric of bounded
geometry. Recall that $A$ is a {\em
Yang-Mills\/} connection if $d_A^*F_A = 0$ (and $d_AF_A=0$) and that if $A$
is anti-self-dual then it satisfies the Yang-Mills equations.

\begin{thm}
\label{thm:Decay}
\cite{DK}, \cite{GroisserParkerDecay}, \cite{Rade}
Let $X$ be a closed, oriented, four-manifold with metric $g$ and let $G$ be
a compact Lie group. Then there exist positive constants $c,\eps$ with the
following significance. If $A$ is a Yang-Mills connection on a $G$
bundle $P$ over an annulus $\Om(x_0;r_0,r_1)$ in $X$ with $0 < 4r_0 < r_1$
and whose curvature $F_A$ satisfies
$$
\|F_A\|_{L^2(\Om(x_0;r_0,r_1))} \le \eps, 
$$
and $r = \dist_g(x,x_0)$, then
$$
|F_A|(x) \le c\left(\frac{r_0^2}{r^4} + \frac{1}{r_1^2}\right)
\|F_A\|_{L^2(\Om(x_0;r_0,r_1))},\qquad 2r_0 \le r\le r_1/2.
$$
\end{thm}

A similar estimate, for anti-self-dual connections and flat metrics, was
proved by Donaldson using certain differential inequalities arising from
the Chern-Simons function \cite[Appendix]{DonApplic}, \cite[Proposition
7.3.3]{DK}.

Though not mentioned explicitly in \cite{DK}, \cite{Rade}, the proof of Theorem
\ref{thm:Decay} extends to give the following more general decay estimate:

\begin{cor}
\label{cor:Decay}
Continue the hypotheses of Theorem \ref{thm:Decay}. Then for any
integer $k\geq 0$ we have
\begin{equation}
\label{eq:GenRadeDecay}
|\cov_A^kF_A|(x) 
\leq 
\frac{c_k}{r^{k}}\left(\frac{r_0^2}{r^4} + \frac{1}{r_1^2}\right)
\|F_A\|_{L^2(\Omega(x_0;r_0,r_1))},\qquad 2r_0 \le r\le r_1/2.
\end{equation}
\end{cor}

For example, the corollary follows immediately by combining R\aa de's Lemma
2.2 and Theorem $1'$ (the version of his Theorem 1 for a cylinder
$(t_0,t_1)\times S^3$ in place of the annulus $\Omega(r_0,r_1)$).  

We shall also need to establish the existence of a universal
estimate for the $L^{\sharp,2}$ norm of the curvature of a
$g$-anti-self-dual connection.

\begin{lem}
\label{lem:UniversalLSharp2BoundFA}
\cite[Lemma 5.9]{FLKM1}
Let $(X,g)$ be a $C^\8$, closed, oriented, Riemannian four-manifold
and let $\kappa\geq 0$ be an integer. Then there is a constant
$c=c(X,g,\kappa)$ such that $\|F_A\|_{L^{\sharp,2}(X,g)} \leq c$ for
all $[A] \in M_\kappa^w(X,g)$.
\end{lem}

For any $N>1$ and $\la>0$, Lemma 7.2.10 in \cite{DK} provides a
clever construction of a cut-off function $\be$ on $\RR^4$ such that $\be(x)
= 0$ for $|x|\le N^{-1}\sqrt{\la}$, $\be(x) = 1$ for $|x|\ge
N\sqrt{\la}$, and satisfying the crucial estimate
$$
\|\cov\be\|_{L^4} \le c(\log N)^{-3/4},
$$
where $c$ is a constant independent of $N$ and $\la$. We will need the
following refinement of this result; for the definition and properties of
the $L^\sharp$ and related families of Sobolev norms, see \cite[\S
4]{FeehanSlice}.

\begin{lem}
\label{lem:dBetaEst}
\cite[Lemma 5.8]{FLKM1}
There is a positive constant $c$ such that the following holds. For any
$N>4$ and 
$\la>0$, there is a $C^\8$ cut-off function $\be=\be_{N,\la}$ on $\RR^4$
such that 
$$
\be(x) = 
\begin{cases}
1 &\text{if }|x|\ge \half\sqrt{\la}, \\
0 &\text{if }|x|\le N^{-1}\sqrt{\la}, 
\end{cases}
$$
and satisfying the following estimates:
\begin{align}
\|\cov\be\|_{L^4} &\le c(\log N)^{-3/4}, \tag{1}\\
\|\cov^2\be\|_{L^2} &\le c(\log N)^{-1/2}, \tag{2}\\
\|\cov\be\|_{L^{2\sharp}} &\le c(\log N)^{-1/2}, \tag{3}
\\
\|\cov\be\|_{L^2} &\le c(\log N)^{-1}\sqrt{\lambda}, \tag{4}
\\
\|\cov^2\be\|_{L^{4/3}} &\le c(\log N)^{-1}\sqrt{\lambda}. \tag{5}
\end{align}
\end{lem}

\subsection{Estimates of \boldmath{$\fS(A,\Phi)$} for spliced pairs}
\label{subsec:EstimateCutoffMonoPair}
An important application of the decay estimates for anti-self-dual
connections in Theorem \ref{thm:Decay} is to bound $\fS(A,\Phi)$ when
$(A,\Phi)$ is an approximate solution to the extended $\PU(2)$ monopole
equations produced by the splicing construction of \S
\ref{sec:Splicing}. As we shall discover in \S \ref{sec:Existence}, a
key measure of the size of $\fS(A,\Phi)$ in
$\Gamma(\Lambda^+\otimes\fg_E)\oplus \Gamma(V^+)$ is given by
\begin{equation}
\label{eq:SpecialNormPUMonMeasure}
\|(v,\psi)\|_{L^{\sharp,2;2}(X)}
=
\|v\|_{L^{\sharp,2}(X)}
+
\|\psi\|_{L^2(X)},
\end{equation}
where $\|v\|_{L^{\sharp,2}(X)} = \|v\|_{L^{\sharp}(X)} +
\|v\|_{L^{2}(X)}$ and the $L^\sharp$ family of Sobolev norms is defined in
\cite{FeehanSlice}, following Taubes' development in \cite{TauPath},
\cite{TauStable}, \cite{TauConf}.

\begin{prop}
\label{prop:CutoffMonoPairEst}
Let $(A,\Phi)$ be a spliced pair on the bundles $(\fg_E,V^+)$ over $(X,g)$,
produced by the splicing construction of \S \ref{sec:Splicing} (see
definitions \eqref{eq:SplicedSpinor} and
\eqref{eq:SplicedConnection}), where the background family
$\sU\Subset\sC(\ft_\ell)$ of pairs is a relatively open, precompact
submanifold. Then, there are constants
 and $c = c(g)$, $C = C(\sU,\kappa,g)$, $\lambda_0 =
\lambda_0(\sU,\kappa,g)$, such that for all $\lambda <\lambda_0$, where
$\lambda := \max_{1\leq i\leq m}\lambda_i$, we have
\begin{equation}
\label{eq:CutoffMonoPairEst}
\begin{aligned}
\|\fS(A,\Phi)\|_{L^{\sharp,2;2}(X,g)} 
&\leq 
C\lambda^{1/6}
+ \|\fS(A_0,\Phi_0)\|_{L^{\sharp,2;2}(X,g)} 
\\
&\quad 
+ \sum_{i=1}^m\|F^{+,g}(A_i)\|_{L^{\sharp,2}(B(x_i,\sqrt{\lambda_i}/4))},
\end{aligned}
\end{equation}
where $\lambda_i=\lambda[A_i]$.
If $F_{A_i}^{+,\delta} = 0$ on $\RR^4$ for all $i$ then
\begin{equation}
\label{eq:CutoffMonoPairWithASDS4Est}
\|\fS(A,\Phi)\|_{L^{\sharp,2;2}(X,g)} 
\leq 
C\lambda^{1/6}
+ \|\fS(A_0,\Phi_0)\|_{L^{\sharp,2;2}(X,g)}. 
\end{equation}
\end{prop}

\begin{pf}
Let $\eps$ be a positive constant which is less than or equal to that of
the decay estimate in Theorem \ref{thm:Decay} for the curvature of a
Yang-Mills connection. Chose a constant $0<\varrho=L^{-1}\leq 1$ small
enough that $\varrho\sqrt{8\pi^2\kappa_i} <
\eps$, for all $i\geq 1$. Choose $\lambda_0$ small enough that
$\sqrt{\lambda_0}\leq\varrho/16$.

For the sake of exposition we shall assume without loss in the proof that the
connections $A_i$ are $\delta$-anti-self-dual, even though this is not
required by our splicing construction or the hypotheses of Proposition
\ref{prop:CutoffMonoPairEst}. Indeed, for our application in \cite{FLConj}
the moduli spaces of mass-centered instantons on $S^4$ will be
replaced by diffeomorphic moduli spaces of spliced, almost-mass-centered,
approximate instantons. However, the changes required to accommodate the
more general situation needed in \cite{FLConj} are elementary.

The Chebychev inequality (see Lemma \ref{lem:Chebychev}) implies that
$$
\|F_{A_i}\|_{L^2(\RR^4-B(0,\lambda_i/\varrho),\delta)}
\leq
\varrho\sqrt{8\pi^2\kappa_i} < \eps,
$$
since $\varrho\sqrt{8\pi^2\kappa_i} < \eps$, by definition of
$\varrho$, and as (Definition \ref{defn:CenterScale})
$$
\lambda[A_i]^2
=
\frac{1}{8\pi^2\kappa}
\int_{\RR^4}|x|^2|F_{A_i}|_{\delta}^2\,d^4x = \lambda_i^2.
$$
Hence, provided $\lambda<\lambda_0$, the $\delta$-anti-self-dual
connections $A_i$ over $S^4$ obey the hypotheses of Theorem \ref{thm:Decay}
when $r_0 = \lambda_i/\varrho$ for small enough $\varrho$ and $r_1=\8$.

Recall from definitions \eqref{eq:SplicedSpinor} and
\eqref{eq:SplicedConnection} that
\begin{equation}
\label{eq:SplicedPair}
(A,\Phi)
:= 
\begin{cases}
  (A_0,\Phi_0) &\text{on $X - \cup_{i=1}^m B(x_i,8\lambda_i^{1/3})$},
  \\
  (A_0, \chi_{x_i,8\lambda_i^{1/3}}\Phi_0) 
&\text{on $\Omega(x_i;4\lambda_i^{1/3},8\lambda_i^{1/3})$},
  \\
  (\Gamma + \chi_{x_i,4\sqrt{\lambda_i}}\sigma_i^*A_0,0) &\text{on
    $\Omega(x_i;2\sqrt{\lambda_i},4\sqrt{\lambda_i})$},
  \\
  (\Gamma,0) &\text{on
    $\Omega(x_i;\sqrt{\lambda_i}/2,2\sqrt{\lambda_i})$},
  \\
    (\Gamma + (1-\chi_{x_i,\sqrt{\lambda_i}/2})
    (\varphi_n\circ\varphi_i^{-1})^*\tau_i^*A_i,0) &\text{on
    $\Omega(x_i;\sqrt{\lambda_i}/4,\sqrt{\lambda_i}/2)$},
  \\
    ((\varphi_n\circ\varphi_i^{-1})^*A_i,0)
    &\text{on $B(x_i,\sqrt{\lambda_i}/4)$},
\end{cases} 
\end{equation}
where the cut-off functions are defined in equation
\eqref{eq:ChiCutoffFunctionDefn}.  
Over the region $X - \cup_{i=1}^m B(x_i,8\lambda_i^{1/3})$ we have
$(A,\Phi) = (A_0,\Phi_0)$, while over the annuli
$\Omega(x_i;\sqrt{\lambda_i}/2,2\sqrt{\lambda_i})$ we have $F_A =
0$, and over the balls $B(x_i,\sqrt{\lambda_i}/4)$ we have $F_A =
F_{A_i}$. Thus, it remains to consider the three annuli
where cutting off occurs and the balls $B(x_i,\sqrt{\lambda_i}/4)$.

Over the region $X - \cup_{i=1}^m B(x_i,8\lambda_i^{1/3})$ the definition
\eqref{eq:SplicedPair} of $(A,\Phi)$ gives $(A,\Phi) = (A_0,\Phi_0)$
and so $\fS(A,\Phi) = \fS(A_0,\Phi_0)$. Our hypothesis on the precompact
family $\sU$ of background pairs $[A_0,\Phi_0]$ is that
$\|\fS(A_0,\Phi_0)\|_{L^{\sharp,2;2}(X)}$ is small:
\begin{equation}
\label{eq:LpEstimateMonoPairBallComp}
\|\fS(A,\Phi)\|_{L^{\sharp,2;2}(X - \cup_{i=1}^m B(x_i,8\lambda_i^{1/3}))}
\leq
\|\fS(A_0,\Phi_0)\|_{L^{\sharp,2;2}(X)}.
\end{equation}
Over the outer annulus $\Omega(x_i;4\lambda_i^{1/3},8\lambda_i^{1/3})$ the
definition \eqref{eq:SplicedPair} of $(A,\Phi)$ gives $(A,\Phi) =
(A_0,\chi_i\Phi_0)$ where we temporarily use the abbreviation $\chi_i = 
\chi_{x_i,8\lambda_i^{1/3}}$. Our expression \eqref{eq:PT} for
$\fS(A,\Phi)$ then yields
$$
\fS(A,\Phi)
=
\begin{pmatrix}
F_{A_0}^+ - \tau\rho^{-1}(\chi_i\Phi_0\otimes\chi_i\Phi_0^*)_{00},
\\
\chi_iD_{A_0}\Phi_0 + \rho(\vartheta+d\chi_i)\Phi_0
\end{pmatrix}.
$$
Given $2< p \leq \8$, choose $4 < q\leq \8$ for the remainder of the proof
by setting $1/p=1/4+1/q$.  Since $[A_0,\Phi_0]$ varies in a precompact
family, using the $L^4$ estimate for $d\chi_i$ in Lemma \ref{lem:dBetaEst}
and the preceding expression for $\fS(A,\Phi)$ yields
\begin{align}
%\begin{equation}
%\begin{aligned}
\notag
&\|\fS(A,\Phi)\|_{L^{\sharp,2;2}(\Omega(x_i;4\lambda_i^{1/3},8\lambda_i^{1/3}))}
\\
\label{eq:LpEstimateMonoPairOuterAnnulus}
&\leq
c\|\fS(A,\Phi)\|_{L^p(\Omega(x_i;4\lambda_i^{1/3},8\lambda_i^{1/3}))}
\quad\text{(by \cite[Lemma 4.1]{FeehanSlice}, for $2<p\leq\8$)} 
\\
\notag
&\leq
C\left(\lambda_i^{4/(3p)} 
+ \|d\chi_i\|_{L^4}
\|\Phi_0\|_{L^q(\Omega(x_i;4\lambda_i^{1/3},8\lambda_i^{1/3}))}\right)
\\
\notag
&\leq
C(\lambda_i^{4/(3p)} + \lambda_i^{4/(3q)})
\leq
C\lambda_i^{4/(3p)-1/3}.
%\end{aligned}
%\end{equation}
\end{align}
This concludes the estimate for $\fS(A,\Phi)$ over the outer annulus.

Over the middle annulus $\Omega(x_i;2\sqrt{\lambda_i},4\sqrt{\lambda_i})$ 
we see that \eqref{eq:SplicedPair} yields
\begin{equation}
\label{eq:ExpandCurvMiddleAnnulus}
F_A = \chi_iF_{A_0} + d\chi_i\wedge\si_i^*A_0 
+ (\chi_i^2-\chi_i)\si_i^*A_0\wedge \si_i^*A_0,
\end{equation}
where we set $\chi_i = \chi_{x_i,4\sqrt{\lambda_i}}$ for convenience.
Since the connection one-form $\sigma_i^*A_0$ is in radial gauge with
respect to the point $x_i\in X$, we have $|\varphi_i^*\sigma_i^*A_0|(x)
\leq C|x|$, where $C=\|F_{A_0}\|_{L^\8(B(x_i,\varrho/2))}$, and thus for
any $1\leq p \leq \8$, 
\begin{equation}
\label{eq:LpEstimateBackroundConnOneFormBall}
\begin{aligned}
\|\sigma_i^*A_0\|_{L^p(B(x_i,4\sqrt{\lambda_i}))}
&\leq
c\lambda_i^{2/p}\|\sigma_i^*A_0\|_{L^\8(B(x_i,4\sqrt{\lambda_i}))}
\\
&\leq
C\lambda_i^{(2/p)+1/2}.
\end{aligned}
\end{equation}
{}From the expression \eqref{eq:ExpandCurvMiddleAnnulus} for $F_A$ over
the annulus $\Omega(x_i;2\sqrt{\lambda_i},4\sqrt{\lambda_i})$ 
--- where $\fS(A,\Phi)=(F_A^{+,g},0)$ --- and the $L^p$
estimates \eqref{eq:LpEstimateBackroundConnOneFormBall} for
$\sigma_i^*A_0$ and that of Lemma \ref{lem:dBetaEst} for $d\chi_i$, we
obtain an estimate for $F_A^{+,g}$,
%\begin{equation}
%\begin{aligned}
\begin{align}
\notag
&\|\fS(A,\Phi)
\|_{L^{\sharp,2;2}(\Omega(x_i;2\sqrt{\lambda_i},4\sqrt{\lambda_i}))}
\\
\notag
&=
\|F_A^{+,g}\|_{L^{\sharp,2}(\Omega(x_i;2\sqrt{\lambda_i},4\sqrt{\lambda_i}))} 
\\
\notag
&\leq 
c\|F_A^{+,g}\|_{L^p(\Omega(x_i;2\sqrt{\lambda_i},4\sqrt{\lambda_i}))}
\quad\text{(by \cite[Lemma 4.1]{FeehanSlice}, for $2<p\leq\8$)} 
\\
\label{eq:LpEstimateMonoPairMiddleAnnulus}
&\leq 
c\left(\|F_{A_0}^{+,g}\|_{L^p(B(x_i,4\sqrt{\lambda_i}))}
+ \|d\chi_i\|_{L^4}
\|\si_i^*A_0\|_{L^q(\Omega(x_i;2\sqrt{\lambda_i},4\sqrt{\lambda_i}))}
\right.
\\
\notag
&\quad
\left.
+\|\si_i^*A_0\|_{L^{2p}(\Omega(x_i;2\sqrt{\lambda_i},4\sqrt{\lambda_i}))}^2
\right)
\\
\notag
&\leq C\left(\lambda_i^{2/p} + \lambda_i^{(2/q)+1/2} 
+ \lambda_i^{2(p+1/2)}\right)
\leq 
C\lambda_i^{2/p}.
\end{align}
%\end{aligned}
%\end{equation}
This completes the required estimates for $\fS(A,\Phi)$ over the annulus
$\Omega(x_i;2\sqrt{\lambda_i},4\sqrt{\lambda_i})$.

Next we turn to the estimate of $\fS(A,\Phi)=(F_A^{+,g},0)$ over the inner
annulus $\Omega(x_i;\sqrt{\lambda_i}/4,\sqrt{\lambda_i}/2)$. The hypotheses
imply that $F_{A_i}^{+,\delta}=0$ on $\RR^4$,
where $\delta$ is the Euclidean metric on $\RR^4$. Since the metric $g$ is
approximately Euclidean in normal coordinates, $\varphi_i^*g \approx
\delta$, we shall assume $\varphi_i^*g = \delta$ without loss and regard
the annulus $\Omega(x_i;\sqrt{\lambda_i}/4,\sqrt{\lambda_i}/2) \subset X$
as identified isometrically with the annulus
$\Omega(0;\sqrt{\lambda_i}/4,\sqrt{\lambda_i}/2) \subset \RR^4$, even
though we shall {\em not\/} assume $F_{A_i}^{+,g}\equiv 0$. This time, our
expression \eqref{eq:SplicedPair} for the connection $A$ yields
\begin{equation}
\label{eq:ExpandCurvInnerAnnulus}
F_A = \chi_iF(\tau_i^*A_i) +
 d\chi_i\wedge \tau_i^*A_i 
+ (\chi_i^2-\chi_i)\tau_i^*A_i 
\wedge \tau_i^*A_i,
\end{equation}
where we now abbreviate $\chi_i = 1-\chi_{x_i,\sqrt{\lambda_i}/2}$ and
$\tau_i^*A_i = (\varphi_n\circ\varphi_i^{-1})^*\tau_i^*A_i$.

Recall from \cite[Equation (7.12)]{TauSelfDual} or \cite{UhlRem} that, as
$\tau_i^*A_i$ is in radial gauge with respect to the south pole of $S^4$,
then the local connection one-form $\varphi_n^*\tau_i^*A_i$ obeys
\cite[p. 146]{FU} $(\varphi_n^*\tau_i^*A_i)_r
\equiv \partial_r\lrcorner\varphi_n^*\tau_i^*A_i = 0$ and
$\varphi_n^*\tau_i^*A_i(\8,\theta) = 0$, where
$r=|x|=\varphi_n^{-1}(\cdot)$, so 
$$
(\varphi_n^*F_{A_i})_{r\theta}
=
\frac{\partial}{\partial r}(\varphi_n^*\tau_i^*A_i)_\theta,
$$
where $(\varphi_n^*\tau_i^*A_i)_\theta
= \partial_\theta\lrcorner\varphi_n^*\tau_i^*A_i$
and thus
$$
\varphi_n^*\tau_i^*A_i(r,\theta)
=
-\int_r^\8 (\varphi_n^*F_{A_i})_{t\theta} \,dt.
$$
Hence, Theorem \ref{thm:Decay} (with $r_1=\8$) yields a bound for
$|F_{A_i}|_\delta$ and $|\varphi_n^*\tau_i^*A_i|_\delta$ 
on $\Omega(0;2r_0,\8) = \RR^4-B(2r_0)$,
\begin{equation}
\label{eq:InvariantPtEstSphereConnOneFormBallComp}
\begin{aligned}
|\varphi_n^*F_{A_i}|_\delta(x)
&\leq
c\frac{r_0^2}{r^4}\|F_{A_i}\|_{L^2(\Omega(0;r_0,\8),\delta)},
\\
|\varphi_n^*\tau_i^*A_i|_\delta(x)
&\leq
c\frac{r_0^2}{r^3}\|F_{A_i}\|_{L^2(\Omega(0;r_0,\8),\delta)},
\end{aligned}
\end{equation}
and therefore, on the annulus 
$\Omega(0;\sqrt{\lambda_i}/4,\sqrt{\lambda_i}/2)$, we have --- setting $r_0 =
L\lambda_i$ and using $\|F_{A_i}\|_{L^2(\Omega(0;r_0,\8),\delta)} < \eps \ll
1$ when $L\gg 1$ (but otherwise fixed),
\begin{equation}
\label{eq:PointEstimateSphereConnOneFormBallComp}
\begin{aligned}
|\varphi_n^*F_{A_i}|_\delta(x)
&\leq
c\frac{\lambda_i^2}{r^4},
\\
|\varphi_n^*\tau_i^*A_i|_\delta(x)
&\leq
c\frac{\lambda_i^2}{r^3},
\quad 
x\in \Omega(0;2L\lambda_i,\8).
\end{aligned}
\end{equation}
(Compare \cite[Proposition 7.7]{TauSelfDual} and \cite[Lemmas 9.1 \&
9.2]{TauFrame}; the above decay estimates for the connection $A_i$ with
scale $\lambda_i$ matches that of \cite[Equations (7.16) \&
(7.17)]{TauSelfDual}, obtained by a slightly different method.)  We can now
integrate the preceding estimates 
\eqref{eq:PointEstimateSphereConnOneFormBallComp}
to get the following $L^p$ bound, for any
$1\leq p\leq \8$,
\begin{equation}
\label{eq:LpEstimateSphereConnOneFormBallComp}
\begin{aligned}
\|F_{A_i}\|_{L^p(\Omega(0;\sqrt{\lambda_i}/4,\sqrt{\lambda_i}/2),\delta)}
&\leq
c\lambda_i^{2/p},
\\
\|\tau_i^*A_i\|_{L^p(\Omega(0;\sqrt{\lambda_i}/4,\sqrt{\lambda_i}/2),\delta)}
&\leq
c\lambda_i^{(2/p)+1/2},
\end{aligned}
\end{equation}
{}From the expression \eqref{eq:ExpandCurvInnerAnnulus} for $F_A$ over the
inner annulus $\Omega(x_i;\sqrt{\lambda_i}/4,\sqrt{\lambda_i}/2)$ and the
$L^p$ estimates \eqref{eq:LpEstimateSphereConnOneFormBallComp} for
$F_{A_i}$, $\tau_i^*A_i$ and that of Lemma \ref{lem:dBetaEst} for
$d\chi_i$, we obtain
\begin{equation}
\begin{aligned}
\label{eq:LpEstimateMonoPairInnerAnnulus}
&\|\fS(A,\Phi)
\|_{L^{\sharp,2;2}(\Omega(x_i;\sqrt{\lambda_i}/4,\sqrt{\lambda_i}/2),g)} 
\\
&\leq
c\|\fS(A,\Phi)
\|_{L^p(\Omega(x_i;\sqrt{\lambda_i}/4,\sqrt{\lambda_i}/2),g)} 
\quad\text{(by \cite[Lemma 4.1]{FeehanSlice}, for $2<p\leq\8$)}
\\
&=
c\|F_A^{+,g}
\|_{L^p(\Omega(x_i;\sqrt{\lambda_i}/4,\sqrt{\lambda_i}/2),g)} 
\\
&\leq 
c\|F^{+,g}_{A_i}
\|_{L^p(\Omega(0;\sqrt{\lambda_i}/4,\sqrt{\lambda_i}/2),\delta)} 
+ c\|d\chi_i\|_{L^4}
\|\tau_i^*A_i
\|_{L^q(\Omega(0;\sqrt{\lambda_i}/4,\sqrt{\lambda_i}/2),\delta)}
\\
&\quad
+ c\|\tau_i^*A_i
\|_{L^{2p}(\Omega(0;\sqrt{\lambda_i}/4,\sqrt{\lambda_i}/2),\delta)}^2
\\
&\leq 
C\left(\lambda_i^{2/p} + \lambda_i^{(2/q) + 1/2} 
+ \lambda_i^{2(p+(1/2))}\right)
\leq 
C\lambda_i^{2/p},
\end{aligned}
\end{equation}
noting that $1/q = 1/p-1/4$.
This completes the required estimate for $F_A^{+,g}$ over the inner
annulus.

Finally, we turn to the balls $B(x_i,\sqrt{\lambda_i}/4)$. We use {\em
normal, but not necessarily normal geodesic\/} coordinates 
(see Remark \ref{rmk:NormalCoordChoice})
centered at $x_i$ to attach the connections $A_i$,
defined on $\RR^4 = (TX)_{x_i}$ with its standard metric $\delta$.
Just as in \cite[Equation (8.20)]{TauSelfDual}, using
$F_{A_i}^{+,\delta}=0$ and
$$
F_{A_i}^{+,g} = \thalf(1+*_g)F_{A_i} = \thalf(*_g-*_{\delta})F_{A_i}
$$
gives, via $(\varphi_i^*g)_{\mu\nu} \equiv g_{\mu\nu} = \delta_{\mu\nu} +
O(r)$ for $x$ near $x_i\in X$, 
\begin{align*}
\|F_{A_i}^{+,g}\|_{L^{\sharp,2}(B(x_i,\sqrt{\lambda_i}/4))}
&\leq 
c\|*_g-*_{\delta}\|\cdot\|F_{A_i}\|_{L^{\sharp,2}(B(x_i,\sqrt{\lambda_i}/4))}
\\
&\leq 
c\lambda_i^{1/2}\|F_{A_i}\|_{L^{\sharp,2}(B(x_i,\sqrt{\lambda_i}/4))}.
\end{align*}
According to Lemma \ref{lem:UniversalLSharp2BoundFA}, we have for the
standard round metric $g_0$ of radius one on $S^4$, 
\begin{equation}
\label{eq:ProofUniversalLSharp2BoundFAonS4}
\|F_{A_i}\|_{L^\sharp(S^4,g_0)}
\leq 
c\|F_{A_i}\|_{L^2(S^4,g_0)},
\end{equation}
for $c=c(\kappa_i)$ and
all anti-self-dual connections $A_i$ over $S^4$ with $c_2(A_i)=\kappa_i$.
Therefore, combining the past two estimates, yields
\begin{equation}
\label{eq:LSharp2EstimateMonoPairBall}
\|\fS(A,\Phi)\|_{L^{\sharp,2;2}(B(x_i,\sqrt{\lambda_i}/4))}
=
\|F_A^{+,g}\|_{L^{\sharp,2}(B(x_i,\sqrt{\lambda_i}/4))}
\leq
C\lambda_i^{1/2}.
\end{equation}
The required bounds for $\fS(A,\Phi)$ now follow by combining the estimates
\eqref{eq:LpEstimateMonoPairOuterAnnulus},
\eqref{eq:LpEstimateMonoPairMiddleAnnulus}, and
\eqref{eq:LpEstimateMonoPairInnerAnnulus} for $\fS(A,\Phi)$ over the annuli
where cutting off occurs, together with the estimates 
\eqref{eq:LpEstimateMonoPairBallComp} for
$\fS(A,\Phi)$ on the complement of the balls in $X$ and the estimates
\eqref{eq:LSharp2EstimateMonoPairBall} over the small balls, and fixing a
value of $2<p\leq \8$. We take $p$ as close to $2$ as possible (in order to
get the best decay rate in terms of powers of $\lambda$), so a
convenient value is $p=8/3$.
\end{pf}

\begin{rmk}
\label{rmk:NormalCoordChoice}
The reasons for allowing any normal coordinates, rather than restricting to
normal geodesic coordinates which would give sharper estimates, is
explained in \cite{FLConj}.
\end{rmk}

\begin{rmk}
\label{rmk:ApproxPairBoundsForExtPUBackground}
Suppose $B'\Subset B\subset X$ are concentric balls.
If the precompact family $\sU\subset\sC(\ft_\ell)$ is a space of solutions
to an extended $\PU(2)$ monopole equation defined by an eigenvalue cutoff
constant $\mu>0$ (see \S \ref{sec:Existence}), 
$$
\sC(A_0,\Phi_0)=\Pi_{A_0,\Phi_0,\mu}\sC(A_0,\Phi_0),
$$
then Corollary \ref{cor:L2_1InhomoRegLocal} implies that --- when
$\|F_{A_0}\|_{L^2(B)} + \|\Phi_0\|_{L^4(B)} < \eps$ --- we have the bound
$$
\|(A_0-\Gamma,\Phi_0)\|_{L^2_{k,\Gamma}(B')} 
\leq 
\left(\|F_{A_0}\|_{L^2(B)} + \|\Phi_0\|_{L^2(B)}^2
+\|\Phi_0\|_{L^2_{1,A_0}(B)}\right),
$$
where $c=c(B',B,k)$ and $d_\Gamma^*(A_0-\Gamma)=0$,
since standard linear elliptic theory gives
$$
\|\Pi_{A_0,\Phi_0,\mu}\sC(A_0,\Phi_0)\|_{L^2_{k,\Gamma}(B')} 
\leq
c\|\sC(A_0,\Phi_0)\|_{L^2(B)}.
$$
The condition $\|F_{A_0}\|_{L^2(B)}<\eps$, for small enough $\eps(B)$,
ensures that Theorem \ref{thm:CoulombBallGauge} yields
$\|A_0-\Gamma\|_{L^2_{1,\Gamma}(B)} \leq c\|F_{A_0}\|_{L^2(B)}$, with
$A_0-\Gamma$ in Coulomb gauge.
\end{rmk}

\begin{rmk}
\label{rmk:ProjPerpfSSpliced}
As we shall see in \S \ref{sec:Existence}, it would be preferable to
estimate the quantity $\Pi_{A,\Phi,\mu}^\perp\fS(A,\Phi)$ directly for the
purposes of solving the extended $\PU(2)$ monopole equations, as
$\Pi_{A,\Phi,\mu}^\perp\fS(A,\Phi)$ will be smaller than
$\fS(A,\Phi)$, because we only assume that the background pairs
$(A_0,\Phi_0)$ obey $\Pi_{A_0,\Phi_0,\mu}^\perp\fS(A_0,\Phi_0)=0$. 
However, as explained in \S \ref{sec:Existence}, our
estimate for $\fS(A,\Phi)$ will suffice.
\end{rmk}

%end of file
%file: global.tex

\section{Global, uniform elliptic estimates for the anti-self-dual and
Dirac operators}
\label{sec:Global}
We recall some estimates \cite{FeehanSlice}, \cite{FeehanKato} for the
operators $d_A^+$ and $D_A$, with the property that all constants depend on
the connection $A$ at most through the $L^2$ norm of the curvature
$F_A$. These estimates are used in \S \ref{sec:Eigenvalue} to allow us to
identify the small-eigenvalue eigenspaces of the Laplacian
$d_{A,\Phi}^1d_{A,\Phi}^{1,*}$ associated to the linearization
$d_{A,\Phi}^1$ of the $\PU(2)$ monopole equations at a pair $(A,\Phi)$ and
in \S \ref{sec:Existence} to give a uniform bound for a partial right
inverse of this linearization.

Ultimately the desired global, uniform elliptic estimates rely on the
following Bochner-Weitzenb\"ock formulas, which we collect here for
convenience: 
\begin{align}
d_Ad_A^* + 2d_A^{+,*}d_A^+ 
&= 
\cov_A^*\cov_A + \{\Ric,\cdot\}
- 2\{F_A^-,\cdot\}, \label{eq:BW1} 
\\
2d_A^+d_A^{+,*}
&= 
\cov_A^*\cov_A - 2\{\sW^+,\cdot\} + \frac{R}{3} +
\{F_A^+,\cdot\}, \label{eq:BW+} 
\\
D_A^*D_A 
&= 
\cov_A^*\cov_A + \frac{R}{4}
+ \rho(F_A^+) + \frac{1}{2}\rho(F^+_{A_d}+F^+_{A_e}), \label{eq:BWDirac+} 
\\
D_AD_A^* 
&= 
\cov_A^*\cov_A + \frac{R}{4}
+ \rho(F_A^-) + \frac{1}{2}\rho(F^+_{A_d}+F^+_{A_e}), \label{eq:BWDirac-} 
\end{align}
for the Laplacians on $\Gamma(\Lambda^1\otimes\fg_E)$,
$\Gamma(\Lambda^+\otimes\fg_E)$, $\Gamma(V^+)$, and
$\Gamma(V^-)$, respectively; see \cite[p. 94]{FU} for the first
two formulas and \cite[Lemma 4.1]{FL1} for the third and fourth. Here, $A$
is viewed as an $\SO(3)$ connection on $\fg_E$ (see the conventions
described in \S \ref{sec:Prelim} and $A_d$, $A_e$ are the
$\U(1)$ connections on $\det W^+$ and $\det E$, respectively.

As we shall see in \S \ref{sec:Eigenvalue} and \S \ref{sec:Existence},
the presence of the $F_A^-$ term in the Bochner formula \eqref{eq:BWDirac-}
for $D_AD_A^*$, unlike the $F_A^+$ term in the Bochner formula
\eqref{eq:BW+}, means that we cannot simply imitate the strategy of
\cite{TauStable} in our construction of gluing solutions of the $\PU(2)$
monopole equations.

\subsection{Elliptic estimates for the anti-self-dual operator}
\label{subsec:EllipticASD}
In this subsection we recall the elliptic estimates for the {\em
anti-self-dual operator\/} $d_A^+$ described in \cite{FeehanSlice}.  The
following Kato-Sobolev inequality is well-known and will be frequently
employed; see Lemma 4.6 in \cite{TauSelfDual} for the case $p=2$.

\begin{lem}
\label{lem:Kato}
Let $X$ be a $C^\8$, closed, oriented, Riemannian, four-manifold 
and let $1\le p<4$. Then there is a positive constant $c$ with the
following significance.  Let $A$ be an orthogonal $L^p_1$ connection
on a Riemannian vector bundle $V$ over $X$. Then, for any $a\in
L^p_1(V)$ and $4/3\leq q<\8$ defined by $1/p=1/q+1/4$, we have
\begin{equation}
\label{eq:Kato}
\|a\|_{L^q(X)} \le c\|a\|_{L^p_{1,A}(X)}.
\end{equation}
\end{lem}

The estimates for $d_A^+$ we require are treated in detail in
\cite{FeehanSlice}, based on earlier work of Taubes, so we just summarize
the results here; the $L^\sharp$ and related families of Sobolev-Taubes
norms are defined in \cite[\S 4]{FeehanSlice}.

\begin{lem}
\label{lem:LinftyL22CovLapEstv}
\cite[Lemma 5.9]{FeehanSlice}
Let $X$ be a $C^\8$, closed, oriented, Riemannian four-manifold. Then
there are positive constants $c$ and $\eps=\eps(c)$ with the following
significance.  Let $E$ be a Hermitian, rank-two vector bundle over $X$ and
let $A$ be an orthogonal $L^2_4$ connection on $\fg_E$ with curvature
$F_A$, such that $\|F_A\|_{L^{\sharp,2}(X)} <\eps$.  Then the following
estimate holds for any $v\in L^{\sharp,2}_2(\Lambda^+\otimes\fg_E)$:
\begin{equation}
\label{eq:LinftyL22CovLapEstv}
\|v\|_{L^2_{2,A}(X)} + \|v\|_{C^0(X)}
\le c(1+\|F_A\|_{L^2(X)})
(\|d_A^+d_A^{+,*}v\|_{L^{\sharp,2}(X)} + \|v\|_{L^2(X)}).
\end{equation}
\end{lem}

The preceding lemma is proved via a (subtle) integration-by-parts argument
(due to Taubes), together with the Bochner formula \eqref{eq:BW+} to
replace the covariant Laplacian $\cov_A^*\cov_A$ by $d_A^+d_A^{+,*}$.
In \S \ref{sec:Existence} we shall use the following consequence of Lemma
\ref{lem:LinftyL22CovLapEstv}:

\begin{cor}
\label{cor:L21AEstdA*v}
\cite[Corollary 5.10]{FeehanSlice}
Continue the hypotheses of Lemma \ref{lem:LinftyL22CovLapEstv}. Then: 
\begin{equation}
\label{eq:L21AEstdA*v}
\|d_A^{+,*}v\|_{L^2_{1,A}(X)}
\le c(1+\|F_A\|_{L^2(X)})
(\|d_A^+d_A^{+,*}v\|_{L^{\sharp,2}(X)} + \|v\|_{L^2(X)}).
\end{equation}
\end{cor}

Note that any $a\in \Gamma(\Lambda^1\otimes\fg_E)$ which is
$L^2$-orthogonal to $\Ker d_A^+$ is given by $a = d_A^{+,*}v$, for some $v
\in \Gamma(\Lambda^+\otimes\fg_E)$.  We shall also use the following
$L^2_{1,A}$ estimates for $v \in
\Gamma(\Lambda^+\otimes\fg_E)$; similar estimates are given in \cite[Lemma
5.2]{TauSelfDual} and in 
\cite[Appendix A]{TauIndef}; they follow from an
elementary application of the Bochner formula \eqref{eq:BW+}.

\begin{lem}
\label{lem:L21AEstv}
\cite[Lemma 6.6]{FLKM1}
Let $X$ be a $C^\8$, closed, oriented, Riemannian four-manifold.  Then
there are positive constants $c$ and $\eps=\eps(c)$
with the following significance. Let $A$
be an $L^2_4$ orthogonal connection on an $\SO(3)$ bundle $\fg_E$ over
$X$ such that $\|F_A^+\|_{L^2}<\eps$.
Then, for all $v\in L^{4/3}_2(\Lambda^+\otimes\fg_E)$,
\begin{align}
\label{eq:L21AEstv}
\|v\|_{L^2_{1,A}(X)} 
&\le 
c(1+\|F_A^+\|_{L^2(X)})^{1/2}(\|d_A^{+,*}v\|_{L^2(X)}+\|v\|_{L^2(X)}),
\\
\label{eq:L2EstdA*v}
\|d_A^{+,*}v\|_{L^2(X)} 
&\leq 
(1+\|F_A^+\|_{L^2(X)})^{1/2}\left(\|d_A^+d_A^{+,*}v\|_{L^{4/3}(X)} 
+ \|v\|_{L^2(X)}\right)
\end{align}
\end{lem}

\subsection{Elliptic estimates for the Dirac operator}
\label{subsec:EllipticDirac}
In this subsection we recall the elliptic estimates for the Dirac
operator $D_A$ proved in \cite{FeehanSlice}. An application of Lemma 5.5 in
\cite{FeehanSlice} and the Bochner formula \eqref{eq:BW+} yields

\begin{lem}
\label{lem:LinftyL22Estphi}
\cite[Lemma 7.2]{FeehanKato}
Let $X$ be a closed, oriented four-manifold with metric $g$. Let $0\leq
M<\8$ be a constant. Then there are positive constants $c(g,M)$ and
$\eps(g)$ with the following significance. Let $(\rho,W^+,W^-)$ be a
\spinc structure over $X$ with \spinc connection determined by the
Levi-Civita connection on $T^*X$ and a complex $L^2_4$ connection $A_d$ on
$\det W^+$ with $\|F_{A_d}\|_{L^\8(X)}\leq M$.  Let $E$ be a Hermitian,
rank-two vector bundle over $X$ with unitary $L^2_4$ connection $A_e$ on
$\det E$ with $\|F_{A_e}\|_{L^\8(X)}\leq M$. Denote $V=W\otimes E$ and
$V^\pm = W^\pm\otimes E$. Suppose $A$ is an orthogonal $L^2_4$ connection
on $\fg_E$ for which $X=U_2^+\cup U_\8^+$, with $U_2^+, U_\8^+
\Subset X$ open subsets such that $\|F_A^+\|_{L^{\sharp,2}(U_2^+)} <\eps\ll
1$ and $\|F_A^+\|_{C^0(U_\8^+)}
\leq M <\8$. If $\phi^+\in L^2_3(X,V^+)$, then
\begin{equation}
\label{eq:LinftyL22Estphi}
\begin{aligned}
\|\phi^+\|_{C^0\cap L^2_{2,A}(X)}
&\leq
c(1+\|F_A\|_{L^2(X)})^2(1+\|F_A^+\|_{C^0(U_\8^+)})^2
\\
&\quad\times (\|D_A^- D_A^+\phi^+\|_{L^{\sharp,2}(X)} + \|\phi^+\|_{L^2(X)}).
\end{aligned}
\end{equation}
The analogous estimate holds for $\phi^-\in L^2_3(X,V^-)$, with $F_A^+$
replaced by $F_A^-$.
\end{lem}

Since $\|D_A^+\phi^+\|_{L^2_{1,A}}\le c\|\phi^+\|_{L^2_{2,A}}$, 
Lemma \ref{lem:LinftyL22Estphi} yields an $L^2_{1,A}$ estimate for
$D_A^+\phi$:

\begin{cor}
\label{cor:L21AEstD_Apmphi}
\cite[Corollary 7.3]{FeehanKato}
Continue the hypotheses of Lemma \ref{lem:LinftyL22Estphi}. Then: 
\begin{equation}
\label{eq:L21EstDAphi}
\begin{aligned}
\|D_A^+\phi^+\|_{L^2_{1,A}(X)}
&\le
c(1+\|F_A\|_{L^2(X)})^2(1+\|F_A^+\|_{C^0(U_\8^+)})^2
\\
&\quad\times (\|D_A^- D_A^+\phi^+\|_{L^{\sharp,2}(X)} + \|\phi^+\|_{L^2(X)}).
\end{aligned}
\end{equation}
The analogous estimate holds for $D_A^-\phi^-\in L^2_2(X,V^+)$, with $F_A^+$
replaced by $F_A^-$.
\end{cor}

Now Lemma \ref{lem:LinftyL22Estphi} and Corollary \ref{cor:L21AEstD_Apmphi}
provide useful $C^0\cap L^2_{2,A}(X)$ elliptic estimates for positive spinors
$\phi^+\in\Gamma(X,V^+)$ and useful $L^2_{1,A}(X)$ elliptic estimates for
$D_A^+\phi^+\in\Gamma(X,V^-)$ even when $A$ bubbles, because we still
have uniform $C^0$ bounds for $F^+_A$ away from the bubble points (on
the set $U_\8$) and small $L^{\sharp,2}$ bounds around the bubble points
(on the set $U_2$). However, this is never the case for $F^-_A$ in such
applications --- because $F^-_A$ is neither $C^0$-bounded nor
$L^{\sharp,2}$-small around the bubble points --- so neither Lemma
\ref{lem:LinftyL22Estphi} nor Corollary \ref{cor:L21AEstD_Apmphi} provide
useful $C^0\cap L^2_{2,A}(X)$ elliptic estimates for $\phi^-\in\Gamma(X,V^-)$
or useful $L^2_{1,A}(X)$ elliptic estimates for
$D_A^-\phi^-\in\Gamma(X,V^+)$. 

To partly address this problem and obtain estimates strong enough to meet
the demands of \S \ref{sec:Eigenvalue} and \S \ref{sec:Existence}, we
recall that eigenspinors (both positive and negative) satisfy a useful decay
estimates \cite{FeehanKato} on the complement of balls where the curvature
is $L^2$-concentrated. These were used in \cite{FeehanKato} to derive
Sobolev estimates for spinors on the complement of such balls:

\begin{prop}
\label{prop:LinftyL22EstphiCompOfBubbles}
\cite[Proposition 7.4]{FeehanKato}
Continue the hypotheses of Lemma \ref{lem:LinftyL22Estphi}.
If $\phi \in L^2_2(X,V^-)$, write $\phi = \phi_0 + \phi_\perp$,
where $\phi_0 \in \Ker D_A^-$ and $\phi_\perp \in (\Ker D_A^-)^\perp$, so
that $\phi_\perp = D_A^+\varphi$. Then
\begin{equation}
\label{eq:LinftyL22EstphiPerpCompOfBubbles}
\begin{aligned}
\|\phi_\perp\|_{L^2_{1,A}(X)}
&\le
c(1+\|F_A\|_{L^2(X)})^2(1+\|F_A^+\|_{C^0(U_\8)})^2
\\
&\quad\times (\|D_A^-\phi\|_{L^{\sharp,2}(X)} + \|\varphi\|_{L^2(X)}).
\end{aligned}
\end{equation}
If $U = X-\cup_{i=1}^m B(x_i,\lambda_i^{1/2})$ and $U' = X-\cup_{i=1}^m
B(x_i,2\lambda_i^{1/3})\Subset U$ and $\max_i\lambda_i\leq 1$, then
\begin{equation}
\label{eq:LinftyL22EstEigenphiCompOfBubbles}
\begin{aligned}
\|\phi_0\|_{C^0\cap L^2_{2,A}(U')}
&\le
c(1+\|F_A\|_{L^2(X)})(1+\|F_A^-\|_{C^0(U)})
\|\phi_0\|_{L^2(U)}.
\end{aligned}
\end{equation}
\end{prop}

In contrast, a simple integration-by-parts argument yields the
corresponding bound for $\phi\in\Gamma(X,V^+)$:

\begin{lem}
\label{lem:L21AEstphi}
\cite[Lemma 7.5]{FeehanKato}
Continue the hypotheses of Lemma \ref{lem:LinftyL22Estphi}.
Then, for any $\phi\in L^2_1(X,V^+)$, we have
\begin{equation}
\label{eq:L21AEstphi}
\|\phi\|_{L^2_{1,A}(X)} 
\le c\left(1 + \|F_A^+\|_{C^0(U_\8)}\right)^{1/2}
(\|D_A^+\phi\|_{L^2(X)} + \|\phi\|_{L^2(X)}).
\end{equation}
\end{lem}

\begin{rmk}
Lemma 5.9 in \cite{FL1} yields the bound $\|F_A\|_{L^\sharp} \leq
c\|F_A\|_{L^2}$ if $A$ is an anti-self-dual connection. The proof yields
the same answer if $A$ satisfies another elliptic equation. In the present
application, $A$ will be obtained by splicing anti-self-dual connections over
$S^4$ onto a background connection over $X$ satisfying the extended
$\PU(2)$ monopole equations. Since the background family is always
precompact, we may assume a bound of the form $\|F_A\|_{L^\sharp} \leq
c\|F_A\|_{L^2}$ without loss.
\end{rmk}

{}From the Hodge theorem, we have an $L^2$-orthogonal splitting
$$
L^2_{k-1}(V^-)
=
\Ker D_A^*\oplus \Ran D_A,
$$
where $D_A^*: L^2_{k+1}(V^-) \to L^2_k(V^+)$ and $D_A:
L^2_k(V^+) \to L^2_{k-1}(V^-)$. It remains to estimate
the harmonic components of sections in $L^2_{k-1}(V^-)$. Because
of the presence of the $F_A^-$ term in the Bochner formula \eqref{eq:BWDirac-}
for $D_AD_A^*$, we cannot expect to obtain global $L^p(X)$ estimates for
$\psi \in \Ker D_A^* = \Ker D_AD_A^*$ when $p>2$. However, in the course of
deriving an estimate for the Laplacian $d_{A,\Phi}^1d_{A,\Phi}^{1,*}$ we
shall see that it is enough to derive an $L^2_{1,A}$ estimate for $\psi
\in \Ker D_A^*$ upon restriction to $\supp\Phi$. For this purpose, it is
enough to have decay estimates for such $\psi$ away from the points of
curvature concentration and in the support of $\Phi$.

%end of file
%file: dirichlet.tex

\section{The Dirichlet problem for the Dirac Laplacian}
\label{sec:Dirichlet}
Our goal in this section is to prove the following elliptic estimate. This
estimate plays a crucial role in the proof of Corollary
\ref{cor:L21AEstPAaphi}, which gives a uniform bound for a partial right
inverse of the linearization of the $\PU(2)$ monopole equations.

\begin{thm}
\label{thm:Lp2SpinorAnnulusEst}
Let $M$ be a positive constant. Let $0<4r_0<r_1<\8$ and consider an annulus
$\tilde\Omega = \Omega(r_0,r_1) = B(r_1)-\bar B(r_0) \subset \RR^4$ centered at
the origin, with $C^\8$ metric $g$ of bounded geometry and injectivity
radius large compared with $r_1$. Let $\Omega =
\Omega(t_0,t_1) \Subset \tilde\Omega$ and $\Omega' =
\Omega(2t_0,t_1/2) \Subset \Omega$, where $r_0 \leq t_0 < t_1 \leq r_1$ and
$r_0t_0^{-3}, r_1^{-2} \leq M$.  
Let $A$ be an $C^\8$ \spinc connection on a Hermitian bundle
$V|_{\tilde\Omega}$, where $V=V^+\oplus V^-$ and $V^\pm \cong
\tilde\Omega\times\CC^4$.  Suppose $F_A$ and the curvature $\Rm_g$ of the
Levi-Civita connection for $g$ on $T\RR^4$ obeys
\begin{gather}
\label{eq:BoundedCurvHypothesisMain}
\|\mathrm{Rm}_g\|_{L^\8(\Omega)} \leq M,
\\
\label{eq:L2CurvHypothesisMain}
\|F_A\|_{L^2(\tilde\Omega)}  < \eps
\quad\text{and}\quad
\|\mathrm{Rm}_g\|_{L^2(\Omega)} < \eps.
\end{gather}
Then there is a positive constant $\eps =
\eps(M,p,r_1)$ such that the following holds.
For any $1<p\leq 2$ there is a constant $C=C(M,p,r_1)$ such that for all
$\psi\in L^p_2(\Omega,V)$ we have
\begin{equation}
\label{eq:Lp2SpinorAnnulusEst}
\|\psi\|_{L^p_{2,A}(\Omega')}
\leq
C(\|D_AD_A^*\psi\|_{L^p(\Omega)} + \|\psi\|_{L^2_{1,A}(\Omega)}).
\end{equation}
\end{thm}

In our applications of Theorem \ref{thm:Lp2SpinorAnnulusEst} in 
\S \ref{sec:Existence}, the
constant $r_1$ will typically be fixed while the constant $r_0$ is allowed
to approach zero: the essential point in the estimate
\eqref{eq:Lp2SpinorAnnulusEst} is that the dependence of the constant $C$
on $\Omega$ is known and, in particular, remains bounded as $r_0\to 0$. The
annulus $\Omega(r_0,r_1)$ will be thought of as surrounding a ball where
the curvature of an $\SO(3)$ connection $A^o$ on $\fg_E$ is allowed to
bubble inside the ball but remains $L^\8$ bounded on $\Omega(r_0,r_1)$,
where $V=W\otimes E$, with \spinc connection on $V$ determined by $A^o$ on
$\fg_E$, $A_d$ on $\det W^+$, $A_e$ on $\det E$, and the Levi-Civita
connection on $T^*X$. Usually, the constant $r_0$ will be a constant
multiple of $\lambda$, the local scale of the
connection $A$. The slightly unusual conditions on the relative sizes of the
radii $r_0$, $r_1$ arises because of the application in the proof of our
decay estimates for eigenspinors \cite[Theorem 1.1]{FeehanKato}. The
dependence of the constant $C$ on $r_1$ arises through
\begin{itemize}
\item
The first Dirichlet eigenvalue of the scalar Laplacian on a ball of radius
$r_1$ in $\RR^4$ (Remark \ref{rmk:EigenvalueRescaling}), and
\item
The dependence of the norm $K=K(\Omega)$ of
the Sobolev embedding \cite[Theorem V.5.4]{Adams} 
\begin{equation}
\label{eq:L21toL4SobolevEmbedding}
L^2_1(\Omega)\subset L^4(\Omega)
\end{equation}
on the geometry of the largest interior cone defining the interior cone
property for $\Omega$: the hypothesis $r_0<r_1/4$ ensures that this
depends on $r_1$ but not $r_0$. The constant $r_1$ determines
the maximum height of the cone with respect to which the annulus
$\Omega(r_0,r_1)$ satisfies an interior cone condition --- see \cite[\S
IV.3 \& Theorem V.5.4]{Adams}.
\end{itemize}
The Sobolev constant $K$ also depends, in general, on the $C^0$
bounds for the curvature $\Rm_g$ of the Riemannian metric $g$ and its
injectivity radius, which we shall always assume to be much larger than $r_1$
\cite[\S II.7]{Aubin}. 

We can assume without loss that $A$ extends to a $C^\8$
connection (not necessarily Yang-Mills) on $\RR^4\cup\{\8\}$ with small
curvature and write 
\begin{equation}
\label{eq:AisProduct+a}
A = \Gamma + a, 
\end{equation}
with $\Gamma$ being the product
connection on $\RR^4\times\CC^4$. Rather than assume that the \spinc
connection $A$ is Yang-Mills, in our application to connections $A$
produced by the splicing construction of
\S \ref{sec:Splicing}, it suffices to assume that the $\SO(3)$ connection
$A^o$ is Yang-Mills, while the fixed $\U(1)$ connections $A_d$ and $A_e$ obey
a curvature constraint of the same shape as
\eqref{eq:L2CurvHypothesisMain}. (Of course, we can always assume without
loss that $A_d$ and $A_e$ are Yang-Mills.) If $a=A-\Gamma$ is a
connection one-form on $(\RR^4\less\{0\})\cup\{\8\}$ in radial gauge with
respect to the point at infinity, whose extension we take to be Yang-Mills
for simplicity (as will be true in our application),
then Corollary \ref{cor:Decay} and the
proof of inequality \eqref{eq:InvariantPtEstSphereConnOneFormBallComp} imply
\begin{equation}
\label{eq:CovDerivRadialGaugeOneFormDecay}
|\cov_\Gamma^k a|_\delta(x)
\leq
c_k\frac{r_0^2}{r^{3+k}}\|F_{A}\|_{L^2(\Omega(r_0,\8),\delta)},
\quad 2r_0 < r < \8,
\end{equation}
provided only that $\|F_{A}\|_{L^2(\Omega(r_0,\8),\delta)} < \eps$, where
$\eps$ is a universal constant and $k\geq 0$ is any integer. 
We denote $r = |x| = \dist_\delta(0,x)$,
and $\delta$ is the standard metric on $\RR^4$. This result, aside from a
consideration of the precise constants involved, is Lemma 9.2 in
\cite{TauFrame}. In our application, we at most need the estimate
\eqref{eq:CovDerivRadialGaugeOneFormDecay} when $k\leq 1$, $A$ has
scale $\lambda$, and $r_0 = L\lambda$, with $L\gg 1$ a fixed, universal
constant. Inequality \eqref{eq:CovDerivRadialGaugeOneFormDecay} yields the
useful bound
\begin{equation}
\label{eq:C1EstConnectionOneFormMain}
\|a\|_{L^\8_{1,\Gamma}(\Omega)} \leq M,
\end{equation}
while the bound
\begin{equation}
\label{eq:BoundedFAHypothesisMain}
\|F_A\|_{L^\8(\Omega)} \leq M
\end{equation}
would follow directly from Theorem \ref{thm:Decay}, where $\Omega$
in inequalities \eqref{eq:C1EstConnectionOneFormMain} and
\eqref{eq:BoundedFAHypothesisMain}
can be taken to be $\Omega(2\sqrt{r_0},\8)$ if $r_0\leq 1$ and simply
$\Omega(2r_0,\8)$ if $r_0\geq 1$.

\subsection{The Dirichlet problem for scalar, second-order elliptic
equations}  
\label{subsec:DirichletScalar}
To begin, we recall some results from \cite{GT} concerning elliptic,
second-order differential equations for functions on domains in Euclidean
space. First, we shall need the following consequence of the Calderon-Zygmund
inequality \cite[Theorem 9.9 \& Corollary 9.10]{GT}, \cite{Stein}:

\begin{lem}
\cite[Lemma 9.17]{GT} 
\label{lem:ScalarPDEHomogeneousBound}
Let $\Omega\subset\RR^n$ be a domain whose boundary
$\partial\Omega$ is a $C^{1,1}$ submanifold. Let
\begin{equation}
\sL = a^{\mu\nu}(x)\frac{\partial^2}{\partial x_\mu\partial x_\nu}
+ b^\mu(x)\frac{\partial}{\partial x_\mu} + c(x)
\end{equation}
be a {\em strictly elliptic\/} second-order differential operator on
$\Omega$, so $a^{\mu\nu}\xi_\mu\xi_\nu\geq \lambda|\xi|^2$ for all
$\xi\in\RR^n$ and $x\in\Omega$, for some positive constant
$\lambda$. Require that the coefficients obey $a^{\mu\nu}\in
C^0(\bar\Omega)$, $b^\mu, c\in L^\8(\Omega)$, and $c\leq 0$, with
$|a^{\mu\nu}|, |b^\mu|, |c| \leq \Lambda$, for some positive constant
$\Lambda$. Let $K=\max\{\|a^{\mu\nu}\|_{C^{0,1}(\bar\Omega)},
\|b^\mu\|_{C^{0,1}(\bar\Omega)}, \|c\|_{L^{\8}(\Omega)}\}$. Then there is a
constant $C=C(K,n,p,\lambda,\Lambda,\Omega)$ such that for all $u\in
L^p_2(\Omega)\cap L^p_1(\Omega;\partial\Omega)$, $1<p<\8$,
\begin{equation}
\label{eq:ScalarPDEHomogeneousBound}
\|u\|_{L^p_2(\Omega)} \leq C\|\sL u\|_{L^p(\Omega)}. 
\end{equation}
\end{lem}

See also Theorem 9.14 in \cite{GT} for a similar statement, under slightly
different hypotheses. We shall also need results concerning the existence
and uniqueness of solutions to the Dirichlet boundary-value problem:

\begin{thm}
\label{thm:ScalarPDEExistenceUniqueness}
\cite[Theorem 9.15 \& Corollary 9.18]{GT}
Continue the hypotheses of Lemma \ref{eq:ScalarPDEHomogeneousBound}. If
$f\in L^p(\Omega)$ and $\varphi\in L^p_2(\Omega)$, with $1<p<\8$, then the
Dirichlet problem
\begin{equation}
\label{eq:SobolevDirichletProblem}
\sL u = f\quad\text{on }\Omega,\quad
u - \varphi \in L^p_1(\Omega;\partial\Omega)
\end{equation}
has a unique solution $u\in L^p_2(\Omega)$. If $p>n/2$ and $\varphi \in
C^0(\partial\Omega)$, then the Dirichlet problem
\begin{equation}
\label{eq:DirichletProblem}
\sL u = f\quad\text{on }\Omega,\quad
u = \varphi \in \partial\Omega
\end{equation}
has a unique solution $u\in L^p_{2,\loc}(\Omega)\cap C^0(\bar\Omega)$.  
\end{thm}

Lastly, we recall the following regularity results for solutions to the
Dirichlet problem:

\begin{thm}
\label{thm:ScalarPDERegularity}
\cite[Theorem 9.19]{GT}
Continue the hypotheses of Lemma
\ref{eq:ScalarPDEHomogeneousBound}. Suppose the coefficients of $\sL$
belong to $C^{k-1,1}(\Omega)$ (respectively, $C^{k-1,\alpha}(\Omega)$), and
$f\in L^q_{k,\loc}(\Omega)$ (respectively, $C^{k-1,\alpha}(\Omega)$), with
$1<p,q<\8$ and $0<\alpha<1$. If $u\in L^p_{2,\loc}(\Omega)$ is a solution
to $\sL u = f$, then $u\in L^q_{k+2,\loc}(\Omega)$ (respectively,
$C^{k+1,\alpha}(\Omega)$).

Furthermore, if $\Omega$ is a $C^{k+1,1}$ domain, $\sL$ has
coefficients in $C^{k-1,1}(\bar\Omega)$ (respectively,
$C^{k-1,\alpha}(\bar\Omega)$), and $f\in L^q_k(\Omega)$ (respectively,
$C^{k-1,\alpha}(\bar\Omega)$), then $u\in L^q_{k+2}(\Omega)$
(respectively, $C^{k+1,\alpha}(\bar\Omega)$).
\end{thm}

\subsection{The Dirichlet problem for harmonic spinors on an annulus}
\label{subsec:DiracDirichletAnnulus}
The main purpose of this subsection is to prove Theorem
\ref{thm:DirichletDiracLaplacian}, the analogue for the Dirac Laplacian of
the existence, uniqueness, and regularity results for the scalar,
second-order elliptic equation in the preceding subsection.

We outline the strategy in a series of lemmas. We first have the following
analogues of Theorems 8.3 and 8.13 and Corollary 8.7 in \cite{GT}, which
solve the Dirichlet problem for second-order, linear elliptic operators on
functions over domains in Euclidean space. While we could
appeal to standard results for elliptic boundary-value problems for
systems, such as those of
\cite{AgmonDouglisNirenberg1},
\cite{AgmonDouglisNirenberg2}, \cite[\S 20.2]{Hormander}, and
\cite[\S\S6.2--6.5]{Morrey}, it is essential for our application in the
remainder of our article that we know the precise dependence, if any, of
the constants appearing in estimates on the geometry of the domain $\Omega$
and the curvature $F_A$ of the connection $A$. Hence, we instead proceed by
modifying the arguments of \cite[\S\S 8.1--8.4 \& 9.5, 9.6]{GT} for the
Dirichlet boundary-value problem for the second-order, linear, elliptic
operator on functions. 

Note that because we are only considering the Dirichlet problem in Theorem
\ref{thm:DirichletDiracLaplacian} for covariant Laplacians or Laplacians
plus zeroth-order terms, we do not encounter the obstructions to
prescribing local boundary conditions for operators of Dirac type
discussed in \cite{APS1}.

We write $L_k^p(\Omega,\partial\Omega)$ for the
closure in $L_k^p(\Omega)$ of the subspace $C^\8_0(\Omega)$, rather than
$W^{k,p}_0(\Omega)$ and $W^{k,p}(\Omega)$, as is customary in \cite{Adams}
or \cite{GT}.

\begin{thm}
\label{thm:DirichletDiracLaplacian}
Let $0<4r_0<r_1<\8$ and let $\Omega=\Omega(r_0,r_1)$ be the open annulus
$B(r_1)-\bar B(r_0)$ in $\RR^4$ with $C^\8$ metric $g$ on $\RR^4$ of bounded
geometry, with injectivity radius $\geq r_1$ and
$\|\Rm_g\|_{L^\8(\Omega)}\leq M <\8$.
Let $A$ be a $C^\8$ \spinc connection on a Hermitian vector
bundle $V=V^+\oplus V^-$, where $V^\pm\cong\Omega\times\CC^4$ and let $D_A$
be the Dirac operator on $\Gamma(\Omega,V)$, with Bochner form
$\cov_A^*\cov_A +\sR_A$. Then there is a positive constant $\eps =
\eps(\min\{1,r_1\},M)$ such that if 
\begin{equation}
\label{eq:L2BochnerFormHypothesis}
\|\sR_A\|_{L^2(\Omega)} < \eps, 
\end{equation}
then the following holds.  Let $1<p<\8$ and $\psi\in L^p_2(\Omega,V)$.
\begin{enumerate}
\item
There is a unique solution $h \in L^p_2(\Omega,V)$ to the Dirichlet problem
\begin{equation}
\label{eq:DirichletDiracProblem}
D_A^2h = 0\quad\text{on $\Omega$ and}\quad h-\psi 
\in L^p_2(\Omega,\partial\Omega;V).
\end{equation}
If $p>2$ and $\psi\in C^0(\partial\Omega)$, then 
there is a unique solution $h \in C^0(\bar\Omega,V)\cap
L^p_{2,\loc}(\Omega,V)$ to the Dirichlet problem 
\begin{equation}
\label{eq:C0DirichletDiracProblem}
D_A^2h = 0\quad\text{on $\Omega$ and}
\quad h = \psi \quad\text{on $\partial\Omega$}.
\end{equation}
\item
The solution $h$ obeys the {\em a priori\/} estimate
\begin{equation}
\label{eq:L21DirichletDiracAPrioriEst}
\|h\|_{L^2_{1,A}(\Omega)} \leq c\|\psi\|_{L^2_{1,A}(\Omega)},
\end{equation}
where $c = c(r_1,M)$.
%\marginpar{ [CHECK DEP ON $\sR_A$]}
\end{enumerate}
The same results holds for the Dirichlet problem for the covariant
Laplacian, $\cov_A^*\cov_A$, as in that case one simply has $\sR_A=0$.
\end{thm}

Note that when $V=W\otimes E$, the exact form of $\sR_A$ is given by the
Bochner formulas \eqref{eq:BWDirac+}, \eqref{eq:BWDirac-}, with ``$A$'' in
those expressions being the $\SO(3)$ connection on $\fg_E$.

We first record a few elementary integration-by-parts formulas which we
shall need for the proof of Theorem \ref {thm:DirichletDiracLaplacian}.
Let $\Omega\subset X$ be an open subset with boundary
$\partial\Omega\subset X$ a compact, smooth submanifold.  {}From
\cite[Equation (II.8.1)]{LM} and the Divergence Theorem, we have
\begin{equation}
\label{eq:LaplacianIntByParts}
\int_\Omega\langle\cov_A^*\cov_A\varphi_1,\varphi_2\rangle\,dV
=
\int_\Omega\langle\cov_A\varphi_1,\cov_A\varphi_2\rangle\,dV
-\int_{\partial\Omega}\langle\cov_{A,\nu}\varphi_1,\varphi_2\rangle\,dS,
\end{equation}
where $\nu$ is the outward-pointing unit-normal vector field on $\partial
\Omega$. Similarly, as 
$$
d\langle\varphi_1,\varphi_2\rangle 
= 
\langle\cov_A\varphi_1,\varphi_2\rangle 
+ \langle\varphi_1,\cov_A\varphi_2\rangle,
$$
we have
\begin{equation}
\label{eq:IntByParts}
\int_\Omega\langle\cov_A\varphi_1,\varphi_2\rangle\,dV
=
-\int_\Omega\langle\varphi_1,\cov_A\varphi_2\rangle\,dV
+\int_{\partial\Omega}\langle\varphi_1,\varphi_2\rangle\,dS.
\end{equation}
We shall also need the corresponding integration-by-parts formulas for the
Dirac operator. Recall from \cite[Proposition II.5.3 \& Equation (II.5.7)]{LM}
that 
\begin{align}
\label{eq:DiracIntByParts}
\int_\Omega\langle D_A\psi_1,\psi_2\rangle\,dV
&=
\int_\Omega\langle \psi_1,D_A\psi_2\rangle\,dV
+ \int_{\partial\Omega}\langle \rho(\nu^*)\psi_1,\psi_2\rangle\,dS,
\\
\label{eq:DiracLaplacianIntByParts}
\int_\Omega\langle D_A^2\psi_1,\psi_2\rangle\,dV
&=
\int_\Omega\langle D_A\psi_1,D_A\psi_2\rangle\,dV
+ \int_{\partial\Omega}\langle \rho(\nu^*)D_A\psi_1,\psi_2\rangle\,dS,
\end{align}
where $\nu^* = g(\cdot,\nu)$ is the one-form on $\partial\Omega$ dual to
the outward-pointing unit-normal vector field $\nu$. We now turn to some
preparatory lemmas.

\begin{lem}
\label{lem:FirstScalarDirichletEigenvalueBall}
Let $\mu_1(\Omega;\Delta)$ denote the first eigenvalue of the Laplacian
$\Delta$ on $C^\8(\Omega,\partial\Omega;\CC)$, for a domain
$\Omega\subset\RR^n$ with the standard metric and Dirichlet boundary
condition on $\partial\Omega$. If $B$ is the unit ball and $\Omega_\eps$
the annulus $B(0,1)-\bar B(0,\eps)$, then $\mu_1(B;\Delta)$ is the first
positive zero of the Bessel function $J_{(n-2)/2}$ and
$\mu_1(\Omega_\eps;\Delta) > \mu_1(B;\Delta)$ for all $0<\eps<1$. In
particular, if $n=4$, then $\mu_1(B;\Delta) \approx 3.83171$, the first
positive zero of the Bessel function $J_1$. 
\end{lem}

\begin{proof}
The calculation of the first (and all higher) eigenvalues of the Dirichlet
Laplacian on $B\subset\RR^n$ can be found, for example, in
\cite[pp. 106--107]{Folland}. The lower bound on
$\mu_1(\Omega_\eps;\Delta)$ is a consequence of domain monotonicity of
eigenvalues for the Laplacian with Dirichlet boundary conditions
\cite[Corollary I.1]{Chavel}.
\end{proof}

Lemma \ref{lem:FirstScalarDirichletEigenvalueBall} leads to a lower bound
on the first eigenvalue of the Dirac Laplacian on the annulus
$\Omega_\eps$. Recall that the Bochner formula for the Dirac Laplacian
gives $D_A^2 = \cov_A^*\cov_A+\sR_A$ on $C^\8(\Omega,V)$ and if we restrict
to sections of $V^+$ or $V^-$, then $\sR_A$ is replaced by $\sR_A^+$ or
$\sR_A^-$.

\begin{lem}
\label{lem:FirstDiracDirichletEigenvalueBall}
Continue the notation and hypotheses of Lemma
\ref{lem:FirstScalarDirichletEigenvalueBall}. Let
$\mu_1(\Omega;\cov_A^*\cov_A)$, respectively $\mu_1(\Omega;D_A^2)$, denote
the first eigenvalue of the covariant Laplacian $\cov_A^*\cov_A$,
respectively the Dirac Laplacian $D_A^2$, on
$C^\8(\Omega,\partial\Omega;V)$ with Dirichlet boundary conditions on
$\partial\Omega$. Let $K=K(\Omega)$ denote the norm of the Sobolev
$L^2_1(\Omega)\to L^4(\Omega)$.  Then
\begin{equation}
\label{eq:FirstCovLapDirichletEigenvalueBall}
\mu_1(\Omega_\eps;\cov_A^*\cov_A) \geq \mu_1(\Omega_\eps;\Delta) 
\geq \mu_1(B;\Delta) > 0,
\end{equation}
and
\begin{equation}
\label{eq:FirstDiracDirichletEigenvalueBall}
\mu_1(\Omega_\eps;D_A^2) 
\geq 
(1-K\|\sR_A\|_{L^2(\Omega)}) \mu_1(B;\Delta)
- K\|\sR_A\|_{L^2(\Omega)},
\end{equation}
which is positive for small enough $\|\sR_A\|_{L^2(\Omega_\eps)}$.
If we restrict the Dirac Laplacian $D_A^2$ to sections of $V^+$ or $V^-$,
then $\sR_A$ is replaced by $\sR_A^+$ or $\sR_A^-$ above.
\end{lem}

\begin{rmk}
\label{rmk:EigenvalueRescaling}
The balls and annuli in Lemmas \ref{lem:FirstDiracDirichletEigenvalueBall}
and \ref{lem:FirstScalarDirichletEigenvalueBall} have outer radii equal to
$1$. If the radii are multiplied by a factor $r_1$, then the eigenvalues
are multiplied by a factor $r_1^{-2}$ \cite[\S XII.7]{Chavel}.
\end{rmk}

\begin{proof}
The first Dirichlet eigenvalue of the covariant Laplacian $\cov_A^*\cov_A$ on
$C^\8(\Omega,\partial\Omega;V)$ is given by
$$
\mu_1(\Omega;\cov_A^*\cov_A) 
=
\inf_{\substack{\|\psi\|_{L^2(\Omega)}=1,\\ \psi=0\text{ on }\partial\Omega}} 
\|\cov_A\psi\|_{L^2(\Omega)}^2,
$$
where the infimum is taken over smooth sections. The Kato inequality
$|\cov|\psi||\leq |\cov_A\psi|$ and the analogous characterization of
$\mu_1(\Omega;\Delta)$ imply that
$$
\|\cov_A\psi\|_{L^2(\Omega)}^2 
\geq
\|\cov|\psi|\|_{L^2(\Omega)}^2
\geq
\mu_1(\Omega;\Delta).
$$
The desired estimate now follows from Lemma
 \ref{lem:FirstScalarDirichletEigenvalueBall}.

Similarly, the first Dirichlet eigenvalue of the Dirac Laplacian $D_A^2$ on
$C^\8(\Omega,\partial\Omega;V)$ is given by
$$
\mu_1(\Omega;D_A^2) 
=
\inf_{\substack{\|\psi\|_{L^2(\Omega)}=1,\\ \psi=0\text{ on }\partial\Omega}} 
\|D_A\psi\|_{L^2(\Omega)}^2,
$$
where the infimum is taken over smooth sections. Applying the Bochner
formula and integration by parts, we see that
\begin{align*}
\|D_A\psi\|_{L^2(\Omega)}^2
&=
(D_A^2\psi,\psi)_{L^2(\Omega)}
\\
&=
(\cov_A^*\cov_A\psi,\psi)_{L^2(\Omega)} + (\sR_A\psi,\psi)_{L^2(\Omega)}
\\
&\geq
\|\cov_A\psi\|_{L^2(\Omega)}^2 
- \|\sR_A\|_{L^2(\Omega)}\|\psi\|_{L^4(\Omega)}^2
\\
&\geq
\|\cov_A\psi\|_{L^2(\Omega)}^2 
- c\|\sR_A\|_{L^2(\Omega)}(\|\psi\|_{L^2(\Omega)}^2
+ \|\cov_A\psi\|_{L^2(\Omega)}^2)
\\
&\geq
(1-c\|\sR_A\|_{L^2(\Omega)})\|\cov_A\psi\|_{L^2(\Omega)}^2 
- c\|\sR_A\|_{L^2(\Omega)}\|\psi\|_{L^2(\Omega)}^2
\\
&\geq
(1-c\|\sR_A\|_{L^2(\Omega)})\mu_1(\Omega;\cov_A^*\cov_A) 
- c\|\sR_A\|_{L^2(\Omega)}.
\end{align*}
The inequality \eqref{eq:FirstDiracDirichletEigenvalueBall} now follows
from the eigenvalue estimate \eqref{eq:FirstCovLapDirichletEigenvalueBall}. 
\end{proof}

We can now turn to our proof of the {\em a priori} $L^2_1$ estimate.

\begin{proof}[Proof of inequality \eqref{eq:L21DirichletDiracAPrioriEst}
in Theorem \ref{thm:DirichletDiracLaplacian}]
Since
\begin{equation}
\label{eq:L21Esth-1}
\|h\|_{L^2_{1,A}(\Omega)}
\leq 
\|h-\psi\|_{L^2_{1,A}(\Omega)} + \|\psi\|_{L^2_{1,A}(\Omega)},  
\end{equation}
it suffices to estimate the $L^2_{1,A}(\Omega)$ norm of $h-\psi$.  Observe
that
\begin{equation}
\label{eq:L21Esth-2}
\begin{aligned}
\|h-\psi\|_{L^2_{1,A}(\Omega)}
&\leq 
\|h-\psi\|_{L^2(\Omega)} + \|\cov_A(h-\psi)\|_{L^2(\Omega)} 
\\
&\leq
\mu_1^{-1/2}\|D_A(h-\psi)\|_{L^2(\Omega)} 
+ \|\cov_A(h-\psi)\|_{L^2(\Omega)}
\quad\text{(by Lemma \ref{lem:FirstDiracDirichletEigenvalueBall})}
\\
&\leq c\|\cov_A(h-\psi)\|_{L^2(\Omega)},
\end{aligned} 
\end{equation}  
where $\mu_1 = \mu_1(\Omega;D_A^2)$.
Using integration by parts \eqref{eq:LaplacianIntByParts} and the facts
that $h-\psi=0$ on $\partial\Omega$ and $D_A^2h = 0$ on $\Omega$,
we see that
\begin{align*}
\|\cov_A(h-\psi)\|_{L^2(\Omega)}^2
&=
(\cov_A^*\cov_A(h-\psi),h-\psi)_{L^2(\Omega)}
\\
&=
(D_A^2(h-\psi),h-\psi)_{L^2(\Omega)} + (\sR_A(h-\psi),h-\psi)_{L^2(\Omega)}
\\
&=
-(D_A^2\psi,h-\psi)_{L^2(\Omega)} + (\sR_A(h-\psi),h-\psi)_{L^2(\Omega)}
\\
&=
-(D_A\psi,D_A(h-\psi))_{L^2(\Omega)} + (\sR_A(h-\psi),h-\psi)_{L^2(\Omega)}
\\
&\leq
c\|\cov_A\psi\|_{L^2(\Omega)}\|\cov_A(h-\psi)\|_{L^2(\Omega)}
+ \|\sR_A\|_{L^2(\Omega)}\|h-\psi\|_{L^4(\Omega)}^2
\\
&\leq
c\|\cov_A\psi\|_{L^2(\Omega)}\|\cov_A(h-\psi)\|_{L^2(\Omega)}
+ K\|\sR_A\|_{L^2(\Omega)}\|h-\psi\|_{L^2_{1,A}(\Omega)}^2,
\end{align*}
where $K=K(\Omega)$ is the norm of the Sobolev embedding $L^2_{1,A}(\Omega)
\to L^4(\Omega)$.
Therefore, using the interpolation inequality $xy\leq \eps x^2+\eps^{-1}
y^2$ and the bound $\|\sR_A\|_{L^2(\Omega)} < \eps\cdot\min\{1,K^{-1}\}$,
(by hypothesis \eqref{eq:L2BochnerFormHypothesis})
the preceding inequality gives
\begin{equation}
\label{eq:L21Esth-3}
\|\cov_A(h-\psi)\|_{L^2(\Omega)}
\leq 
c\|\cov_A\psi\|_{L^2(\Omega)} 
+ \eps\left(\|\cov_A(h-\psi)\|_{L^2(\Omega)} 
+ \|h-\psi\|_{L^2(\Omega)}\right).
\end{equation}
Combining inequalities \eqref{eq:L21Esth-2} and \eqref{eq:L21Esth-3}
and using rearrangement yields
\begin{equation}
\label{eq:L21Esth-4}
\|h-\psi\|_{L^2_{1,A}(\Omega)}
\leq
c\|\cov_A\psi\|_{L^2(\Omega)}.
\end{equation}
Hence, the a priori estimate \eqref{eq:L21DirichletDiracAPrioriEst} follows
from inequalities \eqref{eq:L21Esth-1} and \eqref{eq:L21Esth-4}.
\end{proof}

With the preceding lemmas at hand, we can now complete the proof of 
Theorem \ref{thm:DirichletDiracLaplacian}.

\begin{proof}[Proof of Theorem \ref{thm:DirichletDiracLaplacian}]
Observe that the given Dirichlet problem, $D_A^2h'=0$ and $h'-\psi' \in
L^p_2(\Omega,\partial\Omega;V)$, is equivalent to 
\begin{equation}
\label{eq:EquivDirichletDiracProblem}
D_A^2h = \psi \quad\text{on $\Omega$ and}\quad
h \in L^p_2(\Omega,\partial\Omega;V),
\end{equation}
if we take $h = h'-\psi' \in L^p_2(\Omega,\partial\Omega;V)$ and $\psi =
-D_A^2\psi' \in L^p(\Omega;V)$. Let $\sH = L^2_1(\Omega,\partial\Omega;V)$
and define a bilinear form $\sB:\sH\times\sH\to\CC$ by 
$$
\sB(\varphi_1,\varphi_2) 
= 
(\varphi_1,D_A^2\varphi_2)_{L^2(\Omega)} 
= 
(D_A\varphi_1,D_A\varphi_2)_{L^2(\Omega)}.
$$ 
Similarly, define $\psi^* \in \sH^*$ by setting 
$$
\psi^*(\varphi) = (\varphi,\psi)_{L^2(\Omega)},
\quad\text{for all }\varphi \in \sH,
$$
noting that $(\varphi,\psi)_{L^2(\Omega)} \leq 
\|\varphi\|_{L^2(\Omega)}\cdot \|\psi\|_{L^2(\Omega)}$, so
$\psi^*:\sH\to\CC$ is a continuous map, as claimed. Integrating by parts
and arguing as in the proof of Lemma
\ref{lem:FirstDiracDirichletEigenvalueBall} we see that
\begin{align*}
\|\varphi\|_{L^2_{1,A}(\Omega)}^2
&=
\|\varphi\|_{L^2(\Omega)}^2 + \|\cov_A\varphi\|_{L^2(\Omega)}^2
\\
&\leq
\mu_1^{-1}\|D_A\varphi\|_{L^2(\Omega)}^2
+ \|D_A\varphi\|_{L^2(\Omega)}^2 - (\sR_A\varphi,\varphi)_{L^2(\Omega)}
\\
&\leq
(1+\mu_1^{-1})\|D_A\varphi\|_{L^2(\Omega)}^2
+ \|\sR_A\|_{L^2(\Omega)}\|\varphi\|_{L^4(\Omega)}^2
\\
&\leq
(1+\mu_1^{-1})\|D_A\varphi\|_{L^2(\Omega)}^2
+ K\|\sR_A\|_{L^2(\Omega)}\|\varphi\|_{L^2_{1,A}(\Omega)}^2,
\end{align*}
where the constant $K(\Omega)$ depends only on the geometry of cone
defining the interior cone property for $\Omega$ \cite[\S IV.3]{Adams} by
the Sobolev embedding $L^2_1(\Omega)\subset L^4(\Omega)$ \cite[Theorem
V.5.4]{Adams}. By hypothesis \eqref{eq:L2BochnerFormHypothesis}, 
we may suppose
$\|\sR_A\|_{L^2(\Omega)} \leq \half K^{-2}$ and hence rearrangement in
the last inequality yields
$$
\|\varphi\|_{L^2_{1,A}(\Omega)}^2
\leq
c'\sB(\varphi,\varphi),
\quad\text{for all }\varphi \in \sH,
$$ 
with $c' = 2(1+\mu_1^{-1})$.  Thus, $\sB$ is a bounded, {\rm coercive\/}
bilinear form on $\sH$ (see, for example, \cite[Equation (5.11)]{GT}). The
Lax-Milgram theorem \cite[\S 5.8]{GT} then implies that the map $\sB:\sH\to
\sH^*$, $\varphi\mapsto \sB(\cdot,\varphi)$ is continuous and bijective ---
with continuous inverse by the open mapping theorem --- and consequently
there is a unique solution $h\in\sH$ to the equation
$$
\sB(\varphi,h) = (\varphi,\psi)_{L^2(\Omega)},
\quad\text{for all }\varphi \in \sH.
$$
In other words, there is unique weak solution $h\in
L^2_1(\Omega,\partial\Omega;V)$ to the Dirichlet problem
$D_A^2h=\psi$. Standard regularity arguments 
%\marginpar{[NEXT SECTION; JUST BOUNDS] }
ensure that $h$ is a strong
(or classical) solution, with the asserted regularity properties for a
given $\psi$ \cite{GT}. Let 
$$
G_A:\sH\to \sH, \quad\psi\mapsto \psi^* \mapsto \sB(\cdot,h) \mapsto h
$$
be the Green's operator for the Dirichlet problem $D_A^2h=\psi$, $h\in
L^2_1(\Omega,\partial\Omega;V)$. Since $G_A$ can be written as a
composition of continuous operators,
$$
L^2_1(\Omega;V) \to L^2(\Omega;V) \to L^2(\Omega;V)^* 
\to L^2_1(\Omega,\partial\Omega;V)^*
\to L^2_1(\Omega,\partial\Omega;V),
$$
with the first embedding being compact by Rellich's theorem \cite{Adams}
(for bounded $\Omega$), then $G_A$ is compact. Therefore, $G_A$ has a
discrete spectrum and zero as the only possible limit point; the spectrum
is real since $G_A$ is seen to be self-adjoint in the usual way. Hence,
$D_A^2$ has a discrete, real spectrum with no limit points.
\end{proof}

\subsection{\boldmath{$L^p_2$} elliptic estimates for spinors over annuli}
\label{subsec:Lp2DiracLaplacianAnnulus}
In this subsection we complete the proof of Theorem
\ref{thm:Lp2SpinorAnnulusEst}, the main result of this section.

We shall need elliptic estimates for negative spinors, that is $L^2_{k+1}$
sections of $V^-=W^-\otimes E$, over open subsets of $X$ where the curvature
$F_A$ obeys a $C^0$ bound which is uniform with respect to $A$. The desired
estimates are fairly delicate, as the domain in question is of the form $U
= X - \cup_{i=1}^m B(x_i,4\lambda_i^{1/3})$ and we require estimates
which are uniform with respect to the ball radii. In particular, the
estimates we need do not follow from a simple application of the standard
elliptic estimate for $D_AD_A^*$ over $X$ using cutoff functions equal to
one on the complement of these balls and equal to zero on the smaller balls
$B(x_i,2\lambda_i^{1/3})$.

We shall need the following analogue --- in the case of Sobolev norms of
sections of a Hermitian bundle $V$ defined via covariant derivatives --- of
the interpolation inequality given by \cite[Theorem 7.27]{GT} for Sobolev
norms of functions on domains in $\RR^n$: 

\begin{lem}
\label{lem:Interpolation}
Let $\Omega$ be an open subset of $\RR^n$, not necessarily bounded.
Let $\Gamma$ denote the product connection on $V \cong \Omega\times\CC^r$
and suppose the unitary connection $A$ on $V$ is given by $A =
\Gamma + a$, with $\|a\|_{L^\8_{1,\Gamma}(\Omega)} \leq
\|F_A\|_{L^\8(\Omega)}$. Suppose $\|F_A\|_{L^\8(\Omega)} \leq M$, for
some positive constant $M$. Then for all $1\leq p < \8$, there is a
constant $C = C(n,p,M)$ such the following holds: for any $u\in
L^p_2(\Omega,\partial\Omega;V)$ and $0<\eps< \half M^{-1}$,
\begin{equation}
\label{eq:CovDerivInterpolation}
\|\cov_A u\|_{L^p(\Omega)}
\leq 
\eps\|\cov_A^2 u\|_{L^p(\Omega)} + C\eps^{-1}\|u\|_{L^p(\Omega)}.
\end{equation}
\end{lem}

\begin{proof}
{}From \cite[Theorem 7.27]{GT} (see the last inequality in the proof of
Theorem 7.27 on \cite[p. 172]{GT}) we have the interpolation inequality
\begin{equation}
\label{eq:GTInterpolation}
\|\cov_\Gamma\psi\|_{L^p(\Omega')}
\leq
c(\eps\|\cov_\Gamma^2\psi\|_{L^p(\Omega')} + \eps^{-1}\|\psi\|_{L^p(\Omega')}),
\end{equation}
for some constant $c(n,p)$ and all $\eps>0$. Observe that 
\begin{equation}
\label{eq:SecondCovDerAExpansion}
\cov_A^2u 
= 
\cov_\Gamma^2u + a\otimes \cov_\Gamma u +
(\cov_\Gamma a + a\otimes a)\otimes u.
\end{equation}
The interpolation inequality \eqref{eq:GTInterpolation}
and the identity \eqref{eq:SecondCovDerAExpansion} now yield an
interpolation inequality for $\cov_A u$:
\begin{align*}
\|\cov_A u\|_{L^p(\Omega)}
&\leq
\|\cov_\Gamma u\|_{L^p(\Omega)} + \|a\|_{L^\8(\Omega)}\|u\|_{L^p(\Omega)}
\\
&\leq
\eps\|\cov_\Gamma^2 u\|_{L^p(\Omega)} 
+ (1+\eps^{-1})\|a\|_{L^\8(\Omega)}\|u\|_{L^p(\Omega)}
\quad\text{(by \eqref{eq:GTInterpolation})}
\\
&\leq
\eps\|\cov_A^2 u\|_{L^p(\Omega)} 
+ \eps\|a\|_{L^\8(\Omega)}\|\cov_\Gamma u\|_{L^p(\Omega)} 
\\
&\quad + \left(\eps(\|\cov_\Gamma a\|_{L^\8(\Omega)}) + \|a\|_{L^\8(\Omega)}^2)
+ (1+\eps^{-1})\|a\|_{L^\8(\Omega)}\right)
\|u\|_{L^p(\Omega)}
\\
&\leq
\eps\|\cov_A^2 u\|_{L^p(\Omega)} 
+ \eps\|a\|_{L^\8(\Omega)}\|\cov_A u\|_{L^p(\Omega)} 
\\
&\quad 
+ \left(\eps(\|\cov_\Gamma a\|_{L^\8(\Omega)}) + 2\|a\|_{L^\8(\Omega)}^2)
+ (1+\eps^{-1})\|a\|_{L^\8(\Omega)}\right)
\|u\|_{L^p(\Omega)}.
\end{align*}
Consequently, for $\eps < \half\|a\|_{L^\8(\Omega)}^{-1}$, rearrangement
yields 
$$
\|\cov_A u\|_{L^p(\Omega)}
\leq 
2\eps\|\cov_A^2 u\|_{L^p(\Omega)} + C\eps^{-1}\|u\|_{L^p(\Omega)},
$$
where $C$ depends on $\|a\|_{L^\8_{1,\Gamma}(\Omega)}$. Replacing $\eps$ by
$\eps/2$ and using $\|a\|_{L^\8_{1,\Gamma}(\Omega)} \leq
c\|F_A\|_{L^\8(\Omega)}$ yields the conclusion.
\end{proof}

Modulo Proposition \ref{prop:SystemsPDEHomogeneousBound} and
Corollary \ref{cor:DiracLaplacianHomogeneousBound} below, we can complete
the proof of Theorem \ref{thm:Lp2SpinorAnnulusEst}:

\begin{proof}[Proof of Theorem \ref{thm:Lp2SpinorAnnulusEst}, given 
Proposition \ref{prop:SystemsPDEHomogeneousBound} and
Corollary \ref{cor:DiracLaplacianHomogeneousBound}]
According to Theorem \ref{thm:DirichletDiracLaplacian} there is a unique
solution $h\in L^p_{2,A}(\Omega)$ to the Dirichlet problem
$$
D_AD_A^*h = 0\text{ on }\Omega 
\quad\text{and}\quad h = \psi\text{ on }\partial\Omega.
$$
We can estimate the $C^0\cap L^2_{2,A}(\Omega')$ norm of $h$ in
terms of $\|h\|_{L^2(\Omega)}$
using \cite[Theorem 1.1]{FeehanKato}, since $\Omega'
\subseteq \Omega(2r_0,r_1/2)$, $A$ is Yang-Mills with $L^2$-small
curvature, and $g$ obeys the bounded geometry conditions of 
 \cite[\S 1.4]{FeehanKato}. We can therefore
estimate the $C^0\cap L^2_{2,A}(\Omega')$ norm of $h$ in terms of
$\|\psi\|_{L^2_{1,A}(\Omega)}$, using 
estimate \eqref{eq:L21DirichletDiracAPrioriEst}. Thus
$$
\psi-h \in L^p_2(\Omega,\partial\Omega;V),
$$
and so Corollary \ref{cor:DiracLaplacianHomogeneousBound} implies that
\begin{equation}
\label{eq:SystemsCZLpEst}
\|\cov_A^2(\psi-h)\|_{L^p(\Omega)}
\leq 
C\|D_AD_A^*(\psi-h)\|_{L^p(\Omega)}
=
C\|D_AD_A^*\psi\|_{L^p(\Omega)}.
\end{equation}
Consequently, we have
\begin{align*}
\|\cov_A^2\psi\|_{L^p(\Omega')}
&\leq
\|\cov_A^2(\psi-h)\|_{L^p(\Omega')}
+ \|\cov_A^2h\|_{L^p(\Omega')}
\\
&\leq
\|\cov_A^2(\psi-h)\|_{L^p(\Omega)} + c\|h\|_{L^2(\Omega)}
\quad\text{(by \cite{FeehanKato})}
\\
&\leq
C(\|\cov_A^2(\psi-h)\|_{L^p(\Omega)} + \|\psi\|_{L^2_{1,A}(\Omega)})
\quad\text{(by \eqref{eq:L21DirichletDiracAPrioriEst})}
\\
&\leq
C(\|D_AD_A^*\psi\|_{L^p(\Omega)} + \|\psi\|_{L^2_{1,A}(\Omega)})
\quad\text{(by \eqref{eq:SystemsCZLpEst})}.
\end{align*}
But Lemma \ref{lem:Interpolation} implies that
$$
\|\cov_A(\psi-h)\|_{L^p(\Omega')}
\leq
\eps\|\cov_A^2(\psi-h)\|_{L^p(\Omega')}
+ C\eps^{-1}\|\psi-h\|_{L^p(\Omega')}
$$
and so, because $\|\cov_A\psi\|_{L^p(\Omega')} \leq
\|\cov_A(\psi-h)\|_{L^p(\Omega')} + \|\cov_A h\|_{L^p(\Omega')}$ and $p\leq
2$, we obtain the $L^p_{2,A}$ estimate 
$$
\|\psi\|_{L^p_{2,A}(\Omega')}
\leq
C(\|D_AD_A^*\psi\|_{L^p(\Omega)} + \|\psi\|_{L^2_{1,A}(\Omega)}).
$$
This completes the proof.
\end{proof}

The use of Lemma \ref{lem:Interpolation} in the preceding proof is
unnecessary when $p\leq 2$, as we hypothesized, but it indicates the
structure of the argument if this assumption were dropped (for example, if
we had useful $L^p_2(\Omega')$ estimates for $h$ when $2<p<\8$).

Thus we see that the key remaining ingredient in our argument is to prove
the following systems analogue of Lemma
\ref{lem:ScalarPDEHomogeneousBound}:

\begin{prop}
\label{prop:SystemsPDEHomogeneousBound}
Let $M$ be a positive constant. Let $\Omega\Subset \RR^4$ be a $C^1$
domain, where $\RR^n$ has a $C^\8$ metric $g$ of bounded geometry and
injectivity radius $\varrho$, let $\mu_1(\Omega,g)$ be the first
Dirichlet eigenvalue of the scalar Laplacian on $\Omega$, and let $1<p<\8$.    
Assume $\Omega$
obeys an interior cone condition for a fixed cone $\sK$.  Let $A$ be an
$C^\8(\Omega)$ \spinc connection on a Hermitian bundle $V|_\Omega$, where
$V=V^+\oplus V^-$ and $V^\pm \cong
\Omega\times\CC^4$.  Suppose 
\begin{align}
\label{eq:BoundedCurvHypothesis}
\|F_A\|_{L^\8(\Omega)}, 
\quad \|\mathrm{Rm}_g\|_{L^\8(\Omega)} &\leq M,
\end{align}
and that $A = \Gamma + a$, with 
\begin{equation}
\label{eq:C1EstConnectionOneForm}
\|a\|_{L^\8_{1,\Gamma}(\Omega)} \leq c\|F_A\|_{L^\8(\Omega)},
\end{equation}
where $c=c(M)$. Then there is a positive constant $\eps =
\eps(\sK,M,p,\varrho,\mu_1)$ such that the following holds.
For any $1<p\leq 2$ there is a constant $C=C(\sK,M,p,\varrho,\mu_1)$ such
that for all $u\in L^p_2(\Omega,\partial\Omega;\CC^r)$ we have
\begin{equation}
\label{eq:SystemsCalderonZygmund}
\|\cov_A^2u\|_{L^p(\Omega)} \leq C\|\cov_A^*\cov_A u\|_{L^p(\Omega)}. 
\end{equation}
\end{prop}

\begin{cor}
\label{cor:DiracLaplacianHomogeneousBound}
Continue the hypotheses of Proposition
\ref{prop:SystemsPDEHomogeneousBound} and, in addition require that
\begin{equation}
\label{eq:L2CurvHypothesis}
\|F_A\|_{L^2(\Omega)}, 
\quad \|\mathrm{Rm}_g\|_{L^2(\Omega)} \leq \eps,
\end{equation}
Then there is a constant $C=C(\sK,M,p,\varrho,\mu_1)$
such that for all $\psi\in L^p_2(\Omega,\partial\Omega;V)$,
\begin{equation}
\label{eq:DiracLaplacianCalderonZygmund}
\|\cov_A^2\psi\|_{L^p(\Omega)} \leq C\|D_A^2 u\|_{L^p(\Omega)}. 
\end{equation}
\end{cor}

\begin{rmk}
For the domain $\Omega=\Omega(r_0,r_1)$ of interest, the dependence of
$\sK(\Omega)$ on the geometry of $\Omega$ is explained in the paragraph
following the statement of Theorem \ref{thm:Lp2SpinorAnnulusEst}, while the
dependence of $\mu_1=\mu_1(\Omega)$ on $\Omega$ is explained in Lemma
\ref{lem:FirstDiracDirichletEigenvalueBall} and Remark
\ref{rmk:EigenvalueRescaling}.
\end{rmk}

The refinement we give in Proposition \ref{prop:SystemsPDEHomogeneousBound}
and Corollary \ref{cor:DiracLaplacianHomogeneousBound} 
is crucial for our purposes because the dependence of the constant $C$ in
\cite{GT} on the domain $\Omega$ is often difficult to
discern from the proofs in \cite{GT}. While we might fix $r_0=1$ and
consider $r_1$ large in $\Omega(r_0,r_1)$, there is usually an implicit
assumption in the arguments of \cite[Chapters 7--9]{GT} that the domain
$\Omega$ is bounded. So we cannot necessarily assume that $r_1=\8$ is
allowed and thus infer that the constant $C$ is independent of $r_1\in
(1,\8)$. For this purpose, the simplest strategy is to take the following
corollary of the {\em Calderon-Zygmund inequality\/} \cite[Theorem 9.9]{GT}
as our starting point --- as do Gilbarg and Trudinger, for example, in
their proof of \cite[Theorem 9.14]{GT}:

\begin{thm}[Calderon-Zygmund]
\label{thm:CalderonZygmund}
\cite[Corollary 9.10]{GT}
Let $\Omega\subset\RR^n$ be a domain, not necessarily bounded, where
$\RR^n$ has its standard metric and let $1<p<\8$. Then there is a constant
$c(n,p)$ such that for all $u\in L^p_2(\Omega,\partial\Omega;\CC)$,
\begin{equation}
\label{eq:CalderonZygmund}
\|\cov^2u\|_{L^p(\Omega)} \leq c\|\Delta u\|_{L^p(\Omega)}. 
\end{equation}
If $p=2$, then $\|\cov^2u\|_{L^2(\Omega)} = \|\Delta u\|_{L^2(\Omega)}$. 
\end{thm}

In particular, the constant $c(n,p)$ in inequality
\eqref{eq:CalderonZygmund} is independent of $\Omega$; thus we can
establish the dependence of $C$ on the radii of $\Omega(r_0,r_1)$ by
perturbation from this special case. 
%\marginpar{[CHECK $\Omega$ BOUNDED CONDITIONS.]}

\begin{proof}[Proof of Proposition \ref{prop:SystemsPDEHomogeneousBound}]
First observe that if $u=(u_1,\dots,u_r)\in
L^p_2(\Omega,\partial\Omega;\CC^r)$, for some integer $r\geq 1$, then each
$u_j\in L^p_2(\Omega,\partial\Omega;\CC^r)$, $1\leq j\leq r$, obeys the
estimate \eqref{eq:CalderonZygmund} and so the same holds for $u$ with
$\cov$ now defined by the product connection, $\cov_\Gamma$, on
$\Omega\times\CC^r$ and the Levi-Civita connection on $T^*\RR^n$ for
the standard metric on $\Omega\subset\RR^n$. 

To go from the standard metric $\delta^{ij}$ on
$\RR^n$ to a metric $g^{ij}$ of bounded geometry on $\RR^n$ we can invoke
the standard argument of \cite[Lemma 6.1 \& Theorem 9.11]{GT} to pass from
the metric $\delta^{ij}$ to a constant coefficient metric $a^{ij}$ and then
to the given metric $g^{ij}$ to show that the bound
\eqref{eq:CalderonZygmund} continues to hold when $\delta^{ij}$ is replaced
by $g^{ij}$. Alternatively, we can more simply observe that the proof of
the Calderon-Zygmund inequality \cite[Equation (9.27)]{GT} allows the
Laplacian $\Delta_g$ on $C^\8(\RR^n)$ to be defined by a Riemannian metric
$g^{ij}$ of bounded geometry, that is $C^2(\bar\Omega)$-bounded curvature
and injectivity radius bounded below by a fixed positive constant. (See
\cite[pp. 6--7]{Chavel} for the first Green's identity and \cite{Aubin} for the
existence of the fundamental solution for the Laplacian $\Delta_g$, both of
which are used in \cite{GT} when $g^{ij} = \delta^{ij}$.) The
constant $c(n,p)$ in inequality \eqref{eq:CalderonZygmund} is then simply
replaced by $c(M,n,p)$, where $M$ is an $L^\8(\Omega)$ bound on the
Riemannian curvature of the metric $g$.

Lastly, we need to replace the product connection $\Gamma$ on
$\Omega\times\CC^r$ with the given connection on the Hermitian vector
bundle $V|_\Omega\simeq \Omega\times\CC^r$. Writing $A = \Gamma+a$, observe
that, schematically,
\begin{equation}
\label{eq:LaplacianAExpansion}
\cov_A^*\cov_Au 
= 
\cov_\Gamma^*\cov_\Gamma u + a\otimes \cov_\Gamma u +
(\cov_\Gamma a + a\otimes a)\otimes u.
\end{equation}
The identities \eqref{eq:SecondCovDerAExpansion} and
\eqref{eq:LaplacianAExpansion} yield
\begin{align*}
\|\cov_A^2u\|_{L^p(\Omega)} 
&\leq
\|\cov_\Gamma^2u\|_{L^p(\Omega)} 
+ \|a\|_{L^\8(\Omega)}\|\cov_\Gamma u\|_{L^p(\Omega)}
\\
&\quad + (\|\cov_\Gamma a\|_{L^\8(\Omega)} + \|a\|_{L^\8(\Omega)}^2)
\|u\|_{L^p(\Omega)}
\\
&\leq
C(\|\cov_\Gamma^*\cov_\Gamma u\|_{L^p(\Omega)} 
+ \|\cov_\Gamma u\|_{L^p(\Omega)} + \|u\|_{L^p(\Omega)})
\\
&\leq
C(\|\cov_A^*\cov_A u\|_{L^p(\Omega)} 
+ \|\cov_A u\|_{L^p(\Omega)} + \|u\|_{L^p(\Omega)}).
\end{align*}
So, using Lemma \ref{lem:Interpolation} and rearrangement we have
\begin{equation}
\label{eq:LpSecondCovDeriv-AlmostFinal}
\|\cov_A^2u\|_{L^p(\Omega)} 
\leq 
C(\|\cov_A^*\cov_A u\|_{L^p(\Omega)} + \|u\|_{L^p(\Omega)}).
\end{equation}
Here, $C$ now also depends on an upper bound for
$\|a\|_{L^\8_{1,\Gamma}(\Omega)}$. 

It remains to argue that the term $\|u\|_{L^p(\Omega)}$ on the right-hand
side of inequality \eqref{eq:LpSecondCovDeriv-AlmostFinal} can be omitted.
We first show 

\begin{claim}
\label{claim:LpSecondCovDeriv-L2AlmostFinal}
Continue the hypotheses of Proposition
\ref{prop:SystemsPDEHomogeneousBound}.  Assume $|\Omega| = \vol_g(\Omega)
< \8$ and $n=4$. Then for any $1< p<\8$, there is a constant
$C=(\sK,p,|\Omega|,M)$ such that
\begin{equation}
\label{eq:LpSecondCovDeriv-L2AlmostFinal}
\|\cov_A^2u\|_{L^p(\Omega)} 
\leq 
C(\|\cov_A^*\cov_A u\|_{L^p(\Omega)} + \|u\|_{L^2(\Omega)}),
\end{equation}
\end{claim}

\begin{proof}
We can assume $p\neq 2$ throughout the proof, as the result is immediate
from inequality \eqref{eq:LpSecondCovDeriv-AlmostFinal} in this case. 
If $1< p< 2$ then H\"older's inequality gives
\begin{equation}
\|u\|_{L^p(\Omega)} \leq |\Omega|^{2p/(2-p)}\cdot\|u\|_{L^2(\Omega)},
\end{equation}
and via inequality \eqref{eq:LpSecondCovDeriv-AlmostFinal} we are done in
this case.

Hence, assume $2< p<\8$ for the remainder of the proof and note that
\begin{equation}
\label{eq:FiniteVolSobolevEmbedding}
\|u\|_{L^2(\Omega)} \leq |\Omega|^{2p/(p-2)}\cdot\|u\|_{L^p(\Omega)}.
\end{equation}
{}From the interpolation inequality \cite[Equation (7.10)]{GT},
we have
\begin{equation}
\label{eq:OmitLpuGTInterpolation}
\|u\|_{L^p(\Omega)} 
\leq \eps\|u\|_{L^{2p}(\Omega)} + \eps^{-\eta}\|u\|_{L^2(\Omega)},
\end{equation}
where $\eta = p - 2 > 0$ and $\eps$ is any positive constant.  Now the
Sobolev embedding theorem \cite[Theorem V.5.4]{Adams}, the interpolation
inequality \eqref{eq:CovDerivInterpolation}, and the finite-volume $L^p$
estimate \eqref{eq:FiniteVolSobolevEmbedding} yield 
\begin{equation}
\label{eq:OmitLpuSobolevEmbedding}
\begin{aligned}
\|u\|_{L^{2p}(\Omega)} 
&\leq c(\sK,p)(\|\cov_A^2u\|_{L^2(\Omega)} 
+ \|\cov_Au\|_{L^2(\Omega)} + \|u\|_{L^2(\Omega)})
\\
&\leq c(\sK,M,p)(\|\cov_A^2u\|_{L^2(\Omega)} 
+ \|u\|_{L^2(\Omega)})
\\
&\leq c(\sK,|\Omega|,M,p)(\|\cov_A^2u\|_{L^p(\Omega)} 
+ \|u\|_{L^2(\Omega)}).
\end{aligned}
\end{equation}
(The dependence of $c$ on the cone $\sK$ follows from the Sobolev embedding
theorem; see \cite[\S IV.3 \& Theorem V.5.4]{Adams}.) Choosing $\eps <
\half c^{-1}$ in inequality \eqref{eq:OmitLpuGTInterpolation}, we can now
combine the estimates \eqref{eq:OmitLpuGTInterpolation},
\eqref{eq:OmitLpuSobolevEmbedding}, and
\eqref{eq:LpSecondCovDeriv-AlmostFinal} and use rearrangement to give the 
desired result when $2<p<\8$.
\end{proof}

Our eigenvalue estimate \eqref{eq:FirstCovLapDirichletEigenvalueBall} 
yields
\begin{equation}
\label{eq:L2EigenvalueCovLapEst}
\|u\|_{L^2(\Omega)} \leq \mu_1^{-1}\|\cov_A^*\cov_A u\|_{L^2(\Omega)},
\end{equation}
where $\mu_1 = \mu_1(\Omega;\cov_A^*\cov_A)\geq  \mu_1(\Omega;\Delta)>0$
(for any bounded $C^1$ domain $\Omega$) 
and so, for $2\leq p < \8$, inequalities
\eqref{eq:LpSecondCovDeriv-L2AlmostFinal} and
\eqref{eq:FiniteVolSobolevEmbedding} give the desired estimate
\eqref{eq:SystemsCalderonZygmund} for $2\leq p < \8$, namely
\begin{equation}
\label{eq:LpSecondCovDeriv-pgeq2AlmostFinal}
\|\cov_A^2u\|_{L^p(\Omega)} 
\leq 
C\|\cov_A^*\cov_A u\|_{L^p(\Omega)}, 
\quad\text{for $2\leq p < \8$}.
\end{equation}
Note that for $2\leq p < \8$, inequalities
\eqref{eq:OmitLpuGTInterpolation}, \eqref{eq:OmitLpuSobolevEmbedding},
\eqref{eq:L2EigenvalueCovLapEst}, \eqref{eq:FiniteVolSobolevEmbedding}, and
\eqref{eq:LpSecondCovDeriv-pgeq2AlmostFinal} yield
\begin{equation}
\label{eq:LpuBoundedByLpCovLap}
\|u\|_{L^p(\Omega)} \leq C\|\cov_A^*\cov_A u\|_{L^p(\Omega)},
\quad\text{for $2\leq p < \8$}.
\end{equation}
We use the preceding estimate to
extend the bound \eqref{eq:LpSecondCovDeriv-pgeq2AlmostFinal}
to $1<p\leq 2$ by duality. The 
bound \eqref{eq:LpuBoundedByLpCovLap} is equivalent to
\begin{equation}
\label{eq:LpGreenuBoundedByLpu}
\|G_Au\|_{L^p(\Omega)} \leq C\|u\|_{L^p(\Omega)},
\quad\text{for $2\leq p < \8$},
\end{equation}
where $G_A = (\cov_A^*\cov_A)^{-1}$ is the $L^2$-self-adjoint Dirichlet
Green's operator.   

For $1<p\leq 2$, write $1=1/p+1/p'$, with $2\leq p' < \8$ and observe that
\begin{align*}
\|G_Au\|_{L^p(\Omega)}
&=
\sup\{(G_Au,v)_{L^2(\Omega)}: \|v\|_{L^{p'}(\Omega)} = 1\}
\\
&=
\sup\{(u,G_Av)_{L^2(\Omega)}: \|v\|_{L^{p'}(\Omega)} = 1\}
\\
&\leq
\sup\{\|u\|_{L^p(\Omega)}\|G_Av\|_{L^{p'}(\Omega)}: 
\|v\|_{L^{p'}(\Omega)} = 1\}
\\
&\leq
\|u\|_{L^p(\Omega)}\cdot
\sup\{C\|v\|_{L^{p'}(\Omega)}: \|v\|_{L^{p'}(\Omega)} = 1\}
\quad\text{(by \eqref{eq:LpGreenuBoundedByLpu})}
\\
&= C\|u\|_{L^p(\Omega)}.
\end{align*}
Hence, inequalities \eqref{eq:LpGreenuBoundedByLpu} and thus
\eqref{eq:LpuBoundedByLpCovLap} also hold for $1<p\leq 2$ and so combining
\eqref{eq:LpuBoundedByLpCovLap} with
\eqref{eq:LpSecondCovDeriv-AlmostFinal} yields the desired estimate
\eqref{eq:SystemsCalderonZygmund} for $1<p\leq 2$.
\end{proof}

\begin{proof}[Proof of Corollary \ref{cor:DiracLaplacianHomogeneousBound}]
This follows by the Bochner formula $D_A^2 = \cov_A^*\cov_A + \sR_A$ and
the argument at the conclusion of the proof of Proposition
\ref{prop:SystemsPDEHomogeneousBound}, now using our $L^2$ hypothesis on
$\sR_A$ via \eqref{eq:L2CurvHypothesis} and the eigenvalue estimate 
\eqref{eq:FirstDiracDirichletEigenvalueBall}
for $D_A^2$ in place of the ones
\eqref{eq:FirstCovLapDirichletEigenvalueBall}, 
\eqref{eq:L2EigenvalueCovLapEst} for $\cov_A^*\cov_A$. 
\end{proof}

%end of file
%file: eigenvalue.tex

\section{Eigenvalue estimates for PU(2)-monopole Laplacians}
\label{sec:Eigenvalue}
We need to analyze the behavior of the small eigenvalues of the Laplacian
$d_{A,\Phi}^1d_{A,\Phi}^{1,*}$ on $L^2(X,\Lambda^+\otimes\fg_E)\oplus
L^2(X,V^-)$ as the point $[A,\Phi]$ in $\sC(\ft)$ converges, in the Uhlenbeck
topology, to a point $([A_0,\Phi_0],\bx)$ in
$\sC(\ft_\ell)\times\Sym^\ell(X)$, where $\ft$ is the \spinu structure
$(\rho,W^+,W^-,E)$ on $(X,g)$ and $\ft_\ell = (\rho,W^+,W^-,E_\ell)$. 
Our derivation in \cite{FL4} of the eigenvalue bounds
for the Laplacian $d^{0,*}_{A,\Phi}d_{A,\Phi}^0$ on $L^2(X,\fg_E)$ is
closely modelled on the arguments given here for the Laplacian
$d_{A,\Phi}^1d_{A,\Phi}^{1,*}$ on $L^2(X,\Lambda^+\otimes\fg_E)\oplus
L^2(X,V^-)$.

Perhaps the earliest results of this type and those best known to geometric
analysts are due to C. H. Taubes, for the anti-self-dual Laplacian
$d_A^+d_A^{+,*}$
\cite{TauSelfDual}, \cite{TauIndef}, and our basic method ultimately goes
back to Taubes 
--- though not the details, which are quite different for $\PU(2)$
monopoles and considerably more challenging. In \cite[Lemma 4.8]{TauIndef},
Taubes describes the small eigenvalues of $d_A^+d_A^{+,*}$, where $A$ is an
approximately anti-self-dual connection constructed by splicing $|\bx|\leq
\ell$ one-instantons onto the background product connection
$\Ga$. Analogous arguments are used by Mrowka in the proof of Theorem
6.1.0.0.4 \cite{MrowkaThesis} in order to establish his `main eigenvalue
estimate', though the eigenvalue estimates are more precise here in that we
give both upper and lower bounds.

More recently, there has been work of S. Cappell, R. Lee, and E. Y. Miller
\cite{CappellLeeMiller1}, who considered the effect of cut-and-paste and
stretching constructions on the eigenvalues of a first-order, self-adjoint,
elliptic operator on a closed $n$-manifold containing a long tube. This
generalizes the situation of \cite{MrowkaThesis} somewhat, so it is of
interest to summarize their main result at this point.
To set notation, temporarily
let $X$ be a closed smooth $n$-manifold, $Y\subset X$ a
separating, codimension-one submanifold decomposing $X$ as $X=X_1\cup_Y X_2$.
Let $D$ be some first-order, elliptic self-adjoint operator
on $X$ of ``Atiyah-Patodi-Singer type'' on a neighborhood $[-1,1]\times
Y$ of $Y$; thus (1) $D=\sigma({\partial/ {\partial t}}
+\widehat{D})$ in this neighborhood. Let $X(T)$ denote the manifold obtained
by replacing the neighborhood of $Y$ by $[-T,T]\times Y$. The
operator $D$ has a natural extension to $X(T)$ using $(1)$ to define $D$ on
the neighborhood of $Y$. Define $X_1(T)= X_1\cup ([0,T]\times
Y)$, $X_1(\infty)=X_1\cup ([0,\infty)\times Y)$; note that
$D$ determines an operator on $X_1(T)$ and $X_1(\infty)$ using $(1)$
again. Similarly define $X_2(T)$ and $X_2(\infty)$.

The results of this \cite{CappellLeeMiller1}
concern constructing eigenvectors with low
eigenvalues of $D$ on $X(T)$ in terms of corresponding eigenvectors of $D$
on $X_i(T)$ and $L^2$ eigenvectors of $D$ on $X_i(\infty)$.

The main theorem of \cite{CappellLeeMiller1}
roughly states that for $T$ large enough,
the eigenvectors of $D$ on $X(T)$ whose eigenvalues lie in $(-T^{3/2
+\epsilon},T^{3/2 +\epsilon})$ can be approximated by splicing together
$L^2$ eigenvectors on $X_1(\infty)$ and $X_2(\infty)$. (If the
kernel of $\widehat{D}$ is non-empty, one needs also to ``splice in'' certain
extended $L^2$ solutions to $Du=0$ on $X_i(\infty)$.) Moreover, any
eigenvalue in the range $(-T^{3/2 +\epsilon},T^{3/2 +\epsilon})$
actually lies in the range $(-\exp(cT),\exp(cT))$ for some $c>0$; that is,
``small eigenvalues are exponentially small in $T$''. 
%\marginpar{\tiny preceding 4 paragraphs needed?}

In the present article, our main application of the eigenvalue bounds of
this section will be to deriving an estimate for the partial right inverse
$P_{A,\Phi,\mu} =
d_{A,\Phi}^{1,*}(d_{A,\Phi}^1d_{A,\Phi}^{1,*})^{-1}\Pi_{A,\Phi,\mu}^\perp$
for $d_{A,\Phi}^1$, as given in Corollary \ref{cor:L21AEstPAaphi}. Our
estimates on $L^2$-orthogonal projections of approximate eigenvectors will
play an important role in \cite{FL4} in our proofs of injectivity and
surjectivity of the gluing maps and of Uhlenbeck continuity of the gluing
maps.

In the next definition we collect some of the key technical properties we
shall require of the pairs $(A,\Phi)$ so we can describe the
small-eigenvalue behavior of the Laplacians of interest in this
section. While our main eigenvalue estimate, Theorem \ref{thm:H2SmallEval},
is only applied in this paper to the case where $(A,\Phi)$ is a spliced
pair, it is perhaps useful to try to isolate the minimal properties we
actually need.

\begin{defn}
\label{defn:LambdaClosePair}
Let $(X,g)$ be a closed, oriented four-manifold with Riemannian metric
$g$. Let $E_\ell$ be a Hermitian, rank-two vector bundle over $X$, let
$(\rho,W^+,W^-)$ be a
\spinc structure on $X$, and set $V_\ell = W\otimes E_\ell$, and $V_\ell^\pm =
W^\pm\otimes E_\ell$. Write $\ft_\ell = (\rho,W^+,W^-,E_\ell)$.
Let $\ell\geq 1$ be an integer and suppose
$(A_0,\Phi_0,\bx)$ is an $L^2_k$ representative on
$(\fg_{E_\ell},V_\ell^+)$ of a 
point in $\sC(\ft_\ell)\times\Sym^\ell(X)$, with $k\geq 4$.  Let
$\nu_2[A_0,\Phi_0]$ be the least positive eigenvalue of the
Laplacian $d_{A_0,\Phi_0}^1d_{A_0,\Phi_0}^{1,*}$ on
$L^2(\Lambda^+\otimes\fg_{E_\ell})\oplus L^2(V_\ell^-)$. Suppose $C$ and
$\lambda_0\leq 1$ are positive constants depending at most on $g$ and the
maximum of $\|F_{A_0}\|_{L^2_{k-1,A_0}(X)}$, $\|\Phi_0\|_{L^2_{k,A_0}(X)}$,
$\nu_2[A_0,\Phi_0]$, $\nu_2[A_0,\Phi_0]^{-1}$,
$\dim\Ker d_{A_0,\Phi_0}^{1,*}$, and $\ell$.

For $m = |\bx|$ and $\lambda_p\in (0,\lambda_0]$, $1\leq p\leq m$, define a
precompact open subset $U_0\Subset X\less\bx$ and balls surrounding the
points $x\in\bx$ by setting 
%\marginpar{\tiny Check $U_0$ usage. radii?}
\begin{equation}
\label{eq:UZero}
U_0
=
X-\mathop{\cup}\limits_{p=1}^m\bar B(x_p,4\lambda_p^{1/3})
\quad\text{and}\quad
B_p
= 
B(x_p,\lambda_p^{1/3}/4),\quad 1\leq p\leq m.
\end{equation}
Let $E$ be a Hermitian, rank-two vector bundle over $X$ with $\det E = \det
E_\ell$ and $c_2(E)=c_2(E_\ell)+\ell$. Write $\ft = (\rho,W^+,W^-,E)$.
Denote $V = W\otimes E$ and $V^\pm =
W^\pm\otimes E$. Let $\kappa_p\geq 1$ denote the multiplicities of the
points $x_p\in\bx$, so $\sum_{p=1}^m\kappa_p = \ell$.  Fix an isomorphism
$E|_{X\less\bx} \cong E_\ell|_{X\less\bx}$.  For $\lambda = \max_{1\leq p\leq
m}\lambda_p$, we say that an $L^2_k$ pair $(A,\Phi)$ on $(\fg_E,V^+)$,
representing a point in $\sC(\ft)$, is $\lambda$-{\em close\/} to the ideal
triple $(A_0,\Phi_0,\bx)$ on $(\fg_{E_\ell},V_\ell^+)$ if the following hold:
\begin{align}
\label{eq:APhiA0Phi0L4Close}
\|(A,\Phi)-(A_0,\Phi_0)\|_{L^2_{1,A_0}(U_0)} &\leq C\lambda^{1/3},
\\
\label{eq:FASphereBallIntegerMass}
|\|F_A\|_{L^2(B_p)} - (8\pi^2\kappa_p)^{1/2}| &\leq 1/100, 
\quad 1\leq p\leq m,
\\
\label{eq:FA+SphereBallL2Small}
\|F_A^{+,g}\|_{L^2(B_p)} &\leq C\lambda^{1/3}, \quad 1\leq p\leq m,
\\
\label{eq:PhiLInftyBallSmall}
\|\Phi\|_{L^\8(X)} &\leq C,
\\
\label{eq:PhiIsZeroOnBallSmall}
\Phi &= 0 \quad\text{on }B_p.
\end{align}
\end{defn}

\begin{rmk}
\label{rmk:LambdaClosePairConstants}
Throughout this section, the constant $C$ may be allowed to increase or the
constant $\lambda_0$ may be allowed to decrease from line to line. However,
their dependencies will remain as stated in Definition
\ref{defn:LambdaClosePair}. As always, we use $c$ to denote a constant
which depends at most on the Riemannian metric of $X$, but is otherwise
universal. 
\end{rmk}

We remark that if $[A_\alpha,\Phi_\alpha]$ is a sequence of points in
$\sC(\ft)$ which converges to $[A_0,\Phi_0,\bx]$ in the Uhlenbeck topology
\cite[Definition 4.19]{FL1}, then, after passing to a subsequence, there is
a sequence of $L^2_{k+1,\loc}(X\less\bx)$ bundle isomorphisms
$u_\alpha:E|_{X\less\bx}\to E_\ell|_{X\less\bx}$ and a constant $\alpha_0$
such that for all $\alpha\geq
\alpha_0(\lambda)$, the pairs $u_\alpha(A_\alpha,\Phi_\alpha)$ obey the 
condition \eqref{eq:APhiA0Phi0L4Close}. If $[A_\alpha,\Phi_\alpha]$ is a
sequence of $\PU(2)$ monopoles, then the proof of Theorem 4.20 in
\cite{FL1} implies that the pairs $u_\alpha(A_\alpha,\Phi_\alpha)$ also
obey conditions \eqref{eq:FASphereBallIntegerMass},
\eqref{eq:FA+SphereBallL2Small}, and \eqref{eq:PhiLInftyBallSmall}.

Proposition \ref{prop:CutoffMonoPairEst} (and its proof) ensures that
conditions \eqref{eq:FASphereBallIntegerMass} and
\eqref{eq:FA+SphereBallL2Small} are obeyed by the approximate, extended
$\PU(2)$ monopoles produced by splicing in \S \ref{sec:Splicing}. The condition
\eqref{eq:PhiLInftyBallSmall} is also obeyed by such pairs, as one can
easily see by definition \eqref{eq:SplicedSpinor}. Though we require
$\Phi \equiv 0$ on each ball $B_p$ in
condition \eqref{eq:PhiIsZeroOnBallSmall}, this is really only used at one
point in the argument, namely in the proof of Lemma
\ref{lem:H2DeAEvectorEst} via the definition 
\eqref{eq:FromKernelSphereToSmallEvalVecs} of certain approximate
eigenvectors.
Condition \eqref{eq:APhiA0Phi0L4Close} is trivially obeyed
since, for a slightly different choice of $U_0$ (with balls
$B(x_p,4\lambda_p^{1/3})$ replaced by $B(x_p,8\lambda_p^{1/3}$), we have
$(A,\Phi) = (A_0,\Phi_0)$ on $U_0$.

\begin{thm}
\label{thm:H2SmallEval}
Continue the notation of Definition \ref{defn:LambdaClosePair}. Suppose
$\lambda\in(0,\lambda_0]$ and let $(A,\Phi)$ be an $L^2_k$ pair on
$(\fg_E,V)$ which is $\lambda$-close to the ideal pair $(A_0,\Phi_0,\bx)$.
Assume that the kernel of $d_{A_0,\Phi_0}^1d_{A_0,\Phi_0}^{1,*}$ has
dimension $n$ and let $N = n + 2\ell$. Let $\{\mu_l[A,\Phi]\}_{l\geq 1}$
denote the eigenvalues of the Laplacian $d_{A,\Phi}^1d_{A,\Phi}^{1,*}$ on
$L^2(\Lambda^+\otimes\fg_E)\oplus L^2(V^-)$, repeated according to their
multiplicity and in ascending order. Then
\begin{align*}
\mu_l[A,\Phi] &\leq
\begin{cases}
C(-\log\lambda)^{-3/2}, &\text{if }l=1,\dots,N, \\
K + C(-\log\lambda)^{-3/4}, &\text{if }l=N+1,
\end{cases} 
\\
\mu_l[A,\Phi] &\geq K- C(-\log\lambda)^{-3/4}, \quad \text{if }l\geq N+1,
\end{align*}
where
$$
K = \min\{\nu_2[A_0,\Phi_0],\nu_2[A_1],\dots,\nu_2[A_m]\}.
$$
\end{thm}

\begin{cor}
\label{cor:H2SmallEval}
Let $(X,g)$ be a closed, oriented, Riemannian four-manifold. Suppose
$[A_\alpha,\Phi_\alpha]$ is a sequence of points in $\sC(\ft)$ which
converges, in the Uhlenbeck topology, to a point
$[A_0,\Phi_0,\bx]\in\sC(\ft_\ell)\times\Sym^\ell(X)$ and suppose the
kernel of $d_{A_0,\Phi_0}^1d_{A_0,\Phi_0}^{1,*}$ has dimension $n$.  Then
the first $n+2\ell$ eigenvalues (counted with multiplicity) of
$d_{A_\alpha,\Phi_\alpha}^1d_{A_\alpha,\Phi_\alpha}^{1,*}$ converge to
zero, while the $(n+2\ell+1)$-st eigenvalue of
$d_{A_\alpha,\Phi_\alpha}^1d_{A_\alpha,\Phi_\alpha}^{1,*}$ converges to
$K$, the minimum of the least positive eigenvalues of
$d_{A_0,\Phi_0}^1d_{A_0,\Phi_0}^{1,*}$ and $D_{A_p,g_p}D_{A_p,g_p}^*$,
$1\leq p\leq m$.
\end{cor}

The proof of Theorem \ref{thm:H2SmallEval} takes up the remainder
of this section.  Since the argument is fairly lengthy, it seems useful to
outline the basic strategy. It is helpful to view of the proof as a
linearized version of the proof of the gluing theorem for $\PU(2)$ monopoles:
\begin{itemize}
\item
We establish upper bounds for the first $N+m+1$ eigenvalues of
$d_{A,\Phi}^1d_{A,\Phi}^{1,*}$ by constructing a set of approximately
$L^2$-orthonormal, $N+m+1$ approximate eigenvectors via cutting off the
first $n+1$ eigenvectors for $d_{A_0,\Phi_0}^1d_{A_0,\Phi_0}^{1,*}$ and the
first $2\kappa_p+1$ eigenvectors for each of the Dirac Laplacians
$D_{A_p,g_p}D_{A_p,g_p}^*$, $1\leq p\leq m$. The error between the
approximate eigenvectors and their $L^2$-orthogonal projections onto the
small-eigenvalue eigenspace can also be estimated. This is the content of
\S \ref{subsec:Eigen2UpperBoundSmall}. 
\item
In the reverse direction, we construct an approximately $L^2$-orthonormal
basis for $\Ker(d_{A_0,\Phi_0}^{1,*}\oplus\oplus_{p=1}^m D_{A_p,g_p}^*)$ by
cutting off the eigenvectors corresponding to the first $N$ eigenvalues of
$d_{A,\Phi}^1d_{A,\Phi}^{1,*}$. Again, we estimate the error between the
approximate eigenvectors and their $L^2$-orthogonal projections onto the
kernel. We use the resulting kernel basis to compute a lower bound for the
$(N+1)$-st eigenvalue of $d_{A,\Phi}^1d_{A,\Phi}^{1,*}$. This is the
content of \S \ref{subsec:ApproxOrthoBasisEigenvectors} and 
\S \ref{subsec:Eigen2LowerBoundFirstNonSmallEigenvalue}.
\end{itemize}

\subsection{Technical preliminaries}
\label{subsec:Eigen2TechPrelim}
In this subsection we collect some preparatory results concerning estimates
for eigenvectors, eigenvalue bounds, and decay estimates before proceeding
to the different stages of the proof of Theorem \ref{thm:H2SmallEval}.  

\subsubsection{Induced metrics over the four-sphere}
\label{subsubsec:AlmostRoundFourSphereMetric}
The background metric $g$ on $X$ induces, by cutting off near the points
$x_p\in X$, metrics on $S^4$ and $\RR^4$ which are close to standard. We
record here the precise cut-and-paste scheme, together with estimates of
their deviation from the standard metrics. (See, for example
\cite[Definition 3.11 \& Lemma 3.12]{FeehanGeometry}.) 

Define a cutoff function $\beta_p$ on $X$ so that $\beta_p=1$ on
$B(x_p,\quarter\lambda_p^{1/3})$ and $\beta_p=0$ on
$X-B(x_p,\half\lambda_p^{1/3})$.  We define families of almost flat
metrics $\delta_p$ on $\RR^4=(TX)_{x_p}$, depending on the metric $g$ on
$X$, gluing site $x_p$, and scale $\lambda_p$, by setting
\begin{equation}
\label{eq:DefnAlmostFlatMetric}
\delta_p
:=
\begin{cases}
\beta_p\delta + (1-\beta_p)\exp_{f_p}^*g
&\text{on $B(0,\lambda_p^{1/3}/4)$},
\\
\delta&\text{on $(TX)_{x_p}-B(0,\lambda_p^{1/3}/2)$},
\end{cases}
\end{equation}
where $\beta_p$ is given by definition
\eqref{eq:ChiCutoffFunctionDefn}.  Here, $\delta$ is simply the
Euclidean metric on $(TX)_{x_p}$ defined by the metric $g$ on $TX$, noting
that $(\exp_{f_p}^*g)(0) = \delta$.  We then obtain a family of conformally
equivalent metrics $g_p$ on $S^4 = (TX)_{x_p}
\cup\{\8\}$ and rescaled metrics $\tg_p$
which are close to the standard metric $g_0$ by setting
\begin{align}
\label{eq:RescaledAlmostRoundMetric}
h^*g_p(x) &= \frac{4}{(1+|x|)^2}\delta_p(x),
\quad x\in\RR^4,
\\
h^*\tg_p(x) &= \frac{4}{(1+|x|)^2}\delta_p(\lambda_p x),
\quad x\in\RR^4,
\end{align}
where $h:\RR^4\to S^4-\{s\}$ is inverse to stereographic projection. The
definition of geodesic, normal (or `Gaussian') coordinates implies that on
a small ball around the points $x_p\in X$, the metric $\delta_p$ is
$C^1$-close to $g$ and that the difference between $g$ and $\delta_p$ is
$C^2$-bounded, uniformly with respect to the ball radius \cite[Chapter
1]{Aubin}, \cite[Equations (8.8)--(8.12)]{TauSelfDual}. If we only use
normal, rather than geodesic normal coordinates, then the metric $\delta_p$
is only $C^0$ close to $g$. 

\begin{lem}
\label{lem:FourSphereMetricEstimates}
Continue the above notation. Then $\tg_p$ converges to $g_0$ in
$C^\8(S^4\less s)$ as $\lambda_p$ converges to zero.
\end{lem}
%\marginpar{\tiny Does this lemma really get used? Add proof if so}

\subsubsection{Induced connections over the four-sphere}
Similarly, connections over $X$ induce connections on $S^4$: we shall
need estimates for the self-dual component of the curvatures of these
connections and the second Chern classes of the associated Hermitian vector
bundles.

Fix an index $p\in\{1,\dots,m\}$, let $\varrho$ be the injectivity radius
of $(X,g)$, and fix a constant $\de\in(0,\quarter\varrho)$.  Choose an
orthogonal frame for $\fg_{E_\ell}|_{x_p}$ and use parallel translation via the
connection $A_0$ along radial geodesics from $x_p\in X$ to trivialize
$\fg_{E_\ell}$ over the ball $B(x_p,\varrho)$, so that
$\fg_{E_\ell}|_{B(x_p,\varrho)} \to B(x_p,\varrho)\times \su(2)$ is the
resulting smooth bundle map. Let $\Gamma$ denote the product connection on
$B(x_p,\varrho)\times \su(2)$. We have $\|F_{A_0}\|_{L^\8(X)} \le C$, for
some positive constant $C$, and so $\|F_{A_0}\|_{L^2(B(x_p,\de))} \le
C\de^2$. Thus, we may suppose that $\de$ is fixed small enough that Theorem
\ref{thm:CoulombBallGauge} provides an $\SU(2)$ gauge transformation
of $B(x_p,\de)\times \su(2)$ so that $a_0 = A_0|_{B(x_p,\de)} - \Gamma$
(after relabeling) obeys 
$$
\|a_0\|_{L^4(B(x_p,\de))} + \|\cov_\Ga a_0\|_{L^2(B(x_p,\de))}
\le c\|F_{A_0}\|_{L^2(B(x_p,2\de))}.
$$
Our hypothesis (Definition \ref{defn:LambdaClosePair}) fixed an isomorphism
$E|_{X\less\bx}\cong E_\ell|_{X\less\bx}$ and so the $A_0$-trivialization
$\fg_{E_\ell}|_{B(x_p,\varrho)} \cong B(x_p,\varrho)\times \su(2)$ induces
an isomorphism $\fg_E|_{X\less\bx}\cong (B(x_p,\de)\less\{x_p\})\times
\su(2)$. With respect to this trivialization, write $A = \Ga + a_p$ over
$B(x_p,\de)\less\{x_p\}$, where
$a_p\in\Om^1(B(x_p,\de)\less\{x_p\},\su(2))$. Let $\Omega_p =
\Omega(x_p;\quarter\lambda^{1/3},\half\lambda^{1/3})$
denote the open annulus $B(x_p,\half\lambda^{1/3})-
\barB(x_p,\quarter\lambda^{1/3})$ in $X$ and let $B_p =
B(x_p,\half\lambda^{1/3})$. Then, for any $\lambda\in (0,\lambda_0)$, the
condition \eqref{eq:APhiA0Phi0L4Close} implies that $A$ obeys
$$
\|a_p-a_0\|_{L^4(\Omega_p)}
+ \|\cov_{A_0}(a_p-a_0)\|_{L^2(\Omega_p)}
\leq 
C\lambda^{1/3}.
$$
Since $\|F_{A_0}\|_{L^\8(X)} \le C$, we have
\begin{equation}
\label{eq:a0Est}
\|a_0\|_{L^4(B_p)}
+ \|\cov_{\Gamma}a_0\|_{L^2(B_p)}
\leq 
C\lambda^{1/3},
\end{equation}
and therefore
\begin{equation}
\label{eq:aiSeqEst}
\|a_p\|_{L^4(\Omega_p)}
+ \|\cov_{\Gamma}a_p\|_{L^2(\Omega_p)}
\leq 
C\lambda^{1/3}.
\end{equation}
Let $\beta_p$ be the cutoff function of \S
\ref{subsubsec:AlmostRoundFourSphereMetric} and let $g_p$ be the
Riemannian metric on $S^4$ which coincides with $g$ on $B_p$
(after identifying the point $x_p\in X$ with the north pole $n\in S^4$) and
extends $g_p$ outside $B(x_p,\lambda_p^{1/3})=B(n,\lambda_p^{1/3})$ to
give a smooth metric on $S^4$.  Define an $\SO(3)$ bundle $F^p$ over $S^4$
by setting
$$
F^p := \begin{cases} \fg_{E}&\text{over }B_p, \\
S^4\less\{n\}\times \su(2) &\text{over }S^4\less\{n\},
\end{cases}
$$
where the identification of the $\SO(3)$ bundles $\fg_E$ and
$S^4\less\{n\}\times \su(2)$ over the annulus $B_p\less\{x_p\}$ is induced
from the $\SO(3)$ bundle isomorphism $\fg_E|_{B_p\less\{n\}}\cong
B_p\less\{n\}\times\su(2)$ described in the preceding paragraphs.  
We cut off the connection $A$ on $\fg_E$ over
the annulus $\Om_p$ and thus obtain an $L^2_k$ connection $A_p$ on the
$\SO(3)$ bundle $F^p$ over $S^4$ by setting
$$
A_p
:= 
\begin{cases}
A&\text{on }\fg_E|_{B(n,\quarter\lambda_p^{1/3})}, 
\\
\Ga + \beta_pa_p 
&\text{on }S^4\less\{n\}\times \su(2).
\end{cases}
$$
Our hypothesis \eqref{eq:FA+SphereBallL2Small} on the pair
$(A,\Phi)$ implies that $A$ obeys
\begin{equation}
\label{eq:F+SeqEst}
\|F^+_A\|_{L^2(B_p)} \le C\lambda_p^{1/3}.
\end{equation}
Since
$$
F^+_{A_p}
= \beta_pF^+_A + (d\beta_p\wedge a_p)^+ 
 + \beta_p(\beta_p-1)(a_p\wedge a_p)^+,
$$
we see that
$$
\|F^+_{A_p}\|_{L^2(S^4)}
\le 
\|F^+_A\|_{L^2(B_p)} 
+ \sqrt{2}\|d\beta_p\|_{L^4(X)}\|a_p\|_{L^4(\Omega_p)} 
+ \sqrt{2}\|a_p\|_{L^4(\Omega_p)}^2.
$$
and therefore the estimates \eqref{eq:a0Est}, \eqref{eq:aiSeqEst}, and
\eqref{eq:F+SeqEst} imply that
\begin{equation}
\label{eq:SDBubbleCurvLimit}
\|F^+_{A_p}\|_{L^2(S^4)}
\le 
C\lambda_p^{1/3}.
\end{equation}
Similarly, as
$$
F_{A_p}
= \beta_pF_A + d\beta_p\wedge a_p
 + \beta_p(\beta_p-1)a_p\wedge a_p,
$$
we see that
$$
\|F_{A_p}-F_A\|_{L^2(\Omega_p)} 
\le 
\|\cov_\Ga a_p\|_{L^2(\Omega_p)}
+ c\|a_p\|_{L^4(\Omega_p)}
+ \|a_p\|_{L^4(\Omega_p)}^2,
$$
and therefore, the estimates \eqref{eq:a0Est} and \eqref{eq:aiSeqEst} yield
\begin{equation}
\label{eq:BubbleCurvDiffLimit}
\|F_{A_p}-F_A\|_{L^2(\Omega_p)} 
\le 
C\lambda_p^{1/3}.
\end{equation}
Since by $\|F_A\|_{L^2(B_p)}^2 = 8\pi^2\kappa_p$ by condition
\eqref{eq:FASphereBallIntegerMass}, while $A_p = A$ on $B_p-\Omega_p$ and
$A_p = \Gamma$ on $\RR^4-B_p$, we see that the error estimate
\eqref{eq:BubbleCurvDiffLimit} and the Chern-Weil formula imply that
$-\quarter p_1(F^p) = \kappa_p$. In particular, we can choose a Hermitian
vector bundle $E^p$ over $S^4$ with
\begin{equation}
\label{eq:BubbleCharge}
\fg_{E^p} = F^p \quad\text{and}\quad \det E^p = S^4\times\CC,
\end{equation}
and $c_2(E^p) = \kappa_p$.

\subsubsection{Estimates for small-eigenvalue eigenvectors of the
$\PU(2)$-monopole Laplacian}
Recall from \cite[Equation (2.36)]{FL1} that the first-order operator 
$$
d_{A,\Phi}^{1,*}:L^2_{k+1}(\Lambda^+\otimes\fg_E)\oplus L^2_{k+1}(V^-)
\to L^2_k(\Lambda^+\otimes\fg_E)\oplus L^2_k(V^+)
$$
is the formal adjoint of $d_{A,\Phi}^1$ and is given by
\begin{equation}
\label{eq:LAPhi*}
d_{A,\Phi}^{1,*}\left(\begin{matrix} v \\ \psi \end{matrix}\right) =
\left(\begin{matrix}
d_{A}^*v + (\rho(\cdot)\Phi)^*\psi \\
D_{A,\vartheta}^*\psi 
-(\Phi\otimes(\cdot)^*+(\cdot)\otimes{\Phi}^*)_{00}^*(\tau\rho^{-1})^*v
\end{matrix}\right). 
\end{equation}
The application of cut and paste techniques to the problem of estimating
eigenvalues and eigenvectors is simplest when we have $C^0\cap L^2_{2,A_0}$
estimates for eigenvectors of the Laplacian
$d_{A_0,\Phi_0}^1d_{A_0,\Phi_0}^{1,*}$: this is possible since the pair
$[A_0,\Phi_0]$ varies in a family with strong Sobolev-norm bounds on the
curvature $F_{A_0}$ and the spinor $\Phi_0$. By retracing the proof of
Lemma 5.9 in
\cite{FeehanSlice}, we can obtain an estimate with the mildest possible
bounds on $(F_{A_0},\Phi_0)$ --- though such optimality is not strictly
necessary. That {\em some\/} estimate should exist, with constant a polynomial in
$\|F_{A_0}\|_{L^2_{l,A_0}(X)}$ and $\|\Phi_0\|_{L^2_{l+1,A_0}(X)}$, is a
consequence of standard elliptic theory and Uhlenbeck's local Coulomb-gauge
fixing theorem \cite{UhlLp}. Then, by reworking the proof of \cite[Lemma
6.7]{FLKM1}, we easily obtain:

\begin{lem}
\label{lem:LInftyEstDe2A0Phi0Evec}
Let $X$ be a $C^\8$, closed, oriented, Riemannian four-manifold.  Then
there is a constant $c$ with the following significance.
Let $(A,\Phi)$ be an $L^2_4$ pair on $(\fg_E,V)$
over $X$. Let $(\eta,\psi)\in L^2_3(X,\Lambda^+\otimes\fg_E)\oplus
L^2_3(X,V^-)$ be an eigenvector of $d_{A,\Phi}^1d_{A,\Phi}^{1,*}$ with
eigenvalue $\mu[A,\Phi]$. Then
$$
\|(\eta,\psi)\|_{L^\8(X)} + \|(\eta,\psi)\|_{L^2_{2,A}(X)}  
\le C[A,\Phi]\|\eta\|_{L^2(X)},
$$
where $C[A,\Phi] := C([A,\Phi],\mu[A,\Phi])$ with $\mu=\mu[A,\Phi]$ in
$$
C([A,\Phi],\mu)
:=
c(1+\|F_A\|_{L^\8(X)}+\|\Phi\|_{L^\8_{1,A}})^4(1+\mu)(1+\sqrt{\mu}).
$$
More generally, if $(\eta,\psi)$ lies
in the span of the eigenvectors of $d_{A,\Phi}^1d_{A,\Phi}^{1,*}$ with
eigenvalues less than or equal to a positive constant $\mu$, then
$$
\|(\eta,\psi)\|_{L^\8(X)} + \|(\eta,\psi)\|_{L^2_{2,A}(X)}
\le 
C([A,\Phi],\mu)\|(\eta,\psi)\|_{L^2(X)}.
$$
\end{lem}

We shall also need uniform $L^2_{1,A}(X)$ estimates for $\eta$ and
$L^2_{1,A}(U)$ estimates for $\psi$, when $(\eta,\psi)$ is an eigenvector
of $d_{A,\Phi}^1d_{A,\Phi}^{1,*}$ and $U\Subset X\less\bx$, but in the more
difficult case where we cannot assume $[A,\Phi]$ varies in a family for
which $\|F_A\|_{L^2_{l,A}(X)}$ is uniformly bounded when $l>0$.

\begin{lem}
\label{lem:L21GlobalLocalEstDe2APhiEvec}
Let $X$ be a $C^\8$, closed, oriented, Riemannian four-manifold.  Then
there is a positive constant $c$ with the following significance. Suppose
$[A,\Phi]$ obeys the conditions of Definition \ref{defn:LambdaClosePair}.
Let $(\eta,\psi)\in 
L^2_{k-1}(X,\Lambda^+\otimes\fg_E)\oplus L^2_{k-1}(X,V)$ be an eigenvector
of $d_{A,\Phi}^1d_{A,\Phi}^{1,*}$ with eigenvalue $\mu[A,\Phi]$. Then
\begin{equation}
\label{eq:L21GlobalLocalEstDe2APhiEvec}
\|\eta\|_{L^2_{1,A}(X)} + \|\psi\|_{L^2_{2,A}(U)}  
\le 
C[A,\Phi]\|(\eta,\psi)\|_{L^2(X)},
\end{equation}
$U = X-\cup_{p=1}^m B(x_p,c_p\lambda_p^{1/3})$, 
where the $c_p$ are universal positive constants, $C([A,\Phi],\mu)$ is a
polynomial in $\sqrt{\mu}$, where $\mu=\mu[A,\Phi]$, $\|F_A\|_{L^2(X)}$,
$\|F_{A_0}\|_{L^2_{k-1,A_0}(X)}$, and $\|\Phi_0\|_{L^2_{k-1,A_0}(X)}$.
More generally, if $(\eta,\psi)$ lies in the span of the eigenvectors of
$d_{A,\Phi}^1d_{A,\Phi}^{1,*}$ with eigenvalues less than or equal to a
positive constant $\mu$, then
\begin{equation}
\label{eq:L21GlobalLocalEstDe2APhiEProj}
\|\eta\|_{L^2_{1,A}(X)} + \|\psi\|_{L^2_{2,A}(U)}  
\le 
C([A,\Phi],\mu)\|(\eta,\psi)\|_{L^2(X)}.
\end{equation}
\end{lem}

\begin{proof}
For convenience, write $(a,\phi) = d_{A,\Phi}^{1,*}(\eta,\psi)$. First
consider the $L^2_{1,A}(X)$ estimate for the component $\eta$.  Since
$a=d_A^{+,*}\eta + (\cdot\Phi)^*\psi$ from equation \eqref{eq:LAPhi*}, Lemma
\ref{lem:L21AEstv} implies that
\begin{align*}
\|\eta\|_{L^2_{1,A}} 
&\le 
C(\|d_A^{+,*}\eta\|_{L^2} + \|\eta\|_{L^2}) 
\quad\text{(by Lemma \ref{lem:L21AEstv})}
\\
&\le 
C(\|a\|_{L^2} + \|\Phi\|_{L^\8}\|\psi\|_{L^2} + \|\eta\|_{L^2})
\quad\text{(by Equation \eqref{eq:LAPhi*})}
\\
&\le 
C(\|d_{A,\Phi}^{1,*}(\eta,\psi)\|_{L^2} + \|(\eta,\psi)\|_{L^2})
\\
&\leq C\|(\eta,\psi)\|_{L^2},
\end{align*}
where we use $\|d_{A,\Phi}^{1,*}(\eta,\psi)\|_{L^2} 
= \sqrt{\mu[A,\Phi]}\|(\eta,\psi)\|_{L^2}$ to obtain the last line above.

Second, we consider the $L^2_{1,A}(U)$ estimate for the component $\psi$.
{}From the proof of Theorem \ref{thm:L21AEstPAaphi} we have the estimates
\eqref{eq:Lp1psiU} and \eqref{eq:L21phi-Final}, 
\begin{align*}
\|\psi\|_{L^2_{2,A}(U)} 
&\leq
C(\|D_A\phi\|_{L^2} + \|\eta\|_{L^2_{1,A}} + \|\psi\|_{L^2}),
\\
\|\phi\|_{L^2_{1,A}} 
&\le 
C(\|d_{A,\Phi}^1(a,\phi)\|_{L^2} + \|(a,\phi)\|_{L^2}).
\end{align*}
(Note that as we only seek an estimate for $\psi$ on the subset $U\Subset
X\less\bx$, we do not necessarily require that $\Phi\equiv 0$ on the balls
$B(x_p,c_p\lambda_p^{1/3})$, as we do in the proof of Theorem
\ref{thm:L21AEstPAaphi}, where estimates for solutions to the linearized,
extended $\PU(2)$-monopole equations are needed over all of $X$.)   
Combining the preceding two estimates, noting that
$d_{A,\Phi}^1d_{A,\Phi}^{1,*}(\eta,\psi)
= \mu[A,\Phi]\cdot(\eta,\psi)$, and applying our $L^2_{1,A}(X)$
estimate for $\eta$ gives
\begin{align*}
\|\psi\|_{L^2_{2,A}(U)} 
&\leq
C(\|d_{A,\Phi}^1d_{A,\Phi}^{1,*}(\eta,\psi)\|_{L^2}
+ \|d_{A,\Phi}^{1,*}(\eta,\psi)\|_{L^2} + \|(\eta,\psi)\|_{L^2})
\\
&=
C\|(\eta,\psi)\|_{L^2}.
\end{align*}
This completes the proof. 
\end{proof}

\subsubsection{A positive lower bound for the first eigenvalue of the
anti-self-dual Laplacian on the four-sphere}
We shall need a lower bound for the least positive eigenvalue of the
Laplacian $d_A^+d_A^{+,*}$ on $L^2(S^4,\Lambda^+\otimes\fg_E)$. If $A$ is
an anti-self-dual connection with respect to the standard round metric of
radius one on $S^4$ (and thus scalar curvature $12$), then the Bochner
formula \eqref{eq:BW+} for $d_A^+d_A^{+,*}$ gives
$$
d_A^+d_A^{+,*}
=
\thalf\cov_A^*\cov_A + 2.
$$
More generally, the method of proof of Lemma
\ref{lem:DiracS4VanishingKernel} yields the following estimate for the
least eigenvalue of $d_A^+d_A^{+,*}$ on $S^4$.

\begin{lem}
\label{lem:ASDLaplaceS4VanishingKernel}
Let $K_1$ be the norm of the Sobolev embedding $L^2_1(S^4)\to
L^4(S^4)$. Then for all $0 <\eps < K_1/100$ the following holds. Let $S^4$
have a metric $g'$ which is $L^2_2$-close to the standard round metric $g$
of radius one, so $\|g'-g\|_{L^2_2(S^4,g)}<\eps/100$. Let $E$ be a
Hermitian, rank-two vector bundle over $S^4$ with $c_2(E) \geq 0$ and
orthogonal connection $A$ on $\fg_E$. If $\|F_A^+\|_{L^2(S^4)} <
\eps/100K_1$ then $\Ker d_A^{+,*} = \{0\}$ and the least positive
eigenvalue of $d_A^+d_A^{+,*}$ obeys the lower bound
$$
\mu_1(S^4;d_A^+d_A^{+,*}) \geq 2 - \eps.
$$
\end{lem}

\begin{proof}
First suppose $g'=g$.  The equality in the assertion of the lemma follows
from the fact that the scalar curvature of the round metric on the
$n$-sphere $S^n$ of radius $1$ (and thus constant sectional curvature $1$)
is equal to $n(n-1)$ by
\cite[pp. 37 \& 60]{Chavel} and so the Bochner formula  \eqref{eq:BW+} for
the Laplacian gives
$$
d_A^+d_A^{+,*}
=
\thalf\cov_A^*\cov_A + 2 + \thalf\{F_A^+,\cdot\},
$$
Hence, integrating by parts and applying H\"older's inequality, we see that
$$
\|d_A^{+,*}v\|_{L^2}^2
\geq
\thalf\|\cov_Av\|_{L^2}^2 + 2\|v\|_{L^2}^2 -
\thalf\|F_A^+\|_{L^2}\|v\|_{L^4}^2.
$$
By Lemma \ref{lem:Kato} we have $\|v\|_{L^4} \leq
K_1\|v\|_{L^2_{1,A}}$, where $K_1$ is a universal constant, so
\begin{align*}
\|d_A^{+,*}v\|_{L^2}^2
&\geq
\thalf(1-K_1\|F_A^+\|_{L^2})\|\cov_Av\|_{L^2}^2 
+ (2-\thalf K_1\|F_A^+\|_{L^2})\|v\|_{L^2}^2
\\
&\geq
(2-\thalf K_1)\|F_A^+\|_{L^2})\|v\|_{L^2}^2,
\end{align*}
where the factor in the last inequality is positive if we choose $\eps <
4/K_1$.  This leads to a lower bound for the first eigenvalue
$\mu_1(S^4;d_A^+d_A^{+,*})$,
$$
\mu_1(S^4;d_A^+d_A^{+,*})
=
\inf_{\|v\|_{L^2}=1}
(d_A^+d_A^{+,*}v,v)_{L^2}
\geq
2-\thalf K_1\|F_A^+\|_{L^2} > 0,
$$
and $\Ker d_A^{+,*} = \{0\}$, as desired. The argument is very similar if the
metric on $S^4$ is only $L^2_2$-close to the standard metric, with a slightly
smaller choice of $\eps$, as the Riemann curvature tensor of the
Levi-Civita connection for $g'$ is given by $\Rm_{g'} = \Rm_g +
\Rm_{g'}-\Rm_g$ with $\|\Rm_{g'}-\Rm_g\|_{L^2_2(S^4,g)} \leq \eps/100K_1$.
\end{proof}

\subsubsection{Small eigenvalues of the Dirac Laplacian}
We shall also need estimates for the small eigenvalues of the Dirac
Laplacian $D_AD_A^*$ on $L^2(X,V^-)$. As the spectra of $D_AD_A^*$ and
$D_A^*D_A$ coincide on the complements of their kernels, it suffices to
consider the small (but non-zero) eigenvalues of $D_A^*D_A$ --- for which
we have a useful Bochner formula \eqref{eq:BWDirac+} involving only the
component $F_A^+$ of $F_A$ --- rather than $D_AD_A^*$, whose Bochner formula
\eqref{eq:BWDirac-} involves the component $F_A^-$ of $F_A$. For this case,
the proof of Theorem 6.1 in \cite{FLKM1}, which gives estimates for the
small eigenvalues of $d_A^+d_A^{+,*}$ on
$L^2_{k-1}(X,\Lambda^+\otimes\fg_E)$, yields 

\begin{thm}
\label{thm:DiracSmallEval}
Let $M$ be a positive constant and let $(X,g)$ be a closed, oriented,
Riemannian four-manifold.  Then there are positive constants $c=c(g)$,
$C = C(g,M)$, and small enough $\lambda_0 = \lambda_0(C)$
such that for all $\lambda \in (0,\lambda_0]$, the following holds.  Let
$E$, $E_\ell$ be Hermitian, rank-two vector bundles over $X$ with $\det E =
\det E_\ell$, $c_2(E)=c_2(E_\ell)+\ell$, for some integer $\ell\geq 0$, and
suppose $c_2(E)\leq M$. Suppose $\bx\in\Sym^\ell(X)$ and define a
precompact open subset of $X\less\bx$ by setting $U:= X-\cup_{x\in\bx}\bar
B(x,4\sqrt{\lambda})$.  Let $A$ be a connection representing a point in the
open subset of $\sB_E(X)$ defined by the following constraints:
%\marginpar{\tiny Check KM1 reference and result}
\begin{itemize}
\item $\|A-A_0\|_{L^2(U)} < c\sqrt{\lambda}$ and
$\|A-A_0\|_{L^2_{1,A_0}(U)} < M$, (noting that
$E|_{X\less\bx} \cong E_\ell|_{X\less\bx}$).
\end{itemize}
Let $\nu_2[A_0]$ be the {\em least positive eigenvalue\/} of the Laplacian
$D_{A_0}^*D_{A_0}$ on $L^2(X,V_\ell^+)$ and suppose
$\max\{\nu_2[A_0],\nu_2[A_0]^{-1}\} \leq M$. Assume that the kernel of
$D_{A_0}^*D_{A_0}$ has dimension less than or equal to $M$. Let
$\{\mu_i^+[A]\}_{i\geq 1}$ denote the {\em positive\/} eigenvalues,
repeated according to their multiplicity, of the Laplacian $D_A^*D_A$ on
$L^2(X,V^+)$. Then
$$
\dim_\CC\Ker D_A \leq \dim_\CC\Ker D_{A_0},
$$
and, furthermore,
\begin{equation}
\begin{aligned}
\mu_i^+[A] 
&\leq
\begin{cases}
C\lambda, &\text{if }1\leq i\leq \dim_\CC\Ker D_{A_0},
\\
\nu_2[A_0] + C\sqrt{\lambda}, &\text{if }i=\dim_\CC\Ker D_{A_0}+1,
\end{cases} 
\\
\mu_i^+[A] 
&\geq 
\nu_2[A_0]- C\sqrt{\la}, 
\quad \text{if }i\geq \dim_\CC\Ker D_{A_0}+1.
\end{aligned}
\end{equation}
\end{thm}

\begin{cor}
\label{cor:DiracSmallEval}
Continue the hypotheses and notation of Theorem \ref{thm:DiracSmallEval}.
Let $\{\mu_i^-[A]\}_{i\geq 1}$ denote the {\em positive\/} eigenvalues,
repeated according to their multiplicity, of the Laplacian $D_AD_A^*$ on
$L^2(X,V^-)$. Then, 
$$
\dim_\CC\Ker D_A^* \leq \dim_\CC\Ker D_{A_0}^* + \ell,
$$
where $\ell = c_2(E)-c_2(E_\ell)\geq 0$ and, furthermore, 
\begin{align*}
\mu_i^-[A] 
&\leq
\begin{cases}
C\lambda, &\text{if }1\leq i\leq\dim_\CC\Ker D_{A_0},
\\
\nu_2[A_0] + C\sqrt{\lambda}, &\text{if }i=\dim_\CC\Ker D_{A_0}+1,
\end{cases} 
\\
\mu_i^-[A] &\geq \nu_2[A_0]- C\sqrt{\la},  
\quad \text{if }i\geq \dim_\CC\Ker D_{A_0}+1.
\end{align*}
\end{cor}

\begin{proof}
We have the bound on $\dim_\CC\Ker D_A^*$ because $\dim_\CC\Ker D_A^* =
\dim_\CC\Ker D_A - \Ind_\CC D_A$ and as, due to excision \cite[\S 7.1]{DK}
and Lemma \ref{lem:DiracS4VanishingKernel}, we have $\Ind_\CC D_A =
\Ind_\CC D_{A_0} - \ell$. By the remarks preceding the corollary, we have
$\mu_i^-[A] = \mu_i^+[A]$ for all $i\geq 1$, since these eigenvalues are
positive, so the result follows from Theorem \ref{thm:DiracSmallEval}.
\end{proof}

\begin{rmk}
The small-eigenvalue estimates in Theorem
\ref{thm:DiracSmallEval} and Corollary \ref{thm:DiracSmallEval} are
required by the proof of Theorem \ref{thm:L21AEstPAaphi}. We need the
spectral splitting in order to apply Theorem 1.1 in
\cite{FeehanKato} and estimate negative spinors on $X$.
\end{rmk}

\subsubsection{A lower bound for the first positive eigenvalue of the Dirac
Laplacian on the four-sphere}
We recall from the proof of Proposition 2.28 in \cite{FL1} that the complex
index of the Dirac operator $D_A:\Gamma(S^+\otimes E)\to
\Gamma(S^-\otimes E)$, where $(\rho,S^+,S^-)$ is the standard spin
structure on $S^4$ and $E$ is a Hermitian, rank-two vector bundle with
$c_2(E)=l=-\quarter p_1(\fg_E)\geq 0$, is given by  
\begin{equation}
\label{eq:IndexDiracS4}
\Ind_\CC D_A 
= 
\quarter\left(p_1(\fg_E) + (c_1(W^+)+c_1(E))^2 - \sigma(S^4)\right)
= -l.
\end{equation}
Hence, assuming $\Ker D_A = \{0\}$, we have $\dim_\RR\Ker D_A^* = 2l$; the
vanishing condition is guaranteed by the following lemma.

\begin{lem}
\label{lem:DiracS4VanishingKernel}
Let $K_1$ be the norm of the Sobolev embedding $L^2_1(S^4)\to L^4(S^4)$. Then
for all $0 <\eps < K_1/100$ the following holds. Let $S^4$ have a metric $g'$
which is $L^2_2$-close to the standard metric $g$, so
$\|g'-g\|_{L^2_2(S^4,g)}<\eps/100$. Let $(\rho,S^+,S^-)$ be the standard
spin structure on $S^4$ \cite[\S 6.4]{SalamonSWBook} and let $E$ be a
Hermitian, rank-two vector bundle over $S^4$ with $c_2(E) \geq 0$ and
orthogonal connection $A$ on $\fg_E$. Let $A_d$, $A_e$ denote the unitary
connections on the (trivial) lines bundles $\det S^+$, $\det E$. If
$$
\|F_A^+\|_{L^2(S^4)} + \|F_{A_d}^+\|_{L^2(S^4)} + \|F_{A_e}^+\|_{L^2(S^4)} 
< 
\eps/100K_1, 
$$
then $\Ker D_A = \{0\}$, $\Ker D_A^*$ has
real dimension $2c_2(E)$, and the least positive eigenvalues of $D_A^*D_A$
and $D_AD_A^*$ obey the lower bound
$$
\mu_1(S^4;D_AD_A^*) = \mu_1(S^4;D_A^*D_A) \geq 3 - \eps.
$$
\end{lem}

\begin{proof}
The equality in the assertion of the lemma follows from the fact that aside
from zero eigenvalues, the spectra of the Laplacians $D_AD_A^*$ and
$D_A^*D_A$ coincide \cite[Lemma 1.6.5]{Gilkey}. Thus we need only consider
the least positive eigenvalue of $D_A^*D_A$. For clarity of exposition,
assume first that $g'=g$ and that $F_{A_d}^+ = F_{A_e}^+ = 0$.
 
The scalar curvature of the round metric on the $n$-sphere $S^n$ of radius
$1$ (and thus constant sectional curvature $1$) is equal to $n(n-1)$ by
\cite[pp. 37 \& 60]{Chavel} and so the Bochner formula for the Laplacian
$D_A^*D_A$ \cite[Lemma 4.1(a)]{FL1} gives
$$
D_A^*D_A = \cov_A^*\cov_A + 3 + \rho(F_A^+),
$$
where $A$ is viewed as a unitary connection on $E$ (rather than an
orthogonal connection on $\fg_E$). Hence, integrating by parts and
applying H\"older's inequality, we see
that
$$
\|D_A\phi\|_{L^2}^2
\geq
\|\cov_A\phi\|_{L^2}^2 + 3\|\phi\|_{L^2}^2 -
\|F_A^+\|_{L^2}\|\phi\|_{L^4}^2.
$$
By Lemma \ref{lem:Kato} we have $\|\phi\|_{L^4} \leq
K_1\|\phi\|_{L^2_{1,A}}$, where $K_1$ is a universal constant, so
\begin{align*}
\|D_A\phi\|_{L^2}^2
&\geq
(1-K_1\|F_A^+\|_{L^2})\|\cov_A\phi\|_{L^2}^2 
+ (3-K_1\|F_A^+\|_{L^2})\|\phi\|_{L^2}^2
\\
&\geq
(3-K_1\|F_A^+\|_{L^2})\|\phi\|_{L^2}^2,
\end{align*}
where the last inequality follows if we choose $\eps\leq 1/K_1$.
This leads to a lower bound for the first eigenvalue $\mu_1(S^4;D_A^*D_A)$,
$$
\mu_1(S^4;D_A^*D_A)
=
\inf_{\|\phi\|_{L^2}=1}
(D_A^*D_A\phi,\phi)_{L^2}
\geq
3-K_1\|F_A^+\|_{L^2} > 0,
$$
and $\Ker D_A = \{0\}$, as desired. The argument is very similar if the metric on
$S^4$ is only $L^2_2$-close to the standard metric or if $F_{A_d}^+\neq
0$, $F_{A_e}^+\neq 0$, with a slightly smaller choice of $\eps$.
\end{proof}
 
\subsubsection{Pointwise decay estimates for harmonic spinors on the
four-sphere} 
It will also be useful to estimate the effect of cutting off harmonic
spinors on the complement of small balls around the north pole where the
curvature of the connection $A_p$ is concentrated.
%\marginpar{\tiny stereographic projection: have also used $\varphi_{n,s}$?}

\begin{lem}
\label{lem:S4HarmonicSpinorDecay}
Let $K_2\geq 1$ be a constant.  Then there are constants $\eps(K_2)>0$ and, if
$M>0$ is constant, a constant $c(K_2,M)$ with the following significance.
Let $S^4$ have a metric $g$ equivalent to the standard round metric of
radius one $g_1$, so $K_2^{-1}g_1 \leq g \leq K_2g_1$.  Let $h:\RR^4\cong
S^4\less\{s\}$ be inverse stereographic projection, and let $B(r_0)$ denote
the ball $\{x\in
\RR^4:r<r_0\}$, where $r=|x|$. Let $A$ be an $L^2_k$ unitary connection on
a Hermitian bundle $E$ over 
$S^4$, let $(\rho,S^+,S^-)$ be the standard spin structure over $S^4$, and
let $V=V^+\oplus V^-$ with $V^\pm = S^\pm\otimes E$. Suppose $A$ obeys
$$
\|F_A\|_{L^\8(S^4-h(B(1)),g)} \leq M
\text{ and }
\|F_A\|_{L^2(S^4-h(B(r_0)),g)} \leq \eps.
$$
Let $D_{A,g}:\Gamma(S^4,V)\to\Gamma(S^4,V)$ be the corresponding Dirac
operator. Suppose $\psi\in L^2(S^4,V)$ obeys $D_{A,g}\psi = 0$. Then for
all $x\in \RR^4-B(2r_0)$ we have
\begin{equation}
\label{eq:S4HarmonicSpinorDecay}
|h^*\psi|(x) 
\leq 
c\eps^{3/2}r_0(1+r^{-3})\|\psi\|_{L^2(S^4-h(B(r_0)),g)}.
\end{equation}
Furthermore, if $2r_0 \leq r_2 < \8$, 
\begin{equation}
\label{eq:S4HarmonicSpinorL2Decay}
\|\psi\|_{L^2(S^4-h(B(r_2)),g)} 
\leq 
\begin{cases}
c\eps^{3/2}r_0r_2^{-1}\|\psi\|_{L^2(S^4-h(B(r_0)),g)} &\text{ if }r_2 \leq 1,
\\
c\eps^{3/2}r_0r_2^{-2}\|\psi\|_{L^2(S^4-h(B(r_0)),g)} &\text{ if }r_2 \geq 1.
\end{cases}
\end{equation}
\end{lem}

\begin{proof}
We may suppose without loss that $g$ is the standard round metric of radius
one, so $(h^*g)_x = f^2|dx|^2$, where $f(x) = 2/(1+|x|^2)$,
$x\in\RR^4$, and $\delta_x = |dx|^2$ is the standard metric on $\RR^4$. 
%\marginpar{\tiny Make sure \cite{FeehanKato} can handle approx round metrics}
Let $\tilde\psi = f^{3/2}\cdot h^*\psi$ and recall that $D_{A,\delta} =
f^{-5/2}\circ D_{A,g} \circ f^{3/2}$ \cite[Remark 4.8]{FL1}, 
\cite[Theorem II.5.4]{LM}, and so
$D_{A,\delta}\tilde\psi = 0$. Theorem 1.1 in
\cite{FeehanKato} (with $r_1=\8$) implies that for all $x \in \RR^4-B(2r_0)$,
we have
\begin{equation}
\label{eq:R4HarmonicSpinorDecay}
|\tilde\psi|(x) 
\leq 
c\eps^{3/2}r_0r^{-3}\|\tilde\psi\|_{L^2(\RR^4-B(r_0),\delta)}.
\end{equation}
Our first task is to show that $\|\tilde\psi\|_{L^2(\RR^4-B(r_0),\delta)}$
is bounded by $c\|\psi\|_{L^2(S^4-h(B(r_0)),g)}$.

For $y\in\RR^4 = S^4\less\{n\}$ (say), we recall that $|\partial/\partial
r|_\delta = 1$, while $|\partial/\partial r|_g = 2(1+r^2)^{-1}$. Thus,
$$
\dist_g(h(y),s)
=
\int_0^{|y|} 2(1+t^2)^{-1} \,dt
=
\arctan(|y|)
\leq 
2|y|.
$$
Suppose $r_1>4r_0$, but fixed, and consider the region $\RR^4-B(r_1)$.
Observe that, as $dV_g = \sqrt{\det g}\,d^4x = f^4\,d^4x$ and, when
$|x|\geq 1$ and setting $y=x^{-1}$ (identifying $\RR^4$ with the
quaternions), we have $f(x)^{-1} = \half(1+|x|^2) \leq |x|^2 = |y|^{-2} \leq
\quarter \dist_g(h(y),s)^{-2}$, so 
\begin{align*}
\|\tilde\psi\|_{L^2(\RR^4-B(r_1),\delta)}^2
&=
\int_{\RR^4-B(r_1)}|h^*\psi|^2 f^3\,d^4x
\\
&=
\int_{\RR^4-B(r_1/2)}|h^*\psi|^2 f^{-1}\sqrt{\det g}\,d^4x
\\
&\leq
\sup_{p\in S^4}\int_{S^4-h(B(r_1))}c\dist_g(\cdot,p)^{-2}|\psi|^2\,dV_g
\\
&=
\|\psi\|_{L^{2\sharp}(S^4-h(B(r_1)),g)}^2.
\end{align*}
Hence, the embedding $L^2_1\subset L^{2\sharp}$ \cite[Lemma
4.1]{FeehanSlice} implies that
\begin{equation}
\label{eq:R4LessBallHarmonicSpinorL2Est}
\begin{aligned}
\|\tilde\psi\|_{L^2(\RR^4-B(r_1),\delta)}
&\leq
c\|\psi\|_{L^{2\sharp}(S^4-h(B(r_1)),g)}
\\
&\leq
c\|\psi\|_{L^2_{1,A}(S^4-h(B(r_1/2)),g)}.
\end{aligned}
\end{equation}
On $S^4-h(B(r_1))$, where $\|F_A\|_{L^\8(S^4-h(B(r_1)),g)} \leq M$ by
hypothesis, we have the uniform elliptic estimate
\begin{equation}
\label{eq:S4HarmonicSpinorL21EllipticEst}
\|\psi\|_{L^2_{1,A}(S^4-h(B(r_1/2)),g)}
\leq
c(r_1)\|\psi\|_{L^2(S^4-h(B(r_1/4)),g)},
\end{equation}
via the Bochner formula for $D_{A,g}^2$, 
because $D_{A,g}\psi=0$ on $S^4$. Combining the estimates
\eqref{eq:R4LessBallHarmonicSpinorL2Est} and
\eqref{eq:S4HarmonicSpinorL21EllipticEst} gives
\begin{equation}
\label{eq:CompBigBallR4S4HarmonicSpinorL2Est}
\|\tilde\psi\|_{L^2(\RR^4-B(r_1),\delta)}
\leq
c(r_1)\|\psi\|_{L^2(S^4-h(B(r_1/4)),g)}.
\end{equation}
Now consider the region $\Omega(r_0,r_1)$.  The harmonic spinor $\psi$
on $S^4$ obeys, for universal and fixed $r_1\leq 1$ (without loss), and all
$r\in [r_0,r_1]$, the estimate
$$
|\tilde\psi|(x)
= 
f(x)^{3/2}|h^*\psi|(x) 
\leq 
(2/(1+r_1^2))^{3/2}|h^*\psi|(x) 
\leq
2\sqrt{2}|h^*\psi|(x).
$$
In particular, using $dV_g = f^4\,d^4x$, we have
\begin{align*}
\|\tilde\psi\|_{L^2(\Omega(r_0,r_1),\delta)}^2
&=
\int_{\Omega(r_0,r_1)}|\tilde\psi|^2\,d^4x
\\
&=
\int_{\Omega(r_0,r_1)}|f^{-3/2}\cdot\tilde\psi|^2 f^3\,dx
\\
&=
\int_{\Omega(r_0,r_1)}f^{-1}|h^*\psi|^2f^4\,d^4x
\\
&=
\half\int_{\Omega(r_0,r_1)}(1+r^2)|h^*\psi|^2\,dV_g
\\
&\leq
\|\psi\|_{L^2(h(\Omega(r_0,r_1)),g)}^2,
\end{align*}
and thus, for all $0\leq r_0< r_1\leq 1$,
\begin{equation}
\label{eq:AnnulusR4S4HarmonicSpinorL2Est}
\|\tilde\psi\|_{L^2(\Omega(r_0,r_1),\delta)}
\leq
\|\psi\|_{L^2(h(\Omega(r_0,r_1)),g)}.
\end{equation}
Combining our $L^2(\RR^4,\delta)$ estimate
\eqref{eq:CompBigBallR4S4HarmonicSpinorL2Est} for $\tilde\psi$ over
$\RR^4-B(r_1)$ and our estimate \eqref{eq:AnnulusR4S4HarmonicSpinorL2Est}
over $\Omega(r_0,r_1)$ yields
\begin{equation}
\label{eq:CompSmallBallR4S4HarmonicSpinorL2Est}
\|\tilde\psi\|_{L^2(\RR^4-B(r_0),\delta)}
\leq
c\|\psi\|_{L^2(S^4-h(B(r_0)),g)}.
\end{equation}
Hence, on $S^4$, using $h^*\psi = f^{-3/2}\cdot\tilde\psi$ and $r=|x| \in
[2r_0,\8)$ and our decay estimate \eqref{eq:R4HarmonicSpinorDecay} for
$\tilde\psi$ on $\RR^4-B(2r_0)$, we obtain a decay estimate for $\psi$ on
$S^4-h(B(2r_0))$,
\begin{align*}
|h^*\psi|(x) 
&\leq 
c\eps^{3/2}r_0r^{-3}f^{-3/2}\|\psi\|_{L^2(S^4-h(B(r_0)),g)}
\\
&\leq 
c\eps^{3/2}r_0r^{-3}(1+r^2)^{3/2}\|\psi\|_{L^2(S^4-h(B(r_0)),g)}.
\end{align*}
This yields the pointwise bound \eqref{eq:S4HarmonicSpinorDecay}.

Hence, for all $r_2 > 2r_0$ and setting $E =
\|\psi\|_{L^2(S^4-h(B(r_0)),g)}$ for convenience, we have 
\begin{align*}
\|\psi\|_{L^2(S^4-h(B(r_2)),g)}^2
&=
\int_{\RR^4-B(r_2)}|h^*\psi|^2 f^4\,d^4x
\\
&\leq
c\eps^3r_0^2\int_{r_2}^\8 r^{-6}f^{-3}E^2f^4 r^3\,dr
\\
&=
c\eps^3r_0^2E^2\int_{r_2}^\8 r^{-3}(1+r^2)^{-1}\,dr
\\
&=
c\eps^3r_0^2E^2\left(r_2^{-2} - \log(1+r_2^{-2})\right).
\end{align*}
The result \eqref{eq:S4HarmonicSpinorL2Decay} follows.
\end{proof}

We note that the proof of the preceding lemma is much simpler in the case
of harmonic spinors over $\RR^4$:

\begin{lem}
\label{lem:R4HarmonicSpinorDecay}
There is a constant $\eps>0$ and, if $M>0$, a constant $c(M)$ with the
following significance.  Let $\RR^4$ have an asymptotically flat metric $g$
(in the sense, for example, of \cite[Equation (1.9)]{FeehanKato}).  Let
$h:\RR^4\cong S^4\less\{s\}$ be inverse stereographic projection, and let
$B(r_0)$ denote the ball $\{x\in
\RR^4:r<r_0\}$, where $r=|x|$. Let $A$ be an $L^2_k$ unitary connection on
a Hermitian bundle $E$ over $S^4$ and denote its pullback to $\RR^4$ via
$h$ also by $A$.  Let $(\rho,S^+,S^-)$ be the standard spin structure over
$\RR^4$, and let $V=V^+\oplus V^-$ with $V^\pm = S^\pm\otimes h^*E$. Suppose
$A$ obeys
$$
\|F_A\|_{L^\8(\RR^4-B(1))} \leq M
\text{ and }
\|F_A\|_{L^2(\RR^4-B(r_0))} \leq \eps.
$$
Let $D_{A,g}:\Gamma(\RR^4,V)\to\Gamma(\RR^4,V)$ be the corresponding Dirac
operator. Suppose $\psi\in L^2(\RR^4,V)$ obeys $D_{A,g}\psi = 0$. Then for
all $x\in \RR^4-B(2r_0)$ we have
\begin{equation}
\label{eq:R4HarmonicSpinorDecay-2}
|\psi|(x) 
\leq 
c\eps^{3/2}r_0r^{-3}\|\psi\|_{L^2(\RR^4-h(B(r_0)),g)}.
\end{equation}
Furthermore, if $2r_0 \leq r_2 < \8$, 
\begin{equation}
\label{eq:R4HarmonicSpinorL2Decay}
\|\psi\|_{L^2(\RR^4-B(r_2))} 
\leq 
c\eps^{3/2}r_0r_2^{-1}\|\psi\|_{L^2(S^4-h(B(r_0)),g)}.
\end{equation}
\end{lem}

\begin{proof}
Inequality \eqref{eq:R4HarmonicSpinorDecay-2} follows from Theorem 1.1 in
\cite{FeehanKato}. For the second estimate, observe that, setting $E =
\|\psi\|_{L^2(\RR^4-B(r_0))}$ for convenience, we have 
\begin{align*}
\|\psi\|_{L^2(\RR^4-B(r_2))}^2
&=
\int_{\RR^4-B(r_2)}|\psi|^2 \,d^4x
\\
&\leq
c\eps^3r_0^2E^2\int_{r_2}^\8 r^{-6}r^3\,dr
=
c\eps^3r_0^2r_2^{-2}E^2.
\end{align*}
The result \eqref{eq:R4HarmonicSpinorL2Decay} follows.
\end{proof}

\subsection{An upper bound for the small eigenvalues}
\label{subsec:Eigen2UpperBoundSmall}
We turn to the proof of the assertions of Theorem \ref{thm:H2SmallEval}
concerning the upper bound for the small eigenvalues.

We define a cut-off function $\chi_0$ on $X$ which is equal to zero on the
balls $B(x_p,4\lambda_p^{1/3})$, and setting $\chi_0=1$ on $X-\cup_{p=1}^m
B(x_p,8L_p\lambda_p^{1/3})$ (see Lemma \ref{lem:dBetaEst} for its
definition). We recall that $\|d\chi_0\|_{L^4(X)} \leq c(\log L)^{-3/4}$,
where $L = \max_{1\leq p\leq m}L_p$ and so, if we set $L_p = \lambda_p^{-1/6}$
for $1\leq p\leq m$ we have
$\|d\chi_0\|_{L^4(X)} \leq c(-\log\lambda)^{-3/4}$ with $\lambda \in
(0,\half]$. Recall from Definition \ref{defn:LambdaClosePair} that
$$
U_0 := X -\bigcup_{p=1}^m\bar B(x_p,4\lambda_p^{1/3}),
$$
so $U_0\Subset X\less\bx$ and $U_0$ is the interior of $\supp\chi_0$. 

Similarly, define cut-off functions $\chi_p$, $1\leq p\leq m$, on $X$ which
are equal to $1$ on the balls $B(x_p,\eighth L_p^{-1}\lambda_p^{1/3})$ and
equal to zero on the complement of the balls
$B(x_p,\quarter\lambda_p^{1/3})$. As before, $\|d\chi_p\|_{L^4(X)} \leq
c(\log L)^{-3/4}$ and so, if $L_p = \lambda_p^{-1/6}$, we have
$\|d\chi_p\|_{L^4(X)} \leq c(-\log\lambda)^{-3/4}$ with $\lambda \in
(0,\half]$.

For the remainder of this section, we assume for simplicity that $\chi_0$
is defined as in the preceding paragraph and observe that
$\supp\chi_0\Subset X\less\bx$. Our main task in this subsection is to find
an upper bound for the first $N+m+1$ eigenvalues of
$d_{A,\Phi}^1d_{A,\Phi}^{1,*}$.

\begin{prop}
\label{prop:H2SmallEvalUpperBound}
Continue the hypotheses and notation of Theorem
\ref{thm:H2SmallEval}.  Let $\nu_2[A_0,\Phi_0]$ be the least positive
eigenvalue of $d_{A_0,\Phi_0}^1d_{A_0,\Phi_0}^{1,*}$ and let
$\mu_1[A,\Phi]\le\cdots\le\mu_{N+m+1}[A,\Phi]$ be the first $N+m+1$
eigenvalues of $d_{A,\Phi}^1d_{A,\Phi}^{1,*}$. Then there are positive
constants $C$ and $\lambda_0$ (see Remark \ref{rmk:LambdaClosePairConstants})
such that for all $\lambda \in (0,\lambda_0]$, the following holds:
\begin{equation}
\label{eq:H2DeAEvalEst}
\mu_l[A,\Phi] 
\le
\begin{cases}
C(-\log\lambda)^{-3/2}, &\text{if }l=1,\dots,N, \\
K + C(-\log\lambda)^{-3/4}, &\text{if }l=N+1,
\end{cases} 
\end{equation}
where $K = \min\{\nu_2[A_0,\Phi_0],\nu_2[A_1],\dots,\nu_2[A_m]\}$ and
$\nu[A_p]$ is the least positive eigenvalue of the Dirac Laplacian
$D_{A_p,g_p}D_{A_p,g_p}^*$ on $L^2(S^4,V_p^-)$, with $A_p$ and $g_p$ the
induced connection and Riemannian metric from $A$ and $g$ over $X$ near the
point $p\in X$.
\end{prop}

Before we begin the proof proper, let
$\{(\xi_i,\varphi_i)\}_{i=1}^{n}\subset
L^2(X,\Lambda^+\otimes\fg_{E_\ell})\oplus L^2(X,V_\ell^-)$ be an
$L^2$-orthonormal basis for $\Ker d_{A_0,\Phi_0}^{1,*}$ and let
$(\xi_{n+1},\varphi_{n+1}) \in(\Ker d_{A_0,\Phi_0}^{1,*})^\perp$ be an
$L^2$ unit eigenvector corresponding to the least positive eigenvalue
$\nu_2[A_0,\Phi_0]$ of $d_{A_0,\Phi_0}^1d_{A_0,\Phi_0}^{1,*}$.  Let
$\{\varphi_{pi}\} \subset L^2(S^4,V_p^-)$ be an $L^2$-orthonormal basis for
$\Ker D_{A_p,g_p}^*$, for $1\leq p\leq m$ and $1 \leq i\leq 2\kappa_p$; let
$\varphi_{p,2\kappa_p+1} \in(\Ker D_{A_p,g_p}^*)^\perp$ be $L^2$ unit
eigenvectors corresponding to the least positive eigenvalue $\nu_2[A_p]$ of
$D_{A_p,g_p}D_{A_p,g_p}^*$, $1\leq p\leq m$. (We know from equation
\eqref{eq:BubbleCharge} that $c_2(E^p)=\kappa_p$ and from Lemma
\ref{lem:DiracS4VanishingKernel} that 
$\dim_\RR\Ker D_{A_p,g_p}^* = 2\kappa_p$.) We begin by constructing a set of
approximate eigenvectors, $\{(\xi_{pi}',\varphi_{pi}')\}$, for the first
$N+m+1$ eigenvalues of $d_{A,\Phi}^1d_{A,\Phi}^{1,*}$. The most natural
choice is given by
\begin{align}
\notag
(\xi_{0i}',\varphi_{0i}')
&:=
(\chi_0\xi_i,\chi_0\varphi_i),\quad 1 \leq i \leq n+1,
\\
\label{eq:FromKernelSphereToSmallEvalVecs}
(\xi_{pi}',\varphi_{pi}')
&:=
(0,\chi_p\varphi_{pi}),\quad 1\leq p\leq m\text{ and }1 \leq i\leq 2\kappa_p+1,
\end{align}
recalling that $\sum_{p=1}^m\kappa_p = \ell$ by hypothesis. 

When $p\neq 0$ one can see from equation
\eqref{eq:LAPhi*} for $d_{A,\Phi}^{1,*}$ that this choice yields
\begin{equation}
\label{eq:Expandd1*OnCutoffEvector0phi}
\begin{aligned}
d_{A,\Phi}^{1,*}(\xi_{pi}',\varphi_{pi}') 
&= 
\chi_p((\cdot\Phi)^*\varphi_{pi},D_{A_p,g_p}^*\varphi_{pi}) 
+ (0,\rho(d\chi_p)\varphi_{pi})
\\
&=
(0,\chi_pD_{A_p,g_p}^*\varphi_{pi}) 
+ (0,\rho(d\chi_p)\varphi_{pi})
\quad\text{on }B_p,
\end{aligned}
\end{equation}
since the condition \eqref{eq:PhiIsZeroOnBallSmall} on $\Phi$ and
$B_p$ implies that $\Phi\equiv 0$ on $\supp\chi_p$.

We can now proceed with the proof of Proposition
\ref{prop:H2SmallEvalUpperBound}. 

\begin{lem}
\label{lem:H2DeAEvectorEst}
Continue the above notation. Then there are positive constants $C$ and
$\lambda_0$ such that for all $\lambda\in (0,\lambda_0]$, the following
hold:
\begin{align}
\label{eq:H2DeAEvectorEstX}
\|d_{A,\Phi}^{1,*}(\xi_{0i}',\varphi_{0i}')\|_{L^2(X)} 
&\le 
\begin{cases}
C\lambda^{1/3}, &\text{if }1\le i\le n, \\
\sqrt{\nu_2[A_0,\Phi_0]}+C\lambda^{1/3}, &\text{if }i= n+1,
\end{cases}
\\
\label{eq:H2DeAEvectorEstS4}
\|d_{A,\Phi}^{1,*}(\xi_{pi}',\varphi_{pi}')\|_{L^2(X)} 
&\le 
\begin{cases}
C(-\log\lambda)^{-3/4}, &\text{if }1\le i\le 2\kappa_p, \\
\sqrt{\nu_2[A_p]}+C(-\log\lambda)^{-3/4}, &\text{if }i= 2\kappa_p+1.
\end{cases}
\end{align}
\end{lem}

\begin{proof}
To verify the first assertion in Lemma \ref{lem:H2DeAEvectorEst}, note that
$(\xi_{0i}',\varphi_{0i}') = \chi_0(\xi_i,\varphi_i)$ and so for $1\le i\le
n+1$ (schematically),
\begin{equation}
\begin{aligned}
\label{eq:ExpanddA+*OnCutoffEvectorXi}
d_{A,\Phi}^{1,*}(\xi_{0i}',\varphi_{0i}') 
&= 
\chi_0d_{A,\Phi}^{1,*}(\xi_i,\varphi_i) + d\chi_0\otimes (\xi_i,\varphi_i)
\\
&= 
\chi_0d_{A_0,\Phi_0}^{1,*}(\xi_i,\varphi_i) + d\chi_0\otimes (\xi_i,\varphi_i)
+ \chi_0(A-A_0,\Phi-\Phi_0)\otimes (\xi_i,\varphi_i).
\end{aligned}
\end{equation}
For $1\leq p\leq m$ and $1\leq i\leq 2\kappa_p+1$, equation
\eqref{eq:Expandd1*OnCutoffEvector0phi} yields the corresponding
expressions for $d_{A,\Phi}^{1,*}(\xi_{pi}',\varphi_{pi}')$. Now
\begin{align}
\label{eq:KernelPlus1L2EvectorEstX}
\|d_{A_0,\Phi_0}^{1,*}(\xi_i,\varphi_i)\|_{L^2(X)} 
&=
\begin{cases}
0, &\text{if $1\leq i\leq n$},
\\
\sqrt{\nu_2[A_0,\Phi_0]}\|(\xi_i,\varphi_i)\|_{L^2(X)}, &\text{if $i = n+1$}.
\end{cases}
\\
\label{eq:KernelPlus1L2EvectorEstS4}
\|D_{A_p,g_p}^*\varphi_{pi}\|_{L^2(S^4)} 
&=
\begin{cases}
0, &\text{if $1\leq i\leq 2\kappa_p$},
\\
\sqrt{\nu_2[A_p]}\|\varphi_{pi}\|_{L^2(S^4)}, &\text{if $i = 2\kappa_p+1$}.
\end{cases}
\end{align}
By Lemma \ref{lem:dBetaEst} (yielding $\|d\chi_0\|_{L^2} \leq
c\lambda^{1/3}$), Lemma \ref{lem:LInftyEstDe2A0Phi0Evec}, and equation
\eqref{eq:ExpanddA+*OnCutoffEvectorXi}, we have for $p=0$ and $1\leq i\leq
n+1$, 
\begin{align*}
&\|d_{A,\Phi}^{1,*}(\xi_{0i}',\varphi_{0i}')\|_{L^2} 
\\
&\le (\|d\chi_0\|_{L^2}
+ \|\chi_0(A-A_0,\Phi-\Phi_0)\|_{L^2})\|(\xi_i,\varphi_i)\|_{L^\8} 
\\
&\quad + \|d_{A_0,\Phi_0}^{1,*}(\xi_i,\varphi_i)\|_{L^2} 
\quad\text{(by Equation \eqref{eq:ExpanddA+*OnCutoffEvectorXi})}
\\
&\le c\lambda^{1/3}\|(\xi_i,\varphi_i)\|_{L^\8} 
+ \|d_{A_0,\Phi_0}^{1,*}(\xi_i,\varphi_i)\|_{L^2} 
\quad\text{(by Lemma \ref{lem:dBetaEst} \& Inequality
\eqref{eq:APhiA0Phi0L4Close})} 
\\
&\le C\lambda^{1/3}\|(\xi_i,\varphi_i)\|_{L^2} 
+ \|d_{A_0,\Phi_0}^{1,*}(\xi_i,\varphi_i)\|_{L^2} 
\quad\text{(by Lemma \ref{lem:LInftyEstDe2A0Phi0Evec}).}
\end{align*}
Therefore, the first assertion in Lemma \ref{lem:H2DeAEvectorEst}
follows from the preceding inequality and
the eigenvalue estimates \eqref{eq:KernelPlus1L2EvectorEstX}.

For $1\leq p\leq m$ and $1\leq i\leq 2\kappa_p+1$ and noting that
$\varphi_{pi}$ is an eigenvector of $D_{A_p,g_p}D_{A_p,g_p}^*$, so Theorem
1.2 in \cite{FeehanKato} applies,
\begin{align*}
&\|d_{A,\Phi}^{1,*}(\xi_{pi}',\varphi_{pi}')\|_{L^2(X)} 
\\
&\le \|d\chi_p\|_{L^2}\|\varphi_{pi}\|_{L^\8(\supp d\chi_p)} 
+ \|D_{A_p,g_p}^*\varphi_{pi}\|_{L^2(\supp\chi_p)}
\quad\text{(by Equation \eqref{eq:Expandd1*OnCutoffEvector0phi})}
\\
&\le c\lambda^{1/3}\|\varphi_{pi}\|_{L^\8(\supp d\chi_p)}
+ \|D_{A_p,g_p}^*\varphi_{pi}\|_{L^2(S^4)} 
\quad\text{(by Lemma \ref{lem:dBetaEst})} 
\\
&\le C\lambda^{1/3}\|\varphi_{pi}\|_{L^2(S^4)}
+ \|D_{A_p,g_p}^*\varphi_{pi}\|_{L^2(S^4)}, 
\end{align*}
with the last inequality following from
\cite[Theorem 1.2]{FeehanKato}, giving the estimate for $\varphi_{pi}$.
Therefore, the final assertion in Lemma \ref{lem:H2DeAEvectorEst} follows
from the preceding inequality and the eigenvalue estimates
\eqref{eq:KernelPlus1L2EvectorEstS4}.
\end{proof}

\begin{rmk}
\label{rmk:WhyWeAssumeMetricS4IsAlmostRound}
It is worth noting that if the Dirac operator $D_{A_p,g_p}$ were replaced
by $D_{A_p,g_0}$, where we use $g_0$ --- the standard round metric on
$S^4$ --- rather than the approximately round metric $g_p$ induced from $g$
near $x_p\in X$, then our derivation of the estimate
\eqref{eq:H2DeAEvectorEstS4} would require a uniform $L^2_{1,A_p}(S^4)$
estimate for the negative eigenspinor $\varphi_{pi}$. However, as we
discuss further in \cite{FeehanKato}, such an estimate appears impossible
to obtain due to the unfavorable structure of the Bochner formula
\eqref{eq:BWDirac-}.
\end{rmk}

\begin{lem}
\label{lem:H2DeAEvectorInnerProductEst}
Continue the above notation. Then there are positive constants $C$ and
$\lambda_0$ such that for all $\lambda\in (0,\lambda_0]$, the following
holds:
\begin{equation}
\label{eq:H2DeAEvectorInnerProductEst}
\left|((\xi_{pi}',\varphi_{pi}'),
(\xi_{qj}',\varphi_{qj}'))_{L^2(X)} - \delta_{pq}\delta_{ij}\right|
\le 
\begin{cases}
C\lambda^{4/3}, &\text{if $p=q=0$ and $1\le i,j\le n+1$},
\\
C\lambda^{4/3}, &\text{if $p=q\neq 0$ and $1\le i,j\le 2\kappa_p+1$},
\\
0&\text{if $p \neq q$}.
\end{cases}
\end{equation}
\end{lem}
%\marginpar{\tiny Middle estimate, p=q=0, used to be $\lambda^{2/3}$, so now
%should be improvement where this used.}

\begin{proof}
Suppose $p=q=0$ and $1\le i,j\le n+1$. Using $(\xi_{0i}',\varphi_{0i}') 
= \chi_0(\xi_i,\varphi_i) = (\xi_i,\varphi_i) + (\chi_0-1)(\xi_i,\varphi_i)$, 
\begin{align*}
((\xi_{0i}',\varphi_{0i}'),(\xi_{0j}',\varphi_{0j}'))_{L^2(X)}
&= 
((\xi_i,\varphi_i),(\xi_j,\varphi_j))_{L^2(X)}
+ 2((\chi_0-1)(\xi_i,\varphi_i),(\xi_j,\varphi_j))_{L^2(X)}
\\
&\quad + ((\chi_0-1)(\xi_i,\varphi_i),(\chi_0-1)(\xi_j,\varphi_j))_{L^2(X)}.
\end{align*}
But $((\xi_i,\varphi_i),(\xi_j,\varphi_j))_{L^2(X)}=\de_{ij}$, while
the second two terms on the right above are bounded by
\begin{align*}
&\left|2((\chi_0-1)(\xi_i,\varphi_i),(\xi_j,\varphi_j))_{L^2(X)}
+ ((\chi_0-1)(\xi_i,\varphi_i),(\chi_0-1)(\xi_j,\varphi_j))_{L^2(X)}\right| 
\\
&\qquad\le 
c\lambda^{4/3}\|(\xi_i,\varphi_i)\|_{L^\8}\|(\xi_j,\varphi_j)\|_{L^\8} 
\le 
C\lambda^{4/3},
\end{align*}
courtesy of the definition of $\chi_0$ (which implies that
$\Vol(\supp(1-\chi_0)) \leq c\lambda^{4/3}$) and
Lemma \ref{lem:LInftyEstDe2A0Phi0Evec}.  This gives the
conclusion when $p=q=0$. 

The case $p\neq q$ follows from the fact that the supports of $\chi_p$ and
$\chi_q$ are disjoint when $p\neq q$. 

The assertion for $p=q\neq 0$ and $1\le i,j\le 2\kappa_p+1$ follows by the
much the same argument as for $p=q=0$. Indeed, using
$(\xi_{pi}',\varphi_{pi}')  
= (0,\chi_p\varphi_{pi}) = (0,\varphi_{pi}) +
(0,(\chi_p-1)\varphi_{pi})$,  
\begin{align*}
((\xi_{pi}',\varphi_{pi}'),(\xi_{pj}',\varphi_{pj}'))_{L^2(X)}
&= 
(\varphi_{pi},\varphi_{pj})_{L^2(S^4)}
+ 2((\chi_p-1)\varphi_{pi},\varphi_{pj})_{L^2(S^4)}
\\
&\quad + ((\chi_p-1)\varphi_{pi},(\chi_p-1)\varphi_{pj})_{L^2(S^4)}.
\end{align*}
But $(\varphi_{pi},\varphi_{pj})_{L^2(S^4)}=\de_{ij}$, while
the remaining two terms on the right above are bounded by (see Lemma
\ref{lem:S4HarmonicSpinorDecay}) 
\begin{align*}
&\left|2((\chi_p-1)\varphi_{pi},\varphi_{pj})_{L^2(S^4)}
+ ((\chi_p-1)\varphi_{pi},(\chi_p-1)\varphi_{pj})_{L^2(S^4)}\right| 
\\
&\le 
3\|\varphi_{pi}\|_{L^2(S^4-h(B_p))}\|\varphi_{pj}\|_{L^2(S^4-h(B_p))} 
\\
&\le
c\lambda^{4/3}\|\varphi_{pi}\|_{L^2(S^4)}\|\varphi_{pj}\|_{L^2(S^4)} 
\quad\text{(by Lemma \ref{lem:S4HarmonicSpinorDecay})}
\\
&\leq
c\lambda^{2/3},
\end{align*}
where $B_p = B(\quarter\lambda_p^{1/3})$, noting that $1-\chi_p$ is
supported on $S^4-h(B_p)$ and we applied Lemma
\ref{lem:S4HarmonicSpinorDecay} with $r_0 = L\lambda_p$, for a constant
$L\gg 1$, and $r_2 = \quarter\lambda_p^{1/3}\ll 1$, observing that
$\|F_{A_p}\|_{L^2(S^4-h(B(L\lambda_p)))} \leq \eps$ and
$\|F_{A_p}\|_{L^2(S^4-h(B(1)))} \leq c\kappa_p$, by hypothesis of Theorem
\ref{thm:H2SmallEval} on $(A,\Phi)$ and the definition of $A_p$ on $S^4$.
\end{proof}

\begin{proof}[Proof of Proposition \ref{prop:H2SmallEvalUpperBound}]
Let $\{(\xi_l',\varphi_l')\}_{l=1}^{N+m+1}$ denote an ordering of the set
of approximate eigenvectors $\{(\xi_{pi}',\varphi_{pi}'): 0\leq p\leq m,
1\leq i\leq n_p+1\}$, where $n_0 = n$ and $n_p = 2\kappa_p$ when $p\geq 1$.
Now $\|(\xi_l',\varphi_l')\|_{L^2(X)}\ge 1-C\lambda^{2/3}$ by Lemma
\ref{lem:H2DeAEvectorInnerProductEst} 
and so the first and second assertions of Lemma \ref{lem:H2DeAEvectorEst} yield
\begin{equation}
\label{eq:L2EstApproxEvectorsX}
\begin{aligned}
\|d_{A,\Phi}^{1,*}(\xi_l',\varphi_l')\|_{L^2} 
&\le 
C(-\log\lambda)^{-3/4}(1-C\lambda^{2/3})^{-1}\|(\xi_l',\varphi_l')\|_{L^2}
\\
&\le 
C(-\log\lambda)^{-3/4}\|(\xi_l',\varphi_l')\|_{L^2}, \quad 1\le l\le N,
\\
\|d_{A,\Phi}^{1,*}\xi_{N+1}'\|_{L^2} 
&\le 
\left(\sqrt{K}+C(-\log\lambda)^{-3/4}\right)\|\xi_{n+1}'\|_{L^2},
\end{aligned}
\end{equation}
for small enough $\lambda_0$ and
recalling that $K = \min\{\nu_2[A_0,\Phi_0],\nu_2[A_1],\dots,\nu_2[A_m]\}$. 

Let $F_l':=[(\xi_1',\varphi_1'),\dots,(\xi_l',\varphi_l')]$ denote the span
of the vectors $(\xi_1',\varphi_1'),\dots,(\xi_l',\varphi_l')$ in
$L^2(X,\Lambda^+\otimes\fg_E)\oplus L^2(X,V^-)$ for $l=1,\dots,N+m+1$.
By applying Gram-Schmidt orthonormalization to the basis
$\{(\xi_1',\varphi_1'),\dots,(\xi_l',\varphi_l')\}$ for $F_l'$, the
estimates \eqref{eq:L2EstApproxEvectorsX} imply that
\begin{equation}
\|d_{A,\Phi}^{1,*}(\eta,\psi)\|_{L^2} 
\le 
\begin{cases}
C(-\log\lambda)^{-3/4}\|(\eta,\psi)\|_{L^2}, &\text{if }(\eta,\psi)\in F_{N}', 
\\
\left(\sqrt{K}+C(-\log\lambda)^{-3/4}\right)\|(\eta,\psi)\|_{L^2}, 
&\text{if }(\eta,\psi)\in F_{N+1}'. 
\end{cases}
\label{eq:dA+*SmallEvalEst}
\end{equation}
Recall from the statement of Theorem \ref{thm:H2SmallEval} that
$\mu_1[A,\Phi]\le\cdots\le\mu_{N+1}[A,\Phi]$ are the first $N+1$
eigenvalues of $d_{A,\Phi}^1d_{A,\Phi}^{1,*}$; let
$(\eta_l,\psi_l)$, for $l=1,\dots,N+1$, be the corresponding
$L^2$-orthonormal eigenvectors. Let $F_l'$ denote the span of the
eigenvectors $(\eta_1,\psi_1),\dots,(\eta_l,\psi_l)$ in
$L^2(X,\Lambda^+\otimes\fg_E)\oplus L^2(X,V^-)$ and let $F_l^\perp$ be its
$L^2$-orthogonal complement in $L^2(X,\Lambda^+\otimes\fg_E)\oplus
L^2(X,V^-)$. Note that $\dim (F_l^\perp\cap F'_{N}) \ge N-l$. Then,
for $l=1,\dots,N$, 
\begin{align*}
\mu_l[A,\Phi] 
&= 
\inf_{(\eta,\psi)\in F_{l-1}^\perp}
\frac{\|d_{A,\Phi}^{1,*}(\eta,\psi)\|_{L^2}^2}{\|(\eta,\psi)\|_{L^2}^2} 
\le 
\inf_{(\eta,\psi)\in F_{l-1}^\perp\cap F'_{N}}
\frac{\|d_{A,\Phi}^{1,*}(\eta,\psi)\|_{L^2}^2}{\|(\eta,\psi)\|_{L^2}^2} 
\\
&\le 
C(-\log\lambda)^{-3/2}, 
\end{align*}
where we use the estimates
\eqref{eq:dA+*SmallEvalEst} to obtain the final inequality above.  This
gives the first inequality in the stated eigenvalue bounds
\eqref{eq:H2DeAEvalEst}. Therefore, $d_{A,\Phi}^1d_{A,\Phi}^{1,*}$ has at
least $N$ eigenvalues which are less than or equal to $C(-\log\lambda)^{-3/2}$.
Similarly,
\begin{align*}
\mu_{N+1}[A,\Phi] 
&= 
\inf_{(\eta,\psi)\in F_{N}^\perp}
\frac{\|d_{A,\Phi}^{1,*}(\eta,\psi)\|_{L^2}^2}{\|(\eta,\psi)\|_{L^2}^2} 
\le 
\inf_{(\eta,\psi)\in F_{N}^\perp\cap F'_{N+1}}
\frac{\|d_{A,\Phi}^{1,*}(\eta,\psi)\|_{L^2}^2}{\|(\eta,\psi)\|_{L^2}^2} 
\\
&\le \left(\sqrt{K} + C(-\log\lambda)^{-3/4}\right)^2
\le K + C(-\log\lambda)^{-3/4},
\end{align*}
where the estimates \eqref{eq:dA+*SmallEvalEst} yield the penultimate
inequality above. This completes the proof of the second inequality in the
eigenvalue bounds \eqref{eq:H2DeAEvalEst}.
\end{proof}

\subsection{A map from the small-eigenvalue eigenspace to the direct sum of
kernels} 
\label{subsec:ApproxOrthoBasisEigenvectors}
In this section we decompose the eigenspace spanned by the eigenvectors for
the first $N$ eigenvalues of $d_{A,\Phi}^1d_{A,\Phi}^{1,*}$ and identify
the decomposition with the kernels of the operators $d_{A_0,\Phi_0}^{1,*}$
and $D_{A_p,g_p}^*$, for $p=1,\dots,m$.  For this purpose, we define a
block-diagonal, first-order, linear elliptic operator
\begin{equation}
\label{eq:DiagonalAdjointOperator}
T_{\bA,\bPhi}^*
=
d_{A_0,\Phi_0}^{1,*}\oplus\bigoplus_{p=1}^m D_{A_p,g_p}^*
\end{equation}
with domain and range, respectively, given by
\begin{align*}
&\left(L^2_{k-1}(X,\Lambda^+\otimes\fg_{E_\ell})\oplus L^2_{k-1}(X,V_\ell^-)\right)
\oplus \bigoplus_{p=1}^m
L^2_{k-1}(S^4,V_p^-),
\\
&\left(L^2_{k-2}(X,\Lambda^1\otimes\fg_{E_\ell})\oplus L^2_{k-2}(X,V_\ell^+)\right)
\oplus \bigoplus_{p=1}^m
L^2_{k-2}(S^4,V_p^+).
\end{align*}
We shall construct an approximately $L^2$-orthonormal basis for $\Ker
T_{\bA,\bPhi}^*$ by cutting-off the eigenvectors corresponding to the first
$N$ eigenvalues of $d_{A,\Phi}^1d_{A,\Phi}^{1,*}$, estimating the error
between the cut-off eigenvectors and their $L^2$-orthogonal projections
onto $\Ker T_{\bA,\bPhi}^*$ and then applying Gram-Schmidt
orthonormalization to obtain an $L^2$-orthonormal basis.

Let $\Pi_{A_0,\Phi_0}$ be the $L^2$-orthogonal projection from
$L^2(X,\Lambda^+\otimes\fg_{E_\ell})\oplus L^2(X,V_\ell^-)$ onto $\Ker
d_{A_0,\Phi_0}^{1,*}$ and let $\Pi_{A_0,\Phi_0}^\perp = 1-\Pi_{A_0,\Phi_0}$
be the $L^2$-orthogonal projection from
$L^2(X,\Lambda^+\otimes\fg_{E_\ell})\oplus L^2(X,V_\ell^-)$ onto $(\Ker
d_{A_0,\Phi_0}^{1,*})^\perp = \Ran d_{A_0,\Phi_0}^1$.  Recall that
$\Pi_{A_p}$ is the $L^2$-orthogonal projection from $L^2(S^4,V_p^-)$ onto
$\Ker D_{A_p,g_p}^*$ and let $\Pi_{A_p}^\perp = 1-\Pi_{A_p}$ be the
$L^2$-orthogonal projection from $L^2(X,V_\ell^-)$ onto $(\Ker
D_{A_p,g_p}^*)^\perp = \Ran D_{A_p,g_p}$.

Then the $L^2$-orthogonal projection from the
domain of $T_{\bA,\bPhi}$ onto $\Ker T_{\bA,\bPhi}^*$ is given by
\begin{equation}
\label{eq:DiagonalKernelProj}
\Pi_{\bA,\bPhi}
=
\Pi_{A_0,\Phi_0}\oplus \bigoplus_{p=1}^m \Pi_{A_p}.
\end{equation}
Recall from the proof of Proposition \ref{prop:H2SmallEvalUpperBound} that
$\{(\eta_l,\psi_l): 1\leq l\leq N+1\}$ was defined to be the
set of $L^2$-orthonormal eigenvectors associated with the eigenvalues
$\mu_1[A,\Phi]\le\cdots\le\mu_{N+1}[A,\Phi]$ of
$d_{A,\Phi}^1d_{A,\Phi}^{1,*}$. We then obtain a set of $N+1$ approximate
eigenvectors for $T_{\bA,\bPhi}$ by cutting off:
\begin{equation}
\label{eq:VectorCutoffEigenvectors}
(\bfeta_l',\bpsi_l')
=
\left((\chi_0\eta_l,\chi_0\psi_l),
(0,\chi_1\psi_l),\dots,(0,\chi_m\psi_l)\right),
\quad l = 1, \dots,N+1.
\end{equation}
Our main goal in this subsection is to prove:

\begin{prop}
\label{prop:H2SmallEvalApproxOrthonBasis}
Continue the hypotheses and notation of Theorem
\ref{thm:H2SmallEval}. Then there are positive constants 
$C$ and $\lambda_0$ such that for all $\lambda \in (0,\lambda_0]$, the set
$\{\Pi_{\bA,\bPhi}(\bfeta_l',\bpsi_l')\}_{l=1}^N$ is a basis for $\Ker
T_{\bA,\bPhi}^*$ and the following estimates hold:
\begin{gather}
\tag{1}
\|\Pi_{\bA,\bPhi}^\perp(\bfeta_l',\bpsi_l')\|_{L^2}
\leq C(-\log\lambda)^{-3/4},\quad 1\le l\le N,
\\
\tag{2}
\left|\left(\Pi_{\bA,\bPhi}(\bfeta_l',\bpsi_l'),
\Pi_{\bA,\bPhi}(\bfeta_k',\bpsi_k')\right)_{L^2}
-\de_{kl}\right| \le C(-\log\lambda)^{-3/2}, \quad 1\le k,l\le N,
\\
\tag{3}
\begin{aligned}
\|\Pi_{\bA,\bPhi}(\bfeta_{N+1}',\bpsi_{N+1}')\|_{L^2} 
&\le C(-\log\lambda)^{-3/4},
\\
\|\Pi_{\bA,\bPhi}^\perp(\eta_{N+1}',\psi_{N+1}')\|_{L^2}
&\ge 1-C(-\log\lambda)^{-3/4}. 
\end{aligned}
\end{gather}
\end{prop}

The following estimates are the key to the proofs of Propositions
\ref{prop:H2SmallEvalApproxOrthonBasis} and
\ref{prop:H2SmallEvalLowerBound}.

\begin{lem}
\label{lem:H2DeA0EvectorEst} 
Continue the above notation. Then there are positive constants $C$ and
$\lambda_0$ such that for all $\lambda \in (0,\lambda_0]$, the following
hold:
\begin{align}
\label{eq:H2DeA0EvectorEst-L43} 
\|d_{A_0,\Phi_0}^{1,*}(\chi_0\eta_l,\chi_0\psi_l)\|_{L^{4/3}(X)} 
&\le 
C(-\log\lambda)^{-3/4}\|(\eta_l,\psi_l)\|_{L^2(X)}, \quad 1\le l\le N,
\\
\label{eq:H2DeA0EvectorEst-L2} 
\|d_{A_0,\Phi_0}^{1,*}(\chi_0\eta_{N+1},\chi_0\psi_{N+1})\|_{L^2(X)} 
&\le 
\|d_{A,\Phi}^{1,*}(\eta_{N+1},\psi_{N+1})\|_{L^2(\supp\chi_0)}
\\
\notag
&\quad + C(-\log\lambda)^{-3/4}\|(\eta_{N+1},\psi_{N+1})\|_{L^2(X)}.
\end{align}
\end{lem}

\begin{proof}
Observe that the definition \eqref{eq:LAPhi*} of $d_{A_0,\Phi_0}^{1,*}$ and
writing $(A,\Phi)=(A_0,\Phi_0)+(A-A_0,\Phi-\Phi_0)$ on $X\less\bx$ implies
that, schematically,
$$
d_{A_0,\Phi_0}^{1,*}(v,\psi)
=
d_{A,\Phi}^{1,*}(v,\psi) - \{(A-A_0,\Phi-\Phi_0),(v,\psi)\},
$$
and so (schematically) for $p=0$ and $1\leq l\leq N$ we have
\begin{equation}
\label{eq:ExpanddA+*OnCutoffEvectorEta}
\begin{aligned}
d_{A_0,\Phi_0}^{1,*}(\chi_0\eta_l,\chi_0\psi_l) 
&= 
d\chi_0\otimes(\eta_l,\psi_l) 
+ \chi_0 d_{A,\Phi}^{1,*}(\eta_l,\psi_l)  
\\
&\quad - \chi_0\{(A-A_0,\Phi-\Phi_0),(\eta_l,\psi_l)\}.
\end{aligned}
\end{equation}
Observe that
\begin{equation}
\label{eq:SqRootL43L2EvalueBound}
\|d_{A,\Phi}^{1,*}(\eta_l,\psi_l)\|_{L^{4/3}} 
\leq
c\|d_{A,\Phi}^{1,*}(\eta_l,\psi_l)\|_{L^2} 
=
c\sqrt{\mu_l[A,\Phi]}\|(\eta_l,\psi_l)\|_{L^2}.
\end{equation}
Thus, by the upper bounds for the eigenvalues $\mu_l[A,\Phi]$ in
\eqref{eq:H2DeAEvalEst} (see Proposition \ref{prop:H2SmallEvalUpperBound})
and Lemmas \ref{lem:dBetaEst} and the identity
\eqref{eq:ExpanddA+*OnCutoffEvectorEta}, we have for $1\le l\le N$,
\begin{align*}
\|d_{A_0,\Phi_0}^{1,*}(\chi_0\eta_l,\chi_0\psi_l)\|_{L^{4/3}} 
&\le (\|d\chi_0\|_{L^2}+\|\chi_0(A-A_0,\Phi-\Phi_0)\|_{L^2})
\|(\eta_l,\psi_l)\|_{L^4(\supp\chi_0)} 
\\
&\quad + \|d_{A,\Phi}^{1,*}(\eta_l,\psi_l)\|_{L^{4/3}} 
\quad\text{(by Equation \eqref{eq:ExpanddA+*OnCutoffEvectorEta})}
\\
&\le C(\|d\chi_0\|_{L^2}+\|\chi_0(A-A_0,\Phi-\Phi_0)\|_{L^2} 
+ \sqrt{\mu_l[A,\Phi]})\|(\eta_l,\psi_l)\|_{L^2}
\\
&\quad\text{(by Inequality \eqref{eq:SqRootL43L2EvalueBound} and 
Lemma \ref{lem:L21GlobalLocalEstDe2APhiEvec})}
\\
&\le C(\|d\chi_0\|_{L^2}+\|\chi_0(A-A_0,\Phi-\Phi_0)\|_{L^2} 
+ \sqrt{\mu_l[A,\Phi]})\|(\eta_l,\psi_l)\|_{L^2}
\\
&\le C(-\log\lambda)^{-3/4}\|(\eta_l,\psi_l)\|_{L^2}, 
\end{align*}
where the final bound follows from our hypothesis 
\eqref{eq:APhiA0Phi0L4Close} on $(A-A_0,\Phi-\Phi_0)$, Lemma
\ref{lem:dBetaEst}, and the small-eigenvalue upper bounds
\eqref{eq:H2DeAEvalEst}. This gives the 
$L^{4/3}$ estimate \eqref{eq:H2DeA0EvectorEst-L43}.

For the $L^2$ estimate \eqref{eq:H2DeA0EvectorEst-L2}, the identity
\eqref{eq:ExpanddA+*OnCutoffEvectorEta} gives
\begin{align*}
&\|d_{A_0,\Phi_0}^{1,*}(\chi_0\eta_{N+1},\chi_0\psi_{N+1})\|_{L^2(X)} 
\\
&\le 
\|d\chi_0\|_{L^4}\|(\eta_{N+1},\psi_{N+1})\|_{L^4(\supp\chi_0)} 
\\
&\quad 
+ \|(A-A_0,\Phi-\Phi_0)\|_{L^4(\supp\chi_0)}
\|(\eta_{N+1},\psi_{N+1})\|_{L^4(\supp\chi_0)} 
\\
&\quad 
+ \|d_{A,\Phi}^{1,*}(\eta_{N+1},\psi_{N+1})\|_{L^2(\supp\chi_0)} 
\\
&\le c(\|d\chi_0\|_{L^4} + \|(A-A_0,\Phi-\Phi_0)\|_{L^4(\supp\chi_0)})
\|(\eta_{N+1},\psi_{N+1})\|_{L^2}
\\
&\quad 
+ \|d_{A,\Phi}^{1,*}(\eta_{N+1},\psi_{N+1})\|_{L^2(\supp\chi_0)}
\quad\text{(by Lemma \ref{lem:L21GlobalLocalEstDe2APhiEvec})}
\\
&\le c((-\log\lambda)^{-3/4} + \lambda^{1/3})\|(\eta_{N+1},\psi_{N+1})\|_{L^2}
+ \|d_{A,\Phi}^{1,*}(\eta_{N+1},\psi_{N+1})\|_{L^2(\supp\chi_0)},
\end{align*}
where the final inequality follows from our hypothesis 
\eqref{eq:APhiA0Phi0L4Close} on
$\|(A-A_0,\Phi-\Phi_0)\|_{L^4(\supp\chi_0)}$, and our estimates for
$d\chi_0$. This completes the proof of Lemma \ref{lem:H2DeA0EvectorEst}.
\end{proof}

\begin{lem}
\label{lem:H2DeApEvectorEst} 
Continue the above notation. Then there are positive constants $C$ and
$\lambda_0$ such that for all $\lambda \in (0,\lambda_0]$, the following
hold:
\begin{align}
\label{eq:H2DeApEvectorEst-L43} 
\|D_{A_p,g,p}^*\chi_p\psi_l\|_{L^{4/3}(S^4)} 
&\le 
C\lambda^{1/3}\|(\eta_l,\psi_l)\|_{L^2(S^4)}, \quad 1\le l\le N,
\\
\label{eq:H2DeApEvectorEst-L2} 
\|D_{A_p,g,p}^*\chi_p\psi_{N+1}\|_{L^2(S^4)} 
&\le 
\|d_{A,\Phi}^{1,*}(\eta_{N+1},\psi_{N+1})\|_{L^2(\supp\chi_p)}
\\
\notag
&\quad + C\lambda^{1/6}\|(\eta_{N+1},\psi_{N+1})\|_{L^2(X)}.
\end{align}
\end{lem}

\begin{proof}
Denote $(a_l,\phi_l) = d_{A,\Phi}^{1,*}(\eta_l,\psi_l)$, for $1\leq l\leq N+1$.
We have $(A,\Phi)=(A_p,0)+(0,\Phi)$ on $\supp\chi_p\Subset S^4\less\{s\}$
and so (schematically) the definition \eqref{eq:LAPhi*} of
$d_{A,\Phi}^{1,*}$ implies that for $1\leq l\leq N+1$,
$$
\phi_l = D_{A_p,g_p}^*\psi_l - \Phi\otimes\eta_l,
$$
and so, using 
$$
D_{A_p,g_p}^*(\chi_p\psi_l)
=
d\chi_p\otimes \psi_l + \chi_pD_{A_p,g_p}^*\psi_l,
$$
we see that (viewed equivalently either on $X$ or $S^4$) for $1\leq l\leq N+1$,
\begin{equation}
\label{eq:Expandd1*OnCutoffEvector0psiN+1}
D_{A_p,g_p}^*(\chi_p\psi_l)
= 
d\chi_p\otimes \psi_l + \chi_p\phi_l + \chi_p\Phi\otimes\eta_l.
\end{equation}
Thus, by the identity \eqref{eq:Expandd1*OnCutoffEvector0psiN+1} and
H\"older's inequality, we get for $1\leq l\leq N$,
\begin{align*}
&\|D_{A_p,g_p}^*(\chi_p\psi_l)\|_{L^{4/3}} 
\\
&\le 
\|d\chi_p\|_{L^2}\|\psi_l\|_{L^4(\supp d\chi_p)}
+ \|\chi_p\Phi\|_{L^4}\|\eta_l\|_{L^2(\supp \chi_p)}
+ \|\chi_p\phi_l\|_{L^{4/3}} 
\\
&\le 
(\|d\chi_p\|_{L^2} + \|\chi_p\Phi\|_{L^4})\|(\eta_l,\psi_l)\|_{L^2(X)}
+ \|d_{A,\Phi}^{1,*}(\eta_l,\psi_l)\|_{L^{4/3}}
\quad\text{(by Lemma \ref{lem:L21GlobalLocalEstDe2APhiEvec})}
\\
&\le 
C\left(\|d\chi_p\|_{L^2} + \|\chi_p\Phi\|_{L^4} 
+ c(\Vol B_p)^{1/4}\sqrt{\mu_l[A,\Phi]}\right)\|(\eta_l,\psi_l)\|_{L^2(X)}
\\
&\quad\text{(by Lemma \ref{lem:L21GlobalLocalEstDe2APhiEvec}
\& Inequality \eqref{eq:SqRootL43L2EvalueBound})}
\\
&\le 
C\lambda^{1/3}(1 + (-\log\lambda)^{-3/4})
\|(\eta_l,\psi_l)\|_{L^2},
\end{align*}
where the final inequality follows from the definition of $\chi_p$, our
hypothesis \eqref{eq:APhiA0Phi0L4Close}
on $\Phi|_{B_p}$, and the upper bounds \eqref{eq:H2DeAEvalEst} for
$\mu_l[A,\Phi]$. This establishes the $L^{4/3}$ estimate
\eqref{eq:H2DeApEvectorEst-L43} for $1\leq l\leq N$.

We turn to the proof of the estimate \eqref{eq:H2DeApEvectorEst-L2}. This
time, for $l=N+1$, the identity \eqref{eq:Expandd1*OnCutoffEvector0psiN+1} and
H\"older's inequality and the fact that $(a_l,\phi_l) =
d_{A,\Phi}^{1,*}(\eta_l,\psi_l)$ gives
\begin{align*}
&\|D_{A_p,g_p}^*(\chi_p\psi_l)\|_{L^2} 
\\
&\le 
\|d\chi_p\|_{L^{8/3}}\|\psi_l\|_{L^8(\supp d\chi_p)}
+ \|\chi_p\Phi\|_{L^4}\|\eta_l\|_{L^4(\supp \chi_p)}
+ \|\chi_p\phi_l\|_{L^2} 
\\
&\le 
(\|d\chi_p\|_{L^{8/3}} + \|\chi_p\Phi\|_{L^4})\|(\eta_l,\psi_l)\|_{L^2(X)}
+ \|d_{A,\Phi}^{1,*}(\eta_l,\psi_l)\|_{L^2(\supp \chi_p)}
\quad\text{(by Lemma \ref{lem:L21GlobalLocalEstDe2APhiEvec})}
\\
&\le 
C\left(\lambda^{1/6} + \lambda^{1/3}\right)\|(\eta_l,\psi_l)\|_{L^2}
+ \|d_{A,\Phi}^{1,*}(\eta_l,\psi_l)\|_{L^2(\supp \chi_p)},
\end{align*}
where the final inequality follows from the definition of $\chi_p$ and our
hypothesis \eqref{eq:APhiA0Phi0L4Close}
on $\Phi|_{B_p}$. This establishes the $L^2$ estimate
\eqref{eq:H2DeApEvectorEst-L2} and completes the proof
of Lemma \ref{lem:H2DeApEvectorEst}.
\end{proof}

Note again that the derivation of the estimates
\eqref{eq:H2DeApEvectorEst-L43} and \eqref{eq:H2DeApEvectorEst-L2} would be
difficult or impossible to obtain if the Dirac operator $D_{A_p,g_p}$ were
replaced by $D_{A_p,g_0}$, using $g_0$ --- the standard round metric on
$S^4$ --- rather than the approximately round metric $g_p$ induced from $g$
near $x_p\in X$ (see Remark \ref{rmk:WhyWeAssumeMetricS4IsAlmostRound}).

We consider the problem of estimating the orthogonal projection errors:

\begin{lem}
\label{lem:H2PerpProjA0Phi0EigenvecEst} 
Continue the above notation. Then there are positive constants $C$ and
$\lambda_0$ such that for all $\lambda \in (0,\lambda_0]$, the following
holds:
\begin{equation}
\label{eq:H2PerpProjA0Phi0EigenvecEst} 
\|\Pi_{A_0,\Phi_0}^\perp\chi_0(\eta_l,\psi_l)\|_{L^2(X)}
\leq
C(-\log\lambda)^{-3/4}\|(\eta_l,\psi_l)\|_{L^2(X)},
\quad 1\le l\le N.
\end{equation}
\end{lem}

\begin{proof}
We need to estimate $\Pi_{A_0,\Phi_0}^\perp\chi_0(\eta_l,\psi_l)$, for
$1\leq l\leq N$, where $(\eta_l,\psi_l)$ are orthonormal eigenvectors
corresponding to the first $N+1$ eigenvalues of the Laplacian
$d_{A,\Phi}^1d_{A,\Phi}^{1,*}$ and $\Pi_{A_0,\Phi_0}$ is $L^2$-orthogonal
projection from $L^2(X,\Lambda^+\otimes\fg_{E_\ell})\oplus L^2(X,V_\ell^-)$ onto
$\Ker d_{A_0,\Phi_0}^{1,*}$. Since $\Ker d_{A_0,\Phi_0}^{1,*}$ is a closed
subspace and $(\Ker d_{A_0,\Phi_0}^{1,*})^\perp = \Ran d_{A_0,\Phi_0}^1$,
there is a unique $d_{A_0,\Phi_0}^1(a_l,\phi_l) \in \Ran d_{A_0,\Phi_0}^1
\subset L^2(X,\Lambda^+\otimes\fg_{E_\ell})\oplus L^2(X,V_\ell^-)$ such that
$$
\Pi_{A_0,\Phi_0}^\perp\chi_0(\eta_l,\psi_l)
=
d_{A_0,\Phi_0}^1(a_l,\phi_l)
=
\frac{(\chi_0(\eta_l,\psi_l),d_{A_0,\Phi_0}^1(a_l,\phi_l))_{L^2}}
{\|d_{A_0,\Phi_0}^1(a_l,\phi_l)\|_{L^2}^2}d_{A_0,\Phi_0}^1(a_l,\phi_l).
$$
Hence, integrating by parts, we see that
$$
\|\Pi_{A_0,\Phi_0}^\perp\chi_0(\eta_l,\psi_l)\|_{L^2}
\leq
\|d_{A_0,\Phi_0}^{1,*}(\chi_0(\eta_l,\psi_l))\|_{L^{4/3}}
\frac{\|(a_l,\phi_l)\|_{L^4}}{\|d_{A_0,\Phi_0}^1(a_l,\phi_l)\|_{L^2}}.
$$
We may suppose without loss that $(a_l,\phi_l) \in (\Ker
d_{A_0,\Phi_0}^1)^\perp = \Ran d_{A_0,\Phi_0}^{1,*}$ and therefore
we can take $(a_l,\phi_l) = d_{A_0,\Phi_0}^{1,*}(v_l,s_l)$. Standard
elliptic estimates for the Laplacian
$d_{A_0,\Phi_0}^1d_{A_0,\Phi_0}^{1,*}$ then give
\begin{align*}
\|(a_l,\phi_l)\|_{L^4}
&\leq
c\|(a_l,\phi_l)\|_{L^2_{1,A_0}}
\leq
C\|(v_l,s_l)\|_{L^2_{2,A_0}}
\quad\text{(by Lemma \ref{lem:Kato})}
\\
&\leq
C(\|d_{A_0,\Phi_0}^1d_{A_0,\Phi_0}^{1,*}(v_l,s_l)\|_{L^2}
+ \|(v_l,s_l)\|_{L^2})
\\
&\leq
C(1+\mu_1[A_0,\Phi_0]^{-1})
\|d_{A_0,\Phi_0}^1d_{A_0,\Phi_0}^{1,*}(v_l,s_l)\|_{L^2}
\\
&=
C\|d_{A_0,\Phi_0}^1(a_l,\phi_l)\|_{L^2}.
\end{align*}
Hence, combining the preceding estimates, we see
$$
\|\Pi_{A_0,\Phi_0}^\perp\chi_0(\eta_l,\psi_l)\|_{L^2}
\leq
c\|d_{A_0,\Phi_0}^{1,*}(\chi_0(\eta_l,\psi_l))\|_{L^{4/3}}.
$$
The preceding estimate, together with inequality
\eqref{eq:H2DeA0EvectorEst-L43}, therefore yields
estimate \eqref{eq:H2PerpProjA0Phi0EigenvecEst} for $l=1,\dots, N$, and
completes the proof of Lemma  \ref{lem:H2PerpProjA0Phi0EigenvecEst}.
\end{proof}

\begin{lem}
\label{lem:H2PerpProjAp0EigenvecEst} 
Continue the above notation. Then there are positive constants $C$ and
$\lambda_0$ such that for all $\lambda \in (0,\lambda_0]$, the following
holds:
\begin{equation}
\label{eq:H2PerpProjAp0EigenvecEst} 
\|\Pi_{A_p}^\perp\chi_p\psi_l\|_{L^2(S^4)}
\leq
C\lambda^{1/3}\|(\eta_l,\psi_l)\|_{L^2(X)},
\quad 1\le l\le N.
\end{equation}
\end{lem}

\begin{proof}
We need to estimate $\Pi_{A_p}^\perp\chi_p\psi_l$, for
$1\leq l\leq N$, where $(\eta_l,\psi_l)$ are orthonormal eigenvectors
corresponding to the first $N+1$ eigenvalues of the Laplacian
$d_{A,\Phi}^1d_{A,\Phi}^{1,*}$ and $\Pi_{A_p}$ is $L^2$-orthogonal
projection from $L^2(S^4,V_p^-)$ onto
$\Ker D_{A_p,g_p}^*$. Since $\Ker D_{A_p,g_p}^*$ is a closed
subspace and $(\Ker D_{A_p,g_p}^*)^\perp = \Ran D_{A_p,g_p}$,
there is a unique $D_{A_p,g_p}\phi_{pl} \in \Ran D_{A_p,g_p}
\subset L^2(S^4,V_p^-)$ such that
$$
\Pi_{A_p}^\perp\chi_p\psi_l
= D_{A_p,g_p}\phi_{pl} =
\frac{(\chi_p\psi_l,D_{A_p,g_p}\phi_{pl})_{L^2}}
{\|D_{A_p,g_p}\phi_{pl}\|_{L^2}^2}D_{A_p,g_p}\phi_{pl}.
$$
Hence, integrating by parts, we see that
$$
\|\Pi_{A_p}^\perp\chi_p\psi_l\|_{L^2}
\leq
\|D_{A_p,g_p}^*(\chi_p\psi_l)\|_{L^{4/3}}
\frac{\|\phi_{pl}\|_{L^4}}{\|D_{A_p,g_p}\phi_{pl}\|_{L^2}}.
$$
On the other hand, we can bound $\|\phi_{pl}\|_{L^4}$ via
\begin{align*}
\|\phi_{pl}\|_{L^4}
&\leq
c\|\phi_{pl}\|_{L^2_{1,A_p}}
\quad\text{(by Lemma \ref{lem:Kato})}
\\
&\leq
C(\|D_{A_p,g_p}\phi_{pl}\|_{L^2} + \|\phi_{pl}\|_{L^2})
\quad\text{(by Lemma \ref{lem:L21AEstphi})}
\\
&\leq
C(1 + \nu_2[A_p]^{-1})\|D_{A_p,g_p}\phi_{pl}\|_{L^2}
\leq C\|D_{A_p,g_p}\phi_{pl}\|_{L^2},
\end{align*}
the final inequality following from Lemma \ref{lem:DiracS4VanishingKernel}
(together with our hypothesis that $\|F_{A_p}^{+,g_p}\|_{L^2(S^4)}$ is
small). Hence, combining the preceding estimates, we see
$$
\|\Pi_{A_p}^\perp\chi_p\psi_l\|_{L^2}
\leq
c\|D_{A_p,g_p}^*(\chi_p\psi_l)\|_{L^{4/3}}.
$$
The preceding estimate, together with inequality
\eqref{eq:H2DeApEvectorEst-L43}, yields
estimate \eqref{eq:H2PerpProjAp0EigenvecEst} for $l=1,\dots, N$, and
proves Lemma  \ref{lem:H2PerpProjAp0EigenvecEst}.
\end{proof}

\begin{proof}[Proof of Assertion (1) in Proposition
\ref{prop:H2SmallEvalApproxOrthonBasis}]
Combining Lemmas \ref{lem:H2PerpProjA0Phi0EigenvecEst} and 
\ref{lem:H2PerpProjAp0EigenvecEst}
proves the first assertion of Proposition
\ref{prop:H2SmallEvalApproxOrthonBasis}.
\end{proof}

We now turn to the remaining orthogonality assertions in Proposition
\ref{prop:H2SmallEvalApproxOrthonBasis}.

\begin{lem}
\label{lem:HowApproxOrthogEvectorEst} 
Continue the above notation. Then there are positive constants $C$ and
$\lambda_0$ such that for all $\lambda \in (0,\lambda_0]$, the following
hold:
\begin{equation}
\label{eq:HowApproxOrthogEvectorEst} 
\left|((\bfeta_l',\bpsi_l'),(\bfeta_k',\bpsi_k'))_{L^2}-\de_{kl}\right|
\le C\lambda^{2/3}, \quad 1\le k,l\le N+1.
\end{equation}
\end{lem}

\begin{proof}
{}From \eqref{eq:VectorCutoffEigenvectors} we see that we can write
\begin{align*}
(\bfeta_l',\bpsi_l')
&=
\left((\chi_0\eta_l,\chi_0\psi_l),
(0,\chi_1\psi_l),\dots,(0,\chi_m\psi_l)\right)
\\
&=
\left((\eta_l,\psi_l)_{B_0},
(\eta_l,\psi_l)_{B_1},\dots,(\eta_l,\psi_l)_{B_m}\right)
\\
&\quad + \left((1-\chi_0)(\eta_l,\psi_l)_{B_0},
(-\eta_l,(1-\chi_1)\psi_l)_{B_1},\dots,
(-\eta_l,(1-\chi_m)\psi_l)_{B_m}\right).
\end{align*}
Here, we write $X = U_0\cup\cup_{p=1}^m B_p$ (for simplicity --- the
constants defining $U_0$ in \eqref{eq:UZero} differ only by universal
factors),  with $U_0 = \supp \chi_0$ and
$(\eta_l,\psi_l)_{B_p}$ denotes the restriction of $(\eta_l,\psi_l)$ from
$X$ to $B_p$.  Remark that the cutoff functions $\chi_p$ are chosen so that
$d\chi_p$ is supported on annuli of the form
$(c_-\lambda^{1/3},c_+\lambda^p)$, where $0 < p\leq 1/3$ and the annulus
surrounds a curvature bubble with scale $\lambda\in (0,1]$, in order to
apply the local estimates for eigenspinors given in \cite{FeehanKato};
there is no loss in taking $p=1/3$ and fixed constants
$0<c_-<c_+<\8$. Hence, for $1\le k,l\le N+1$,
%\marginpar{\tiny consistent choice of $\chi_p$? Where is \cite{FeehanKato}
%being applied in this prop?}
\begin{align*}
((\bfeta_l',\bpsi_l'),(\bfeta_k',\bpsi_k'))_{L^2}
&=
\left( (\eta_l,\psi_l),(\eta_k,\psi_k) \right)_{L^2(U_0)}
+ \sum_{p=1}^m\left( (\eta_l,\psi_l),(\eta_k,\psi_k) \right)_{L^2(B_p)}
\quad\text{(Term 1)}
\\
&\quad 
+ 2\left( (\eta_l,\psi_l),
(1-\chi_0)(\eta_k,\psi_k) \right)_{L^2(U_0)}\quad\text{(Term 2)}
\\
&\quad 
+ 2\sum_{p=1}^m\left( (\eta_l,\psi_l),
(-\eta_k,(1-\chi_p)\psi_k) \right)_{L^2(B_p)}\quad\text{(Term 3)}
\\
&\quad 
+ \left( (1-\chi_0)(\eta_l,\psi_l), 
(1-\chi_0)(\eta_k,\psi_k) \right)_{L^2(U_0)}\quad\text{(Term 4)}
\\
&\quad 
+ \sum_{p=1}^m\left( (-\eta_l,(1-\chi_p)\psi_l), 
(-\eta_k,(1-\chi_p)\psi_k) \right)_{L^2(B_p)} \quad\text{(Term 5)}.
\end{align*}
Now, the first inner product (Term 1) is equal to
$$
\left( (\eta_l,\psi_l),(\eta_k,\psi_k) \right)_{L^2(U_0)}
+
\sum_{p=1}^m\left( (\eta_l,\psi_l),
(\eta_k,\psi_k) \right)_{L^2(B_p)}
=
((\eta_l,\psi_l),(\eta_k,\psi_k))_{L^2(X)}
=
\de_{kl}.
$$
The second inner product (Term 2) is bounded by twice 
\begin{align*}
&\|(\eta_l,\psi_l)\|_{L^4(U_0)}\|(1-\chi_0)(\eta_k,\psi_k)\|_{L^{4/3}(U_0)}
\\
&\leq
\|1-\chi_0\|_{L^2(U_0)}\|(\eta_k,\psi_k)\|_{L^4(U_0)}^2
\\
&\leq C\|1-\chi_0\|_{L^2(U_0)}\|(\eta_k,\psi_k)\|_{L^2(X)}^2 
\quad\text{(by Lemmas \ref{lem:Kato} and \ref{lem:L21GlobalLocalEstDe2APhiEvec})}
\\
&\leq C\lambda^{2/3},
\end{align*}
where the final inequality follows from the definition of $\chi_0$.
For $1\leq p\leq m$, the $p$th terms in the third inner product 
(Term 3) are bounded by twice
\begin{align*}
&(\eta_l,\eta_k)_{L^2(B_p)} 
+ ((1-\chi_p)\psi_l,(1-\chi_p)\psi_k)_{L^2(B_p)}
\\
&\leq
\|\eta_l\|_{L^4(X)}\|\eta_k\|_{L^4(X)}\cdot(\Vol B_p)^{1/2}
\\
&\quad + \|\psi_l\|_{L^4(B_p-\supp\chi_p)}\|\psi_k\|_{L^4(B_p-\supp\chi_p)}
\|1-\chi_p\|_{L^2(B_p)}
\\
&\leq
C\left((\Vol B_p)^{1/2} + \|1-\chi_p\|_{L^2(B_p)}\right)
\|(\eta,\psi)\|_{L^2(X)}
\quad\text{(by Lemmas \ref{lem:Kato} and \ref{lem:L21GlobalLocalEstDe2APhiEvec})}
\\
&\leq C\lambda^{2/3},
\end{align*}
the final inequality following from the definitions of $B_p$ and $\chi_p$.
An upper bound for the small eigenvalues, required by the preceding
applications of Lemma \ref{lem:L21GlobalLocalEstDe2APhiEvec}, is provided
by Proposition \ref{prop:H2SmallEvalUpperBound}.  The fourth and fifth
inner products are also bounded by $C\lambda^{2/3}$.  This completes the
proof of Lemma \ref{lem:HowApproxOrthogEvectorEst}.
\end{proof}

\begin{proof}[Proof of Assertion (2) in Proposition
\ref{prop:H2SmallEvalApproxOrthonBasis}]
We have, for $1\le k,l\le N$,
\begin{align*}
&\left|\left(\Pi_{\bA,\bPhi}(\bfeta_l',\bpsi_l'),
\Pi_{\bA,\bPhi}(\bfeta_k',\bpsi_k')\right)_{L^2}
-\de_{kl}\right|
\\
&\le \left|\left(\Pi_{\bA,\bPhi}(\bfeta_l',\bpsi_l'),
\Pi_{\bA,\bPhi}(\bfeta_k',\bpsi_k')\right)_{L^2}
-\left((\bfeta_l',\bpsi_l'),(\bfeta_k',\bpsi_k')\right)_{L^2}\right|
\\
&\quad + \left|\left((\bfeta_l',\bpsi_l'),(\bfeta_k',\bpsi_k')\right)_{L^2}
-\de_{kl}\right| 
\\
&\le \left|\left(\Pi_{\bA,\bPhi}^\perp(\bfeta_l',\bpsi_l'),
\Pi_{\bA,\bPhi}^\perp(\bfeta_k',\bpsi_k')\right)_{L^2}\right| + C\lambda^{2/3}
\quad\text{(by Lemma  \ref{lem:HowApproxOrthogEvectorEst})}
\\
&\le \|\Pi_{\bA,\bPhi}^\perp(\bfeta_l',\bpsi_l')\|_{L^2}
\|\Pi_{\bA,\bPhi}^\perp(\bfeta_k',\bpsi_k')\|_{L^2} + C\lambda^{2/3} 
\\
&\le C\left((-\log\lambda)^{-3/2} + \lambda^{2/3}\right)
\quad\text{(by Assertion (1) of Proposition \ref{prop:H2SmallEvalUpperBound})}
\end{align*}
and so
\begin{equation}
\label{eq:H2KerDeA0ApproxBasis}
\left|\left(\Pi_{\bA,\bPhi}(\bfeta_l',\bpsi_l'),
\Pi_{\bA,\bPhi}(\bfeta_k',\bpsi_k')\right)_{L^2}
-\de_{kl}\right| \le C(-\log\lambda)^{-3/2}, \quad 1\le k,l\le N.
\end{equation}
Thus, for small enough $\la$, the vectors
$\{\Pi_{\bA,\bPhi}(\bfeta_l',\bpsi_l')\}_{l=1}^{N}$ form an approximately
$L^2$-orthonormal basis for $\Ker T_{\bA,\bPhi}^*$. This proves the second
assertion of Proposition \ref{prop:H2SmallEvalApproxOrthonBasis}.
\end{proof}

\begin{proof}[Proof of Assertion (3) in Proposition
\ref{prop:H2SmallEvalApproxOrthonBasis}]
Finally, we turn to the problem of estimating the orthogonal projections of
$(\bfeta_{N+1}',\bpsi_{N+1}')$. Let
$\{(\bxi_l,\bvarphi_l)\}_{l=1}^{N}$ be an ordered, $L^2$-orthonormal basis
for $\Ker T_{\bA,\bPhi}^*$ obtained by applying Gram-Schmidt
orthonormalization to the basis
$\{\Pi_{\bA,\bPhi}(\bfeta_l',\bpsi_l')\}_{l=1}^{N}$:
\begin{equation}
\label{eq:H2GramSchmidt}
(\tilde\bxi_l,\tilde\bvarphi_l) 
= \Pi_{\bA,\bPhi}(\bfeta_l',\bpsi_l') -
\sum_{k=1}^{l-1}(\Pi_{\bA,\bPhi}(\bfeta_k',\bpsi_k'),
(\bxi_k,\bvarphi_k))_{L^2}\cdot
(\bxi_k,\bvarphi_k)
\end{equation}
with
$$
(\bxi_l,\bvarphi_l)
= 
\frac{(\tilde\bxi_l,\tilde\bvarphi_l)}
{\|(\tilde\bxi_l,\tilde\bvarphi_l)\|_{L^2}},
\quad 1\leq l \leq N.
$$
By the inner-product estimates in Lemma
\ref{lem:HowApproxOrthogEvectorEst} and the orthogonal-complement
estimates of Assertion (1) in Proposition \ref{prop:H2SmallEvalApproxOrthonBasis}
we have, for $1\leq k\leq N$,
\begin{equation}
\begin{aligned}
\label{eq:H2PerpPropjLastEigenvectorEst}
&|((\bfeta_{N+1}',\bpsi_{N+1}'),\Pi_{\bA,\bPhi}(\bfeta_k',\bpsi_k'))_{L^2}|
\\
&=
|((\bfeta_{N+1}',\bpsi_{N+1}'),(\bfeta_k',\bpsi_k'))_{L^2}
-((\bfeta_{N+1}',\bpsi_{N+1}'),
\Pi_{\bA,\bPhi}^\perp(\bfeta_k',\bpsi_k'))_{L^2}|
\\
&\leq
C\lambda^{2/3} + C(1+\lambda^{1/3})(-\log\lambda)^{-3/4}
\leq 
C(-\log\lambda)^{-3/4}.
\end{aligned}
\end{equation}
Consequently, by the Gram-Schmidt equations \eqref{eq:H2GramSchmidt}
together with the inner-product estimates
\eqref{eq:H2PerpPropjLastEigenvectorEst}, and the fact that
$$
\Pi_{\bA,\bPhi}(\bfeta_{N+1}',\bpsi_{N+1}') 
= 
\sum_{k=1}^{N}\left((\bfeta_{N+1}',\bpsi_{N+1}'),
(\bxi_k,\bvarphi_k)\right)_{L^2}\cdot(\bxi_k,\bvarphi_k),
$$
we have
\begin{equation}
\label{eq:H2KerDeA0ApproxPerpVec}
\begin{aligned}
\|\Pi_{\bA,\bPhi}(\bfeta_{N+1}',\bpsi_{N+1}')\|_{L^2} 
&\le C(-\log\lambda)^{-3/4},
\\
\|\Pi_{\bA,\bPhi}^\perp(\bfeta_{N+1}',\bpsi_{N+1}')\|_{L^2} 
&\ge 1-C(-\log\lambda)^{-3/4}.
\end{aligned}
\end{equation}
This proves the third and final assertion of Proposition 
\ref{prop:H2SmallEvalApproxOrthonBasis}.
\end{proof}

\subsection{A lower bound for the first non-small eigenvalue}
\label{subsec:Eigen2LowerBoundFirstNonSmallEigenvalue}
Lastly, we compute a lower bound for the $(N+1)$-st eigenvalue of
$d_{A,\Phi}^1d_{A,\Phi}^{1,*}$.
Recall that $(\eta_{N+1},\psi_{N+1})$ is the $L^2$-unit eigenvector
associated to the $(N+1)$-st eigenvalue $\mu_{N+1}[A,\Phi]$ of
$d_{A,\Phi}^1d_{A,\Phi}^{1,*}$.

The quantity $\nu_2[A_0,\Phi_0]$ is the least positive eigenvalue of the
Laplacian $d_{A_0,\Phi_0}^1d_{A_0,\Phi_0}^{1,*}$ over $X$, whereas
$\nu_2[A_p]$ is the least positive eigenvalue of
$D_{A_p,g_p}D_{A_p,g_p}^*$ over $S^4$. Hence, the least positive eigenvalue
of $T_{\bA,\bPhi}T_{\bA,\bPhi}^*$ is given by
$$
K =
\min\{\nu_2[A_0,\Phi_0],\nu_2[A_1],\dots,\nu_2[A_m]\}.
$$
We have:

\begin{prop}
\label{prop:H2SmallEvalLowerBound}
Continue the above notation. Then there are positive constants 
$C$ and $\lambda_0$ such that for all $\lambda \in (0,\lambda_0]$, the
following holds:
$$
\mu_{N+1}[A,\Phi] \ge K - C(-\log\lambda)^{-3/4}.
$$
\end{prop}

\begin{proof}[Proof of Theorem \ref{thm:H2SmallEval}, given
Proposition \ref{prop:H2SmallEvalLowerBound}]
  Combine Propositions \ref{prop:H2SmallEvalUpperBound} and
  \ref{prop:H2SmallEvalLowerBound}.
\end{proof}

\begin{proof}[Proof of Proposition \ref{prop:H2SmallEvalLowerBound}]
Recall that 
$$
(\bfeta_{N+1}',\bpsi_{N+1}') 
=
\left((\chi_0\eta_{N+1},\chi_0\psi_{N+1}),
\chi_1\psi_{N+1},\dots,\chi_m\psi_{N+1}\right).
$$
Hence,
$$
T_{\bA,\bPhi}^*(\bfeta_{N+1}',\bpsi_{N+1}')
=
d_{A_0,\Phi_0}^{1,*}(\chi_0\eta_{N+1},\chi_0\psi_{N+1})
\oplus
\mathop{\oplus}\limits_{p=1}^m
D_{A_p,g_p}^*\chi_p\psi_{N+1}.
$$
Therefore, because
$$
\|d_{A,\Phi}^{1,*}(\eta_{N+1},\psi_{N+1})\|_{L^2} 
=
\sqrt{\mu_{N+1}[A,\Phi]}\cdot \|(\eta_{N+1},\psi_{N+1})\|_{L^2}, 
$$
and $\|(\eta_{N+1},\psi_{N+1})\|_{L^2} = 1$, we get
\begin{align*}
&\|T_{\bA,\bPhi}^*(\bfeta_{N+1}',\bpsi_{N+1}')\|_{L^2}^2 
\\
&=
\|d_{A_0,\Phi_0}^{1,*}(\chi_0\eta_{N+1},\chi_0\psi_{N+1})\|_{L^2}^2
+ \sum_{p=1}^m\|D_{A_p,g_p}^*\chi_p\psi_{N+1}\|_{L^2}^2
\\
&\leq
\sum_{p=0}^m\|d_{A,\Phi}^{1,*}(\eta_{N+1},\psi_{N+1})\|_{L^2(\supp\chi_p)}^2 
+ C((-\log\lambda)^{-3/2}+\lambda^{1/3})
\\
&\quad + C(-\log\lambda)^{-3/4}
\|d_{A,\Phi}^{1,*}(\eta_{N+1},\psi_{N+1})\|_{L^2(\supp\chi_0)}^2 
\\
&\quad + \sum_{p=1}^m C\lambda^{1/6}
\|d_{A,\Phi}^{1,*}(\eta_{N+1},\psi_{N+1})\|_{L^2(\supp\chi_p)}^2,
\end{align*}
where the last inequality follows from
the estimates \eqref{eq:H2DeA0EvectorEst-L2} and
\eqref{eq:H2DeApEvectorEst-L2}. Therefore,
\begin{align*}
&\|T_{\bA,\bPhi}^*(\bfeta_{N+1}',\bpsi_{N+1}')\|_{L^2}^2 
\\
&\leq
(1+C(-\log\lambda)^{-3/4})\|d_{A,\Phi}^{1,*}(\eta_{N+1},\psi_{N+1})\|_{L^2(X)}^2 
+ C(-\log\lambda)^{-3/2}
\\
&= \mu_{N+1}[A,\Phi](1+C(-\log\lambda)^{-3/4}) + C(-\log\lambda)^{-3/2}.
\end{align*}
Hence, using the upper bound \eqref{eq:H2DeAEvalEst} 
for $\mu_{N+1}[A,\Phi]$, we obtain
$$
\|T_{\bA,\bPhi}^*(\bfeta_{N+1}',\bpsi_{N+1}')\|_{L^2} 
\leq \left(\sqrt{\mu_{N+1}[A,\Phi]} + C(-\log\lambda)^{-3/4}\right).
$$
In the other direction, we have the eigenvalue estimate
$$
\|T_{\bA,\bPhi}^*(\bfeta,\bpsi)\|_{L^2} 
\geq 
\sqrt{K}\cdot \|(\bfeta,\bpsi)\|_{L^2}
$$ 
for all $(\bfeta,\bpsi)\in(\Ker T_{\bA,\bPhi}^*)^\perp$
and so by the orthogonal projection estimates
\eqref{eq:H2KerDeA0ApproxPerpVec},   
\begin{align*}
\|T_{\bA,\bPhi}^*(\bfeta_{N+1}',\bpsi_{N+1}')\|_{L^2}
&\ge 
\sqrt{K}
\cdot\|\Pi_{\bA,\bPhi}^\perp(\bfeta_{N+1}',\bpsi_{N+1}')\|_{L^2} 
\\
&\ge 
\sqrt{K}(1 - C(-\log\lambda)^{-3/4}).
\end{align*}
Combining these upper and lower $L^2$ bounds for
$T_{\bA,\bPhi}^*(\bfeta_{N+1}',\bpsi_{N+1}')$ yields
$$
\sqrt{K}\left(1 - C(-\log\lambda)^{-3/4}\right)
\leq 
\sqrt{\mu_{N+1}[A,\Phi]} + C(-\log\lambda)^{-3/4},
$$
and this gives the lower bound on $\mu_{N+1}[A,\Phi]$ for small enough
$\la_0$ and all $\lambda \in (0,\lambda_0]$.
\end{proof}

\subsection{A map from the direct sum of kernels to the small-eigenvalue
eigenspace} 
\label{subsec:MapKernelToSmallEigenvalueSpace}
Though not necessary for the proof of Theorem \ref{thm:H2SmallEval} (whose
main conclusion is the lower-bound on the first non-small eigenvalue), it
is still convenient at this point to define a map from the direct sum of
the kernels of the operators $d_{A_0,\Phi_0}^{1,*}$ and $D_{A_p,g_p}^*$,
for $p=1,\dots,m$ (as given by the kernel of the diagonal operator
$T_{\bA,\bPhi}$), to the small-eigenvalue eigenspace of the operator
$d_{A,\Phi}^{1,*}$. This map will be important in \cite{FL4}, where we
prove the injectivity and surjectivity of the gluing maps to the extended
$\PU(2)$-monopole moduli space. It is essentially an inverse, modulo a
small error, to the map --- implicit in in the conclusion of Proposition
\ref{prop:H2SmallEvalApproxOrthonBasis} --- from the small-eigenvalue
eigenspace of the operator 
$d_{A,\Phi}^{1,*}$ to the direct sum of the kernels of the operators
$d_{A_0,\Phi_0}^{1,*}$ and $D_{A_p,g_p}^*$ (again as given by the kernel of
the diagonal operator $T_{\bA,\bPhi}$).

Most of the work has already been carried out in \S
\ref{subsec:Eigen2UpperBoundSmall}, where we used cut-and-paste techniques
to construct a set $\{(\xi_l',\varphi_l')\}_{l=1}^N$ of approximate
eigenvectors for the operator $d_{A,\Phi}^{1,*}$, which were approximately
$L^2$-orthonormal. It suffices here to estimate the error between these
approximate eigenvectors and their $L^2$-orthogonal projections onto the
small-eigenvalue eigenspace $\Ran\Pi_{A,\Phi,\mu}$ of $d_{A,\Phi}^{1,*}$
and then estimate their $L^2$ inner products.

\begin{prop}
\label{prop:SmallEvalPerpProjAPhiEigenvecEst} 
Continue the above notation. Then there are positive constants $C$ and
$\lambda_0$ such that for all $\lambda \in (0,\lambda_0]$, the following
holds:
\begin{align}
\label{eq:SmallEvalPerpProjAPhiEigenvecEst} 
\|\Pi_{A,\Phi,\mu}^\perp(\xi_l',\varphi_l')\|_{L^2(X)}
&\leq
C(-\log\lambda)^{-3/4},
\quad 1\le l\le N,
\\
\label{eq:SmallEvalProjAPhiEigenvecEst} 
\left|(\Pi_{A,\Phi,\mu}(\xi_k',\varphi_k'),
\Pi_{A,\Phi,\mu}(\xi_l',\varphi_l'))_{L^2(X)} - \delta_{kl}\right|
&\leq
C(-\log\lambda)^{-3/2},
\quad 1\le k,l\le N.
\end{align}
\end{prop}

\begin{proof}
We need to estimate $\Pi_{A,\Phi,\mu}^\perp(\xi_l',\varphi_l')$, for
$1\leq l\leq N$, where $(\xi_l',\varphi_l')$ are the ordered, approximately
$L^2$-orthonormal, approximate eigenvectors 
corresponding to the first $N$ eigenvalues of the Laplacian
$d_{A,\Phi}^1d_{A,\Phi}^{1,*}$ (see their construction in \S
\ref{subsec:Eigen2UpperBoundSmall}) 
and $\Pi_{A,\Phi,\mu}$ is $L^2$-orthogonal
projection from $L^2(X,\Lambda^+\otimes\fg_E)\oplus L^2(X,V^-)$ onto the
small-eigenvalue eigenspace of $d_{A,\Phi}^1d_{A,\Phi}^{1,*}$.

Since $\Ker d_{A,\Phi}^{1,*}$ is a closed subspace (contained in
$\Ran\Pi_{A,\Phi,\mu}$) and $(\Ker d_{A,\Phi}^{1,*})^\perp = \Ran
d_{A,\Phi}^1$, there is a unique
$$
d_{A,\Phi}^1(a_l,\phi_l) 
\in 
\Ran\Pi_{A,\Phi,\mu}^\perp\subset \Ran d_{A,\Phi}^1
\subset 
L^2(X,\Lambda^+\otimes\fg_E)\oplus L^2(X,V^-)
$$
such that
$$
\Pi_{A,\Phi,\mu}^\perp(\xi_l',\varphi_l')
=
d_{A,\Phi}^1(a_l,\phi_l)
=
\frac{((\xi_l',\varphi_l'),d_{A,\Phi}^1(a_l,\phi_l))_{L^2}}
{\|d_{A,\Phi}^1(a_l,\phi_l)\|_{L^2}^2}d_{A,\Phi}^1(a_l,\phi_l).
$$
Hence, integrating by parts, we see that
$$
\|\Pi_{A,\Phi,\mu}^\perp(\xi_l',\varphi_l')\|_{L^2}
\leq
\|d_{A,\Phi}^{1,*}(\xi_l',\varphi_l')\|_{L^2}
\frac{\|(a_l,\phi_l)\|_{L^2}}{\|d_{A,\Phi}^1(a_l,\phi_l)\|_{L^2}}.
$$
We may suppose without loss that 
$$
(a_l,\phi_l) 
\in 
(\Ker d_{A,\Phi}^1\circ\Pi_{A,\Phi,\mu}^\perp)^\perp 
= 
\Ran \Pi_{A,\Phi,\mu}^\perp\circ d_{A,\Phi}^{1,*},
$$ 
and therefore 
$$
\|(a_l,\phi_l)\|_{L^2} \leq \mu^{-1/2}\|d_{A,\Phi}^1(a_l,\phi_l)\|_{L^2}.
$$
Hence, combining the preceding estimates, we see
$$
\|\Pi_{A,\Phi,\mu}^\perp(\xi_l',\varphi_l')\|_{L^2}
\leq
c\mu^{-1/2}\|d_{A,\Phi}^{1,*}(\xi_l',\varphi_l')\|_{L^2}.
$$
This estimate, together with inequalities \eqref{eq:H2DeAEvectorEstX} (for
$p=0$ and $1\leq i\leq n$) and \eqref{eq:H2DeAEvectorEstS4} (for $1\leq
p\leq m$ and $1\leq i\leq 2\kappa_p$, recalling that
$N=n+2\sum_{p=1}^m\kappa_p$), therefore yields estimate
\eqref{eq:SmallEvalPerpProjAPhiEigenvecEst} for $l=1,\dots, N$.

The $L^2$ inner-product estimates \eqref{eq:SmallEvalProjAPhiEigenvecEst} 
now follow from our $L^2$ inner-product
estimates \eqref{eq:H2DeAEvectorInnerProductEst} 
for $((\xi_k',\varphi_k'),(\xi_l',\varphi_l'))_{L^2(X)}$
and from the $L^2$ orthogonal-projection estimates
\eqref{eq:SmallEvalPerpProjAPhiEigenvecEst}, using
$\Pi_{A,\Phi,\mu}(\xi_l',\varphi_l') =
(\xi_l',\varphi_l') - \Pi_{A,\Phi,\mu}^\perp(\xi_l',\varphi_l')$.
\end{proof}

%end of file
%file: existence.tex

\section{Existence of solutions to the extended PU(2)-monopole equations} 
\label{sec:Existence}
In this section we come to the heart of the matter, namely the
completion of the proof of Theorem \ref{thm:GluingTheorem1}. While our
proof builds on the results of the preceding sections --- by way of our
estimates in \S \ref{sec:Decay} for approximate solutions and our
eigenvalue estimates in \S \ref{sec:Eigenvalue} --- there remains the
substantial task of producing a suitable estimate for a right inverse
$P_{A,\Phi,\mu}$ of the linearization, $d_{A,\Phi}^1$, of the $\PU(2)$
monopole equations \eqref{eq:PT}. This estimate, given in Corollary
\ref{cor:L21AEstPAaphi}, is the principal technical result of our paper and
its proof draws heavily on the global estimates of \S \ref{sec:Global}, as
well as the decay estimates for harmonic spinors discussed in
\cite{FeehanKato} and some delicate applications in \S
\ref{sec:Dirichlet}. Rather 
than attempting to solve the $\PU(2)$ monopole equations \eqref{eq:PT}
directly, for which there are obstructions due to the presence of small
eigenvalues --- as discussed in \S \ref{sec:Eigenvalue} --- we instead
pursue a strategy pioneered by Taubes in \cite{TauIndef} and solve a weaker
version of these equations, the {\em extended $\PU(2)$ monopole
equations\/} (see \eqref{eq:QuickExtPUMonEqnForvpsi}). We set up these
extended equations in \S
\ref{subsec:MonoEqns}. Unlike in \cite{TauIndef}, we shall not attempt to
also solve the {\em obstruction equation\/} for $\PU(2)$ monopoles,
\eqref{eq:QuickExtPUMonEqnForvpsi}, as this will not be necessary for our
topological applications 
\cite{FLConj}, though this problem is nonetheless an interesting one in its
own right \cite{MrowkaPrincetonMorseTalk}.  Once a suitable estimate for
$P_{A,\Phi,\mu}$ is in place, which is the task of \S
\ref{subsec:GlobalEstLinear}, the solution to the extended $\PU(2)$
monopole equations then follows in a standard way by a fixed-point
argument, as we see in \S
\ref{subsec:ExistSolnExtASD}. Finally, we assemble the preceding
ingredients in \S \ref{subsec:ProofGluingTheorem} and finish the proof of 
Theorem \ref{thm:GluingTheorem1}.

\subsection{The extended PU(2) monopole equations}
\label{subsec:MonoEqns}
Given an $L^2_k$ approximate $\PU(2)$ monopole $(A,\Phi)$, our goal in this
section is to solve the quasi-linear system
\begin{equation}
\label{eq:BasicPTDef}
\fS(A+a,\Phi+\phi) 
=
\left(\begin{matrix}
(F^+(A+a))_0 - \tau\rho^{-1}((\Phi+\phi)\otimes(\Phi+\phi)^*)_{00} \\
D_{A+a}(\Phi+\phi) 
\end{matrix}\right)
= 0, 
\end{equation}
for a deformation 
$$
(a,\phi) \in L^2_k(\Lambda^1\otimes\fg_E)\oplus L^2_k(V^+)
$$
such that $(A+a,\Phi+\phi)$ is a $\PU(2)$ monopole.
The strategy, adapted from \cite[\S 7.2]{DK},
\cite{TauSelfDual}, and \cite[\S5]{TauStable},  
is to first compute an {\em a priori} $L^2_{1,A}$ estimate for a solution
$(a,\phi)$ to the linearization of \eqref{eq:BasicPTDef}
and prove existence of an $L^2_1$ solution
$(a,\phi)$ to the full non-linear equation via a Banach-space fixed-point
lemma.  In order for this estimate to be useful we need to
know the precise dependence of the constant on the pair $(A,\Phi)$. 
We then use a separate regularity argument to show that the solution
$(a,\phi)$ is $L^2_k$. From \cite[Equation (3.2)]{FL1}, we recall that
equation \eqref{eq:BasicPTDef} can be rewritten as a first-order, quasi-linear,
partial differential equation 
\begin{equation}
\label{eq:FirstPTDef}
d_{A,\Phi}^1(a,\phi) + \{(a,\phi),(a,\phi)\} = -\fS(A,\Phi),
\end{equation}
where the linearization 
$$
d_{A,\Phi}^1 = (D\fS)_{A,\Phi}:
\begin{matrix}
&L^2_k(\Lambda^1\otimes\fg_E)\\
&\oplus\\ 
&L^2_k(V^+)
\end{matrix} 
\to
\begin{matrix}
&L^2_{k-1}(\Lambda^+\otimes\fg_E)\\
&\oplus\\
&L^2_{k-1}(V^-)
\end{matrix} 
$$ 
is given by \cite[Equation (2.36)]{FL1}
\begin{equation}
\label{eq:LAPhi}
d_{A,\Phi}^1(a,\phi)
= 
\begin{pmatrix} 
d_{A}^+a - \tau\rho^{-1}(\Phi\otimes\phi^*+\phi\otimes{\Phi}^*)_{00} 
\\
D_{A}\phi + \rho(a)\Phi 
\end{pmatrix},
\end{equation}
and the quadratic term $\{(a,\phi),(a,\phi)\}$ is defined by the bilinear
form 
\begin{equation}
\label{eq:PUMonopoleBilinearTerm}
\{(a,\phi),(b,\varphi)\}
=
\begin{pmatrix}
(a\wedge b)^+ - \tau\rho^{-1}(\phi\otimes\varphi^*)_{00} \\
\rho(a)\varphi 
\end{pmatrix},
\end{equation}
for $(b,\varphi) \in L^2_k(\Lambda^1\otimes\fg_E)\oplus L^2_k(V^+)$.

We shall see in \cite{FL4}, using an argument similar to that of
\S \ref{sec:Eigenvalue}, that the composition of the
linearization of the splicing map with $L^2$-orthogonal projection induces
an isomorphism from tangent spaces to the bundle of gluing data onto the
kernel of $d_{A,\Phi}^{0,*} + d_{A,\Phi}^1$, where $d_{A,\Phi}^0$ is the
linearization \cite[Proposition 2.1 \& Equation (2.38)]{FL1} of the gauge
group action at a pair $(A,\Phi)$. Thus, it is natural to seek a
deformation $(a,\phi)$ which is $L^2$-orthogonal to the kernel of
$d_{A,\Phi}^1$ or, equivalently, a deformation $(a,\phi)$ lying the image
of $d_{A,\Phi}^{1,*}$, so
\begin{equation}
\label{eq:DeformationDefn}
(a,\phi) = d_{A,\Phi}^{1,*}(v,\psi),
\quad\text{for } 
(v,\psi) \in L^2_{k+1}(\Lambda^+\otimes\fg_E)\oplus L^2_{k+1}(V^-),
\end{equation}
where the first-order operator $d_{A,\Phi}^{1,*}$ is given by 
equation \eqref{eq:LAPhi*}.
Substituting \eqref{eq:DeformationDefn} into \eqref{eq:FirstPTDef} yields a
second-order, quasi-linear, elliptic partial differential equation
\begin{equation}
\label{eq:SecondPTDef}
d_{A,\Phi}^1d_{A,\Phi}^{1,*}(v,\psi) 
+ 
\{d_{A,\Phi}^{1,*}(v,\psi),d_{A,\Phi}^{1,*}(v,\psi)\} = -\fS(A,\Phi).
\end{equation}
Equation \eqref{eq:SecondPTDef} also serves as a gauge-fixing condition for
a solution $(a,\phi)$ to equation \eqref{eq:FirstPTDef} since
$$
d_{A,\Phi}^{0,*}\circ d_{A,\Phi}^{1,*}(v,\psi)
=
(d_{A,\Phi}^1\circ d_{A,\Phi}^0)^*(v,\psi)
=
(\fS(A,\Phi)\cdot\,)^*(v,\psi),
$$
and thus, for example, $d_{A,\Phi}^{0,*}(a,\phi) = 0$ if $(A,\Phi)$ is a
$\PU(2)$ monopole (see the remarks following equation (2.36) in
\cite{FL1}).

The Laplacian $d_{A,\Phi}^1d_{A,\Phi}^{1,*}$, as we have seen in \S
\ref{sec:Eigenvalue}, will in general have {\em small eigenvalues\/} and so
--- even if the Laplacian has no kernel --- it will not be uniformly
invertible, in the sense that the small eigenvalues tend to zero when we
try to make $\fS(A,\Phi)$ small enough to solve the non-linear equation
\eqref{eq:SecondPTDef} by allowing $[A,\Phi]$ to move closer to an
ideal pair $[A_0,\Phi_0,\bx] \in\sC(\ft_\ell)\times\Sym^{\ell}$. In a
gauge-theory context, such small-eigenvalue problems were first addressed
by Taubes in \cite{TauIndef}, \cite{TauStable} in the case of the
anti-self-dual equation. Thus, in place of the $\PU(2)$ monopole equation
\begin{equation}
\label{eq:QuickPUMonEqnForvpsi}
\fS((A,\Phi)+(a,\phi)) = 0
\end{equation}
for a solution $(a,\phi) = d_{A,\Phi}^{1,*}(v,\psi)$, we follow
\cite{TauIndef}, \cite{TauStable} and consider the pair of equations
\begin{align}
\label{eq:QuickExtPUMonEqnForvpsi}
\Pi_{A,\Phi,\mu}^\perp\fS((A,\Phi)+(a,\phi)) &= 0,
\\
\label{eq:QuickFinitePartPUMonEqnForvpsi}
\Pi_{A,\Phi,\mu}\fS((A,\Phi)+(a,\phi)) &= 0,
\end{align}
for a solution $(a,\phi) = d_{A,\Phi}^{1,*}(v,\psi)$, with
$$
(v,\psi) \in L^2_{k+1}(\Lambda^+\otimes\fg_E)\oplus L^2_{k+1}(V^-)
$$ 
satisfying the constraint 
\begin{equation}
\label{eq:vpsiFiniteRankConstraint}
\Pi_{A,\Phi,\mu}(v,\psi)=0.
\end{equation}
Equation \eqref{eq:QuickExtPUMonEqnForvpsi} is called the {\em extended
$\PU(2)$ monopole equation\/}. Here, $\Pi_{A,\Phi,\mu}$ is the
$L^2$-orthogonal projection from $L^2(\Lambda^+\otimes\fg_E)\oplus
L^2(V^-)$ onto the finite-dimensional subspace spanned by the
eigenvectors of $d_{A,\Phi}^1d_{A,\Phi}^{1,*}$ with eigenvalues in
$[0,\mu]$, where $\mu$ is a positive constant (an upper bound for the small
eigenvalues), and $\Pi_{A,\Phi,\mu}^\perp := \id - \Pi_{A,\Phi,\mu}$.

Because of the small-eigenvalue problem, we shall only explicitly solve
equation \eqref{eq:QuickExtPUMonEqnForvpsi}; the true solutions to
equation 
\eqref{eq:QuickPUMonEqnForvpsi} are then simply the solutions to equation
\eqref{eq:QuickExtPUMonEqnForvpsi} which also satisfy equation
\eqref{eq:QuickFinitePartPUMonEqnForvpsi}.

Equation \eqref{eq:QuickExtPUMonEqnForvpsi} may be rewritten as
\begin{equation}
\label{eq:LongSecOrderExtPUMonEqnForvpsi}
\Pi_{A,\Phi,\mu}^\perp\left(d_{A,\Phi}^1d_{A,\Phi}^{1,*}(v,\psi) 
+ \{d_{A,\Phi}^{1,*}(v,\psi), d_{A,\Phi}^{1,*}(v,\psi)\}\right)
= - \Pi_{A,\Phi,\mu}^\perp\fS(A,\Phi),
\end{equation}
where again $(v,\psi)$ obeys the constraint
\eqref{eq:vpsiFiniteRankConstraint}.  For the
purposes of the regularity theory of \S \ref{sec:Regularity}, it is useful
to rewrite equation
\eqref{eq:LongSecOrderExtPUMonEqnForvpsi} as
%\begin{equation}
\begin{align}
\label{eq:RegLongSecOrderExtPUMonEqnForvpsi}
&d_{A,\Phi}^1d_{A,\Phi}^{1,*}(v,\psi) 
+ \{d_{A,\Phi}^{1,*}(v,\psi), d_{A,\Phi}^{1,*}(v,\psi)\}
\\
\notag
&\quad =
\Pi_{A,\Phi,\mu}
\left(\{d_{A,\Phi}^{1,*}(v,\psi), d_{A,\Phi}^{1,*}(v,\psi)\}\right)
- \Pi_{A,\Phi,\mu}^\perp \fS(A,\Phi).
\end{align}
%\end{equation}
In view of equation \eqref{eq:RegLongSecOrderExtPUMonEqnForvpsi}, it will
be useful to consider estimates for solutions to the following
second-order, quasi-linear elliptic equation,
\begin{equation}
\label{eq:RegLongSecOrderQuasiLinEqnForvpsi}
d_{A,\Phi}^1d_{A,\Phi}^{1,*}(v,\psi) 
+ \{d_{A,\Phi}^{1,*}(v,\psi), d_{A,\Phi}^{1,*}(v,\psi)\}
=
(w,s).
\end{equation}
Of course, equation \eqref{eq:RegLongSecOrderQuasiLinEqnForvpsi} reduces to
\eqref{eq:RegLongSecOrderExtPUMonEqnForvpsi} when $(w,s)\in
\Gamma(\Lambda^+\otimes\fg_E)\oplus\Gamma(V^-)$ is given by the right-hand
side of equation \eqref{eq:RegLongSecOrderExtPUMonEqnForvpsi}. 

Equation \eqref{eq:QuickFinitePartPUMonEqnForvpsi}
is called the {\em $\PU(2)$ monopole obstruction
equation\/} and takes the form
\begin{equation}
\label{eq:LongSecOrderFinitePartPUMonEqnForvpsi}
\Pi_{A,\Phi,\mu}
\left(\{d_{A,\Phi}^{1,*}(v,\psi), d_{A,\Phi}^{1,*}(v,\psi)\}
- \fS(A,\Phi)\right)
=
0,
\end{equation}
as we can easily see.

On the other hand, for the purposes of applying Banach-space fixed-point
theory to solve the non-linear equation
\eqref{eq:LongSecOrderExtPUMonEqnForvpsi}, it is convenient to write
$$
(v,\psi) = G_{A,\Phi,\mu}(\xi,\varphi) 
\quad\text{for some } 
(\xi,\varphi) \in 
L^2_{k-1}(\Lambda^+\otimes\fg_E)\oplus L^2_{k-1}(V^-),
$$
where 
$$
G_{A,\Phi,\mu}:
\begin{matrix}
&L^2_{k-1}(\Lambda^+\otimes\fg_E)\\
&\oplus\\ 
&L^2_{k-1}(V^-)
\end{matrix} 
\to
\begin{matrix}
&L^2_{k+1}(\Lambda^+\otimes\fg_E)\\
&\oplus\\
&L^2_{k+1}(V^-)
\end{matrix} 
$$
is the Green's operator for the Laplacian
$$
d_{A,\Phi}^1d_{A,\Phi}^{1,*}\circ\Pi_{A,\Phi,\mu}^\perp: 
\begin{matrix}
&L^2_{k+1}(\Lambda^+\otimes\fg_E)\\
&\oplus\\ 
&L^2_{k+1}(V^-)
\end{matrix} 
\to
\begin{matrix}
&L^2_{k-1}(\Lambda^+\otimes\fg_E)\\
&\oplus\\
&L^2_{k-1}(V^-)
\end{matrix} 
$$
Thus, $d_{A,\Phi}^1d_{A,\Phi}^{1,*}G_{A,\Phi,\mu} =
\Pi_{A,\Phi,\mu}^\perp$ and $\Pi_{A,\Phi,\mu}(v,\psi) =
\Pi_{A,\Phi,\mu}G_{A,\Phi,\mu}(\xi,\varphi) = 0$ for any $(\xi,\varphi)
\in L^2_{k-1}(\Lambda^+\otimes\fg_E)\oplus L^2_{k-1}(V^-)$, so the
constraint \eqref{eq:vpsiFiniteRankConstraint} is obeyed. Therefore, setting
$$
P_{A,\Phi,\mu}
=
d_{A,\Phi}^{1,*}G_{A,\Phi,\mu},
$$
it suffices to solve for 
$$
(\xi,\varphi)
\in
L^2_{k-1}(\Lambda^+\otimes\fg_E)\oplus L^2_{k-1}(V^-),
$$ 
such that
\begin{equation}
\label{eq:IntegralPUMonEqnForxipsi}
(\xi,\varphi) 
+ \{P_{A,\Phi,\mu}(\xi,\varphi), P_{A,\Phi,\mu}(\xi,\varphi)\}
= - \fS(A,\Phi),
\end{equation}
because, noting that $\Pi_{A,\Phi,\mu}^\perp(\xi,\varphi) =
\Pi_{A,\Phi,\mu}^\perp d_{A,\Phi}^1P_{A,\Phi,\mu}
(\xi,\varphi)$, 
the resulting $(v,\psi) = G_{A,\Phi,\mu}(\xi,\varphi)$ solves
equation \eqref{eq:LongSecOrderExtPUMonEqnForvpsi}, namely
$$ 
\Pi_{A,\Phi,\mu}^\perp
\left((\xi,\varphi) 
+ \{P_{A,\Phi,\mu}(\xi,\varphi), P_{A,\Phi,\mu}(\xi,\varphi)\}
+ \fS(A,\Phi)\right) = 0,
$$
as one can see by just applying $P_{A,\Phi,\mu}^\perp$ to 
equation \eqref{eq:IntegralPUMonEqnForxipsi}, 
{\em as well as} the constraint $\Pi_{A,\Phi,\mu}(v,\psi) =
\Pi_{A,\Phi,\mu}G_{A,\Phi,\mu}(\xi,\varphi) = 0$.
 
Noting that $\Pi_{A,\Phi,\mu}P_{A,\Phi,\mu}(\xi,\varphi) = 0$, the
solution $(\xi,\varphi)$ to \eqref{eq:IntegralPUMonEqnForxipsi} solves equation
\eqref{eq:QuickFinitePartPUMonEqnForvpsi} or equivalently
equation \eqref{eq:LongSecOrderFinitePartPUMonEqnForvpsi}, which now takes
the shape  
$$
\Pi_{A,\Phi,\mu}
\left(\{P_{A,\Phi,\mu}(\xi,\varphi), P_{A,\Phi,\mu}(\xi,\varphi)\} 
+ \fS(A,\Phi)\right) = 0,
$$
{\em if and only if\/}
\begin{equation}
\label{eq:IntegralFinitePUMonEqnForxipsi}
\Pi_{A,\Phi,\mu}(\xi,\varphi) = 0.
\end{equation}
For applications of Banach-space fixed-point theory, the integral form
\eqref{eq:IntegralPUMonEqnForxipsi} of the extended anti-self-dual equation 
is a little more convenient, while for the purposes of regularity theory
and localization, we use the ``differential'' form
\eqref{eq:RegLongSecOrderExtPUMonEqnForvpsi}. 

\begin{rmk}
\label{rmk:DifferentEigenvalueCutoffs}
If $(a,\phi)$ is a solution to the extended
$\PU(2)$-monopole equations \eqref{eq:QuickExtPUMonEqnForvpsi}
for the small-eigenvalue cutoff parameter
$\mu$, then it is also a solution to the (weaker) equations defined
by a parameter $\mu' \geq \mu$. Indeed, if 
$$
\Pi_{A,\Phi,\mu}^\perp\fS((A,\Phi)+(a,\phi)) = 0,
$$
then we also have 
\begin{align*}
0 
&=
\Pi_{A,\Phi,\mu'}^\perp\Pi_{A,\Phi,\mu}^\perp\fS((A,\Phi)+(a,\phi))
\\
&=
\Pi_{A,\Phi,\mu'}^\perp(1-\Pi_{A,\Phi,\mu})\fS((A,\Phi)+(a,\phi))
\\
&=
\Pi_{A,\Phi,\mu'}^\perp\fS((A,\Phi)+(a,\phi)),
\end{align*}
where the final equality follows 
because $\Ran\Pi_{A,\Phi,\mu}\subseteq \Ran\Pi_{A,\Phi,\mu'}$.
This point will be important in \cite{FLConj}.
\end{rmk}

\subsection{Comments on intrinsic, extended PU(2) monopole equations}\
\label{subsec:CommentIntrinsicExtendedPU2Eqn}
One difficulty with using the extended $\PU(2)$ monopole equations, as
we shall explain further in \cite{FLConj}, is that the equations
\eqref{eq:QuickExtPUMonEqnForvpsi} and 
\eqref{eq:QuickFinitePartPUMonEqnForvpsi}
are not intrinsic to the configuration space $\sC(\ft)$, unlike their
counterpart \eqref{eq:BasicPTDef} in the unobstructed case. Hence, solutions
$(A+a,\Phi+\phi)$ to \eqref{eq:QuickExtPUMonEqnForvpsi} corresponding to
bundles of gluing data defined by the different strata $\Sigma$ of the
symmetric product $\Sym^\ell(X)$ will not fit together (even aside from
the jumping-line or spectral flow problem) to form a smooth submanifold of
$\sC^{*,0}(\ft)$. 

To circumvent this problem we could try to use a different version of
equations \eqref{eq:QuickExtPUMonEqnForvpsi} and 
\eqref{eq:QuickFinitePartPUMonEqnForvpsi}, where the spectral projections
are replaced by ones which are intrinsic to the configuration space $\sC(\ft)$
rather than a bundle of gluing data defined by a single stratum
$\Sigma\subset\Sym^\ell(X)$. Indeed, to describe solutions to the $\SU(2)$
anti-self-dual equation near the product connection, an ``intrinsic,
extended anti-self-dual equation'' was introduced in \cite{DonPoly} and it
is asserted (without proof) in \cite{FrM} that Taubes' approach to gluing
for his version of the (non-intrinsic) extended anti-self-dual equation can
be adapted to the intrinsic, extended equation. Thus inspired by
\cite{DonPoly} and \cite{FrM}, we could consider the pair of equations
\begin{align}
\label{eq:QuickIntrinsicExtPUMonEqnForvpsi}
\Pi_{A,\Phi,\mu}^\perp\fS(A,\Phi) &= 0,
\\
\label{eq:QuickIntrinsicFinitePartPUMonEqnForvpsi}
\Pi_{A,\Phi,\mu}\fS(A,\Phi) &= 0.
\end{align}
Equation \eqref{eq:QuickIntrinsicExtPUMonEqnForvpsi} is called the {\em
intrinsic, extended $\PU(2)$ monopole equation\/}. Here,
$\Pi_{A,\Phi,\mu}$ is the $L^2$-orthogonal projection from
$L^2(\Lambda^+\otimes\fg_E)\oplus L^2(V^-)$ onto the
finite-dimensional subspace spanned by the eigenvectors of
$d_{A,\Phi}^1d_{A,\Phi}^{1,*}$ with eigenvalues in $[0,\mu]$, where
$\mu$ is a positive constant (an upper bound for the small eigenvalues),
and $\Pi_{A,\Phi,\mu}^\perp := \id - \Pi_{A,\Phi,\mu}$.  

An unsatisfactory feature of equation
\eqref{eq:QuickIntrinsicExtPUMonEqnForvpsi} is that the projections
$\Pi_{A,\Phi,\mu} = \Pi_{A'+a,\Phi'+\phi,\mu}$ depend on the deformation
sought from the initial approximate solution, $(A',\Phi')$, and this
complicates the analysis considerably in comparison with the simpler
equation \eqref{eq:QuickExtPUMonEqnForvpsi}, where the projection
$\Pi_{A',\Phi',\mu}$ depends only on the data, $(A',\Phi')$, and not the
deformation sought, $(a,\phi)$. While it may be possible to work around
this problem for the anti-self-dual equation, it appears difficult for
$\PU(2)$ monopoles. For example, if we tried to use the right-inverse
$P_{A,\Phi,\mu}$ for $d_{A,\Phi}^1$ then we would encounter a fairly
serious problem in adapting the proof of Theorem \ref{thm:L21AEstPAaphi}
(and Corollary \ref{cor:L21AEstPAaphi}) --- where we
currently derive an estimate for $P_{A',\Phi',\mu}$ --- as we assumed that
$\Phi'\equiv 0$ where the connection $A'$ bubbles: this assumption is used
in a crucial way in equation \eqref{eq:LpBoundWorstTerm}. On the other
hand, if we use the right-inverse $P_{A',\Phi',\mu}$ for $d_{A',\Phi'}^1$
then we would have to compare the projections $\Pi_{A,\Phi,\mu}$ and
$\Pi_{A',\Phi',\mu}$. Because of the nature of the Bochner formulas involved,
namely \eqref{eq:BW+} and \eqref{eq:BWDirac+}, the two projections appear
surprisingly difficult to compare when $(A,\Phi) = (A'+a,\Phi'+\phi)$ and
we know at most that $\|(a,\phi)\|_{L^2_{1,A'}}$ is small.   

The strategy we use instead in \cite{FLConj} is to construct a space of global
gluing data, $\bGl(\sU_{\mu,\ell},\lambda_0)$ over
$\sU_{\mu,\ell}\times\Sym^{\ell}(X)$, a global splicing map
$\bgamma_{\mu}'$, and then deform $\bgamma_{\mu}'$ to a global gluing map
$\bgamma_{\mu}$ using the results of the present article.

\subsection{Comments on the gluing methods of Taubes, Donaldson, and Mrowka}
\label{sec:GreenOp}
It is convenient at this point to indicate why the method employed by Taubes in
\cite{TauSelfDual},  \cite{TauIndef}, \cite{TauStable} for solving the
anti-self-dual equation, 
$$
d_A^+a + (a\wedge a)^+ = - F_A^+,
$$
does not apply to the $\PU(2)$ monopole equation with merely routine
modifications. His method requires the Bochner-Weitzenb\"ock formula
\eqref{eq:BW+} for the Laplacian $d_A^+d_A^{+,*}$ in order to obtain the
necessary {\em a priori} $L^2_{1,A}$ estimate for $a$. However, as we can
see from the Bochner formulas in \S \ref{sec:Global}, the corresponding
Laplacian in the case of the $\PU(2)$ monopole equation is given by
$$
d_{A,\Phi}^1d_{A,\Phi}^{1,*}
=\left(\begin{matrix}
d_A^+d_A^{+,*} & 0 \\
0 & D_AD_A^*
\end{matrix}\right)
+ 
\begin{pmatrix}
r^{11}_{A,\Phi} & r^{12}_{A,\Phi} \\
r^{21}_{A,\Phi} & r^{22}_{A,\Phi} 
\end{pmatrix},
$$
where the remainder term $R_{A,\Phi} = (r^{ij}_{A,\Phi})$ is a first-order
linear differential operator, with coefficients containing up through
first-order derivatives of $A$ and $\Phi$, and is diagonal only if
$\Phi\equiv 0$.

Thus, while the term $F_A^+$ in the Bochner formula for $d_A^+d_A^{+,*}$
will be bounded in $L^\8$ if $g$ is conformally flat near the gluing sites
and, more generally, $L^{\sharp,2}$-small as desired (for small enough
scales), the term $F_A^-$ in the formula for $D_AD_A^*$ will only be
bounded in $L^{\sharp,2}$ if $(A,\Phi)$ is an approximate $\PU(2)$ monopole
obtained by the splicing construction. In particular, the term
$F_A^-$ is not $L^{\sharp,2}$-small near the gluing sites and so the
rearrangement arguments used to obtain uniform, global elliptic estimates
for $d_A^+d_A^{+,*}$ in \cite[\S 5]{TauStable} or \cite[\S 5]{FeehanSlice}
no longer work.

On the other hand, the method for solving the anti-self-dual equation
employed by Donaldson in \cite{DonConn}, \cite{DK} --- which avoids the use
of Bochner-Weitzenb\"ock formulas by exploiting the conformal invariance of
the anti-self-dual equation---does not apply without significant
modification to the $\PU(2)$ monopole equation since the latter is not
conformally invariant \cite{FL1}. Much the same problem would arise if we
sought to replace the small annuli surrounding points of curvature
concentration with long cylinders --- by analogy with the method described
by Mrowka in
\cite{MrowkaThesis}. Of course, we could (with either conformal model)
transform the equations on $(X,g)$ to an equivalent set of equations, where
the metrics now vary with the tube lengths $T_i$ or neck radii parameters
$\lambda_i$. While such a strategy might circumvent some of the analytical
difficulties associated with the unwanted appearance of $F_A^-$ in the
Bochner formulas, it would lead to models for moduli space ends where the
metric varies (albeit conformally) with the location of the gluing centers
$\bx \in \Sym^\ell(X)$. As will be clearer from our discussion of the
requisite intersection theory in \cite{FLConj}, such behavior seems
undesirable because it makes it more difficult to patch together the
local gluing models (for moduli space ends) and thus form a global model
for the end of the moduli space $\barM(\ft)$ near a level
$M^{\sw}(\fs)\times\Sym^\ell(X)$. 

\subsection{Global estimates: linear theory}
\label{subsec:GlobalEstLinear}
Our next task is to parlay the global elliptic estimates for $d_A^+$ and
$D_A$ in the preceding section into global elliptic estimates for
$d_{A,\Phi}^1$ and its partial right inverse $P_{A,\Phi,\mu}$. The
following elliptic estimates, coupled with our bounds on the small
eigenvalues, form the technical heart of our main gluing result, Theorem
\ref{thm:GluingTheorem1}. We begin with an $L^2_{1,A}$ elliptic estimate for
$d_{A,\Phi}^{1,*}(v,\psi)$:
%\marginpar{\tiny Is this constant dependence statement below correct?} 

\begin{thm}
\label{thm:L21AEstPAaphi}
Let $X$ be a closed, oriented, smooth four-manifold with Riemannian metric $g$.
Given the data of the splicing construction in \S \ref{sec:Splicing}, there
are constants  
\begin{equation}
\label{eq:L21AEstPAaphiConstant}
C 
=
C(\sU_{\ell,\mu},\kappa,\|F_{A_d}\|_{L^\8},\|F_{A_e}\|_{L^\8},
\|\vartheta\|_{L^\8},\mu) < \8,
\end{equation} 
and $\lambda_0 = \lambda_0(C)>0$ such that the following holds.
Let $(A,\Phi)$ be an $L^2_4$ pair on $(\fg_E,V^+)$ produced by the splicing
construction of \S \ref{sec:Splicing}. Then
for all $(v,\psi)\in L^2_3(\Lambda^+\otimes \fg_E)\oplus
L^2_3(V^-)$, 
\begin{equation}
\label{eq:L21AEstPAaphiBound}
\|d_{A,\Phi}^{1,*}(v,\psi)\|_{L^2_{1,A}(X)} 
\le 
C(\|d_{A,\Phi}^1d_{A,\Phi}^{1,*}(v,\psi)\|_{L^{\sharp,2;2}(X)} 
+ \|(v,\psi)\|_{L^2(X)}). 
\end{equation}
\end{thm}

\begin{proof}
Let $U'\subset X$ be the open subset where $|\Phi|>0$, so $\bar U' =
\supp\Phi$.  Denote $(a,\phi) = d_{A,\Phi}^{1,*}(v,\psi)$.  Since $a =
d_A^{+,*}v + (\cdot\Phi)^*\psi$ and, schematically,
$$
\cov_A((\cdot\Phi)^*\psi) 
= 
\cov_A\Phi\otimes\psi + \Phi\otimes\cov_A\psi, 
$$
then H\"older's inequality yields
\begin{align*}
\|a\|_{L^2_{1,A}} 
&\leq 
\|d_A^{+,*}v\|_{L^2_{1,A}} 
+ \|\cov_A((\cdot\Phi)^*\psi)\|_{L^2}
+ \|\Phi\|_{L^\8}\|\psi\|_{L^2(U')}
\\
&\leq
\|d_A^{+,*}v\|_{L^2_{1,A}} 
+ \|\cov_A\Phi\|_{L^4}\|\psi\|_{L^4(U')}
+ \|\Phi\|_{L^\8}\|\cov_A\psi\|_{L^2(U')}
+ \|\Phi\|_{L^\8}\|\psi\|_{L^2(U')}.
\end{align*}
Thus, by Lemma \ref{lem:Kato}, 
\begin{equation}
\label{eq:IntermedL21a-1}
\|a\|_{L^2_{1,A}} 
\leq
\|d_A^{+,*}v\|_{L^2_{1,A}} 
+ c\|\Phi\|_{L^\8\cap L^4_{1,A}}\|\psi\|_{L^2_{1,A}(U')}.
\end{equation}
Using integration by parts, we see that
\begin{align*}
\|d_A^{+,*}v\|_{L^2}^2 
&= 
(d_A^{+,*}v,d_A^{+,*}v)_{L^2} 
= 
(d_A^+d_A^{+,*}v,v)_{L^2} 
\\
&\le 
\|d_A^+d_A^{+,*}v\|_{L^2}\|v\|_{L^2} 
\le 
\frac{1}{2}(\|d_A^+d_A^{+,*}v\|_{L^2}^2 + \|v\|_{L^2}^2),
\end{align*}
and so
\begin{equation}
\label{eq:EasyEstdA*v}
\|d_A^{+,*}v\|_{L^2} 
= 
\frac{1}{\sqrt{2}}(\|d_A^+d_A^{+,*}v\|_{L^2} + \|v\|_{L^2}). 
\end{equation}
Corollary \ref{cor:L21AEstdA*v} gives the following $L^2_{1,A}$ estimate for
$d_A^{+,*}v$: 
$$
\|d_A^{+,*}v\|_{L^2_{1,A}} 
\le  
C(\|d_A^+d_A^{+,*}v\|_{L^{\sharp,2}} + \|v\|_{L^2}).
$$
Recall that 
\begin{equation}
\label{eq:LapvInTermsOfapsi}
d_A^+d_A^{+,*}v = d_A^+a - d_A^+(\cdot\Phi)^*\psi,
\end{equation}
where, schematically,
$$
d_A^+((\cdot\Phi)^*\psi) 
= 
\cov_A\Phi\otimes\psi + \Phi\otimes\cov_A\psi.
$$
This gives, for any $2<q\leq \8$, using the embeddings $L^2_1\subset
L^{2\sharp,4}$ and $L^q\subset L^\sharp$ of \cite[Lemma 4.1]{FeehanSlice},
\begin{equation}
\label{eq:LpBoundWorstTerm}
\begin{aligned}
\|d_A^+(\cdot\Phi)^*\psi\|_{L^{\sharp,2}}
&\leq
\|\cov_A\Phi\|_{L^{2\sharp,4}}\|\psi\|_{L^{2\sharp,4}(U')} 
+ \|\Phi\|_{L^\8}\|\cov_A\psi\|_{L^{\sharp,2}(U')}
\\
&\leq
c\|\cov_A\Phi\|_{L^2_{1,A}}\|\psi\|_{L^2_{1,A}(U')} 
+ c\|\Phi\|_{L^\8}\|\cov_A\psi\|_{L^q(U')}
\\
&\leq
c\|\Phi\|_{C^0\cap L^2_{2,A}}\|\psi\|_{L^q_{1,A}(U')}.
\end{aligned}
\end{equation}
Substituting the estimate \eqref{eq:LpBoundWorstTerm} into the inequality
\eqref{eq:EasyEstdA*v}, via the expression \eqref{eq:LapvInTermsOfapsi} for
$d_A^+d_A^{+,*}v$, then yields  
\begin{equation}
\label{eq:L21dA*v-Final}
\begin{aligned}
\|d_A^{+,*}v\|_{L^2_{1,A}} 
&\le  
C(\|d_A^+a\|_{L^{\sharp,2}} + 
\|d_A^+(\cdot\Phi)^*\psi\|_{L^{\sharp,2}} + \|v\|_{L^2})
\\
&\le  
C(\|d_A^+a\|_{L^{\sharp,2}} + 
\|\Phi\|_{C^0\cap L^2_{2,A}}\|\psi\|_{L^q_{1,A}(U')} + \|v\|_{L^2}).
\end{aligned}
\end{equation}
Therefore, inequalities \eqref{eq:IntermedL21a-1} and \eqref{eq:L21dA*v-Final}
imply that
\begin{equation}
\label{eq:IntermedL21a-2}
\|a\|_{L^2_{1,A}} 
\le
C(\|d_A^+a\|_{L^{\sharp,2}} + \|\psi\|_{L^q_{1,A}(U')} + \|v\|_{L^2}).
\end{equation}
To estimate the $L^{\sharp,2}$ norm of $d_A^+a$ in the right-hand side of
inequality \eqref{eq:IntermedL21a-2}, note that
from the definition \eqref{eq:LAPhi} of $d_{A,\Phi}^1$ we have
$$
d_{A,\Phi}^1\left(\begin{matrix} a \\ \phi \end{matrix}\right)
= 
\left(\begin{matrix} d_A^+a - (\Phi\otimes\phi^*+\phi\otimes\Phi^*)_{00} \\
a\cdot\Phi + D_A\phi \end{matrix}\right),
$$
and thus,
\begin{equation}
\label{eq:L2sharpdAplusa-Final}
\|d_A^+a\|_{L^{\sharp,2}} 
\le 
\|d_{A,\Phi}^1(a,\phi)\|_{L^{\sharp,2;2}} 
+ c\|\Phi\|_{L^\8}\|\phi\|_{L^{\sharp,2}}.
\end{equation}
Choosing $4/3 \leq p \leq 2$ via $1/p = 1/4+1/q$ when $2\leq q\leq 4$, and
noting that
\begin{align*}
\phi
&=
-(\Phi\otimes(\cdot)^*+(\cdot)\otimes\Phi^*)_{00}^*v + D_A^*\psi,
\\
&= \Phi\otimes v + D_A^*\psi \quad\text{(schematically)},
\end{align*}
we have (schematically)
$$
D_AD_A^*\psi
=
\cov_A\Phi\otimes v + \Phi\otimes\cov_A v + D_A\phi,
$$
and so
\begin{equation}
\label{eq:LpDiracLappsi}
\begin{aligned}
\|D_AD_A^*\psi\|_{L^p(U)}
&\leq
\|\cov_A\Phi\|_{L^4}\|v\|_{L^q} + \|\Phi\|_{L^\8}\|\cov_Av\|_{L^p}
+ \|D_A\phi\|_{L^p}
\\
&\leq
c\|\Phi\|_{C^0\cap L^2_{2,A}}\|v\|_{L^2_{1,A}} + \|D_A\phi\|_{L^p},
\end{aligned}
\end{equation}
and thus, recalling that $U'$ obeys an interior cone condition with
geometric cone independent of the (small) size of the ball radii and
confining $q$ to $2<q\leq 4$,
\begin{align*}
\|\psi\|_{L^q_{1,A}(U')} 
&\le 
c\|\psi\|_{L^p_{2,A}(U')} \quad\text{(by \cite[Theorem V.5.4]{Adams})}
\\
&\le 
C(\|D_AD_A^*\psi\|_{L^p(U)} + \|\psi\|_{L^2(U)}) 
\quad \text{(by Theorem \ref{thm:Lp2SpinorAnnulusEst} 
and \cite{FeehanKato})}
\\
&\le C(\|D_A\phi\|_{L^p} + \|v\|_{L^2_{1,A}} + \|\psi\|_{L^2})
\quad\text{(by \eqref{eq:LpDiracLappsi})},
\end{align*}
and consequently
\begin{equation}
\label{eq:Lp1psiU}
\|\psi\|_{L^q_{1,A}(U')} 
\leq
C(\|D_A\phi\|_{L^p} + \|v\|_{L^2_{1,A}} + \|\psi\|_{L^2}).
\end{equation}
Therefore, combining inequalities
\eqref{eq:IntermedL21a-2}, \eqref{eq:L2sharpdAplusa-Final}, and
\eqref{eq:Lp1psiU}, we obtain
\begin{equation}
\label{eq:IntermedL21a-Final}
\|a\|_{L^2_{1,A}} 
\le 
C(\|d_{A,\Phi}^1(a,\phi)\|_{L^{\sharp,2;2}} 
+ \|\phi\|_{L^{\sharp,2}} + \|D_A\phi\|_{L^p}
+ \|v\|_{L^2_{1,A}} + \|\psi\|_{L^2}),
\end{equation}
Lemma \ref{lem:L21AEstphi} and the expression for $d_{A,\Phi}^1(a,\phi)$ gives
\begin{align*}
\|\phi\|_{L^2_{1,A}} 
&\le 
C(\|D_A\phi\|_{L^2} + \|\phi\|_{L^2}) 
\\
&\le 
C(\|d_{A,\Phi}^1(a,\phi)\|_{L^2} + \|\Phi\|_{L^\8}\|a\|_{L^2} 
+ \|\phi\|_{L^2}). 
\end{align*}
Thus,
\begin{equation}
\label{eq:L21phi-Final}
\|\phi\|_{L^2_{1,A}} 
\le 
C(\|d_{A,\Phi}^1(a,\phi)\|_{L^2} + \|(a,\phi)\|_{L^2}).
\end{equation}
Therefore, from the above $L^2_{1,A}$ estimates for $a$ 
in \eqref{eq:IntermedL21a-Final} and for $\phi$ in \eqref{eq:L21phi-Final},  
we see that
\begin{equation}
\label{eq:L21aphi-AlmostFinal}
\|(a,\phi)\|_{L^2_{1,A}} 
\le 
C(\|d_{A,\Phi}^1(a,\phi)\|_{L^{\sharp,2;2}} + \|(a,\phi)\|_{L^2}
+ \|(v,\psi)\|_{L^2} 
+\|\phi\|_{L^{\sharp,2}} + \|v\|_{L^2_{1,A}}). 
\end{equation}
Now $\|\phi\|_{L^{\sharp,2}} \leq c\|\phi\|_{L^2_{1,A}}$ and so the term
$\|\phi\|_{L^{\sharp,2}}$ on the right-hand side of inequality
\eqref{eq:L21aphi-AlmostFinal} 
can be replaced by $\|\phi\|_{L^2}$ via the bound
\eqref{eq:L21phi-Final}. Similarly, $\|v\|_{L^2_{1,A}} \leq
C(\|d_A^{+,*}v\|_{L^2} + \|v\|_{L^2})$ by Lemma \ref{lem:L21AEstv} and so the
interpolation inequality,
$$
\|d_A^{+,*}v\|_{L^2} \leq \eps\|d_A^+d_A^{+,*}v\|_{L^2} + \eps^{-1}\|v\|_{L^2},
$$
together with inequality \eqref{eq:L21dA*v-Final} for
$\|d_A^{+,*}v\|_{L^2_{1,A}}$, inequality \eqref{eq:L2sharpdAplusa-Final} for
$\|d_A^+a\|_{L^{\sharp,2}}$, inequality \eqref{eq:Lp1psiU} for
$\|\psi\|_{L^p_{1,A}(U')}$ , rearrangement for $\eps$ small, and the
preceding bound for $\|\phi\|_{L^{\sharp,2}}$, allow us to replace
$\|v\|_{L^2_{1,A}}$ by $\|v\|_{L^2}$ on the right-hand side of inequality
\eqref{eq:L21aphi-AlmostFinal}. Hence, the estimate
\eqref{eq:L21aphi-AlmostFinal} takes the form
$$
\|(a,\phi)\|_{L^2_{1,A}} 
\le 
C(\|d_{A,\Phi}^1(a,\phi)\|_{L^{\sharp,2;2}} + \|(a,\phi)\|_{L^2}
+ \|(v,\psi)\|_{L^2}). 
$$
Just as in equation \eqref{eq:EasyEstdA*v} we have
$$
\|(a,\phi)\|_{L^2}
\leq 
\frac{1}{2}(\|d_{A,\Phi}^1d_{A,\Phi}^{1,*}(v,\psi)\|_{L^2}
+ \|(v,\psi)\|_{L^2}),
$$
and combining the last two inequalities completes the proof.
\end{proof}

\begin{cor}
\label{cor:L21AEstPAaphi}
Continue the hypotheses and notation of Theorem \ref{thm:L21AEstPAaphi}.
Then for all $(\xi,\varphi)\in L^{\sharp,2}(\Lambda^+\otimes \fg_E)\oplus
L^2(V^-)$, 
\begin{equation}
\label{eq:L21AEstPAaphi}
\|P_{A,\Phi,\mu}(\xi,\varphi)\|_{L^2_{1,A}(X)} 
\le 
C\|(\xi,\varphi)\|_{L^{\sharp,2;2}(X)},
\end{equation}
where $C$ now also depends continuously on $0<\mu<\8$.
\end{cor}

\begin{proof}
We set $(v,\psi)=G_{A,\Phi,\mu}(\xi,\varphi)$, so that 
$$
d_{A,\Phi}^1d_{A,\Phi}^{1,*}(v,\psi)=\Pi_{A,\Phi,\mu}^\perp(\xi,\varphi) 
\quad\text{and}\quad 
(a,\phi)=d_{A,\Phi}^{1,*}(v,\psi)=P_{A,\Phi,\mu}(\xi,\varphi).
$$ 
For the $L^2_{1,A}$ estimate, note that
\begin{align*}
\|(a,\phi)\|_{L^2_{1,A}} 
&\le 
C(\|d_{A,\Phi}^1d_{A,\Phi}^{1,*}(v,\psi)\|_{L^{\sharp,2;2}} 
+ \|(v,\psi)\|_{L^2})
\quad\text{(by Theorem \ref{thm:L21AEstPAaphi})}
\\
&\leq
C(\|d_{A,\Phi}^1d_{A,\Phi}^{1,*}(v,\psi)\|_{L^{\sharp,2;2}} 
+ \mu^{-1}\|d_{A,\Phi}^1d_{A,\Phi}^{1,*}(v,\psi)\|_{L^2}) 
\\
&= 
C\|\Pi_{A,\Phi,\mu}^\perp(\xi,\varphi)\|_{L^{\sharp,2;2}} 
\\
&\leq
C\left(\|(\xi,\varphi)\|_{L^{\sharp,2;2}} 
+ \|\Pi_{A,\Phi,\mu}(\xi,0)\|_{L^4}
+ \|\Pi_{A,\Phi,\mu}(0,\varphi)\|_{L^2}\right) 
\quad\text{(by \cite[Lemma 4.1]{FeehanSlice})}
\\
&\leq
C\left(\|(\xi,\varphi)\|_{L^{\sharp,2;2}} 
+ \|\Pi_{A,\Phi,\mu}(\xi,0)\|_{L^2}
+ \|(0,\varphi)\|_{L^2}\right) 
\quad\text{(by Lemma \ref{lem:L21GlobalLocalEstDe2APhiEvec})}
\\
&\le 
C\|(\xi,\varphi)\|_{L^{\sharp,2;2}}
\quad\text{(by Inequality \eqref{eq:L21GlobalLocalEstDe2APhiEProj})}.
\end{align*}
This completes the proof.
\end{proof}

\subsection{Existence of solutions to the extended PU(2) monopole equations}
\label{subsec:ExistSolnExtASD}
We now turn to the question of existence and uniqueness of solutions
$(v,\psi)$ to equation \eqref{eq:QuickExtPUMonEqnForvpsi} or one of its
equivalent forms, namely equations
\eqref{eq:LongSecOrderExtPUMonEqnForvpsi} or
\eqref{eq:RegLongSecOrderExtPUMonEqnForvpsi}.
We first recall the following elementary fixed-point result \cite[Lemma
7.2.23]{DK}.

\begin{lem}
\label{lem:FixedPoint}
Let $q:\fB\to\fB$ be a continuous map on a Banach space $\fB$ with $q(0)=0$
and  
$$
\|q(x_1)-q(x_2)\|\le K\left(\|x_1\|+\|x_2\|\right)\|x_1-x_2\|
$$
for some positive constant $K$ and all $x_1,x_2$ in $\fB$. Then for each
$y$ in $\fB$ with $\|y\|<1/(10K)$ there is a unique $x$ in $\fB$ such that
$\|x\|\le 1/(5K)$ and $x+q(x)=y$.
\end{lem}

With the aid of Lemma \ref{lem:FixedPoint}, we obtain existence and
uniqueness for solutions to equation \eqref{eq:QuickExtPUMonEqnForvpsi}. 

\begin{thm}
\label{thm:ExtPUMonExist}
Let $(X,g)$ be a closed, oriented, $C^\8$ Riemannian four-manifold.
There are constants $C$ and $\eps(C)$, with the dependence indicated in
condition \eqref{eq:L21AEstPAaphiConstant}, such that the following holds.
Let $(A,\Phi)$ be an $L^2_4$ pair on $(\fg_E,V^+)$ produced by the splicing
construction of \S \ref{sec:Splicing}, with 
\begin{equation}
\label{eq:ExtPUMonExistCondition}
\fS(A,\Phi)\|_{L^{\sharp,2}(X)} < \eps.
\end{equation}
Then there is a unique solution
$(v,\psi)\in C^0\cap L^2_2(\Lambda^+\otimes\fg_E)\oplus L^2_2(V^-)$ to
the extended $\PU(2)$-monopole equation 
\eqref{eq:QuickExtPUMonEqnForvpsi} such that $\Pi_{A,\Phi,\mu}(v,\psi)=0$ and
\begin{equation}
\|d_{A,\Phi}^{1,*}(v,\psi)\|_{L^2_{1,A}(X)}
\le 
C\|\fS(A,\Phi)\|_{L^{\sharp,2}(X)}. 
\end{equation}
Moreover, if $(A,\Phi)$ is an $L^2_l$ pair, for any $l\geq 3$, then the
solution $(v,\psi)$ is contained in $L^2_{l+1}(\Lambda^+\otimes\fg_E)\oplus
L^2_{l+1}(V^-)$.
\end{thm}

\begin{rmk}
\label{rmk:ExtPUMonExistCondition}
The condition \eqref{eq:ExtPUMonExistCondition} is met by choosing a small
enough relatively open submanifold $\sU_{\ell,\mu}\subset \sC(\ft_\ell)$ and
small enough constant $\lambda_0(C)$ (Proposition
\ref{prop:CutoffMonoPairEst}). 
\end{rmk}

\begin{pf}
We try to solve equation \eqref{eq:LongSecOrderExtPUMonEqnForvpsi} for
solutions of the form $(v,\psi)=G_{A,\Phi,\mu}(\xi,\varphi)$, where
$(\xi,\varphi)\in L^{\sharp,2}(\Lambda^+\otimes\fg_E)$ and $G_{A,\Phi,\mu}$
is the Green's operator for the Laplacian $d_{A,\Phi}^1d_{A,\Phi}^{1,*}$ on
$\Ran\Pi_{A,\Phi,\mu}^\perp$. Therefore, as explained in \S
\ref{subsec:MonoEqns}, we seek solutions $(\xi,\varphi)$ to equation
\eqref{eq:IntegralPUMonEqnForxipsi}.  We now apply Lemma
\ref{lem:FixedPoint} to equation \eqref{eq:IntegralPUMonEqnForxipsi}, with
$$
\fB=L^{\sharp,2}(\Lambda^+\otimes\fg_E)\oplus L^2(V^-),
$$
choosing $y=-\fS(A,\Phi)$ and
$$
q((\xi,\varphi)) 
= 
\{d_{A,\Phi}^{1,*}G_{A,\Phi,\mu}(\xi,\varphi), 
d_{A,\Phi}^{1,*}G_{A,\Phi,\mu}(\xi,\varphi)\}, 
$$
so our goal is to solve
$$
(\xi,\varphi) + q((\xi,\varphi)) 
= 
-\fS(A,\Phi)
\quad\text{on }L^{\sharp,2}(\Lambda^+\otimes\fg_E)\oplus L^2(V^-). 
$$
For any $(\xi_1,\varphi_1),(\xi_2,\varphi_2)\in 
L^{\sharp,2}(\Lambda^+\otimes\fg_E)\oplus L^2(V^-)$, 
the estimate in Theorem \ref{thm:L21AEstPAaphi} yields
$$
\|d_{A,\Phi}^{1,*}G_{A,\Phi,\mu}(\xi_i,\varphi_i)\|_{L^2_{1,A}} 
\le 
C\|(\xi_i,\varphi_i)\|_{L^{\sharp,2}}.
$$
Therefore, applying the preceding estimates, H\"older's inequality, Lemma
\ref{lem:Kato}, and the $L^\sharp$-family of embedding and multiplication
results of \cite[Lemmas 4.1 \& 4.3]{FeehanSlice}, we obtain
$$
\|q((\xi_1,\varphi_1))-q((\xi_2,\varphi_2))\|_{L^{\sharp,2}}
\le 
K_2\left(\|(\xi_1,\varphi_1)\|_{L^{\sharp,2}} 
+ \|(\xi_2,\varphi_2)\|_{L^{\sharp,2}}\right) 
\|(\xi_1,\varphi_1)-(\xi_2,\varphi_2)\|_{L^{\sharp,2}},
$$
where $K_2 = cZ_2$ with $c$ a universal constant and $Z_2$ the constant of 
Theorem \ref{thm:L21AEstPAaphi}. 

Thus, provided $\eps \le (10K_2)^{-1}$, 
we have $\|\fS(A,\Phi)\|_{L^{\sharp,2}} <
(10K_2)^{-1}$, and Lemma \ref{lem:FixedPoint} implies that there is a
unique solution 
$(\xi,\varphi)\in L^{\sharp,2}(\Lambda^+\otimes\fg_E)\oplus L^2(V^-)$ 
to equation \eqref{eq:IntegralPUMonEqnForxipsi}
such that $\|(\xi,\varphi)\|_{L^{\sharp,2}}\le (5K_2)^{-1}$.
The $L^{\sharp,2}$ estimate for $(\xi,\varphi)$,
\begin{equation}
\label{eq:LSharp2EstForXiVphi}
\|(\xi,\varphi)\|_{L^{\sharp,2}} \le 2\|\fS(A,\Phi)\|_{L^{\sharp,2}},
\end{equation}
follows from Theorem \ref{thm:L21AEstPAaphi} and a
standard rearrangement argument.
To obtain the stated $L^2_{1,A}$ estimates for $d_{A,\Phi}^{1,*}(v,\psi)$,
observe that 
\begin{equation}
\label{eq:LaplacianOfvpsiWithEigenvalCutoff}
d_{A,\Phi}^{1,*}(v,\psi) 
= 
d_{A,\Phi}^{1,*}G_{A,\Phi,\mu}(\xi,\varphi) 
= 
P_{A,\Phi,\mu}(\xi,\varphi),
\end{equation}
so that the $L^2_{1,A}$ bound is obtained from
\begin{align*}
\|d_{A,\Phi}^{1,*}(v,\psi)\|_{L^2_{1,A}} 
&=
\|P_{A,\Phi,\mu}(\xi,\varphi)\|_{L^2_{1,A}} 
\leq 
C\|(\xi,\varphi)\|_{L^{\sharp,2}}
\quad\text{(by Corollary \ref{cor:L21AEstPAaphi})}
\\
&\le C\|\fS(A,\Phi)\|_{L^{\sharp,2}} 
\quad\text{(by estimate \eqref{eq:LSharp2EstForXiVphi}).}
\end{align*}
The final regularity assertion follows from
Proposition \ref{prop:RegularityTaubesSolution}, our regularity result for
gluing solutions to the (extended) $\PU(2)$-monopole equations. 
\end{pf}

\subsection{Completion of proof of main theorem}
\label{subsec:ProofGluingTheorem}
At this stage, we are ready to complete the proof of 
Theorem \ref{thm:GluingTheorem1}, namely the existence of gluing and
obstruction maps. The conclusion of our proof of Theorem
\ref{thm:GluingTheorem1} is an almost immediate
consequence of Theorem \ref{thm:ExtPUMonExist}.

\begin{proof}[Completion of proof of Theorem \ref{thm:GluingTheorem1}]
According to Theorem \ref{thm:ExtPUMonExist} there are well-defined,
$\sG_E$-equivariant maps $\bdelta_\mu$ and $\bvarphi_\mu$ on the open
subset of $\tsC(\ft)$,
\begin{equation}
\label{eq:IntrinsicDomainTaubesGluing}
\tsC(\ft,\mu,\eps)
:=
\{(A,\Phi)\in\tsC(\ft): \text{\eqref{eq:DeformCondn} holds}\},
\end{equation}
where the pairs $(A,\Phi)$ are required to satisfy the open constraints
\begin{equation}
\label{eq:DeformCondn}
\begin{gathered}
\Spec(d_{A,\Phi}^1d_{A,\Phi}^{1,*}) \subset [0,\thalf\mu)\cup (\mu,\8),
\\
\|\fS(A,\Phi)\|_{L^{\sharp,2;2}(X)} < \eps,
\\
\|P_{A,\Phi,\mu}\|_A < C < \8,
\end{gathered}
\end{equation}
where $C=C(\sU_{\ell,\mu},\kappa,\mu,g)$ is a uniform positive constant and
$\|\cdot\|_A$ is the operator norm for $\Hom(L^{\sharp,2;2}(X), L^2_{1,A}(X))$.

On the open, $\sG_E$-invariant subset $\tsC(\ft,\mu,\eps)\subset
\tsC(\ft)$, the maps $\bdelta_\mu$ and $\bvarphi_\mu$ are given by
\begin{align*}
&\bdelta_\mu:(A,\Phi)\mapsto (A,\Phi) + d_{A,\Phi}^{1,*}(v,\psi),
\\
&\bvarphi_\mu:(A,\Phi)\mapsto \Pi_{A,\Phi,\mu}\fS(\bdelta_\mu(A,\Phi)),
\end{align*}
taking values in $\tsC(\ft)$ and $L^2(\Lambda^+\otimes\fg_E)\oplus
L^2(V^-)$, respectively, where $(v,\psi)\in\Ran \Pi_{A,\Phi,\mu}^\perp$ is
the unique solution to the extended $\PU(2)$ monopole equations 
\eqref{eq:LongSecOrderExtPUMonEqnForvpsi} produced by Theorem
\ref{thm:ExtPUMonExist}, given the approximate $\PU(2)$ monopole $(A,\Phi)$
in the image of $\bgamma'_{\mu,\Sigma}$.

The $\sG_E$-equivariant maps $\bdelta_\mu$ and $\bvarphi_\mu$ descend to an
$S^1$-equivariant map and section, 
\begin{equation}
\begin{aligned}
\label{eq:TaubesDeformationMapnSection}
&\bdelta_\mu: \sC(\ft,\mu,\eps) \to \sC(\ft), 
\\
&\bvarphi_\mu: \sC(\ft,\mu,\eps) \to \fV_\mu
\end{aligned}
\end{equation}
where $\sC(\ft,\mu,\eps) = \tsC(\ft,\mu,\eps)/\sG_E$ and
$\fV_\mu\to \sC(\ft,\mu)$ is the 
continuous, finite-rank, pseudovector bundle 
$$
\fV_\mu
:= 
\{((A,\Phi),(w,s)): (A,\Phi)\in\tsC(\ft,\mu)
\text{ and }(w,s)\in\Ran\Pi_{A,\Phi,\mu}\}/\sG_E,
$$
where we recall that $\Ran\Pi_{A,\Phi,\mu}\subset
L^2_{k-1}(X,\Lambda^+\otimes\fg_E)\oplus L^2_{k-1}(V^-)$. Here, the spaces
$\tsC(\ft,\mu)$ and $\sC(\ft,\mu)$ are defined by simply
omitting the conditions in definition
\eqref{eq:IntrinsicDomainTaubesGluing} involving $\eps$ and the constant $C$.
The map $\bdelta_\mu$ and section $\bvarphi_\mu$ have the desired property that
$$
\bdelta_\mu(\bvarphi_\mu^{-1}(0)) \subset M(\ft),
$$
and are smooth respect to the quotient $L^2_k$ topology upon restriction to
$\sC^0(\ft)$: smoothness follows by the same argument
as used in the proof of smoothness in the case of the Banach-space inverse
function theorem
\cite{AMR}, where the existence of an inverse is established via a
contraction-mapping argument, just as we employ here.  

The gluing map, gluing-data obstruction bundle, and gluing-data obstruction
section are then defined by
\begin{equation}
\label{eq:GluingMapObstructionDefnBundle}
\begin{aligned}
\bga_{\mu,\Sigma} 
&= \bdelta_\mu\circ\bga'_{\mu,\Sigma}:
\bGl^+(\sU_{\ell,\mu},\Sigma,\lambda_0) \to \sC(\ft),
\\
\bchi_{\mu,\Sigma} &= \bvarphi_\mu\circ\bga_{\mu,\Sigma}: 
\bGl^+(\sU_{\ell,\mu},\Sigma,\lambda_0) \to \Xi_\mu,
\\
\Xi_\mu &= \bga^*_{\mu,\Sigma}\fV_\mu,
\end{aligned}
\end{equation}
where the precompact subset $\sU_{\ell,\mu}\subset \sC(\ft_\ell)$ is chosen
according to the hypotheses of Theorem \ref{thm:GluingTheorem1}.
The estimates in Proposition \ref{prop:CutoffMonoPairEst} for
$\fS(A,\Phi)$, the eigenvalue estimates in Theorem \ref{thm:H2SmallEval}
for the Laplacian $d_{A,\Phi}^1d_{A,\Phi}^{1,*}$,
and the operator bounds for $P_{A,\Phi,\mu}$ in Corollary
\ref{cor:L21AEstPAaphi} ensure that pairs $(A,\Phi)$
in the image of $\bga'_{\mu,\Sigma}$ obey the open constraints
\eqref{eq:DeformCondn} for a suitable choice of positive constants $\eps$,
$\mu$ and $C$. Hence, the image of $\bga'_{\mu,\Sigma}$, when
$\bga'_{\mu,\Sigma}$ is restricted to
$\bGl^+(\sU_{\ell,\mu},\Sigma,\lambda_0)$, 
is contained in $\sC(\ft,\mu,\eps)$ for small enough $\lambda_0$ and
$\eps_\ell$ defining the relatively open, precompact submanifold
$\sU_{\ell,\mu} \subset \sC(\ft_\ell,\mu,\eps_\ell)$.
\end{proof}

The rank of $\fV_\mu$ varies over $\sC(\ft,\mu,\eps)$, so that
is why we only refer to it as being a `pseudovector bundle'. However, the
construction of $\bGl^+(\sU_{\ell,\mu},\Sigma,\lambda_0)$ (via the
choice of family $\sU_{\ell,\mu}$ of background pairs $(A_0,\Phi_0)$) in \S
\ref{sec:Splicing} and Theorem
\ref{thm:H2SmallEval} ensure that its image under $\bga'_{\mu,\Sigma}$ is
contained in an open subset of $\sC(\ft,\mu,\eps)$ over which
$\fV_\mu$ has constant rank and is thus a vector bundle when restricted to
this open subset.

%end of file

\ifx\undefined\bysame
\newcommand{\bysame}{\leavevmode\hbox to3em{\hrulefill}\,}
\fi

%end of file
\end{document}